\documentclass[twoside, 10pt]{article}

\usepackage[hmargin = 0.75in, height =9.5in]{geometry}
\usepackage{amssymb,amsthm, amsfonts,amsmath,mathrsfs, amsxtra, url, hyperref}

\usepackage{fancyhdr}
\fancypagestyle{plain}{
 \fancyhf{}

}

\renewcommand{\theequation}{\thesection.\arabic{equation}}
 \numberwithin{equation}{section}

\newtheorem {thm}{Theorem}[section]
\newtheorem {prop}{Proposition}[section]
\newtheorem {lemm}{Lemma}[section]
\newtheorem {deff}{Definition}[section]
\newtheorem {cor}{Corollary}[section]
\newtheorem {rem}{Remark}[section]

\newtheorem {eg}{Example}[section]

\makeatletter
\newenvironment{highlightequation}{%
  \def\tagform@##1{\maketag@@@{(\ignorespaces##1\unskip\@@italiccorr*)}}%
  \ignorespaces
}{%
  \def\tagform@##1{\maketag@@@{(\ignorespaces##1\unskip\@@italiccorr)}}%
  \ignorespacesafterend
}
\makeatother

\def\ba{\begin{array}}
\def\ea{\end{array}}
\def\bea{\begin{eqnarray}}
\def\eea{\end{eqnarray}}
\def\beas{\begin{eqnarray*}}
\def\eeas{\end{eqnarray*}}
\def\bi{\begin{itemize}}
\def\ei{\end{itemize}}
\def\bc{\begin{cases}}
\def\ec{\end{cases}}

\def\bhe{\begin{highlightequation}  }
\def\ehe{\end{highlightequation}  }

\def\a{\alpha}

\def\d{\delta}
\def\e{\varepsilon}
\def\eps{\epsilon}
\def\z{\zeta}
\def\k{\kappa}
\def\l{\lambda}

\def\si{\sigma}
\def\vs{\varsigma}

\def\o{\omega}
\def\vf{\varphi}

\def\D{\Delta}
\def\G{\Gamma}
\def\L{\Lambda}
\def\O{\Omega}

\def\Th{\Theta}
\def\U{\Upsilon}

\def\bF{{\bf F}}

\def\cA{{\cal A}}

\def\cD{{\cal D}}
\def\cE{{\cal E}}
\def\cF{{\cal F}}
\def\cG{{\cal G}}
\def\cH{{\cal H}}
\def\cI{{\cal I}}

\def\cK{{\cal K}}
\def\cL{{\cal L}}

\def\cN{{\cal N}}
\def\cO{{\cal O}}

\def\cR{{\cal R}}

\def\cT{{\cal T}}
\def\cU{{\cal U}}
\def\cV{{\cal V}}

\def\cX{{\cal X}}
\def\cY{{\cal Y}}

\def\hC{\mathbb{C}}

\def\hE{\mathbb{E}}

\def\hH{\mathbb{H}}

\def\hK{\mathbb{K}}

\def\hM{\mathbb{M}}
\def\hN{\mathbb{N}}

\def\hQ{\mathbb{Q}}
\def\hR{\mathbb{R}}
\def\hS{\mathbb{S}}

\def\sB{\mathscr{B}}
\def\sC{\mathscr{C}}
\def\sD{\mathscr{D}}

\def\sN{\mathscr{N}}
\def\sO{\mathscr{O}}

\def\sX{\mathscr{X}}

\def\fA{\mathfrak{A}}

\def\fC{\mathfrak{C}}

\def\fM{\mathfrak{M}}
\def\fN{\mathfrak{N}}

\def\fX{\mathfrak{X}}
\def\fY{\mathfrak{Y}}

\def\fq{\mathfrak{q}}

\def\fm{\mathfrak{m}}

\def\fy{\mathfrak{y}}

\def\fx{\mathfrak{x}}

\def\fg{\mathfrak{g}}

\def\fn{\mathfrak{n}}

\def\ft{\mathfrak{t}}
\def\fb{\mathfrak{b}}
\def\fj{\mathfrak{j}}
\def\fc{\mathfrak{c}}
\def\fs{\mathfrak{s}}

\def\({\textnormal{(}}
\def\){\textnormal{)}}
\def\[{[\n[}
\def\]{]\n]}

\def\ran{\rangle}

\def\no{\noindent}

\def\ss{\smallskip}

\def\q{\quad}
\def\qq{\qquad}

\def\n{\negthinspace}
\def\dn{\n \n}
\def\tn{\n \n \n}

\def\ol{\overline}
\def\ul{\underline}
\def\ua{\mathop{\uparrow}}
\def\da{\mathop{\downarrow}}

\def\wt{\widetilde}
\def\wh{\widehat}
\def\fra{\mathfrak{a}}

\def\pas{{\hbox{$P_t-$a.s.}}}

\def\hb{\hbox}

\def\cd{\cdot}
\def\cds{\cdots}

\def\fa{\,\forall \,}
\def\pa{\partial}
\def\es{\emptyset}

\def\b1{{\bf 1}}
\def\qed{\hfill $\Box$ \medskip}

\def\ti{\n \times \n}
\def\oti{\n \otimes \n}
\def\otii{\n \otimes_{t_i} \n}
\def\otis{\n \otimes_s \n}

\def\df{\n := \n}
\def\ls{\n \le \n}
\def\gs{\n \ge \n}
\def\={\n = \n}
\def\+{\n + \n}
\def\-{\n - \n}
\def\ins{\n \in \n}
\def\ld{\n \land \n}
\def\ve{\n \vee \n}
\def\sb{\n \subset \n}
\def\>{\n > \n}
\def\<{\n < \n}
\def\Cp{\n \cap \n}
\def\cp{\n \cup \n}

\def\limsup{\mathop{\ol{\rm lim}}}

\newcommand{\lsup}[1]{\underset{#1}{\limsup}}

\newcommand{\lmt}[1]{\underset{#1}{\lim}}
\newcommand{\lmtu}[1]{\underset{#1}{\lim} \n \ua \,}
\newcommand{\lmtd}[1]{\underset{#1}{\lim} \n \da \,}

\begin{document}

 \title{\bf   Dynamic Programming Principles for \\
  Optimal   Stopping   with Expectation Constraint
    }

\author{
 Erhan Bayraktar\thanks{ \noindent Department of
  Mathematics, University of Michigan, Ann Arbor, MI 48109; email:
{\tt erhan@umich.edu}.}  \thanks{E. Bayraktar is supported in part by the National Science Foundation
 under  DMS-1613170,
 and in part by the Susan M. Smith Professorship.  Any opinions, findings, and conclusions
 or recommendations expressed in this material are
those of the authors and do not necessarily reflect the views of the National Science Foundation.} $\,\,$,
$~~$Song Yao\thanks{
\noindent Department of
  Mathematics, University of Pittsburgh, Pittsburgh, PA 15260; email: {\tt songyao@pitt.edu}. }
  \thanks{S. Yao is supported   by the National Science
Foundation under DMS-1613208.
} }

\date{}

\maketitle

 \begin{abstract}

We analyze  an   optimal stopping problem with a constraint on the expected cost.
When the reward function and cost function are Lipschitz continuous in state variable,
we show that the value     of such an optimal stopping problem is a continuous function
in current state  and in budget level. Then we  derive  a dynamic programming principle (DPP)
for the value function   in  which  the conditional expected cost acts as an additional state process.
As the optimal stopping problem with expectation constraint can be transformed to
a stochastic optimization problem with supermartingale controls,
we explore   a  {\it second}  DPP  of the value function
and  thus resolve an open question recently raised in [S. Ankirchner, M. Klein, and T. Kruse,
A verification theorem for optimal stopping problems with expectation constraints, Appl. Math. Optim., 2017, pp. 1-33].
Based on these two DPPs, we characterize the value function   as a viscosity solution to the related fully non-linear parabolic Hamilton-Jacobi-Bellman equation.

 \ss \no {\bf MSC:}\;  60G40,  49L20, 93E20,  49L25.

 \ss \no   {\bf Keywords:}\; Optimal stopping  with expectation constraint,
 dynamic programming principle,  shifted processes, 
 shifted stochastic differential equations, flow property,
 stochastic optimization with supermartingale controls,
 fully non-linear parabolic Hamilton-Jacobi-Bellman   equation,
 viscosity solution, Monge-Amp\`ere type equation.

\end{abstract}

  \section{Introduction}

  In this article, we analyze  a continuous-time 
  optimal stopping problem with expectation constraint on  the accumulated cost.
  Suppose that the game begins at time $t$ over the canonical space $\O^t$ of continuous paths.
  Under the Wiener measure $P_t$,
  the coordinator process $W^t \= \{W^t_s\}_{s \in [t,\infty)}$ of $\O^t$ is a Brownian motion.
     Let $\ol{\bF}^t   \= \big\{ \ol{\cF}^t_s \big\}_{s \in [t,\infty)}$
  be the $P_t-$augmentation of the filtration  generated by  $W^t$,
  and let the $\hR^l -$valued state flow $\cX^{t,x}$ evolve 
  from position $x \ins \hR^l $  according to a  stochastic differential equation
    \bea \label{FSDE2}
     \cX_s   \=  x \+ \int_t^s b  ( r,  \cX_r  ) dr \+ \int_t^s  \si (  r, \cX_r  )  \, dW^t_r ,  \q s \ins [t, \infty) .
     \eea

  We aim to maximize the sum $\cR(t,x,\tau)$ of a running reward $ \int_t^\tau   \n f (r,\cX^{t,x}_r) dr    $
  and a terminal reward $\pi  (  \tau   , \cX^{t,x}_\tau     )$
  by choosing an    $\ol{\bF}^t-$stopping time $\tau$,
  which, however, has to satisfy   a  budget constraint   $E_t  [ \int_t^\tau g(r, \cX^{t,x}_r) dr ] \ls y$.
  So the value of such a optimal stopping problem  with expectation constraint is in form of
  \bea \label{def_value_fcn}
 \cV(t,x,y)    \df   \underset{\tau \in \cT^t_x (y)   }{\sup} \,  E_t  \big[ \cR(t,x,\tau) \big]  ,
 \eea
 with $ \cT^t_x (y) \df \big\{\tau  \n : E_t  [ \int_t^\tau  g(r, \cX^{t,x}_r) dr ] \ls y \big\} $
 and $E_t[\cd] \= E_{P_t} [\cd]$.
 In particular, when the   cost rate $g (r,x) $ is a power function of $r$,
 the budget constraint specifies as a moment constraint on   stopping times.

 Kennedy \cite{Kennedy_1982} initiated  the study of optimal stopping problem  with expectation constraint.
 The author used a {\it Lagrange multiplier} method
  to reduce a discrete-time optimal stopping problem with first-moment constraint
  to an unconstrained optimal stopping problem
  and   showed that the optimal value of the dual problem
  is   equal to that of the primal problem.
  Since then, the Lagrangian technique has been prevailing
  in  research of   optimal stopping problems with expectation constraints.

 In the present paper, we develop  a new approach to analyze the optimal stopping problem with expectation constraint
 \eqref{def_value_fcn}. Our main contributions are obtaining the continuity of the value function $\cV$
 and establishing two dynamic programming principles (DPPs) for $\cV$.

 When reward/cost  functions $f,\pi,g $   are Lipschitz continuous in state variable $x$ and
 the cost function $g$ is non-degenerate in sense of (g3),
 we first demonstrate  over a general   probability setting that
 the value function   is continuous  in $(t,x,y)$
  by utilizing a priori estimates of the state process $\cX^{t,x}$
  and   delicately constructing   approximate stopping strategies (see Theorem \ref{thm_continuity}).
 This continuity result together with the properties of shifted processes then allow us to
  derive in Theorem \ref{thm_DPP}  a  DPP  for the value function $\cV$ over the canonical space:
\bea \label{eqn_DPP}
  \cV(t,x,y)  \=   \underset{\tau \in \cT^t_x (y)   }{\sup} \,
E_t  \bigg[ \b1_{\{\tau \le \z (\tau) \}}   \cR(t,x,\tau)
\+  \b1_{\{\tau > \z (\tau) \}} \Big(  \cV \big( \z (\tau), \cX^{t,x}_{\z (\tau)} , \cY^{t,x,\tau}_{\z (\tau)}  \big)
\+ \int_t^{\z (\tau)} \n f(r,\cX^{t,x}_r) dr  \Big)   \bigg] .
\eea
 Here the conditional expected cost $\cY^{t,x,\tau}_s \df E_t  \big[ \int_{\tau \land s}^\tau   g \big(r, \cX^{t,x}_r \big) dr \big|\ol{\cF}^t_s \big]$
 acts as an additional state process   and the intermediate horizon $\z$ can be a general $\ol{\bF}^t-$stopping time
 depending on the stopping rule $\tau$ we select.
 For the ``$\le$" part of \eqref{eqn_DPP}, we exploit the flow property of shifted stochastic differential equations (Proposition \ref{prop_FSDE_shift}) as well as the regular conditional probability distribution
 due to \cite{Stroock_Varadhan};
 while in the ``$\ge$" part, we carefully paste together local $\e-$optimal stopping strategies
 and utilize   the continuity of value function $\cV$.

Also, we can   transform   the optimal stopping problem with expectation constraint to
an unconstrained stochastic optimization problem whose controls are  supermartingales starting from budget level $y$:
Let $\fA_t (y)$ denote all uniformly integrable continuous supermartingales  $\a\=\{\a_s\}_{s \in [t,\infty)}$
with $\a_t \= y$.
As   shown in  Proposition  \ref{prop_surm}, for each nontrivial  $\tau \ins \cT^t_x (y)$
there exists $\a \ins \fA_t (y)$ such that
$\tau$ coincides with  the first hitting time  $\tau(t,x,\a)$  of the process
$ Y^{t,x,\a}_s \df \a_s \- \int_t^s g(r,\cX^{t,x}_r)dr  $, $s \in [t,\infty)$ to   $0$
\big(If $E_t [\int_t^\tau g(r,\cX^{t,x}_r)dr]\=y$,   $\a$ is indeed a true martingale\big).
So   the value function $\cV$ can be alternatively expressed as
$\cV( t, x,y )   \=   \underset{\a \in \fA_t (y) }{\sup} \,   E_t \big[ \cR \big(t, x, \tau (t,x,\a) \big) \big] $.
Correspondingly, we establish a second DPP for the value function $\cV$
over the canonical space  (Theorem \ref{thm_DPP2})
 \bea
  \cV (t,x,y) & \dn \= & \dn   \underset{\a   \in \fA_t (y)   }{\sup} \,
E_t  \bigg[ \b1_{\{\tau (t,x,\a) \le \z (\a) \}}  \cR \big(t, x, \tau (t,x,\a) \big) \nonumber \\
& \dn & \dn \q \qq \qq +  \b1_{\{\tau (t,x,\a) > \z (\a)\}} \Big(  \cV  \big( \z (\a), \cX^{t,x}_{\z (\a)} ,
 Y^{t,x,\a}_{\z (\a)}     \big)
\+ \int_t^{\z (\a)} \n f(r,\cX^{t,x}_r) dr \Big)   \bigg]   , \qq  \qq  \label{eqn_DPP2}
\eea
 and thus justify a postulate recently made by \cite{AKK_2015}  (see Remark 3.3 therein).
  Although the ``$\le$" part of \eqref{eqn_DPP2} can be easily deduced from \eqref{eqn_DPP}, the ``$\ge$" part
   entails an intricate pasting of approximately optimal supermartingale controls.

 In light of these two DPPs, we then show that  the  value function $\cV$
   of the optimal stopping problem with expectation constraint is a  viscosity solution  to
   a related   fully non-linear parabolic Hamilton-Jacobi-Bellman (HJB)  equation
 \bea \label{eq:d223}
 \q
 \begin{cases}
   -  \pa_t u (t,x,y)   \- \frac12 trace \big(  \si (t,x )   \n  \cd  \n   \si^T  (t,x )
   \n \cd \n  D^2_x u(t,x,y)    \big)   \- b^T (t,x )   \n  \cd  \n   D_x u (t,x,y)  \ss \\
   \qq  + g(t,x)    \pa_y u (t,x,y) \- \cH  u (t,x,y)  \= f(t,x) ,
   \q \fa (t,x,y) \ins (0,\infty)  \ti  \hR^l \ti (0,\infty) ,   \ss   \\
   u (t,x,0) \= \pi (t,x) , \q \fa (t,x ) \ins [0,\infty)  \ti  \hR^l ,
 \end{cases}
 \eea
 with the Hamiltonian $\cH  u (t,x,y) \df \underset{a  \in  \hR^d}{\sup}  \, \big\{ \frac12 |a|^2   \pa^2_y u (t,x,y)
 \+  ( D_x ( \pa_y u (t,x,y) )  )^T \dn \cd \n \si (t,x) \n \cd \n  a    \big\}$.
 As pointed out in \cite{Miller_C_2017a},
 the non-linear  HJB   equation \eqref{eq:d223} is a Monge-Amp\`ere type equation.

\no {\bf Relevant Literature.}
  Since  Arrow et al. \cite{ABG_1949} and  Snell \cite{Snell_1952},
  the general theory of (unconstrained) optimal stopping   has been plentifully developed  over decades.
  Expositions of this theory are presented in the monographs \cite{CRS_1971,Neveu_1975,Shiryayev_1978,El_Karoui_1981,Kara_Shr_MF,Peskir_Shiryaev_2006},
    which contain extensive bibliographies and references  to the literature.
  For the recent development of the  optimal stopping under  model uncertainty/non-linear expectations and the closely related controller-stopper-games, see  \cite{Karatzas_Sudderth_2001,Kara_Zam_2005, Follmer_Schied_2004, CDK-2006, Delbaen_2006, Kara_Zam_2008, Riedel_2009, OSNE1,OSNE2,OS_CRM,riedel2012,Bayraktar_Huang_2013,ETZ_2014,ROSVU,NZ_2015,RDOSRT,RDG}   among others.

   As to the optimal stopping with expectation constraint,
   the Lagrange multiplier method 
   introduced in \cite{Kennedy_1982}   was later  developed by many researches   (see e.g. \cite{Pontier_Szpirglas_1984,Moustakides_1986,LSMS_1995,Dokuchaev_1996,Balzer_Jansen_2002,Urusov_2005,Makasu_2009}),
  and has been applied to  various economic and financial problems
  such as  Markov decision processes with constrained stopping times \cite{Horiguchi_2001c,Horiguchi_2001b},
  non-exponential discounting and mean-variance portfolio optimization
   \cite{Pedersen_Peskir_2016,Pedersen_Peskir_2017}  and quickest detection problem  \cite{Peskir_2012}.

 Our stochastic control approach in deriving the second DPP  resembles
 those of two recent papers \cite{AKK_2015}, \cite{Miller_C_2017a}.
 By applying   the martingale representation to the conditional expected cost,
  Ankirchner et al. \cite{AKK_2015}   transformed
  the   optimal stopping problem  with expectation constraint
  to a   stochastic optimization problem in which the stochastic integral  of locally square-integrable controls
  is regarded as an additional state process.
  Miller  \cite{Miller_C_2017a} independently employed the same   method
  to address  the optimal stopping problem with first-moment constraint that
  is embedded in a time-inconsistent optimal stopping problem.
  The idea of expanding the state space by  the conditional probability/expectation process has also appeared
  in the literature dealing with stochastic target problems, see e.g. \cite{BEI_2009,BET_2009,BMN_2014,BER_2015,BPZ_2016}.

 Our paper is distinct from \cite{AKK_2015}, \cite{Miller_C_2017a} in four aspects:
 First,  we first obtain the continuity of the value function $\cV$,
 and using this establish the two DPPs \eqref{eqn_DPP} and \eqref{eqn_DPP2}, which were not addressed by them.
  Second, our value function $\cV$ takes the starting moment $t$ of the game as an input,
 so the related non-linear HJB equation \eqref{eq:d223}  is of parabolic type rather than elliptic type.
 Third,  we need the constraint $E_t [\int_t^\tau g(r,\cX^{t,x}_r) dr] \ls y$
 for the continuity and  the DPPs of the value function,
 although the auxiliary optimal stopping problem   considered in \cite{Miller_C_2017a}
   is subject to    constraint $E[\tau] \=y $ and
 the dynamic programming equation studied by  \cite{AKK_2015}  is for the value function $U$  
 of the optimal stopping with constraint $E[\int_0^\tau g(  X^x_r) dr] \= y$.
 See Remark \ref{rem_constraint} for a comparison of these two types of   constraints.
 Fourth, our discussion of related non-linear HJB equations seems different from theirs. 
 Our Theorem 5.1 obtains that the value function
 $\cV$ is a viscosity supersolution of \eqref{eq:d223},
 and is only a viscosity subsolution of \eqref{eq:d223} with the upper semi-continuous envelope $\ol{\cH} u $ of $\cH u$. 
 By assuming that the  value   $U$ is a smooth function satisfying the DPP,
 Proposition 3.4 of  \cite{AKK_2015}  showed that $U$ is a  supersolution to a similar non-linear HJB equation to \eqref{eq:d223},   and is further a  subsolution if the Hamiltonian is continuous (see Subsection 6.1 of \cite{AKK_2015}
   for an example of discontinuous Hamiltonian). 
   However, possible discontinuity of the Hamiltonian
  was not discussed in  \cite{Miller_C_2017a}.   

 Lately, the optimal stopping with constraint on the distribution of stopping time
 has   attracted a lot of research interests.
 Bayraktar and Miller \cite{Bayraktar_Miller_2016}
 studied the optimal stopping of a Brownian motion  with the restriction that
 the distribution of the stopping time must equal to
 a given measure consisting of finitely-many atoms.
 The applications of such a  distribution-constrained optimal stopping problem
 in mathematical finance include model-free superhedging with an outlook on volatility
 and inverse first-passage-time problem.
 Within a weak formulation on the canonical path space,
  Kallblad \cite{Kallblad_2017}   extended the distribution-constrained optimal stopping problem
  for a general target measure and for path-dependent cost functions.
  From the perspective of  mass transport,
  Beiglboeck et al. \cite{BEES_2016}  obtained a monotonicity principle for the  optimal stopping
  of a Brownian motion under distribution constraint,
  and thus  characterized the constrained optimal stopping rule
  as the first hitting time of a barrier in a suitable  phase space.
 Very recently, Ankirchner et al.  \cite{AKKK_2017} showed that
 for optimally stopping a one-dimensional  Markov process with first-moment constraint on stopping times,
 one only needs to  consider those stopping times at which
 the law of the Markov process   is a weighted sum of three Dirac measures.
There are also some other types of optimal stopping problems with constraints:
see \cite{DTX_2012} for an optimal stopping problem with a reward constraint;
see \cite{Lempa_J_2012,Liang_G_2015,Menaldi_Robin_2016,Liang_Wei_2016}
for optimal stopping with information constraint.

      The rest of the paper is organized as follows:
   In Subsection \ref{subsec:preliminary}, we introduce  notations and make  standing assumptions
   on drift/diffusion coefficients and reward/cost functions.
   In Section \ref{sec:OSEC_cont}, we set up the optimal stopping problem with expectation constraint
   over a general probability space and show the continuity of its value function
   in current state and   budget constraint level.
  Section \ref{sec:shift_prob} explores the measurability/integrability properties of shifted processes
  and the flow property of shifted stochastic differential equations
  as technical preparation for proving our main result, two types of DPPs.
  Then in Subsection \ref{subsec:DPP1}, we derive over   the canonical space  a DPP  for the value function $\cV$
  of the optimal stopping   with expectation constraint
  in which the conditional expected cost acts as an additional state process.
  In subsection \ref{subsec:DPP2}, we transform the the optimal stopping problem with expectation constraint
  to a stochastic optimization problem with supermartingale controls and establish a second DPP for $\cV$.
  Based on   two DPPs, we characterize $\cV$ as the viscosity solution to the related fully nonlinear  parabolic HJB equation
  in Section \ref{sec:PDE}.
  Section \ref{sec:proofs} contains proofs of our results while the demonstration of some auxiliary
  statements with starred labels 
  in these proofs are  relegated to  the Appendix.
  We also include  some technical lemmata in the appendix.

\subsection{Notation and Preliminaries} \label{subsec:preliminary}


 For a   generic Euclidian space $\hE$, we denote its Borel sigma$-$field by  $\sB(\hE)$.
   For any   $x \ins \hE $ and $\d \ins (0,\infty)$,
  $O_\d(x) \df \{x' \ins \hE \n : |x \- x'|   \< \d \}$ denotes
  the open   ball    centered at $ x   $     with radius $\d  $ and its closure is
  $\ol{O}_\d(x) \df \{x' \ins \hE \n : |x \- x'|    \ls  \d \}$.

  Fix  $  l  \ins  \hN$ and $p \ins [1,\infty)$.
  Let  $c(t) \n :  [0,\infty) \n \to \n (0,\infty) $ be a continuous function with $\int_0^\infty \n c(t) dt \< \infty$,
  and let $\fC $ be a constant with  $\fC \ge 1 \+ \int_0^\infty \n c(t) dt   $.
  As $\lmt{t \to \infty} c(t) \= 0$, the continuity of $c(\cd)$ implies that
  $   \|c(\cd)\| \df \underset{t \in [0,\infty)}{\sup} c(t) \< \infty$.
   Also,    let $\rho$ be a  modulus of continuity function and
       denote its inverse function by $\rho^{-1}$.

   We shall consider the following drift/diffusion coefficients and reward/cost functions throughout the paper.

 \no $\bullet$    Let     $b      \n :  (0,\infty) \ti \hR^l \n \to \n  \hR^l $ be
 a $\sB (0,\infty) \oti \sB(\hR^l) \big/ \sB(\hR^l) -$measurable function
  and let  $\si    \n :  (0,\infty) \ti \hR^l \n \to \n \hR^{l \times d}$ be
 a  $\sB (0,\infty) \oti \sB(\hR^l) \big/ \sB(\hR^{l \times d}) -$measurable function such
 that  for any  $t \ins (0,\infty)$ and $ x_1, x_2 \ins \hR^l$
 \bea
 \big|b(t,x_1)\-b(t,x_2) \big|  & \dn \dn \ls  & \dn \dn c(t) |x_1 \-x_2| , \q  |b(t,0)|    \ls c(t), \label{b_cond} \\
 \hb{and} \q  \big|\si(t,x_1)\-\si(t,x_2) \big|  & \dn \dn\ls & \dn \dn \sqrt{c(t)} |x_1 \-x_2| , \q |\si (t,0)|  \ls \sqrt{ c(t) } \,  .  \qq  \label{si_cond}
 \eea

  \no $\bullet$  The running reward function   $ f  \n : (0,\infty) \ti \hR^l \n \to \n \hR  $
 is a $\sB (0,\infty) \oti \sB(\hR^l) \big/ \sB(\hR) -$measurable function  such that
 for any $t \ins (0,\infty)$ and $ x_1, x_2 \ins \hR^l$
 \bea \label{eq:c343}
 \big|f(t,x_1) \- f(t,x_2)\big| \ls c(t) \big( |x_1\-x_2| \ve |x_1\-x_2|^p \big)
 \q   \hb{and}  \q  \big| f  (t, 0 ) \big| \ls c(t) .
 \eea

  \no $\bullet$ The terminal reward function $\pi  \n : [0,\infty) \ti \hR^l \n \to \n \hR $
 is a continuous function   such that for any $t ,t' \ins [0,\infty)$ and $ x,x' \ins \hR^l$
 \bea
 |\pi(t,x) \- \pi(t',x')|   \ls \rho \big( |t \- t'|   \big)
 \+ \fC \big( |x \-x'| \ve |x \-x'|^p \big)
 \q   \hb{and}  \q  \big|\pi (t,0)\big|   \ls \fC   .     \q \label{eq:c321}
\eea

  \no $\bullet$   The cost rate function $g     \n : (0,\infty) \ti \hR^l \n \to \n  (0,\infty)$
 is a $\sB (0,\infty) \oti \sB(\hR^l) \big/ \sB (0,\infty) -$measurable function  satisfying

  \no (g1) $ \big|g(t,x_1) \- g(t,x_2)\big|
 \ls c(t) \big( |x_1\-x_2| \ve |x_1\-x_2|^p \big) $,
   $ \fa  t \ins (0,\infty)$, $ \fa x_1,x_2 \ins \hR^l $; 

  \no (g2) $\int_0^t   g(t,0)   dr \< \infty $, $ \fa t \ins (1,\infty)$;

  \no (g3)   For any $R \ins  (0,\infty)$, there exists $\k_{\n \overset{}{R}} \ins (0,\infty)$ such that
  $  g(t,x)  \gs \k_{\n \overset{}{R}}  $, $ \fa t \ins (0,\infty)$, $ \fa x \ins \hR^l $ with $ |x| \ls R $.
  The constant $\k_{\n \overset{}{R}}$ can be   regarded as the  basic cost rate when the long-term state radius is $R$.

     Moreover,    we will  use the convention $ \inf \es \df \infty $ as well as  the inequality
   \bea
   \( 1 \land n^{q-1} \) \sum_{i=1}^n a_i^q \ls
\bigg(\sum_{i=1}^n a_i \bigg)^q \ls \( 1 \vee n^{q-1} \) \sum_{i=1}^n a_i^q     \label{eqn-d011}
 \eea
 for any $q \ins (0,\infty)$ and any finite subset  $\{a_1, \cds, a_n\} $ of $(0,\infty)$.



\section{Continuity of Value Functions for General Optimal Stopping with Expectation Constraint}

\label{sec:OSEC_cont}

 For an optimal stopping problem with expectation constraint,
 we first discuss the continuity of its value function  over a general complete probability space $(\O,\cF, P)$.

 Let $B$ be a $d-$dimensional standard Brownian motion   on    $(\O,\cF, P)$.
 The $P-$augmentation of its natural filtration
 $ \bF \= \big\{\cF_t \df \si \big( \si ( B_s; s \ins [0,t] ) \cp \sN \big) \big\}_{t  \in  [0, \infty) }$
 satisfies the {\it usual hypothesis},  where $\sN \df \big\{ \cN \sb \O  \n : \cN \sb A \hb{ for some  } A \ins \cF  \hb{ with } P  (A ) \=0   \} $ collects all $P-$null sets in $\cF$.
 Let $\cT$ stand for all $\bF-$stopping times  $\tau$ with $\tau \< \infty$, $P-$a.s.
 For any $\bF-$adapted continuous process $X$, we set $X_* \df \underset{s \in [0,\infty)}{\sup} |X_s|$.

\subsection{Reward Processes}

   Let $(t,x) \ins [0,\infty) \ti \hR^l $. It is well-known that under \eqref{b_cond} and \eqref{si_cond},
   the following stochastic differential equation (SDE) on $\O $
    \bea \label{FSDE}
     X_s   \=  x \+ \int_0^s \n b  ( t\+r,  X_r  )  dr \+ \int_0^s \n \si (  t\+r, X_r  )  dB_r ,  \q s \ins [0, \infty)
     \eea
     admits a unique solution $X^{t,x} \= \{ X^{t,x}_s \}_{s \in [0,\infty)}   $,  
     which is an $\hR^l-$valued, $ \bF-$adapted continuous process  satisfying

 \begin{lemm} \label{lem_X_Lp_estimate}
  Let $q \ins [1,\infty)$ and $(t,x) \ins [0, \infty) \ti \hR^l $.

  \no \(1\)  For some constant $ C_q \gs 1 $   depending on $q$ and $ \int_0^\infty \n c(s) ds$, we have
    \bea
    &&   E\bigg[ \, \underset{s \in [0,\infty)}{\sup}  \big| X^{t,x}_s   \big|^q \bigg]
         \ls  C_q  \big(1\+|x|^q \big)  ;
       \q E\bigg[ \, \underset{s \in [0,\infty)}{\sup}  \big| X^{t,x'}_s \- X^{t,x}_s  \big|^q \bigg]
           \ls  C_q  |x'\-x|^q , ~ \fa x' \ins  \hR^l; \q \hb{and}   \qq \qq   \label{eq:esti_X_1}   \\
  &&   E\bigg[ \, \underset{ \l \in (0,\d]}{\sup}  \big| X^{t,x}_{\tau+\l} \- X^{t,x}_\tau  \big|^q \bigg]
 \ls  C_q \big(1\+|x|^q \big)  \big( \|c(\cd)\|^q \d^q \+ \|c(\cd)\|^{\frac{q}{2}}  \d^{\frac{q}{2}} \big)  ,
 ~ \fa \d \ins (0,\infty), ~ \fa   \tau \ins \cT .    \label{eq:esti_X_2}
    \eea

 \no \(2\) Given $\varpi \ins [1,\infty)$,
 assume   functions $b$ and $\si $ additionally satisfy that for any $ 0 \ls t_1 \< t_2 \< \infty $ and $ x' \ins \hR^l $
 \bea \label{eq:c601}
 \big|b(t_2,x')\-b(t_1,x') \big|    \ls   c(t_1) \rho (t_2\-t_1)  \big( 1\+|x'|^\varpi  \big)
 \q \hb{and} \q     \big|\si(t_2,x') \- \si (t_1,x') \big|  \ls   \sqrt{c(t_1)   } \rho (t_2\-t_1)  \big( 1\+|x'|^\varpi  \big) .
 \eea
 Then it holds for any $t' \ins (t,\infty)$ that
 \bea  \label{eq:c609}
    E \bigg[ \, \underset{s \in [0,\infty)}{\sup}  \big| X^{t' \n,x}_s  \n \- \n  X^{t,x}_s  \big|^q \bigg]
           \ls  C_{q,\varpi} \big(1\+|x|^{q \varpi} \big)  \big( \rho (t' \n \-t) \big)^q  ,
 \eea
 where   $C_{q,\varpi} \gs 1$ is some constant  depending on $q$, $\varpi$ and $ \int_0^\infty \n c(s) ds$.

  \end{lemm}


 Given $t \ins [0,\infty)$, let  the   state process  evolve  from  
 position $x \ins \hR^l$ according to   SDE \eqref{FSDE}.
  If the player chooses to exercise at time   $\tau \ins \cT$,   she will receive
 a running reward $\int_0^\tau   \n f \big(t\+s, X^{t,x}_s \big) \, ds$ and
 a terminal reward $ \pi \big( t\+\tau   , X^{t,x}_\tau   \big) $,
 whose totality is
 \bea \label{eq:c621}
   R (t,x,\tau) \df \int_0^\tau   \n f \big(t\+s, X^{t,x}_s \big) \, ds
   \+ \pi \big( t\+\tau   , X^{t,x}_\tau   \big) .
   \eea
 One can deduce from  \eqref{eq:c343}, \eqref{eq:c321} and  the first inequality in \eqref{eq:esti_X_1} that
 \bea \label{eq:c483}
  E  \big[ | R (t,x,\tau) | \big] \ls 2 \fC \big( 2\+   C_p (1\+|x|^p) \big)
   \n := \n \Psi (x)   .
 \eea
 Given another initial position $x' \ins \hR^l$,   \eqref{eq:c343}, \eqref{eq:c321}, H\"older's inequality and
 the second inequality in \eqref{eq:esti_X_1} imply    that
\bea
 \q \qq  && \hspace{-2cm} E  \big[ \big| R (t, x, \tau ) \- R (t, x', \tau ) \big| \big]
 \ls E  \bigg[ \int_0^\tau     \big| f  (t\+r, X^{t,x}_r  ) \- f  (t\+r, X^{t,x'}_r  ) \big|   dr
 \+     \big| \pi  ( t\+\tau ,   X^{t,x}_{\tau }  ) \- \pi  ( t\+\tau ,   X^{t,x'}_{\tau }  ) \big|   \bigg] \nonumber \\
 && \ls \Big( \int_0^\infty c(t\+r) dr \+ \fC \Big) E  \Big[
 \big( X^{t,x} \n \- \n  X^{t,x'}  \big)_* \+ \big( X^{t,x}  \n \- \n  X^{t,x'}  \big)^p_* \Big]
   \ls 2 \fC \big(   (C_p)^{\frac{1}{p}}  | x \- x'  | \+  C_p  | x \- x'  |^p   \big) . \qq   \label{eq:c467}
\eea

\subsection{Expectation Constraints}

 Let $(t,x) \ins [0,\infty)   \ti  \hR^l$.  As the first inequality in  \eqref{eq:esti_X_1} shows that
    $  ( X^{t,x}_*  )^p 
    \< \infty$,  $P-$a.s.,     (g1)$-$(g3) imply that $P-$a.s.
\bea \label{gh_con}
\int_0^s \n g  \big(t\+r, X^{t,x}_r \big) dr \ls \int_0^s \n g (t\+r,0) dr
\+ \fC \big( X^{t,x}_* \+ (X^{t,x}_*)^p \big) \< \infty, ~ \fa s \ins (0,\infty)   \q \hb{and} \q
    \int_0^\infty \n g  \big(t\+r, X^{t,x}_r \big) dr   \=   \infty  .   \qq
\eea
 Given   $  y  \ins   [0,\infty) $,
 we  try to maximize   the player's expected total wealth $R(t,x,\tau)$ 
 when her expected cost is subject to the following constraint:
  \bea \label{def_constraint}
   E \bigg[ \int_0^\tau  g(t\+r, X^{t,x}_r) dr \bigg]  \ls y  .
 \eea
 Like reward processes
 $ \big\{ \int_0^s \n f \big(t\+r, X^{t,x}_r \big) dr \big\}_{s \in [0,\infty)} $ and
 $ \big\{\pi  (t\+s, X^{t,x}_s  ) \big\}_{s \in [0,\infty)}$, this expectation constraint is  also state-related.
  Hence, starting from the initial state $x \ins     \hR^l$,
  the   value  of the general optimal stopping problem with    expectation constraint $y$ is
 \bea \label{eq:c623}
 V(t,x,y)    \df   \underset{\tau \in \cT_{t,x} (y)   }{\sup} \, E  \big[ R (t,x,\tau) \big]     ,
 \eea
 where  $ \cT_{t,x} (y) \df \big\{ \hb{$\bF-$stopping time } \tau  \n : E \big[ \int_0^\tau  g(t\+r, X^{t,x}_r) dr \big]  \ls y \big\} $.

 For any $\tau \ins \cT_{t,x} (y) $, as $E  \big[ \int_0^\tau  g(t\+r, X^{t,x}_r) dr \big]    \ls  y   \< \infty$,
 one has $ \int_0^\tau  g(t\+r, X^{t,x}_r) dr  \< \infty $, $P-$a.s.
 The second part of \eqref{gh_con} then implies that  $\tau \< \infty$, $P-$a.s.
 So $ \cT_{t,x} (y) \= \big\{\tau \ins  \cT  \n : E \big[ \int_0^\tau  g(t\+r, X^{t,x}_r) dr \big]  \ls y \big\} $.

\begin{eg} \(Moment Constraints\)
For $q \ins (1,\infty) $, $ a  \ins [0,\infty)$ and $ b  \ins (0,\infty)$,
 take $g(\ft,\fx) \df a q \ft^{q-1} \+ b $,  $(\ft,\fx) \ins (0,\infty) \ti \hR^l$.
Then the constraint  \eqref{def_constraint} for $t \= 0 $
specify as the moment constraint   $E[ a\tau^q \+ b \tau ] \ls y$.
\end{eg}

Let $(t,x ) \ins [0,\infty)   \ti  \hR^l   $. It is clear that
\bea \label{eq:c405}
V(t,x,y) \hb{ is   increasing in $y$.}
\eea
 As $ \cT_{t,x} (0) \=   \{0\} $, we see from  \eqref{eq:c483}   that
\bea \label{eq:c433}
\Psi (x) \gs V(t,x,y)   \gs V(  t,x,0  ) \=   E  \big[ \pi  ( t , X^{t,x}_0  ) \big]  \= \pi (t,x) ,
\q \fa (t,x,y) \ins [0,\infty) \ti \hR^l \ti [0,\infty) .
\eea

 When $y \ins (0,\infty)$, we even have the following update of \eqref{eq:c623}.
 \begin{lemm} \label{lem_general_value}
 It holds for any $(t,x,y ) \ins [0,\infty)   \ti  \hR^l  \ti (0,\infty)$ that
 $ V(t,x,y)    \=   \underset{\tau \in \wh{\cT}_{t,x} (y)   }{\sup} \, E  \big[ R (t,x,\tau) \big] $,
 where $\wh{\cT}_{t,x} (y)  \df \big\{\tau \ins \cT_{t,x} (y) \n :  \tau \> 0, \hb{$P-$a.s.}  \big\}$.
 \end{lemm}

   The  value function $V(t,x,y)$ of  the  general  optimal stopping problem with expectation constraint
   is continuous in the following way:

\begin{thm} \label{thm_continuity}

\no \(1\) Given $ t  \ins [0,\infty)  $,    $ V(t,x,y)$  is  continuous in $(x, y ) \ins \hR^l \ti  [0,\infty)  $
 in the sense that  for any $(x,\e) \ins \hR^l \ti (0,1) $,
 there exists $ \d \= \d(t,x, \e) \ins (0,1)  $ such that for any $y \ins   [0,\infty) $
 \beas
 \big| V(t,\fx,\fy ) \-  V(t,x,y)\big|  \ls \e, \q  \fa (\fx,\fy ) \ins \ol{O}_\d (x) \ti \big[(y\-\d)^+, y\+\d\big]  .
 \eeas

\no \(2\) Given $\varpi \ins [1,\infty)$, assume $b,\si$ additionally satisfy \eqref{eq:c601} and
 $ f,g$ additionally satisfy that for any $ 0 \ls t_1  \< t_2 \< \infty $ and $ x' \ins \hR^l $
 \bea \label{eq:c611}
   \big|f(t_2,x')\-f(t_1,x') \big| \vee  \big|g(t_2,x')\-g(t_1,x') \big|
   \ls   c(t_1) \rho (t_2\-t_1)  \big( 1\+|x'|^\varpi  \big)   ,
 \eea
 then   $ V(t,x,y)$  is  continuous in $(t, x, y ) \ins [0,\infty) \ti \hR^l \ti  [0,\infty)  $
 in the sense that  for any $(t,x,\e) \ins  [0,\infty) \ti \hR^l \ti (0,1) $,
 there exists $ \d' \= \d'(t,x,\e) \ins (0,1)  $ such that for any $ y \ins    [0,\infty) $
 \beas
    \big| V(\ft,\fx,\fy ) \-  V(t,x,y)\big|  \ls \e , \q  \fa (\ft,\fx,\fy ) \ins
    \big[(t\-\d')^+, t\+\d'\big] \ti \ol{O}_{\d'}  (x) \ti \big[(y\-\d')^+, y\+\d'\big]  .
\eeas

\end{thm}

 \section{Shifted Processes}

 \label{sec:shift_prob}

  Let us   review   the properties of shifted processes on the canonical space so that
  we can   study two types of   dynamic programming principles  of  
  the optimal stopping problem with expectation constraint over the canonical space.

 Fix $d   \ins  \hN$ and let $t \ins [0,\infty)$. From now on, we consider
 the canonical space  $\O^t   \df  \big\{\o  \ins  \hC \big( [t,\infty)  ; \hR^d \big) \n :  \o(t)  \=  0 \big\}$
 of continuous paths   over period $[t,\infty)$, which   is a separable complete metric space
   under the {\it uniform} norm $ \|\o\|_t   \df  \underset{s \in [t,\infty)}{\sup} |\o(s)|  $.
  Let $\cF^t \df \sB (\O^t)$ be the Borel sigma field of $\O^t$ under $\| \cd \|_t $.
  The canonical process $ W^t \=\{W^t_s\}_{s \in  [t,\infty)} $   of  $\O^t $
  is a  $d-$dimensional standard Brownian motion on $\big(\O^t,  \cF^t \big)$ under  the    Wiener measure $ P_t $.
  Let $  \sN^t $ collect  all $P_t-$null sets, i.e.,
     $  \sN^t \df \big\{ \cN \sb \O^t \n : \cN \sb A \hb{ for some  } A \ins \cF^t \hb{ with } P_t (A ) \=0   \}    $,
          and set $\ol{\cF}^t \df \si (\cF^t \cp \sN^t)$.
         The completion of  $\big(\O^t,  \cF^t, P_t\big)$ is the probability space
        $\big(\O^t,  \ol{\cF}^t, \ol{P}_t \big)$ with $\ol{P}_t \Big|_{\cF^t}    \= P_t $.
        For simplicity,   we   still write $P_t$ for   $\ol{P}_t$ and
        denote   the expectation under  $\ol{P}_t$   by $E_t[\cd]$.
    For any   sub sigma$-$field  $\cG$  of $\ol{\cF}^t$, let
     $L^1  (\cG  ) $  be  the space of all  real-valued,
$  \cG-$measurable random variables $\xi$ with $  E_t  \big[ |\xi|  \big]   \< \infty$.

  We denote the  natural filtration of $ W^t  $ by
  $\bF^t   \=  \big\{ \cF^t_s  \df  \si \big(W^t_r; r  \ins  [t,s]\big)  \big\}_{s \in [t,\infty)}$.
  Its $P_t-$augmentation $ \ol{\bF}^t$   consists of
     $ \ol{\cF}^t_s   \df \si \big(  \cF^t_s  \cp  \sN^t  \big)  $,   $ s \ins [t,\infty)$
     and satisfies the {\it usual hypothesis}.
     Let  $\ol{\cT}^t$  stand for  all  stopping times $\tau$
     with respect to   the filtration  $\ol{\bF}^t$   such that $\tau \< \infty$, \pas, ~
     and set $\ol{\cT}^t_\sharp \df \big\{\tau \ins \ol{\cT}^t \n :
     \tau \hb{ takes countably many values in }[t,\infty) \big\}$.
    For easy reference, 
    we set   $ \cF^t_\infty  \df  \cF^t  $ and
    $  \ol{\cF}^t_\infty   \df   \ol{\cF}^t   $.     

   The following spaces   will be   used in the sequel.

  \no $\bullet$  For any  $q \ins [1,\infty)$,
  let  $\hC^q_t (\hE) \=  \hC^q_{\ol{\bF}^t}   \big( [t,\infty), \hE\big)$ be the space of
  all $\hE-$valued, $\ol{\bF}^t-$adapted processes   $\{X_s\}_{s \in [t,\infty)}$ with    \pas ~ continuous paths
  such that   $ E_t \big[ X^q_*    \big] \< \infty $  with $X_* \df \underset{s \in [t,\infty)}{\sup} |X_s|$.

 \no $\bullet$ Let $\hH^{2,{\rm loc}}_t $ denote all   $ \hR^d - $valued,  $\ol{\cF}^t-$predictable processes
 $\{X_s\}_{s \in [t,\infty)}$ with  $ P_t \big\{ \int_t^s |X_r|^2 dr \< \infty, \, \fa s \ins [t,\infty) \big\} \= 1 $.

  \no $\bullet$ Let $\hM_t $ denote all real-valued, uniformly integrable continuous martingales
 with respect to $\big(\ol{\bF}^t,P_t\big)$.

  \no $\bullet$ Set  $\hK_t \df \big\{K \ins \hC^1_t(\hR) \n : $ for $\pas ~ \o \ins \O^t$, $K_\cd (\o) $ is an
   continuous increasing path
  starting from $0$\big\}. 

 \subsection{Concatenation of Sample Paths}

 Let $0 \ls t \ls s \< \infty  $. We  define  a  {\it translation}  operator $ \Pi^t_s  $ from $ \O^t  $
   to $\O^s$   by
 \beas
 \big( \Pi^t_s  (\o)\big)(r) \df \o (r) \-  \o (s)  ,
 \q   \fa (r, \o )  \ins   [s, \infty)  \ti  \O^t       .
 \eeas
    On the other hand, one can concatenate  $\o \ins \O^t$ and   $ \wt{\o} \ins \O^s$ at   time $s  $ by:
   \beas
 \big(\o \otis    \wt{\o}\big)(r)  \df   \o(r) \, \b1_{\{r \in  [t,s)\}}
   \+ \big(\o(s) \+ \wt{\o}(r) \big) \, \b1_{\{r \in  [s,\infty)\}} , \q \fa  r \ins [t,\infty) ,
 \eeas
  which is   still  of $\O^t$.

 Given $\o \ins \O^t$,
 we set $A^{s,\o} \df \{ \wt{\o} \ins \O^s \n : \o \otis \wt{\o} \ins  A  \} $ for any $A \sb \O^t $;
 and  set    $ \o \otis \wt{A} \df \big\{ \o  \otis \wt{\o} \n : \wt{\o} \ins  \wt{A} \big\} $ for any $\wt{A} \sb \O^s$.
 In particular, $\es^{s,\o} \df  \es$ and $\o \otis \es  \df \es $.

 The next result shows that each   $ A \ins \cF^t_s $  consists of all branches $\o \otis \O^s $ with $\o \ins  A$.

\begin{lemm}  \label{lem_element}
    Let $0 \ls t \ls s \< \infty  $ and  $ A \ins \cF^t_s$.
    It holds for any   $\o \ins A  $ that  $  \o \otis \O^s   \sb A  $ or  $A^{s,\o} \= \O^s $.
  \end{lemm}

  Let $\o \ins \O^t$.   For any    $ \cF^t_s-$measurable random variable $\xi$,
  since the set $   \{ \o'  \ins  \O^t  \n : \xi (\o') \=  \xi (\o) \} \= \xi^{-1} \big(\{\xi(\o)\}\big) $
  belongs to $ \cF^t_s$,   Lemma \ref{lem_element}  implies that
 \bea   \label{eq:bb421}
  \o  \otis \O^s \sb \{\o'  \ins  \O^t  \n : \xi (\o') \=  \xi (\o) \}
 \q \hb{i.e.,} \q
    \xi(  \o \otis \wt{\o} )  \=  \xi(\o), \q
  \fa \wt{\o}  \ins  \O^s .
  \eea
  To wit, the value $ \xi(\o)  $  depends only on   $\o|_{[t,s]}$.

   For any $r \ins [s,\infty]$, the operation $( )^{s,\o}$ projects an $\cF^t_r-$measurable set to
   an $\cF^s_r-$measurable set while the
   operation $ \o \otis \cd$ transforms an $\cF^s_r-$measurable set into  an $\cF^t_r-$measurable set.

 \begin{lemm} \label{lem_concatenation}

 Let $0 \ls t \ls s \< \infty  $, $\o \ins \O^t$ and $r \ins [s,\infty]$.
 We have  $A^{s,\o}  \ins    \cF^s_r$
 for any $A  \ins  \cF^t_r$ and        $\o \otis \wt{A} \ins \cF^t_r $ for any  $ \wt{A} \ins \cF^s_r $.

 \end{lemm}

\subsection{Measurability and Integrability of Shifted Processes}

\label{subsec:shifting-processes}

  Let $0 \ls t \ls s \< \infty  $, let $\xi$ be an   $\hE-$valued random variable on $\O^t$
  and let $X \= \{X_r\}_{r \in [t,\infty)}$ be an  $\hE-$valued process on $\O^t$.
  For any   $\o \ins \O^t$,   we define the shifted random variable $\xi^{s,\o}$
    and the shifted process $X^{s,\o}$ by
      \beas
      \xi^{s,\o}(\wt{\o}) \df \xi ( \o \otis   \wt{\o} ) \q \hb{and} \q
       X^{s,\o}  (r, \wt{\o}) \df   X  (r ,  \o \otis  \wt{\o}) ,
      \q  \fa   ( r, \wt{\o} )  \ins [s,\infty) \ti  \O^s   .
      \eeas

  By Lemma \ref{lem_concatenation},   shifted random variables and shifted processes
  inherit the   measurability of original ones.

 \begin{prop}  \label{prop_shift_meas}
 Let $0 \ls t \ls s \< \infty  $   and   let $\o \ins \O^t$.

 \no \(1\) Let $\xi $ be an $\hE-$valued random variable  on $\O^t$.
  If  $\xi $ is   $\cF^t_r-$measurable for some $r \ins [s,\infty]$,
     the shifted random variable $\xi^{s,\o} $ is $  \cF^s_r -$measurable.

  \no \(2\) Let $ X \= \{X_r \}_{r \in [t,\infty)}$ be an   $\hE-$valued   process on $\O^t$.
   If   $ X $ is   $\bF^t-$adapted,
   the shifted  process     $  X^{s,\o} \= \big\{X^{s,\o}_r \big\}_{  r \in [s,\infty)}$ is   $\bF^s-$adapted.

 \end{prop}

 In virtue of regular conditional probability distribution by \cite{Stroock_Varadhan},
 the shifted random variables carry on the integrability as follows:

   \begin{prop}  \label{prop_shift_integr}
   Let   $0 \ls t \ls s \< \infty  $.
   If   $\xi \ins  L^1 (\cF^t)$,   then it holds for  \pas ~ $\o  \ins  \O^t$  that
   $\xi^{s,\o}  \ins  L^1  ( \cF^s  ) $ and
  \bea   \label{eq:f475}
  E_t   \big[\xi\big| \cF^t_s \big](\o) \= E_s    \big[ \xi^{s,\o} \big] \ins \hR  .
    \eea
    \end{prop}

 Consequently, the shift of a $P_t-$null set still has zero $P_s-$probability.

 \begin{prop}  \label{prop_null_set}
  Let   $0 \ls t \ls s \< \infty  $.

  \no \(1\)  For any  $P_t-$null set    $  \cN   \ins  \sN^t      $,
  it holds  for   \pas ~ $\o \ins \O^t  $ that $\cN^{s,\o}  \ins  \sN^s  $.
 Then for any two real-valued random variables $\xi_1$ and $ \xi_2$ on $\O^t$ with  $\xi_1 \ls \xi_2$,   $P_t-$a.s.,
   it holds for   \pas ~ $\o \ins \O^t  $ that  $\xi^{s,\o}_1 \ls \xi^{s,\o}_2$, $P_s-$a.s.

   \no \(2\) For any $\tau \ins \ol{\cT}^t  $ with $\tau \gs s$,
   it holds for $P_t-$a.s. $\o \ins \O^t$ that   $\tau^{s,\o} \ins \ol{\cT}^s $.

 \end{prop}

 Based on  Proposition \ref{prop_null_set} (1) and Lemma \ref{lem_F_version},
 we can extend Proposition \ref{prop_shift_integr} from raw filtration  $\bF^t$ to augmented filtration  $\ol{\bF}^t$,
 and can show that the shifted    processes  inherit  the integrability  of   original ones.

   \begin{prop}  \label{prop_shift_FP}
   Let   $0 \ls t \ls s \< \infty    $.

\no \(1\)  For any $\ol{\cF}^t_s-$measurable random variable $\xi$, it holds for  \pas ~ $\o \ins \O^t$  that
$\xi^{s,\o} \= \xi (\o) $, $P_s-$a.s.

\no \(2\)     For any $r \ins [s,\infty]$ and $\ol{\cF}^t_r-$measurable random variable $\xi$,
    it holds for  \pas ~ $\o \ins \O^t$  that $\xi^{s,\o} $ is $ \ol{\cF}^s_r  -$measurable.
   If $\xi$ is integrable,  then    it holds for  \pas ~ $\o \ins \O^t$  that
   $\xi^{s,\o}  $ is integrable and
   $E_t  \big[ \xi \big| \ol{\cF}^t_s  \big](\o)   \= E_s    \big[ \xi^{s,\o} \big] \ins \hR$.

\no \(3\) Let $X \= \{X_r\}_{r \in [t,\infty)}  $ be an $\ol{\bF}^t-$adapted process with \pas ~ continuous paths. It holds for \pas ~ $\o \ins \O^t$ that the shifted process
 $ X^{s,\o} \=  \big\{X^{s,\o}_r \big\}_{  r \in [s,\infty)}$ is $\ol{\bF}^s-$adapted   with $P_s-$a.s.  continuous paths.
If $X \ins  \hC^q_t (\hE)$ for some $q \ins [1,\infty)$, then
 $ X^{s,\o}   \ins  \hC^q_s (\hE)  $ for  \pas ~ $\o \ins \O^t$.

 \end{prop}

 Moreover, the shift of a uniformly integrable  martingale are still   uniformly integrable  martingales
 under the augmented filtrations.

  \begin{prop} \label{prop_shift_mart}
  Let  $0 \ls t \ls s \< \infty    $.   For any $M     \= \{M_r\}_{r \in [t,\infty)}
  \ins \hM_t$,   it holds for \pas ~ $\o \ins \O^t$ that
  $ M^{s,\o} \= \{M^{s,\o}_r\}_{r \in  [s,\infty)} $ is of $ \hM_s $.

  \end{prop}

\subsection{Shifted Stochastic Differential Equations}

\label{subsec:shift_SDE}

 Let $(t,x) \ins [0,\infty) \ti \hR^l $. The   SDE \eqref{FSDE2}
 has  a   unique solution $\cX^{t,x} \= \{\cX^{t,x}_s \}_{s \in [t,\infty)} $,
 which    is an $\hR^l-$valued, $\ol{\bF}^t-$adapted continuous process.
 As it holds \pas ~ that
 \beas
 \q \cX_{t+s}  \=  x \+ \int_t^{t+s} b  ( r,  \cX_r  ) dr \+ \int_t^{t+s}  \si (  r, \cX_r  )  \, dW^t_r
 \=  x \+ \int_0^s b  ( t\+r,  \cX_{t+r}  ) dr \+ \int_{r \in [0,s]}  \si ( t\+r,  \cX_{t+r} )  \, dW^t_{t+r} ,
 \q s \ins [0, \infty) ,
 \eeas
 we see that $\big\{\cX^{t,x}_{t+s}\big\}_{s \in [0,\infty)}$ is exactly the unique solution of \eqref{FSDE}
 with the probabilistic specification
 \bea \label{eq:c619}
 \big(\O,\cF,P,\sN,\{B_s\}_{s \in [0,\infty)}, \{\cF_s\}_{s \in [0,\infty)}\big)
 \= \Big(\O^t,  \cF^t, P_t,\sN^t, \{W^t_{t+s}\}_{s \in [0,\infty)},   \big\{\ol{\cF}^t_{t+s}\big\}_{s \in [0,\infty)} \Big) .
 \eea
  Clearly,  $\tau$ is an $\ol{\bF}^t-$stopping time if and only if $\wt{\tau} \df \tau \- t$
  is a stopping time  with respect to the filtration $\big\{\ol{\cF}^t_{t+s}\big\}_{s \in [0,\infty)}$.
  So the corresponding $\cT$ under   setting \eqref{eq:c619} is $\cT \= \big\{\wt{\tau} \= \tau\-t\n : \tau \ins \ol{\cT}^t \big\}$.
  It then follows from    Lemma \ref{lem_X_Lp_estimate} that
  \begin{cor} \label{cor_X_Lp_estimate}
  Let $q \ins [1,\infty)$ and $(t,x) \ins [0, \infty) \ti \hR^l $. 
  For the same constant  $ C_q  $ as in Lemma \ref{lem_X_Lp_estimate},
    \bea
  &&  E_t\bigg[ \, \underset{s \in [t,\infty)}{\sup}  \big| \cX^{t,x}_s   \big|^q \bigg]
         \ls  C_q  \big(1\+|x|^q \big);  \q
      E_t\bigg[ \, \underset{s \in [t,\infty)}{\sup}  \big| \cX^{t,x'}_s \- \cX^{t,x}_s  \big|^q \bigg]
           \ls  C_q  |x'\-x|^q , ~ \fa x' \ins \hR^l ; \q  \hb{and} \qq \qq   \label{eq:esti_X_1b}  \\
  &&  E_t\bigg[ \, \underset{ \l \in (0,\d]}{\sup}  \big| \cX^{t,x}_{\tau+\l} \- \cX^{t,x}_\tau  \big|^q \bigg]
 \ls  C_q \big(1\+|x|^q \big)  \big( \|c(\cd)\|^q \d^q \+ \|c(\cd)\|^{\frac{q}{2}}  \d^{\frac{q}{2}} \big)  ,
 ~ \fa \d \ins (0,\infty) , ~ \fa \tau \ins \ol{\cT}^t .   \;  \q \qq  \label{eq:esti_X_2b}
    \eea


 \end{cor}

      The shift of $\cX^{t,x}$  given path $\o|_{[t,s]}$ turns out to be
 the solution of the shifted stochastic differential equation  \eqref{FSDE} over period $[s,\infty)$ with  initial state $ \cX^{t,x}_s (\o) $:

     \begin{prop} \label{prop_FSDE_shift}
 \(Flow Property\)
 Let $0 \ls t \ls s \< \infty $, $x \ins \hR^l$  and set  $\fX \df \cX^{t,x} $.
 It holds  for \pas ~ $\o \ins  \O^t$ that
 $ P_s \big\{ \wt{\o} \ins \O^s \n :   \fX_r (\o \otis \wt{\o})
 \=  \cX^{s, \fX_s (\o)}_r (\wt{\o}), \; \fa r \ins [s,\infty) \big\} \= 1   $.

\end{prop}

  The proof of Proposition \ref{prop_FSDE_shift}  depends on the following result
  about the convergence of shifted random variables in probability.

 \begin{lemm} \label{lem_shift_converge_proba}
  For any   $ \{\xi_i\}_{i \in \hN} \subset L^1 \big(\ol{\cF}^t\big)$ that converges to 0 in probability $P_t$,
 we can find a subsequence $ \big\{  \wh{\xi}_{\,i} \big\}_{i \in \hN}  $ of it such that for $P_t-$a.s. $\o \in \O^t$,
 $ \big\{ \wh{\xi}^{\,s,\o}_{\,i} \big\}_{i \in \hN}  $   converges to 0 in probability $P_s$.

 \end{lemm}

 \section{Two Dynamic Programming Principle of Optimal Stopping with Expectation Constraint}

\label{sec:DPP}

 In this section, we exploit the flow property of shifted stochastic differential equations
 to establish two types of  dynamic programming principles (DPPs)
 of the optimal stopping problem with expectation constraint over the canonical space.

\subsection{The First Dynamic Programming Principle for $\cV$}

\label{subsec:DPP1}


 Let  the   state process now   evolve  from   time $t  \n \in \n  [0,\infty)$
 and position $x \ins \hR^l$ according to   SDE \eqref{FSDE2}.
  If the player selects to exercise at time $\tau \ins \ol{\cT}^t$, she will receive
 a running reward $\int_t^\tau   \n f (r,\cX^{t,x}_r) dr$
 and a terminal reward $ \pi \big(  \tau   , \cX^{t,x}_\tau   \big) $.
 So the player's   total wealth is
  \beas
  \cR(t,x,\tau) \df  \int_t^\tau   \n f (s,\cX^{t,x}_s) ds \+ \pi \big(  \tau   , \cX^{t,x}_\tau   \big)
  \= \int_0^{\wt{\tau}}  \n f (t\+s,\cX^{t,x}_{t+s}) ds \+ \pi \big(  t\+\wt{\tau}  , \cX^{t,x}_{t+\wt{\tau}}  \big) ,
  \eeas
  which is the payment $R\big(t,x,\wt{\tau}\big)$ in \eqref{eq:c621}   under the specification \eqref{eq:c619}.
  By \eqref{eq:c483} and \eqref{eq:c467}, one has
  \bea \label{eq:c483b}
  E_t \big[ | \cR(t,x,\tau) | \big] \ls   \Psi (x)
  \q \hb{and} \q E_t \big[ \big| \cR(t, x, \tau ) \- \cR(t, x', \tau ) \big| \big]
 \ls   2 \fC \big(   (C_p)^{\frac{1}{p}}  | x \- x'  | \+   C_p  | x \- x'  |^p   \big), ~ \fa x' \ins \hR^l .
  \eea

  Given   $  y  \ins   [0,\infty) $, set $ \cT^t_x (y) \df \big\{\tau \ins \ol{\cT}^t \n : E_t [ \int_t^\tau  g(r, \cX^{t,x}_r) dr ]  \ls y \big\} $.  As $  E_t \big[ \int_t^\tau  g(r, \cX^{t,x}_r) dr \big]
  \= E_t \big[ \int_0^{\wt{\tau}}  g(t\+r, \cX^{t,x}_{t+r}) dr \big] $,  we see that $\big\{\wt{\tau} \= \tau \- t \n : \tau \ins \cT^t_x (y)\big\}$ is the corresponding $T_{t,x}(y)$ under  setting \eqref{eq:c619}.
   Then the maximum of the player's expected   wealth   subject to the budget   constraint
 $   E_t \big[ \int_t^\tau  g(r, \cX^{t,x}_r)   dr \big]  \ls y  $, i.e.,
 \bea \label{eq:d065}
 \cV(t,x,y)    \df   \underset{\tau \in \cT^t_x (y)   }{\sup} \,  E_t \big[  \cR(t,x,\tau) \big]
 \=  \underset{\wt{\tau} \in T_{t,x} (y)   }{\sup} \,  E_t \big[   R \big(t,x,\wt{\tau}\big) \big]
 \eea
 is exactly the value function \eqref{eq:c623}
 of the constrained optimal stopping problem under the specification \eqref{eq:c619}. Then
 \eqref{eq:c433} and Lemma \ref{lem_general_value} show that
\bea
&& \Psi (x)   \gs  \cV(t,x,y)  \gs   \cV (  t,x,0  ) \= \pi (t,x) ,
\q \fa (t,x,y) \ins [0,\infty) \ti \hR^l \ti [0,\infty) \qq \label{eq:c433b} \\
\hb{and} &&  \cV(t,x,y)   \=   \underset{\tau \in \wh{\cT}_{t,x} (y)   }{\sup} \, E  \big[ R (t,x,\tau) \big] ,
 \q \fa (t,x,y ) \ins [0,\infty)   \ti  \hR^l  \ti (0,\infty) , \label{eq:d187}
\eea
  where $\wh{\cT}^t_x (y)  \df \big\{\tau \ins \cT^t_x (y) \n :  \tau \> t, \pas  \big\}$.
  Also,
 \bea \label{eq:c629}
  \hb{Theorem \ref{thm_continuity} still holds for the value function $ \cV $.}
 \eea

 Now, let $(t,x ) \ins [0,\infty)   \ti  \hR^l   $ and let $\tau \ins \ol{\cT}^t $ with
    $  E_t \big[ \int_t^\tau \n g(r, \cX^{t,x}_r) dr \big]  \< \infty$.
    We   define a real-valued, $\ol{\bF}^t-$adapted continuous process:
 \beas
 \cY^{t,x,\tau}_s \df E_t \bigg[ \int_t^\tau \n g \big(r, \cX^{t,x}_r \big) dr \Big|\ol{\cF}^t_s \bigg]
 \- \int_t^{\tau \land s} \n g \big(r, \cX^{t,x}_r \big) dr , \q   s \ins [t,\infty) .
 \eeas
  Since it holds   for any $s \ins [t,\infty)$    that
 \bea \label{eq:c471}
  \cY^{t,x,\tau}_s \= E_t \bigg[ \int_{\tau \land s}^\tau   g \big(r, \cX^{t,x}_r \big) dr \Big|\ol{\cF}^t_s \bigg]
 \= E_t \bigg[ \int_s^{\tau \vee s}   g \big(r, \cX^{t,x}_r \big) dr \Big|\ol{\cF}^t_s \bigg] \ins [0,\infty) , \q P_t-\hb{a.s.,}
 \eea
 the continuity of $\cY^{t,x,\tau}$ implies that
 \bea \label{eq:c485}
 \cN_{t,x,\tau} \df \big\{ \cY^{t,x,\tau}_s \n \notin \n [0,\infty)  \hb{ for some }  s \ins [t,\infty) \big\} \ins \sN^t .
 \eea

 Then we have the first dynamic programming principle for  the value function $\cV$
 in which the conditional expected cost $\cY^{t,x,\tau}$   acts as an additional state process.

\begin{thm} \label{thm_DPP}
 Let $\,t \ins [0,\infty) $.

\no  \(1\) For any $( x, y ) \ins \hR^l \ti [0,\infty)$,
 let  $\{\z (\tau) \}_{\tau \in \cT^t_x (y)} $ be a family of $  \ol{\cT}^t_\sharp -$stopping times. Then
 we have the   DPP \eqref{eqn_DPP},
 where $\underset{\tau \in \cT^t_x (y)   }{\sup} E_t[\cd]$   can be replaced by
 $ \underset{\tau \in \wh{\cT}^t_x (y)   }{\sup} E_t[\cd]$ if $y \> 0$.

\no  \(2\) If $\cV(s,x,y)$ is continuous in $(s,x,y) \ins [t,\infty) \ti \hR^l \ti (0,\infty)$, then
\eqref{eqn_DPP} holds for any $( x, y ) \ins \hR^l \ti [0,\infty)$
and any family $\{\z (\tau) \}_{\tau \in \cT^t_x (y)} $ of $  \ol{\cT}^t -$stopping times.

\end{thm}

\subsection{An Alternative Stochastic Control Problem
 and the Second Dynamic Programming Principle for $\cV$}

\label{subsec:DPP2}

  Fix $t \ins [0,\infty)$ and set  $\fA_t   \df \{ \a \= M \- K \n : (M,K) \ins \hM_t \ti \hK_t  \}  $.
  Clearly, each $ \a \ins \fA_t$ is a uniformly integrable continuous supermartingales
  with respect to $\big(\ol{\bF}^t,P_t\big)$.

    Let  $   x   \ins   \hR^l   $ and $ \a   \ins \fA_t  $.
    We  define  a continuous supermartingale with respect to $\big(\ol{\bF}^t,P_t\big)$
\beas
 Y^{t,x,\a}_s \df  \a_s \- \int_t^s \n g   \big(r, \cX^{t,x}_r \big) dr  ,     \q s \ins [t,\infty) ,
\eeas
 and define an $\ol{\bF}^t-$stopping time
\bea \label{eq:d057}
 \tau (t,x,\a)    \df   \inf\big\{s \ins [t,\infty)  : Y^{t,x,\a}_s \= 0 \big\}   .
\eea

   The uniform integrability of $\a $
   implies that    the   limit  $   \lmt{s \to \infty} \a_s   $   exists in $\hR$, \pas ~
  \if{0}
  Set $ M_s \df \int_t^s \n \a_r dW^t_r $, $s \ins [t,\infty)$.
  For any $j \ins \hN$, define an $\ol{\cT}^t-$stopping time
  $\z_j \df \inf\big\{s \ins [t,\infty) \n : \int_t^s \n |\a_r|^2 dr \> j  \big\} $.
  Since $ \big\{ M_{\z_j \land s} \big\}_{s \in  [t,\infty)}$
  is a uniformly integrable martingale,
  $ \xi_j \df \lmt{s \to \infty} M_{\z_j \land s}   $ exists and belongs to $ L^1 (\ol{\cF}^t, \hR^l)$.
  In particular,   $ \xi_j \ins \hR^l$ except on a $P_t-$null set $\cN_j$.
  For all $\o \ins \O^t $ except on a $P_t-$null set $\wh{\cN}$,
  $\int_t^\infty \n |\a_r (\o)|^2 dr \< \infty$, so we can find a $\fj \= \fj(\o) \ins \hN$
  such that $\z_\fj (\o) \= \infty$. Then for any $\o \in \wh{\cN}^c \cap \Big( \underset{j \in \hN}{\cap} \cN^c_j \Big)$,
  one can deduce that
 \beas
 \lmt{s \to \infty} M   (s,\o) \= \lmt{s \to \infty} M   (z_\fj (\o) \land s,\o)
 \= \lmt{s \to \infty} M_{\z_\fj \land s} (\o) \= \xi_j (\o) \ins \hR .
 \eeas
 \fi
 Since $ \int_t^\infty \n g  \big(r, \cX^{t,x}_r \big) dr \=
 \int_0^\infty \n g  \big(t\+r, \cX^{t,x}_{t+r} \big) dr   \=   \infty $, \pas ~ by
     \eqref{gh_con}, one can deduce   that
  \bea \label{eq:ax041}
  \tau (t,x,\a) \< \infty    , \q  \pas
  \eea
 Namely, $\tau (t,x,\a)   \ins \ol{\cT}^t$.

 Given $\a \ins \fA_t$, the expected wealth $ E_t \big[   \cR \big(t, x', \tau( t, x',   \a    ) \big) \big] $
 is continuous in $x \ins \hR^l$,
 which will play an important role in the demonstration of the second DPP for $\cV$
 (Theorem \ref{thm_DPP2}).

\begin{prop} \label{prop_control_cont}
Let $(t,x) \ins [0,\infty) \ti \hR^l$ and let $\a \ins \fA_t$. For any $\e \ins (0,1)$, there exists
$\d \= \d (t,x,\e) \ins (0,1)$ such that
\beas
  E_t \Big[ \big| \cR \big(t, x', \tau( t, x',   \a    ) \big) \- \cR \big(t, x, \tau( t, x,  \a  ) \big)   \big| \Big]
 \ls \e , \q \fa x' \ins  \ol{O}_\d (x) .
\eeas
\end{prop}

 For any $y \ins (0,\infty)$, we set  $\fA_t(y) \df \big\{\a \ins \fA_t: \a_t \= y, \, \pas \big\}$.

\begin{prop} \label{prop_surm}
 Given $(t,x,y) \ins [0,\infty) \ti  \hR^l \ti  (0,\infty)   $,
 $\a \n \to \n \tau (t,x,\a)$   is   a surjective mapping from $  \fA_t (y)   $ to
   $ \wh{\cT}^t_x (y) $.
\end{prop}

\begin{rem} \label{rem_constraint}
Let $(t,x,y) \ins [0,\infty) \ti  \hR^l \ti  (0,\infty)   $.

\no 1\) Let $\tau \ins \wh{\cT}^t_x (y)$.  Proposition \ref{prop_surm} shows that $\tau \= \tau (t,x,\a)$
for some $\a \ins \fA_t (y)$. In particular, we see from \eqref{eq:d225} of its proof that
$\a$ is a martingale \(resp. supermartingale\) if $   E_t \big[   \int_t^\tau  g(r, \cX^{t,x}_r) dr   \big] \- y \= 0$
\(resp. $\ls 0$\).  To wit, the constraint $E_t \big[\int_t^\tau g(r,\cX^{t,x}_r) dr\big] \= y $ \(resp. $\ls y$\)
corresponds to martingale \(resp. supermartingale\) controls   in the alternative stochastic optimization problem.

In case that  $\a$ is a martingale, we know from the martingale representation theorem that
$\a_s \= y+\int_t^s \fq_r dW^t_r$, $s \ins [t,\infty) $ for some $\fq \ins \hH^{2,{\rm loc}}_t$.
However reversely, for a $ \wt{\fq} \ins \hH^{2,{\rm loc}}_t$,
$ \wt{\a}_s \df y+\int_t^s \wt{\fq}_r dW^t_r $, $s \ins [t,\infty) $ could be a strict local martingale
with  $  E_t \big[\int_t^{\tau (t,x,\wt{\a})} g(r,\cX^{t,x}_r) dr\big] \< y $,
see Example \ref{eg_slm} in the appendix.
This is the reason why \cite{Miller_C_2017a} requires $E[\tau^2] \< \infty$ (see line -4 in page 3 therein)
for the one-to-one correspondence between constrained stopping rules and squarely-integrable controls.

\no 2\) Define the value of the optimal stopping
under the constraint $E_t \big[\int_t^\tau g(r,\cX^{t,x}_r) dr\big] \= y $  by
\beas
\cU(t,x,y) \df \sup\bigg\{  E_t \big[  \cR(t,x,\tau) \big] \n : \tau \ins \ol{\cT}^t \hb{ with } E_t \Big[ \int_t^\tau  g(r, \cX^{t,x}_r) dr \Big]  \= y \bigg\} .
\eeas
 Clearly, $\cU(t,x,y) \ls \cV(t,x,y)$. However, we do not know whether they are equal since
 $ \cU(t,x,y) $ may not be increasing in $y$ \(cf.  line 5 of Lemma 1.1 of \cite{AKK_2015}\).

\no 3\) The constraint $E_t \big[\int_t^\tau g(r,\cX^{t,x}_r) dr\big]  \ls y$
  is necessary for proving the continuity and the first DPP of the   value function $\cV$:
  Even if $\tau_1$ in \eqref{eq:c387} has $ E  \big[\int_0^{\tau_1} g(t\+r, X^{t,x}_r) dr\big] \= y $,
  the approximately optimal stopping time $\wh{\tau}_1$ constructed in the case \eqref{eq:d227}
  may   satisfy   $  E  \big[\int_0^{\wh{\tau}_1} g(t\+r, X^{t,\fx}_r) dr\big] \< (y \- \d)^+ $
  rather than $  E  \big[\int_0^{\wh{\tau}_1} g(t\+r, X^{t,\fx}_r) dr\big] \= (y \- \d)^+ $.
  Even if the $\tau \ins \cT^t_x(y)$ given in Lemma \ref{lem_DPP} reaches
  $ E_t \big[\int_t^\tau g(r,\cX^{t,x}_r) dr\big]  \= y $,
  the pasting $\ol{\tau}$ of $\tau$ with the locally $\e-$optimal stopping times $\tau^i_n$'s in \eqref{eq:d067}
    satisfies $ E_t \big[\int_t^{\ol{\tau}} g(r,\cX^{t,x}_r) dr\big] \< y\+\e$
    but $ E_t \big[\int_t^{\ol{\tau}} g(r,\cX^{t,x}_r) dr\big] \= y\+\e$
    after a series of estimations in \eqref{eq:d229}.

\end{rem}

 By Proposition \ref{prop_surm} and \eqref{eq:d187},  we can alternatively express the   optimal stopping problem
 with expectation constraints \eqref{def_constraint}   as a     stochastic control problem:
 \bea \label{eq:d189}
 \cV( t, x,y )   \=   \underset{\a \in \fA_t (y) }{\sup} \,   E_t \big[ \cR \big(t, x, \tau (t,x,\a) \big) \big] ,
 \q \fa (t,x,y) \ins [0,\infty) \ti \hR^l \ti (0,\infty) .
 \eea

  Moreover, we have the  second dynamic programming principle for the value function $\cV$ in which
  the controlled supermartingale $Y^{t,x,\a} $ serves as an additional state process.

\begin{thm} \label{thm_DPP2}

 Let $\,t \ins [0,\infty) $.

 \no  \(1\) For any $( x, y ) \ins \hR^l \ti [0,\infty)$,
 let  $\{\z (\a) \}_{\a   \in \fA_t (y)} $ be a family of $  \ol{\cT}^t_\sharp -$stopping times.
 Then  we have the   DPP \eqref{eqn_DPP2}.

\no  \(2\) If $\cV(s,x,y)$ is continuous in $(s,x,y) \ins [t,\infty) \ti \hR^l \ti (0,\infty)$, then
\eqref{eqn_DPP2} holds for any $( x, y ) \ins \hR^l \ti [0,\infty)$ and any
   family  $\{\z (\a) \}_{\a   \in \fA_t (y)} $ of $  \ol{\cT}^t  -$stopping times.

\end{thm}

   \section{Related Fully Non-linear Parabolic HJB Equations}

   \label{sec:PDE}

 In this section, 
 we show that  the  value function
   of the optimal stopping problem with expectation constraint is the   viscosity solution  to
   a related   fully non-linear parabolic Hamilton-Jacobi-Bellman (HJB) equation.

 For any  $ \phi (t,x,y) \ins C^{1,2,2}\big([0,\infty) \ti \hR^l \ti [0,\infty)\big)$, we set
 \beas
 \sD  \phi (t,x,y) \df \big(D_x  \phi   ,  D^2_x  \phi   , \pa_y  \phi  , \pa^2_y  \phi  ,  D_x  ( \pa_y  \phi  )    \big)
 (t,x,y) \ins \hR^l \ti \hS_l \ti \hR \ti \hR \ti \hR^l, \q \fa (t,x,y) \ins [0,\infty) \ti \hR^l \ti [0,\infty) ,
 \eeas
 where $\hS_l$   denotes    the set of all  $\hR^{l \times l}-$valued symmetric matrices.

 Recall the  definition of viscosity solutions to a
 parabolic equation with a general (non-linear) Hamiltonian $ H \n : [0,\infty) \ti   \hR^l \ti   \hR  \ti   \hR^l
  \ti   \hS_l \ti \hR \ti \hR \ti \hR^l  \n \to \n   [-\infty, \infty]  $.

\begin{deff}  \label{def:viscosity_solution}

 An upper \(resp.\;lower\) semi-continuous function $u  \n :  [0,\infty)  \ti  \hR^l \ti [0,\infty) \n \to \n  \hR$    is called  a viscosity subsolution \(resp. supersolution\) of
 \beas
 \begin{cases}
  -  \pa_t u (t,x,y)  -   H \big(t,x,   u (t,x,y),  \sD u (t,x,y)  \big)   \= 0 , \q \fa (t,x,y) \ins (0,\infty)  \ti  \hR^l \ti (0,\infty) , \q \ss   \\
   u (t,x,0) \= \pi (t,x) , \q \fa (t,x ) \ins [0,\infty)  \ti  \hR^l
 \end{cases}
 \eeas
  if   $u (t,x,0) \ls $ \(resp.\;$\ge$\) $  \pi (t,x) $,   $\fa (t,x ) \ins [0,\infty)  \ti  \hR^l $,
  and if for any $(t_o,x_o,y_o) \ins (0,\infty)  \ti  \hR^l \ti (0,\infty)$ and
  $ \phi  \ins C^{1,2,2}\big([0,\infty) \ti \hR^l \ti [0,\infty)\big)$
  such that 
  $u  \n - \n  \phi$ attains a strict local maximum $0$ \(resp.\;strict local minimum $0$\)  at $(t_o,x_o,y_o)$, one has
 \beas
        -   \pa_t \phi (t_o,x_o,y_o)
     \-   H \big(t_o,x_o,  \phi (t_o,x_o,y_o), \sD \phi (t_o,x_o,y_o) \big)   \ls   (\hb{resp.\,}   \ge )  \;  0.
 \eeas

  \end{deff}

 For any  $ \phi   \ins C^{1,2,2}\big([0,\infty) \ti \hR^l \ti [0,\infty)\big)$,  we also define
 \beas
  \cL_x\phi(t,x,y)  & \tn \df  & \tn  \frac12 trace \big(  \si (t,x )   \n  \cd  \n   \si^T  (t,x )
   \n \cd \n  D^2_x\phi(t,x,y)    \big)   \+ b^T (t,x )   \n  \cd  \n   D_x\phi(t,x,y) ,     \\
   \cH\phi(t,x,y)  & \tn \df & \tn  \underset{a  \in  \hR^d}{\sup}  \, \Big\{ \frac12 |a|^2   \pa^2_y\phi(t,x,y)
 \+ \big(D_x (\pa_y\phi(t,x,y))  \big)^T \dn \cd \n \si (t,x) \n \cd \n  a    \Big\} \gs 0 ,
 ~   (t,x,y) \ins [0,\infty) \ti \hR^l \ti [0,\infty) ,
 \eeas
 as well as the upper  semi-continuous envelope of $\cH \phi$
  (the smallest upper semi-continuous function above $\cH \phi$)
 \bea \label{eq:d191}
    \ol{\cH} \phi(t,x,y)   \df   \lsup{(t',x',y') \to (t,x,y)} \cH\phi(t',x',y')
    \= \lmtd{\d \to 0} \underset{(t',x',y') \in \sO_\d (t,x,y)}{\sup} \cH\phi(t',x',y')   ,
   ~   (t,x,y) \ins [0,\infty) \ti \hR^l \ti [0,\infty) ,  \q
 \eea
 where $ \sO_\d (t,x,y) \df \big[(t\-\d)^+,t\+\d\big]   \ti  \ol{O}_\d (x) \ti  \big[(y\-\d)^+,y\+\d\big]$.

   \begin{thm} \label{thm_visc_exist}

 Assume that $b,\si$ additionally satisfy \eqref{eq:c601} and  $ f,g$ additionally satisfy   \eqref{eq:c611}.
 Then   the value function  $\cV$ in \eqref{eq:d065} is a  viscosity supersolution of
 \bea
 \begin{cases}
   -  \pa_t u (t,x,y)   \- \cL_x u (t,x,y) \+ g(t,x)    \pa_y u (t,x,y)
   -  \cH u (t,x,y) \- f(t,x) \= 0 , \q \fa (t,x,y) \ins (0,\infty)  \ti  \hR^l \ti (0,\infty) , \q \ss   \\
   u (t,x,0) \= \pi (t,x) , \q \fa (t,x ) \ins [0,\infty)  \ti  \hR^l ,
 \end{cases}
 \label{eq:PDE}
 \eea
 and is   a  viscosity subsolution of
 \bea
 \begin{cases}
   -  \pa_t u (t,x,y)   \- \cL_x u (t,x,y) \+ g(t,x)    \pa_y u (t,x,y)
   -  \ol{\cH} u (t,x,y) \- f(t,x) \= 0 , \q \fa (t,x,y) \ins (0,\infty)  \ti  \hR^l \ti (0,\infty) , \q \ss   \\
   u (t,x,0) \= \pi (t,x) , \q \fa (t,x ) \ins [0,\infty)  \ti  \hR^l .
 \end{cases}
 \label{eq:PDE2}
 \eea

\end{thm}

\begin{rem}
See Section 5.2 of \cite{Miller_C_2017a} for the  connection between the fully non-linear parabolic HJB equation
\eqref{eq:PDE} and generalized Monge-Amp\`ere equations.
\end{rem}

  \section{Proofs}

  \label{sec:proofs}

\subsection{Proofs of Section \ref{sec:OSEC_cont}}

  \no {\bf Proof of Lemma \ref{lem_X_Lp_estimate}:} In this proof,
  we set $ \ol{c} \df \int_0^\infty \n c(s) ds $ and
  let $\fc_q$ denote a generic constant depending only on $q$, whose form may vary from line to line.

 \no {\bf 1)}  Let $T \ins (0,\infty)$ and  set $\wt{q} \df q \ve 2$.
  Given $s \in [0,T]$, we set $\Phi_s \df \underset{r \in [0,s]}{\sup}  \big| X^{t,x}_r    \big|$,  \eqref{FSDE} and \eqref{b_cond} show that
  \bea
   \Phi_s  & \tn \ls & \tn      |x| \+     \int_0^s \n
      \big(  | b  (t\+r, 0 ) | \+ \big| b \big(t\+r, X^{t,x}_r\big) \n - \n  b  (t\+r, 0  )  \big|    \big) dr
 \+   \underset {s' \in [0, s]}{\sup}   \Big|\int_0^{s'} \n   \si(t\+r,    X^{t,x}_r)  dB_r \Big|    \nonumber \\
  &   \tn \ls &  \tn  |x|  \+ \int_0^s \n c(t\+r)  dr  \+ \int_0^s \n c(t\+r)   | X^{t,x}_r    |    dr
       + \n \underset {s' \in [0, s]}{\sup}   \Big|\int_0^{s'} \n   \si(t\+r,    X^{t,x}_r)  dB_r \Big|  , \q
      \hb{$P-$a.s.}      \label{eq:c317}
        \eea
     Taking $\wt{q}-$th power of \eqref{eq:c317}, we can deduce from
     H\"older's inequality, the Burkholder-Davis-Gundy inequality,  \eqref{si_cond} and  Fubini's Theorem   that
    \beas
     && \hspace{-0.7cm} E \big[ \Phi_s^{\wt{q}} \, \big]
        \ls    4^{ \wt{q}-1}   |x|^{\wt{q}} \+  4^{ \wt{q}-1}\ol{c}\,^{\wt{q}} \+
       4^{ \wt{q}-1}     \Big( \int_0^s \n c^{\frac{\wt{q}}{\wt{q}-1}}(t\+r) dr \Big)^{\wt{q}-1} E \bigg[  \int_0^s \n     \big| X^{t,x}_r   \big|^{\wt{q}} dr \bigg]  \+ \fc_q E\bigg[ \Big( \int_0^s \n c(t\+r)  \big(1
       \+ | X^{t,x}_r | \big)^2  d r \Big)^{\frac{\wt{q}}{2}}\bigg]  \nonumber \\
     &&   \ls   4^{ \wt{q}-1}   |x|^{\wt{q}} \+  4^{ \wt{q}-1}\ol{c}\,^{\wt{q}} \+
       4^{ \wt{q}-1}   \Big( \int_0^T \n c^{\frac{\wt{q}}{\wt{q}-1}}(t\+r) dr \Big)^{\wt{q}-1}  \n   \int_0^s \n
          E \big[ \Phi_r^{\wt{q}} \, \big]  dr
    \+   \fc_q \Big( \int_0^T \n c^{\frac{\wt{q}}{\wt{q}-2}}(t\+r) dr \Big)^{\frac{\wt{q}}{2}-1}
    \n   \int_0^s \n E\big[ (1 \+ \Phi_r)^{\wt{q}} \,  \big] d r   \nonumber \\
     & & \ls   \fc_q \bigg[  |x|^{\wt{q}} \+  \ol{c}\,^{\wt{q}}
     \+ T \Big( \int_0^T \n c^{\frac{\wt{q}}{\wt{q}-2}}(t\+r) dr \Big)^{\frac{\wt{q}}{2}-1} \bigg]
      \+ \bigg[ 4^{ \wt{q}-1}   \Big( \int_0^T \n c^{\frac{\wt{q}}{\wt{q}-1}}(t\+r) dr \Big)^{\wt{q}-1}
       \+  \fc_q \Big( \int_0^T \n c^{\frac{\wt{q}}{\wt{q}-2}}(t\+r) dr \Big)^{\frac{\wt{q}}{2}-1} \bigg] \n   \int_0^s \n
          E \big[ \Phi_r^{\wt{q}} \, \big]  dr .   
       \eeas
     An application of  Gronwall's inequality then gives that
       \bea
 E\big[  \Phi_s^q \, \big]   & \tn  \ls  & \tn   1   \+ E\big[ \Phi_s^{\wt{q}} \, \big]
         \ls      1 \+ \fc_q   \bigg[  |x|^{\wt{q}} \+  \ol{c}\,^{\wt{q}}
        \+  T  \Big( \int_0^T \n c^{\frac{\wt{q}}{\wt{q}-2}}(t\+r) dr \Big)^{\frac{\wt{q}}{2}-1}  \bigg]    \nonumber \\
   & \tn     & \tn
   \times    \exp\bigg\{ 4^{\wt{q}-1} \Big( \int_0^T \n c^{\frac{\wt{q}}{\wt{q}-1}}(t\+r) dr \Big)^{\wt{q}-1} s
   \+   \fc_q \Big( \int_0^T \n c^{\frac{\wt{q}}{\wt{q}-2}}(t\+r) dr \Big)^{\frac{\wt{q}}{2}-1} s \bigg\}
    \< \infty ,  \q  \fa s \ins [0,T]. \qq  \qq   \label{eq:c111}
       \eea

 Let $s \in [0,T]$. Since the Burkholder-Davis-Gundy inequality and  \eqref{b_cond} also show   that
   \beas
  && \hspace{-0.7cm} E\bigg[ \underset {s' \in [0, s]}{\sup}   \Big|\int_0^{s'} \n   \si(t\+r,    X^{t,x}_r)  dB_r \Big|^q \bigg]
   \ls \fc_q E\bigg[ \Big( \int_0^s \n  \big| \si(t\+r,    X^{t,x}_r) \big|^2 d r   \Big)^{\frac{q}{2}} \bigg]
   \ls    \fc_q E\bigg[ \Big( \int_0^s \n c(t\+r)  \big(1
       \+    | X^{t,x}_r    | \big)^2  d r \Big)^{\frac{q}{2}}\bigg] \\
  && \ls  \fc_q E\bigg[ (1  \+  \Phi_s  )^{\frac{q}{2}}
       \Big( \int_0^s \n c(t\+r)  \big(1  \+  | X^{t,x}_r    | \big)   d r \big)^{\frac{q}{2}}\bigg]
      \ls  E\bigg[ \frac12 8^{1-q}  (1  \+  \Phi_s  )^q
       \+ \fc_q   \Big( \int_0^s \n c(t\+r)  \big(1   \+ | X^{t,x}_r  |  \big)   d r \Big)^q \bigg]\\
 &&      \ls \frac12 4^{1-q} \big( 1 \+ E \big[ \Phi_s^q \, \big] \big) \+ \fc_q   \Big( \int_0^s \n c(t\+r) dr \Big)^q
       \+ \fc_q E\bigg[ \Big( \int_0^s \n c(t\+r)   \Phi_r   d r \Big)^q \bigg] ,
   \eeas
   taking $q-$th power of \eqref{eq:c317} and using  Fubini's Theorem yield   that
    \bea
   4^{1-q} E\big[  \Phi_s^q \, \big]
       & \tn \ls   & \tn |x|^q \+ \ol{c}^q
          \+ E\bigg[ \Big( \int_0^s \n c(t\+r)  \Phi_r  dr \Big)^q \bigg]
          \+ E\bigg[ \underset {s' \in [0, s]}{\sup}   \Big|\int_0^{s'} \n   \si(t\+r,    X^{t,x}_r)  dB_r \Big|^q \bigg]    \nonumber \\
       & \tn  \ls   & \tn   |x|^q \+   \fc_q  \ol{c}^q
     \+ \frac12 4^{1-q} \big( 1 \+   E\big[  \Phi_s^q \, \big] \big)
     \+  \fc_q \Big( \int_0^s \n c(t\+r) dr \Big)^{q-1} E \bigg[ \int_0^s \n c(t\+r) \Phi_r^q dr \bigg]    .
     \qq  \label{eq:c315}
       \eea
       Here, we applied H\"older's inequality  $\big|\int_0^s \n \fra_r\fb_r dr \big|\ls
       \big( \int_0^s \n |\fra_r|^q dr \big)^{\frac{1}{q}}
       \big( \int_0^s \n |\fb_r|^{\frac{q}{q-1}} dr \big)^{\frac{q-1}{q}}$
       with $\big(\fra_r,\fb_r\big) \= \Big( c^{\frac{1}{q}}(t+r)  \Phi_r ,   c^{\frac{q-1}{q}}(t+r) \Big)$.
       As $ E\bigg[  \underset{r \in [0,s]}{\sup}  \big| X^{t,x}_r   \big|^q \bigg] \< \infty $ by \eqref{eq:c111},
       it follows from \eqref{eq:c315} that for any $s \ins [0,T]$
       \bea
     E \big[  \Phi_s^q \, \big]   \ls 1 \+  2 \ti  4^{q-1} |x|^q \+ \fc_q \ol{c}^q
     \+  \fc_q \ol{c}^{q-1}  \n \int_0^s \n c(t\+r) E \big[  \Phi_r^q \, \big]  dr   .
 \label{eq:c319}
       \eea
  Applying  Gronwall's inequality again  yields that
  $   E \big[  \Phi_s^q \, \big]   \ls   \big( 1 \+  2 \ti  4^{q-1} |x|^q \+ \fc_q \ol{c}^q    \big)
   \exp \big\{\fc_q \ol{c}^{q-1} \n \int_0^s \n c(t\+r) dr \big\}  $, $\fa  s \ins [0,T]$.
    In particular,  taking $s \= T$ and then letting $T \n \to \n \infty$,
    one can deduce    from  the monotone convergence theorem that
    \bea \label{eq:c411}
    E\bigg[ \, \underset{r \in [0,\infty)}{\sup}  \big| X^{t,x}_r   \big|^q \bigg]
       \ls   \big( 1 \+  2 \ti  4^{q-1} |x|^q \+ \fc_q \ol{c}^q    \big)
       \exp\big\{ \fc_q \ol{c}^q  \big\} \, . \qq
    \eea

  \no {\bf 2)} Let $\sX_s \df X^{t,x}_s \n - \n X^{t,x'}_s $, $\fa s \ins [0,\infty)$.
  Given $s \ins [0,\infty)$, we set $\wt{\Phi}_s \df  \underset{r \in [0,s]}{\sup} \big| \sX_r \big|$.
  Since an analogy to \eqref{eq:c317} shows that
  \beas
    \wt{\Phi}_s  \n \le  \n    |x' \- x|
   \n + \n    \int_0^s  \n c(t\+r) | \sX_r | dr
    \n + \n  \underset{s' \in [0,s]}{\sup} \Big|  \n \int_0^{s'}  \dn   \big( \si   (t\+r,X^{t,x}_r  )
    \n - \n  \si   (t\+r , X^{t,x'}_r   ) \big) d B_r \Big| , \q  \hb{$P-$a.s.,}
  \eeas
  the Burkholder-Davis-Gundy inequality and \eqref{si_cond}   imply that
    \beas
  3^{1-q}  E \big[  \wt{\Phi}_s^q   \big]
   & \tn \le&  \tn    |x' \- x|^q   \+   E \bigg[ \Big( \n \int_0^s c(t\+r)| \sX_r | dr \Big)^q \bigg]
   \+ \fc_q E \bigg[ \Big(   \int_0^s  \n  c(t\+r) | \sX_r  |^2  d  r \Big)^{\frac{q}{2}} \bigg]  \nonumber  \\
   & \tn \le& \tn    |x' \- x|^q
   \+   E \bigg[ \Big( \n \int_0^s c(t\+r)| \sX_r | dr \Big)^q \bigg] \+
   \fc_q  E \bigg[ \wt{\Phi}_s^{q/2}  \bigg( \int_0^s  \n c(t\+r)  | \sX_r  |   d  r \bigg)^{\frac{q}{2}} \Bigg]  \nonumber \\
   & \tn \le& \tn    |x' \- x|^q
   \+ \frac12 3^{1-q} E \big[   \wt{\Phi}_s^q \,   \big]
   \+  \fc_q  E \bigg[ \Big( \n \int_0^s c(t\+r)| \sX_r | dr \Big)^q \bigg] .
  \eeas
    Since $  E \big[  \wt{\Phi}_s^q \, \big] \ls 2^{q-1} E \Big[     \big( X^{t,x}_* \big)^q
     \+     \big( X^{t,x'}_* \big)^q \Big] \< \infty$      by Part 1, 
     an analogy to \eqref{eq:c319} shows that
     \beas
      E \big[  \wt{\Phi}_s^q   \big]
     \ls      2 \ti 3^{q-1} |x' \- x|^q
    \+  \fc_q \ol{c}^{q-1} \n  \int_0^s \n c(t\+r)  E \big[ \wt{\Phi}_r^q \big]  dr  ,
    \q   \fa s  \n \in \n  [0, \infty)  .
    \eeas
  Then we see from  Gronwall's inequality  that
   $   E \big[  \wt{\Phi}_s^q   \big]
   \ls  2 \ti 3^{q-1} |x' \- x|^q \exp \big\{\fc_q \ol{c}^{q-1} \int_0^s \n c(t\+r) dr \big\}$,
     $ \fa  s  \ins  [0,\infty) $.
   As $s \n \to \n \infty$,   the monotone convergence theorem implies  that
   $    E\bigg[ \, \underset{r \in [0,\infty)}{\sup}  \big| X^{t,x'}_r \n \- \n  X^{t,x}_r  \big|^q \bigg]
    \ls  2 \ti 3^{q-1} |x' \- x|^q  \exp\big\{ \fc_q \ol{c}^q  \big\}  $.

  \no {\bf 3)}  Let   $ \d \ins (0,\infty) $ and $ \tau \ins \cT $.
  For any $\l \ins (0,\d]$, since it holds $P-$a.s.  that
  \beas
  \hspace{-3mm}
    X^{t,x}_{\tau+\l} \- X^{t,x}_\tau
    \=    \int_\tau^{\tau+\l} b(t\+r,X^{t,x}_r) dr \+ \int_\tau^{\tau+\l} \si (t\+r,X^{t,x}_r) dB_r
    \=    \int_\tau^{\tau+\l} b(t\+r,X^{t,x}_r) dr \+ \int_0^{\tau+\l} \b1_{\{\tau < r < \tau+\d \}} \si (t\+r,X^{t,x}_r) dB_r,
  \eeas
  taking $q-$th power and using \eqref{b_cond} yield that
  \beas
    \big| X^{t,x}_{\tau+\l} \- X^{t,x}_\tau \big|^q
  \ls 2^{q-1} \Big( \int_\tau^{\tau+\l} c(t\+r) \big( 1\+ |X^{t,x}_r| \big) dr \Big)^q
  \+ 2^{q-1} \underset{s \in [0,\infty)}{\sup} \Big| \int_0^s \n \b1_{\{\tau < r < \tau+\d \}} \si (t\+r,X^{t,x}_r) dB_r \Big|^q , \q \hb{$P-$a.s.}
  \eeas
  Then the Burkholder-Davis-Gundy inequality shows that
  \beas
  E \bigg[ \underset{\l \in (0,\d]}{\sup} \big| X^{t,x}_{\tau+\l} \- X^{t,x}_\tau \big|^q  \bigg]
  & \tn \ls  & \tn  2^{q-1} E \bigg[  \Big( \int_\tau^{\tau+\d} c(t\+r) \big( 1\+ |X^{t,x}_r| \big) dr \Big)^q \bigg]
  \+  \fc_q E  \bigg[ \Big( \int_\tau^{\tau+\d} \n \big|\si (t\+r,X^{t,x}_r)\big|^2 dr \Big)^{\frac{q}{2}} \bigg]  \\
   & \tn \ls & \tn  2^{q-1} \d^q \|c(\cd)\|^q
    E \bigg[ \Big( 1\+ \underset{r \in [0,\infty)}{\sup} |X^{t,x}_r| \Big)^q \bigg]
  \+ \fc_q E  \bigg[ \Big( \int_\tau^{\tau+\d} \n c(t\+r) \big( 1\+ |X^{t,x}_r| \big)^2 dr \Big)^{\frac{q}{2}} \bigg] \\
   & \tn \ls & \tn  \fc_q \big( \d^q \|c(\cd)\|^q  \+   \d^{\frac{q}{2}} \|c(\cd)\|^{\frac{q}{2}} \big)
  \bigg( 1 \+ E \bigg[ \,  \underset{r \in [0,\infty)}{\sup} |X^{t,x}_r |^q \bigg] \bigg) ,
  \eeas
 which together with   \eqref{eq:c411} leads to \eqref{eq:esti_X_2}.

   \no {\bf 4)} Now, we assume functions $b$ and $\si$ satisfy \eqref{eq:c601} for some $\varpi \ins [1,\infty)$.
   Let $t' \ins (t,\infty)$  and define    $\wh{\sX}_s \df X^{t' \n,x}_s \- X^{t,x}_s   $, $\fa s \ins [0,T] $.
  By \eqref{eq:c601}, it holds $P-$a.s. that
  \bea
  \big| b(t'\+r,X^{t' \n,x}_r) \- b(t\+r,X^{t,x}_r) \big|
  & \tn \ls  & \tn  \big| b(t'\+r,X^{t' \n,x}_r) \- b(t'\+r,X^{t,x}_r) \big|
  \+ \big| b(t'\+r,X^{t,x}_r) \- b(t\+r,X^{t,x}_r) \big| \qq  \nonumber \\
   & \tn \ls & \tn  c(t'\+r) \big| \wh{\sX}_r \big|
  \+ c(t\+r) \rho(t'\-t) \big( 1\+|X^{t,x}_r|^\varpi \big) , \q \fa r \ins [0,\infty) ,  \label{eq:c603}
  \eea
  and similarly that
  \bea
  \big| \si (t'\+r,X^{t' \n,x}_r) \- \si (t\+r,X^{t,x}_r) \big|
   \ls  \sqrt{ c(t'\+r) } \big| \wh{\sX}_r \big|
  \+ \sqrt{ c(t\+r) } \rho(t'\-t) \big( 1\+|X^{t,x}_r|^\varpi \big) , \q \fa r \ins [0,\infty) .   \label{eq:c605}
  \eea

   Given $s \ins [0,\infty)$, we set
  $\wh{\Phi}_s \df  \underset{r \in [0,s]}{\sup} \big| \wh{\sX}_r \big|$, \eqref{eq:c603} shows  that $P-$a.s.
  \bea \label{eq:c607}
  \wh{\Phi}_s   \ls   \int_0^s  \n c(t'\+r) \big| \wh{\sX}_r \big| dr
  \+ \rho(t'\-t) \int_0^s c(t\+r)   \big( 1\+|X^{t,x}_r|^\varpi \big) dr
    \n + \n  \underset{s' \in [0,s]}{\sup} \Big|  \n \int_0^{s'}  \dn   \big( \si   (t'\+r, X^{t' \n,x}_r  )
    \n - \n  \si   (t\+r , X^{t,x}_r   ) \big) d B_r \Big| .
  \eea
    The Burkholder-Davis-Gundy inequality, \eqref{eqn-d011} and \eqref{eq:c605} imply that
   \beas
 && \hspace{-1.5cm}
  E \bigg[   \underset{s' \in [0,s]}{\sup} \Big|  \n \int_0^{s'}  \dn   \big( \si   (t'\+r, X^{t' \n,x}_r  )
    \n - \n  \si   (t\+r , X^{t,x}_r   ) \big) d B_r \Big|^q \bigg]
    \ls \fc_q E\bigg[ \Big( \int_0^s \big| \si (t'\+r,X^{t' \n,x}_r) \- \si (t\+r,X^{t,x}_r) \big|^2 dr \Big)^{\frac{q}{2}} \bigg]  \\
 && \ls   \fc_q  E\bigg[ \wh{\Phi}^{q/2}_s \Big( \int_0^s  \n  c(t'\+r) |\wh{\sX}_r | dr \Big)^{\frac{q}{2}} \bigg]
 \+ \fc_q \big(\rho (t'\-t)\big)^q  E\bigg[ \Big( \int_0^s  c(t\+r)   \big( 1\+|X^{t,x}_r|^\varpi \big)^2 dr \Big)^{\frac{q}{2}} \bigg] \\
 && \ls   \frac12 3^{1-q}  E \big[  \wh{\Phi}^q_s \big]
 \+ \fc_q E \bigg[   \Big( \int_0^s    c(t'\+r) |\wh{\sX}_r | dr \Big)^q \bigg]
 \+  \fc_q \ol{c}^{\frac{q}{2}}  \big(\rho (t'\-t)\big)^q  E\bigg[ \Big(1\+\underset{r \in [0,s]}{\sup} |X^{t,x}_r|^\varpi \Big)^q \bigg] .
  \eeas
   Taking $q-$th power in \eqref{eq:c607} and using an analogy to \eqref{eq:c319} yield that
  \beas
   3^{1-q} E \big[ \wh{\Phi}^q_s \big]
  & \tn \ls & \tn  E \bigg[ \Big( \int_0^s  \n c(t'\+r) | \wh{\sX}_r | dr \Big)^q \bigg]
  \+ \big(\rho(t'\-t)\big)^q  E \bigg[ \Big( \int_0^s c(t\+r)   \big( 1\+|X^{t,x}_r|^\varpi \big) dr \Big)^q \bigg] \\
  & \tn  & \tn   +  E \bigg[   \underset{s' \in [0,s]}{\sup} \Big|  \n \int_0^{s'}  \dn   \big( \si   (t'\+r, X^{t' \n,x}_r  )
    \n - \n  \si   (t\+r , X^{t,x}_r   ) \big) d B_r \Big|^q \bigg] \\
    & \tn  \ls  & \tn    \frac12 3^{1-q}  E \big[  \wh{\Phi}^q_s \big]
 \+ \fc_q \ol{c}^{q-1}   \int_0^s    c(t'\+r) E \big[ \wh{\Phi}^q_r  \big] dr
 \+ \fc_q \big(  \ol{c}^{\frac{q}{2}} \+  \ol{c}^q \big)  \big(\rho (t'\-t)\big)^q  E\bigg[  1\+\underset{r \in [0,s]}{\sup} |X^{t,x}_r|^{q \varpi}   \bigg] . \q \qq
  \eeas

  As $  E \big[  \wh{\Phi}_s^q \, \big] \ls 2^{q-1} E \Big[     \big( X^{t,x}_* \big)^q
     \+     \big( X^{t' \n,x}_* \big)^q \Big] \< \infty$      by Part 1,
  it then follows from Gronwall's inequality that
     \beas
     E \big[ \wh{\Phi}^q_s \big] \ls \fc_q \big(  \ol{c}^{\frac{q}{2}} \+  \ol{c}^q \big)  \big(\rho (t'\-t)\big)^q
 E\bigg[  1\+\underset{r \in [0,\infty)}{\sup} |X^{t,x}_r|^{q \varpi}   \bigg]  \exp \Big\{\fc_q \ol{c}^{q-1} \n \int_0^s \n c(t'\+r) dr \Big\}   , \q  \fa  s \ins [0,\infty) .
     \eeas
   Letting $s \n \to \n \infty$,
    we can deduce    from  the monotone convergence theorem that
    \beas
         E \bigg[ \underset{r \in [0,\infty)}{\sup} \big|  X^{t' \n,x}_r \- X^{t,x}_r  \big|^q \bigg] \ls \fc_q \big(  \ol{c}^{\frac{q}{2}} \+  \ol{c}^q \big)  \big(\rho (t'\-t)\big)^q
 E\bigg[  1\+\underset{r \in [0,\infty)}{\sup} |X^{t,x}_r|^{q \varpi}   \bigg]  \exp \big\{\fc_q \ol{c}^q \big\} ,
    \eeas
    which together with \eqref{eq:c411} proves \eqref{eq:c609}.  \qed

\no {\bf Proof of \eqref{eq:c483}:}
 We see from   \eqref{eq:c343}   that
   \bea \label{eq:c363}
   |f(t',x')| \ls  |f(t',x')\-f(t',0)|  \+  |f(t',0)| \ls c(t') \big(1\+|x'|\ve |x'|^p \big) \ls c(t') \big(2\+  |x'|^p \big) , \q \fa (t',x') \ins (0,\infty) \ti \hR^l . \q
   \eea
   Similarly,   \eqref{eq:c321}  shows that
 \bea \label{eq:c373}
|\pi  (t',x')|   \ls \fC \big(2\+|x'|^p\big), \q  \fa  (t',x') \ins (0,\infty) \ti \hR^l .
\eea

 Given $\tau \ins  \cT $, Since \eqref{eq:c363}, \eqref{eq:c373} show that
 \bea \label{eq:d117}
  \big|  R (t,x,\tau) \big| \ls \big( 2\+ (X^{t,x}_*)^p \big) \n  \int_0^\infty c(t\+r) dr
   \+ \fC \big( 2\+ (X^{t,x}_*)^p \big)
 \ls 2 \fC \big( 2\+ (X^{t,x}_*)^p \big)  ,
 \eea
   the first inequality in \eqref{eq:esti_X_1} implies   that
 $  E  \big[ \big| R (t,x,\tau) \big| \big] 
 \ls 2 \fC \big( 2\+   C_p (1\+|x|^p) \big) \= \Psi (x)  $.  \qed

\no   {\bf Proof of Lemma \ref{lem_general_value}:}
Let $(t,x,y) \ins  [0,\infty) \ti  \hR^l \ti (0,\infty)$.
 Since $   \cF_0$ consist of $\cF-$measurable sets $A$ with
 $P(A) \= 0$ or $P(A) \= 1$, it holds for any $\tau \ins \cT$ that
  \bea \label{eq:d063}
  P \{\tau \= 0\} \= 1 \hb{ or } P \{\tau \> 0\} \= 1   .
  \eea
  It follows that
 $  V( t, x,y )   \= E \big[ R  (t, x, 0) \big]
 \ve \Big( \underset{\tau \in \wh{\cT}_{t,x}(y) }{\sup} \,   E \big[ R \big(t, x, \tau \big) \big] \Big)
 $. So it suffices to show that
   $E \big[ R  (t, x, 0) \big] \ls \underset{\tau \in \wh{\cT}_{t,x}(y)}{\sup} \, E \big[ R \big(t, x, \tau \big) \big] $.

  We arbitrarily pick up $\tau$ from $ \wh{\cT}_{t,x}(y) $.
  Given $n \ins \hN$, it is clear that   $ \tau_n \df \tau \ld (1/n)  $ also belongs to $ \wh{\cT}_{t,x}(y) $, so
  \bea \label{eq:d185}
  \underset{\tau \in \wh{\cT}_{t,x}(y)}{\sup} \, E \big[ R \big(t, x, \tau \big) \big]
  \gs 
  E \big[ R \big(t, x, \tau_n \big) \big]
  \= E \bigg[ \int_0^{\tau_n} f(t\+s,X^{t,x}_s) ds \+ \pi \big(t\+\tau_n, X^{t,x}_{\tau_n} \big) \bigg] .
  \eea
  An analogy to \eqref{eq:d117} shows that
  $   \big|  R (t,x,\tau_n) \big| \ls \big( 2\+ (X^{t,x}_*)^p \big) \n  \int_t^\infty c(r) dr
   \+ \fC \big( 2\+ (X^{t,x}_*)^p \big)
 \ls 2 \fC \big( 2\+ (X^{t,x}_*)^p \big)$, whose $E-$expectation equals to
 $ 2 \fC \big( 2\+   C_p (1\+|x|^p) \big) \= \Psi (x) $ by  the first inequality in \eqref{eq:esti_X_1b}.
   Then letting $n \to \infty$ in \eqref{eq:d185},
   we can deduce from \eqref{eq:c321}, the continuity of process $X^{t,x}$
   and the dominated convergence theorem that
   \beas
  \q \underset{\tau \in \wh{\cT}_{t,x}(y)}{\sup} \, E \big[ R \big(t, x, \tau \big) \big]
  \gs \lmt{n \to \infty} E \bigg[ \int_0^{\tau_n} f(t\+s,X^{t,x}_s) ds \+ \pi \big(t\+\tau_n, X^{t,x}_{\tau_n} \big) \bigg]
  \= E \big[   \pi  (t, X^{t,x}_0  ) \big] \= \pi  (t, x) \= E \big[ R  (t, x, 0) \big] . \q \;  \hb{\qed}
   \eeas

\no {\bf Proof of Theorem \ref{thm_continuity}: 1)}  Fix $t \ins [0,\infty)$.
 We let $(x,\e) \ins   \hR^l \ti (0,1) $  and set $\e_o \df (5  \+  10 \fC  )^{-1} \e $.
  Since $\fM \df E \big[   (X^{t,x}_*)^p \big] \< \infty$
 by the first inequality in \eqref{eq:esti_X_1}, we can find   $\l_o \= \l_o (t,x, \e)
 \ins  \big(0,\e_o\big) $ such that
\bea \label{eq:c353}
E \big[\b1_A \, (X^{t,x}_*)^p \big]
\< \e_o \; \hb{ for any $A \ins \cF $ with $P (A) \< \l_o $\,. }
\eea
 There   exists   $R \= R (t,x, \e) \ins (0,\infty)$ such that the set
$A_R \df \big\{  X^{t,x}_*  \> R \big\} \ins \cF$
satisfies  $P (A_R) \< \l_o / 2 $.

Let   $ \l \= \l(t,x, \e) \ins   (0, 1  )  $  satisfy that
\bea
 && \hspace{3cm}  \sqrt{\l} \ls \Big(\frac16  \l_o   \k_{\n \overset{}{R}}\Big)  \ld
    \frac{\e_o}{  (2\+ \fM)  \|c(\cd)\|  }   \n \land \rho^{-1}(  \e_o ) \q \hb{and} \q \label{eq:c357} \\
 &&  (C_p)^{\frac{1}{p}} \big(1\+|x| \big)
 \big( \|c(\cd)\| \l^{\frac{1}{2}} \+ \|c(\cd)\|^{\frac{1}{2}}  \l^{\frac{1}{4}} \big)
    \+    C_p \big(1\+|x|^p \big) \big( \|c(\cd)\|^p \l^{\frac{p}{2}} \+ \|c(\cd)\|^{\frac{p}{2}}  \l^{\frac{p}{4}} \big) \ls \e_o \, .  \qq  \qq   \label{eq:c361}
\eea
 We pick up $ \d \= \d(t,x, \e) \ins (0,1)  $ such that
  \bea \label{eq:c351}
  \fC(C_p)^{\frac{1}{p}} \d \+   \fC  C_p \d^p   \ls \l \ld  \e_o  \,  ,
  \eea
  and   fix $ y  \ins   [0,\infty) $.


 \no {\bf 1a)} {\it We first demonstrate that
$V(t,\fx,\fy) \gs   V(t,x,y) \- \e $, $ \fa (\fx,\fy) \ins \ol{O}_\d (x) \ti [(y\-\d)^+,\infty)   $. }

Let $\tau_1  \= \tau_1(t,x,y,\e)  \ins \cT_{t,x} (y)$ such that
\bea \label{eq:c387}
E \big[  R(t,x, \tau_1  )   \big]     \gs    V(t,x,y) \- \e_o  ,
\eea
 and let $\fx \ins \ol{O}_\d (x)$.

     We claim that there exists a  stopping time
  $  \wh{\tau}_1  \=  \wh{\tau}_1  (t,x,\fx, y,  \e) \ins  \cT_{t,\fx} \big( (y\-\d)^+ \big)$  satisfying
 \bea \label{eq:c333}
  {\wh{\tau}_1} \ls \tau_1 \q \hb{and} \q P \big( A^c_R \Cp  \{ \tau_1 \> {\wh{\tau}_1} \+ \sqrt{\l} \, \}\big)
   \< \l_o / 2 .
 \eea
 Set $\d_y \df \d \ld y$, which satisfies   $ y\-\d_y \= (y\-\d) \ve (y \- y) \= (y\-\d)^+$.
 \beas \label{eq:d227}
 \hb{If $ E \big[ \int_0^{\tau_1} \n g (t\+r, X^{t,\fx}_r) dr \big] \ls y \- \d_y  $
 \big(i.e. $\tau_1  \ins  \cT_{t,\fx}  ( (y\-\d)^+  )$\big),   we directly set ${\wh{\tau}_1} \df \tau_1 $.}
 \eeas
  Otherwise,     set $\fra \df E \big[ \int_0^{\tau_1} \n g (t\+r, X^{t,\fx}_r) dr \big]  \- y  \+ \d_y \> 0 $
  (In this case, one must have $y \> 0$).
  Since both  $ \big\{ E \big[\int_0^{\tau_1} \n g (t\+r, X^{t,\fx}_r) dr \big| \cF_s \big] \big\}_{s \in [0,\infty)}$
and $ \big\{ \int_0^s \n g (t\+r, X^{t,\fx}_r) dr \big\}_{s \in [0,\infty)}$ are $\bF-$adapted continuous processes,
\bea \label{eq:c335}
\wh{\tau}_1 \= \wh{\tau}_1 (t,x,\fx,y, \e)
 \df \inf\bigg\{s \ins [0,\infty) \n : E \Big[\int_0^{\tau_1} \n g (t\+r, X^{t,\fx}_r) dr \Big| \cF_s \Big]
\- \int_0^s \n g (t\+r, X^{t,\fx}_r) dr  \ls  \fra  \bigg\}
\eea
defines an   $\bF-$stopping time   which  satisfies
$ E \big[\int_0^{\tau_1} \n g (t\+r, X^{t,\fx}_r) dr \big| \cF_{\wh{\tau}_1} \big]
\- \int_0^{\wh{\tau}_1} \n g (t\+r, X^{t,\fx}_r) dr  \=   \fra $.
  Taking expectation $E [\cd]$ yields that
 \bea \label{eq:c101}
 E \Big[ \int_0^{\wh{\tau}_1} \n g (t\+r, X^{t,\fx}_r) dr \Big]
 \= E \Big[ \int_0^{\tau_1} \n g (t\+r, X^{t,\fx}_r) dr   \Big]  \- \fra \= y \- \d_y \= (y\-\d)^+ , \q
 \hb{so } \, {\wh{\tau}_1} \ins \cT_{t,\fx} \big( (y\-\d)^+ \big) .
\eea
 As $E \big[\int_0^{\tau_1} \n g (t\+r, X^{t,\fx}_r) dr \big| \cF_{\tau_1} \big]
 \- \int_0^{\tau_1} \n g (t\+r, X^{t,\fx}_r) dr \= 0 \<  \fra $,
 we also see that $   \wh{\tau}_1 \ls \tau_1 $.

  The condition (g1), H\"older's inequality,  the second inequality in \eqref{eq:esti_X_1}
   and \eqref{eq:c351} show that
  \bea
  \q  && \hspace{-1.5cm}
       E \Big[ \int_0^\infty \n \big| g(t\+r, X^{t,\fx}_r) \- g(t\+r, X^{t,x}_r) \big| dr  \Big] \nonumber \\
   & & \ls  E \bigg[ \big(   (X^{t,\fx} \n \- \n X^{t,x})_*
  \+  (X^{t,\fx} \n \- \n X^{t,x} )^p_* \big) \n  \int_0^\infty \n c(t\+r) dr  \bigg]
  \ls \fC(C_p)^{\frac{1}{p}} |\fx\-x|\+ \fC C_p |\fx\-x|^p   \ls   \l  . \qq  \label{eq:c331}
  \eea
    Since $ E \big[ \int_0^{\tau_1} \n     g (t\+r, X^{t,x}_r)    dr   \big] \ls y  $ and since
    $\l \gs  \l \ld \e_o  \> \fC (C_p)^{\frac{1}{p}} \d
  \gs \d \gs \d_y $ by \eqref{eq:c351}, one has
 \beas
  \fra  & \tn  \= & \tn  E \Big[ \int_0^{\tau_1} \n g (t\+r, X^{t,\fx}_r) dr \Big] \- y \+ \d_y
  \< E \Big[\int_0^{\tau_1} \n \big( g (t\+r, X^{t,\fx}_r)   \-   g (t\+r, X^{t,x}_r) \big)  dr   \Big] \+ \l \\
   & \tn \ls  & \tn  E \Big[ \int_0^\infty \n \big| g(t\+r, X^{t,\fx}_r) \- g(t\+r, X^{t,x}_r)    \big| dr  \Big] \+ \l
  \ls    2 \l .
 \eeas
 Using \eqref{eq:c331} again, we can deduce from   \eqref{eq:c101} that
   \bea
    2\l & \tn \> & \tn \fra
   \=    E \Big[ \int_{\wh{\tau}_1}^{\tau_1} \n g(t\+r, X^{t,\fx}_r) dr  \Big]
   \gs E \Big[     \int_{\wh{\tau}_1}^{\tau_1} \n g(t\+r, X^{t, x}_r) dr \Big]
   \- E \Big[ \int_0^\infty \n \big| g(t\+r, X^{t,\fx}_r) \- g(t\+r, X^{t,x}_r) \big| dr  \Big] \nonumber  \\
   & \tn \gs  & \tn  E \bigg[    \b1_{A^c_R \cap \{ \tau_1 >  {\wh{\tau}_1}  + \sqrt{\l} \}}
   \int_{\wh{\tau}_1}^{\tau_1} \n g(t\+r, X^{t, x}_r) dr \bigg]   \- \l
  \gs  \k_{\n \overset{}{R}} \, \sqrt{\l} \, P \big( A^c_R \Cp  \{ \tau_1 \> {\wh{\tau}_1} \+ \sqrt{\l} \, \}\big)  \- \l .
  \label{eq:c341}
  \eea
  It follows from \eqref{eq:c357} that
 $ P \big( A^c_R \Cp  \{ \tau_1 \> {\wh{\tau}_1} \+ \sqrt{\l} \, \}\big)
 \< 
   \frac{ 3  \sqrt{\l}  }{ \k_{\n \overset{}{R}}} \ls  \l_o / 2  $, proving the claim \eqref{eq:c333}.

  Set $\cA \df \{ \tau_1 \ls {\wh{\tau}_1} \+    \sqrt{\l}  \}
  \=  \{ {\wh{\tau}_1} \ls \tau_1 \ls {\wh{\tau}_1} \+    \sqrt{\l}  \}  $.
 Since  \eqref{eq:c333} shows that
\beas 
 P (\cA^c ) \= P \{ \tau_1 \> {\wh{\tau}_1} \+    \sqrt{\l}  \}
  \ls P(A_R)  \+ P \big( A^c_R \cap  \big\{ \tau_1 \> {\wh{\tau}_1} \+ \sqrt{\l} \, \big\}\big) \< \l_o \< \e_o  ,
\eeas
   \eqref{eq:c363}$-$\eqref{eq:c357} imply that
   \bea
   E \bigg[ \Big| \int_0^{\wh{\tau}_1 } \n  f  (t\+r, X^{t,x}_r  )  dr \- \int_0^{\tau_1 } \n f  (t\+r,X^{t,x}_r) dr  \Big| \bigg]
    & \tn  \ls   & \tn    E \bigg[   \big(  2 \+ (X^{t,x}_*)^p \big)  \Big( \b1_{\cA^c} \n \int_0^\infty \n c(t\+r) dr
  \+  \b1_\cA \|c(\cd)\|  ( \tau_1  \-  \wh{\tau}_1 ) \Big)    \bigg] \qq  \nonumber \\
      & \tn   \<  & \tn   \fC \big( 2  P(\cA^c) \+   \e_o \big)   \+  \sqrt{\l} (2\+ \fM) \|c(\cd)\|
     \< (1\+3\fC) \e_o \, , \label{eq:c381} \\
 \hb{and} \q E \big[  \b1_{\cA^c}  \big| \pi \big( \wh{\tau}_1  , X^{t,x}_{\wh{\tau}_1 } \big)
 \- \pi \big( \tau_1  , X^{t,x}_{\tau_1 } \big) \big| \big]  & \tn \ls  & \tn
   2 \fC E \big[  \b1_{\cA^c}  \big(2 \+   (X^{t, x}_*   )^p \big)  \big]
  \< 2  \fC \big( 2  P(\cA^c) \+   \e_o \big)  \< 6 \fC   \e_o \, .  \label{eq:c367}
\eea
Also, we can deduce from    \eqref{eq:c321}, \eqref{eq:c357}, H\"older's inequality,
   \eqref{eq:esti_X_2}  and \eqref{eq:c361}   that
 \bea
  && \hspace{-1cm} E \big[  \b1_\cA   \big| \pi \big( \wh{\tau}_1  , X^{t,x}_{\wh{\tau}_1 } \big)
 \- \pi \big( \tau_1  , X^{t,x}_{\tau_1 } \big) \big| \big]
   \ls      E \big[  \b1_\cA  \rho \big( \tau_1 \- \wh{\tau}_1  \big) \big]
\+ \fC E \Big[ \b1_\cA \Big( \big| X^{t,x}_{\wh{\tau}_1 } \n \-    X^{t,x}_{\tau_1 } \big| \+
 \big| X^{t,x}_{\wh{\tau}_1 }  \n \-    X^{t,x}_{\tau_1 } \big|^p \Big)  \Big] \nonumber \\
  & &  \ls     \rho \big(\sqrt{\l}\big)
 \+  \fC  \bigg\{ E \bigg[ \b1_\cA  \underset{r \in (0,\sqrt{\l}\,]}{\sup} \big| X^{t, x}_{\wh{\tau}_1 +r} \n \-   X^{t,x}_{\wh{\tau}_1 } \big|^p \bigg] \bigg\}^{\frac{1}{p}}
 \+ \fC    E \bigg[ \b1_\cA  \underset{r \in (0,\sqrt{\l}\,]}{\sup} \big| X^{t, x}_{\wh{\tau}_1 +r}  \n \-   X^{t,x}_{\wh{\tau}_1 } \big|^p \bigg]   \nonumber  \\
   & &  \ls    \e_o     \+  \fC(C_p)^{\frac{1}{p}} \big(1\+|x| \big)
 \big( \|c(\cd)\| \l^{\frac{1}{2}} \+ \|c(\cd)\|^{\frac{1}{2}}  \l^{\frac{1}{4}} \big)
  \+    \fC  C_p \big(1\+|x|^p \big) \big( \|c(\cd)\|^p \l^{\frac{p}{2}} \+ \|c(\cd)\|^{\frac{p}{2}}  \l^{\frac{p}{4}} \big) \ls (1\+\fC) \e_o \, . \q  \qq  \label{eq:c365}
\eea

  Combining \eqref{eq:c381}, \eqref{eq:c367} and \eqref{eq:c365}    yields that
  \bea \label{eq:c614}
   E \big[ \big|  R(t,x,  \wh{\tau}_1   )\-  R(t,x, \tau_1  ) \big| \big] \< (2   \+ 10 \fC) \e_o  ,
\eea
 which together with   \eqref{eq:c467} and \eqref{eq:c351} show   that
 \beas
E \big[ \big|  R(t,\fx,  \wh{\tau}_1   )\-  R(t,x, \tau_1  ) \big| \big]
\ls E \big[ \big|  R(t,\fx,  \wh{\tau}_1   )\-  R(t,x, \wh{\tau}_1  ) \big| \big]
\+ E \big[ \big|  R(t,x,  \wh{\tau}_1   )\-  R(t,x, \tau_1  ) \big| \big]
\< (4   \+ 10 \fC) \e_o \=  \e \- \e_o  .
\eeas
  Then it follows from \eqref{eq:c405} and \eqref{eq:c387} that
  for any $(\fx, \fy) \ins \ol{O}_\d (x) \ti [(y\-\d)^+,\infty)$,
\bea \label{eq:c395}
V(t,\fx,\fy) \gs V\big( t,\fx,(y\-\d)^+  \big) \gs   E \big[ R(t,\fx, {\wh{\tau}_1} ) \big]
 \> E \big[  R(t,x, \tau_1  )   \big] \- \e \+ \e_o    \gs    V(t,x,y) \- \e    .
\eea

\no {\bf 1b)} {\it To show
$V(t,\fx,\fy) \ls   V(t,x,y) \+ \e $,
$ \fa (\fx,\fy) \ins \ol{O}_\d (x) \ti [0,y\+\d]  $,
 we let $\fx \ins \ol{O}_\d (x)$.  }

 There exists   $\tau_2 \= \tau_2 (t,\fx,y, \e) \ins \cT_{t,\fx} (y\+\d)$ such that
 \bea \label{eq:c399}
 E \big[ R(t,\fx,   \tau_2  ) \big]  \gs  V\big(t,\fx,y\+\d \big)  \- \e_o   .
 \eea
 \if{0}
  Similar to Part (1b),    we can   construct  a  stopping time  $  \wh{\tau}_2 \= \wh{\tau}_2 (t,x,\fx, y, \e)
    \ins \cT_{t,x} (y) $ such that
 $ \wh{\tau}_2 \ls   \tau_2  $ and $  P \big( A^c_R \cap  \{ \tau_2 \> {\wh{\tau}_2} \+ \sqrt{\l} \, \}\big) \< \l_o / 2 $.
  Then   like Part (1d), we can deduce that
 $ E \big[ \big|  R(t,x,  \wh{\tau}_2  )\-  R(t,\fx, \tau_2 ) \big| \big]
\< 
 \e \- \e_o  $.
 \fi

 We claim that we can also construct  a  stopping time  $  \wh{\tau}_2 \= \wh{\tau}_2 (t,x,\fx, y, \e)
  \ins \cT_{t,x} (y) $ satisfying
 \bea \label{eq:c339b}
 \wh{\tau}_2 \ls   \tau_2   \q \hb{and} \q
  P \big( A^c_R \Cp  \{ \tau_2 \> {\wh{\tau}_2} \+ \sqrt{\l} \, \}\big) \< \l_o / 2 .
 \eea
  If $ E \big[ \int_0^{\tau_2} \n g (t\+r, X^{t,x}_r) dr \big] \ls y $ \big(i.e. $\tau_2 \ins \cT_{t,x} (y)$\big),
we directly set $\wh{\tau}_2 \df \tau_2$.
  Otherwise,   set   $ \fb \df E \big[ \int_0^{\tau_2} \n g (t\+r, X^{t,x}_r) dr \big] \- y \> 0  $.
 Similar to \eqref{eq:c335}, $  \wh{\tau}_2  \=  \wh{\tau}_2  (t,x,\fx, y, \e)
 \df \inf\big\{s \ins [0,\infty) \n : E \big[ \int_0^{\tau_2}  \n  g(t\+r, X^{t,x}_r) dr \big|\cF_s \big]
 \- \int_0^s  \n g(t\+r, X^{t,x}_r) dr  \ls  \fb     \big\} $  is   an  $\bF-$stopping time satisfying
$  E \big[ \int_0^{\tau_2} \n g(t\+r, X^{t,x }_r) dr \big| \cF_{\wh{\tau}_2} \big]
 \- \int_0^{\wh{\tau}_2} \n g(t\+r, X^{t,x }_r) dr  \=   \fb  $.
Taking expectation $E[\cd]$ yields that
 \bea \label{eq:c337b}
  E \Big[ \int_0^{\wh{\tau}_2} \n g(t\+r, X^{t,x }_r) dr \Big]
 \=   E \Big[ \int_0^{\tau_2} \n g(t\+r, X^{t,x }_r) dr \Big] \-  \fb
   \= y ,   \q \hb{so } \, {\wh{\tau}_2} \ins \cT_{t,x} (y) .
 \eea
 As $E \big[\int_0^{\tau_2} \n g (t\+r, X^{t,x}_r) dr \big| \cF_{\tau_2} \big]
 \- \int_0^{\tau_2} \n g (t\+r, X^{t,x}_r) dr \= 0 \<  \fb $,
 we also see that $   \wh{\tau}_2 \ls \tau_2 $.

     Since   $ E \big[ \int_0^{\tau_2} \n g (t\+r, X^{t,\fx}_r)    dr   \big] \ls y \+ \d \< y \+ \l $,
 we can deduce from \eqref{eq:c331} and  \eqref{eq:c337b} that
\beas
2\l & \tn \gs  & \tn  E \Big[ \int_0^\infty \n \big| g(t\+r, X^{t,\fx}_r) \- g(t\+r, X^{t,x}_r) \big| dr  \Big] \+\l
\gs E \Big[    \int_0^{\tau_2} \n \big( g(t\+r, X^{t,x}_r)     \-    g(t\+r, X^{t,\fx}_r) \big) dr \Big] \+ \l \\
 & \tn \> & \tn  E \Big[ \int_0^{\tau_2} \n g (t\+r, X^{t,x}_r) dr \Big] \- y \= \fb
 \=     E \Big[ \int_{\wh{\tau}_2}^{\tau_2} \n g(t\+r, X^{t,x}_r) dr  \Big]
  \gs    E \bigg[    \b1_{A^c_R \cap \{ \tau_2 >  {\wh{\tau}_2}  + \sqrt{\l} \}}
   \int_{\wh{\tau}_2}^{\tau_2} \n g(t\+r, X^{t,x}_r) dr \bigg] \\
  & \tn  \gs  & \tn  \k_{\n \overset{}{R}} \, \sqrt{\l} \,
  P \big( A^c_R \Cp  \{ \tau_2 \> {\wh{\tau}_2} \+ \sqrt{\l} \, \}\big) .
\eeas
 By \eqref{eq:c357},  $ P \big( A^c_R \Cp  \{ \tau_2 \> {\wh{\tau}_2} \+ \sqrt{\l} \, \}\big)
 \< 
 \frac{2 \sqrt{\l}  }{ \k_{\n \overset{}{R}}} \< \l_o / 2  $, proving the claim \eqref{eq:c339b}.

   An analogy to \eqref{eq:c381}$-$\eqref{eq:c365} yields that
   $ E \big[ \big|  R(t,x,  \wh{\tau}_2   )\-  R(t,x, \tau_2  ) \big| \big]
\< (2   \+ 10 \fC) \e_o $, so we see from \eqref{eq:c467} and \eqref{eq:c351}
 \beas
E \big[ \big|  R(t,x,  \wh{\tau}_2  )\-  R(t,\fx, \tau_2 ) \big| \big]
\ls E \big[ \big|  R(t,x,  \wh{\tau}_2  ) \-  R(t,x, \tau_2  ) \big| \big]
\+ E \big[ \big|  R(t,x, \tau_2  ) \-  R(t, \fx, \tau_2  ) \big| \big]
\< (4   \+ 10 \fC) \e_o \=  \e \- \e_o  .
\eeas
It then follows from \eqref{eq:c405} and  \eqref{eq:c399}    that
for any $(\fx, \fy) \ins \ol{O}_\d (x) \ti [0,y\+\d]$,
\beas
 V(t,\fx,\fy) \ls V\big(t,\fx,y\+\d \big)   \ls     E \big[ R(t,\fx,  \tau_2 ) \big]    \+ \e_o
 \<  E \big[ R(t,x,  \wh{\tau}_2 ) \big] \+ \e  \ls V(t,x,y) \+ \e ,
\eeas
  which together with \eqref{eq:c395} leads to that
 $\big|V(t,\fx,\fy) \-  V(t,x,y)\big| \ls  \e$, $ \fa (\fx,\fy) \ins \ol{O}_\d (x) \ti \big[(y\-\d)^+, y\+\d\big] $.

\no {\bf 2)} Next, let $\varpi \ins [1,\infty)$, we further assume that $b,\si$ additionally satisfy \eqref{eq:c601} and
 $ f,g$ additionally satisfy \eqref{eq:c611}.

  Fix $(t,x,\e) \ins [0,\infty) \ti  \hR^l \ti (0,1) $.
Given $\ft \ins [0,\infty)$ and $\z \ins \cT$,
  \eqref{eq:c343}, \eqref{eq:c321}, \eqref{eq:c611},   H\"older's inequality, \eqref{eq:c609},
 \eqref{eqn-d011} and  the first inequality in \eqref{eq:esti_X_1}   imply that
\bea
\q && \hspace{-1cm} E \big[ | R(\ft,x,  \z   )\-  R(t,x, \z  )  | \big] \nonumber  \\
&& \ls  E \bigg[\int_0^{\z} \n \Big( \big| f  (\ft\+r, X^{\ft,x}_r  ) \- f  (\ft\+r, X^{t,x}_r  ) \big| \+
\big| f  (\ft\+r, X^{t,x}_r  ) \- f  (t\+r, X^{t,x}_r  ) \big| \Big) dr
\+ \Big| \pi \big( \ft\+\z  , X^{\ft,x}_{\z } \big)
 \- \pi \big( t\+\z  , X^{t,x}_{\z} \big) \Big|  \bigg]  \nonumber  \\
&& \ls  E  \bigg[
 \big( ( X^{\ft,x} \n \- \n  X^{t,x}   )_* \+  ( X^{\ft,x}  \n \- \n  X^{t,x}   )^p_* \big) \Big( \int_0^\infty \n  c(\ft\+r) dr \+ \fC \Big) \bigg]  \+ \rho (|\ft \- t|)
 \+ \rho (|\ft \- t|) E  \Big[ \big(1\+ |X^{t,x}_*|^\varpi \big) \n \int_0^\infty \n c(\ft \ld t\+r)  dr \Big]  \nonumber  \\
&& \ls 2 \fC   C^{1/p}_{p,\varpi} (1\+|x|^\varpi) \rho (|\ft\-t|)
 \+ 2 \fC C_{p,\varpi} (1\+|x|^{p \varpi}) \big(\rho (|\ft\-t|)\big)^p  \+ \rho (|\ft \- t|)
 \+ \fC \rho (|\ft \- t|) \big( 1\+ C_\varpi  (1\+|x|^\varpi ) \big)  . \label{eq:c617}
\eea

  Let us still set $ \e_o $, $\fM$ and take  $ \l_o \= \l_o (t,x, \e)   $,
  $R \= R (t,x, \e)  $, $ \l \= \l (t,x, \e)    $ as in Part 1.
 We now choose  $ \d' \= \d'(t,x, \e) \ins (0,1)  $ such that
  \bea \label{eq:c351b}
  (C_p)^{\frac{1}{p}} \d'  \+  C_p (\d')^p    \+   C^{1/p}_{p,\varpi} \big(1\+|x|^\varpi\big)   \rho (\d')
   \+     C_{p,\varpi} \big(1\+|x|^{p \varpi}\big) \big( \rho (\d') \big)^p
  \+  \rho(\d') \+  \rho(\d') \big( 1\+ C_\varpi (1\+|x|^\varpi)\big)   \ls \frac{\l\ld \e_o}{\fC}   ,
  \eea
  and   fix $ y  \ins   [0,\infty) $.

 \no {\bf 2a)} {\it To show   that
$V(\ft,\fx,\fy) \gs   V(t,x,y) \- \e $, $ \fa (\ft,\fx,\fy) \ins
[(t\-\d')^+,t\+\d' ] \ti \ol{O}_{\d'} (x) \ti [(y\-\d')^+,\infty)   $, we let $(\ft,\fx) \ins [(t\-\d')^+,t\+\d' ] \ti \ol{O}_{\d'} (x)$. }

 The condition (g1), \eqref{eq:c611}, H\"older's inequality,  \eqref{eq:esti_X_1}, \eqref{eq:c609}, \eqref{eqn-d011}
   and \eqref{eq:c351b} show that
  \bea
   && \hspace{-0.7cm}
       E \Big[ \int_0^\infty \n \big| g(\ft\+r, X^{\ft,\fx}_r) \- g(t\+r, X^{t,x}_r) \big| dr  \Big] \nonumber \\
   & &   \ls  E \bigg[ \int_0^\infty \n \Big( \big| g(\ft\+r, X^{\ft,\fx}_r) \- g(\ft\+r, X^{\ft,x}_r) \big|
       \+ \big| g(\ft\+r, X^{\ft,x}_r) \- g(\ft\+r, X^{t,x}_r) \big|
       \+ \big|  g(\ft\+r, X^{t,x}_r) \- g(t\+r, X^{t,x}_r) \big| \Big) dr  \bigg] \nonumber \\
   & & \ls  E \bigg[ \big(   (X^{\ft,\fx}   \- \n X^{\ft,x})_*   \+  (X^{\ft,\fx}   \- \n X^{\ft,x})^p_*
  \+ (X^{\ft,x} \n \- \n X^{t,x})_*   \+  (X^{\ft,x} \n \- \n X^{t,x})^p_* \big)
  \n  \int_0^\infty \n c(\ft\+r) dr  \bigg] \nonumber \\
  && \q + \rho(|\ft\-t|) E\bigg[ \big( 1 \+ (X^{t,x}_*)^\varpi \big) \int_0^\infty c(\ft \ld t \+ r) dr  \bigg] \nonumber \\
  && \ls \fC (C_p)^{\frac{1}{p}} |x\-\fx| \+   \fC  C_p |x\-\fx|^p
  \+ \fC C^{1/p}_{p,\varpi} \big(1\+|x|^\varpi\big)   \rho (|\ft\-t|)
  \+   \fC  C_{p,\varpi} \big(1\+|x|^{p \varpi}\big) \big( \rho (|\ft\-t|) \big)^p
  \+ \fC \rho(|\ft\-t|) \big( 1\+ C_\varpi (1\+|x|^\varpi)\big)   \nonumber  \\
  && \ls \fC (C_p)^{\frac{1}{p}} \d' \+   \fC  C_p (\d')^p
  \+ \fC C^{1/p}_{p,\varpi} \big(1\+|x|^\varpi\big)   \rho (\d')
  \+   \fC  C_{p,\varpi} \big(1\+|x|^{p \varpi}\big) \big( \rho (\d') \big)^p
  \+ \fC \rho(\d') \big( 1\+ C_\varpi (1\+|x|^\varpi)\big)  \ls  \l. \qq  \label{eq:c331b}
  \eea

Let $\tau_3  \= \tau_3(t,x,y,\e)  \ins \cT_{t,x} (y)$ such that
\bea \label{eq:c387b}
E \big[  R(t,x, \tau_3  )   \big]     \gs    V(t,x,y) \- \e_o  .
\eea
 If $ E \big[ \int_0^{\tau_3} \n g (\ft\+r, X^{\ft,\fx}_r) dr \big] \ls   (y\-\d')^+ $,
  we directly set ${\wh{\tau}_3} \df \tau_3 $.
  Otherwise,  we define
  $  \wh{\tau}_3 \= \wh{\tau}_3 (t,\ft,x,\fx,y, \e)
 \df \inf\big\{s \ins [0,\infty) \n : E \big[\int_0^{\tau_3} \n g (\ft\+r, X^{\ft,\fx}_r) dr \big| \cF_s \big]
\- \int_0^s \n g (\ft\+r, X^{\ft,\fx}_r) dr  \ls  \fra'  \big\}$ with
$\fra' \df E \big[ \int_0^{\tau_3} \n g (\ft\+r, X^{\ft,\fx}_r) dr \big]  \- (y\-\d')^+ \> 0 $.
Similar to \eqref{eq:c333},  one can deduce from \eqref{eq:c331b} that
$  \wh{\tau}_3 $ is a   $ \cT_{\ft,\fx} \big( (y\-\d')^+ \big)-$stopping time  satisfying
 \beas 
  {\wh{\tau}_3} \ls \tau_3 \q \hb{and} \q P \big( A^c_{R}
  \Cp  \big\{ \tau_3 \> {\wh{\tau}_3} \+ \sqrt{\l} \, \big\}\big)  \< \l_o / 2 .
 \eeas

 \if{0}
  \no   (1a)
     We claim that there exists a  stopping time
  $  \wh{\tau}_3  \=  \wh{\tau}_3  (t,\ft,x,\fx, y,  \e) \ins  \cT_{\ft,\fx} \big( (y\-\d')^+ \big)$  satisfying
 \bea \label{eq:c333}
  {\wh{\tau}_3} \ls \tau_3 \q \hb{and} \q P \big( A^c_{R} \Cp  \{ \tau_3 \> {\wh{\tau}_3} \+ \sqrt{\l} \, \}\big)
   \< \l_o / 2 .
 \eea
 Set $\d'_y \df \d' \ld y$. 
 If $ E \big[ \int_0^{\tau_3} \n g (\ft\+r, X^{\ft,\fx}_r) dr \big] \ls y \- \d'_y \= (y\-\d')^+ $,
  we directly set ${\wh{\tau}_3} \df \tau_3 $.
  Otherwise,     set $\fra' \df E \big[ \int_0^{\tau_3} \n g (\ft\+r, X^{\ft,\fx}_r) dr \big]  \- y  \+ \d'_y \> 0 $
    (In this case, one must have $y \> 0$).
  Since both  $ \big\{ E \big[\int_0^{\tau_3} \n g (\ft\+r, X^{\ft,\fx}_r) dr \big| \cF_s \big] \big\}_{s \in [0,\infty)}$
and $ \big\{ \int_0^s \n g (\ft\+r, X^{\ft,\fx}_r) dr \big\}_{s \in [0,\infty)}$ are $\bF-$adapted continuous processes,
\bea \label{eq:c335}
\wh{\tau}_3 \= \wh{\tau}_3 (t,\ft,x,\fx,y, \e)
 \df \inf\bigg\{s \ins [0,\infty) \n : E \Big[\int_0^{\tau_3} \n g (\ft\+r, X^{\ft,\fx}_r) dr \Big| \cF_s \Big]
\- \int_0^s \n g (\ft\+r, X^{\ft,\fx}_r) dr  \ls  \fra'  \bigg\}
\eea
defines an  $\bF-$stopping time   which  satisfies
$ E \big[\int_0^{\tau_3} \n g (\ft\+r, X^{\ft,\fx}_r) dr \big| \cF_{\wh{\tau}_3} \big]
\- \int_0^{\wh{\tau}_3} \n g (\ft\+r, X^{\ft,\fx}_r) dr  \=   \fra' $.
  Taking expectation $E [\cd]$ yields that
 \bea \label{eq:c101}
 E \Big[ \int_0^{\wh{\tau}_3} \n g (\ft\+r, X^{\ft,\fx}_r) dr \Big]
 \= E \Big[ \int_0^{\tau_3} \n g (\ft\+r, X^{\ft,\fx}_r) dr   \Big]  \- \fra' \= y \- \d'_y \= (y\-\d')^+ , \q
 \hb{so } \, {\wh{\tau}_3} \ins \cT_{\ft,\fx} \big( (y\-\d')^+ \big) .
\eea
 As $E \big[\int_0^{\tau_3} \n g (\ft\+r, X^{\ft,\fx}_r) dr \big| \cF_{\tau_3} \big]
 \- \int_0^{\tau_3} \n g (\ft\+r, X^{\ft,\fx}_r) dr \= 0 \<  \fra' $,
 one also has  $   \wh{\tau}_3 \ls \tau_3 $.

    Since $ E \big[ \int_0^{\tau_3} \n     g (t\+r, X^{t,x}_r)    dr   \big] \ls y  $ and since
    $\l\gs \l \ld \e_o \> \fC (C_p)^{\frac{1}{p}} \d'   \gs \d' \gs \d'_y $ by \eqref{eq:c351b}, one has
 \beas
  \fra'  & \tn  \= & \tn  E \Big[ \int_0^{\tau_3} \n g (\ft\+r, X^{\ft,\fx}_r) dr \Big] \- y \+ \d'_y
  \< E \Big[\int_0^{\tau_3} \n \big( g (\ft\+r, X^{\ft,\fx}_r)   \-   g (t\+r, X^{t,x}_r) \big)  dr   \Big] \+ \l \\
   & \tn \ls  & \tn  E \Big[ \int_0^\infty \n \big| g(\ft\+r, X^{\ft,\fx}_r) \- g(t\+r, X^{t,x}_r)   \big| dr  \Big] \+ \l\ls    2 \l.
 \eeas
 Using \eqref{eq:c331} again, we can deduce from   \eqref{eq:c101} that
   \bea
    2\l& \tn \> & \tn \fra'
   \=    E \bigg[ \int_{\wh{\tau}_3}^{\tau_3} \n g(\ft\+r, X^{\ft,\fx}_r) dr  \bigg]
  \gs E \bigg[     \int_{\wh{\tau}_3}^{\tau_3} \n g(t\+r, X^{t, x}_r) dr \bigg]
  \- E \Big[ \int_0^\infty \n \big| g(\ft\+r, X^{\ft,\fx}_r) \- g(t\+r, X^{t,x}_r)   \big| dr  \Big] \nonumber  \\
   & \tn \gs  & \tn  E \bigg[    \b1_{A^c_{R} \cap \{ \tau_3 >  {\wh{\tau}_3}  + \sqrt{\l} \}}
   \int_{\wh{\tau}_3}^{\tau_3} \n g(t\+r, X^{t, x}_r) dr \bigg]   \- \l\gs  \k_{\n \overset{}{R}} \, \sqrt{\l} \, P \big( A^c_{R} \Cp  \{ \tau_3 \> {\wh{\tau}_3} \+ \sqrt{\l} \, \}\big)  \- \l.
  \label{eq:c341}
  \eea
  It follows from \eqref{eq:c357} that
 $ P \big( A^c_{R} \Cp  \{ \tau_3 \> {\wh{\tau}_3} \+ \sqrt{\l} \, \}\big)
 \< 
   \frac{ 3  \sqrt{\l}  }{ \k_{\n \overset{}{R}}} \ls  \l_o / 2  $, proving the claim \eqref{eq:c333}.
 \fi

 Using similar arguments to those that lead to \eqref{eq:c614}, one can deduce from
    \eqref{eq:c353}$-$\eqref{eq:c361} that 
   $ E \big[ \big|  R(t,x,  \wh{\tau}_3   )\-  R(t,x, \tau_3  ) \big| \big] \< (2   \+ 10 \fC) \e_o $.
 \if{0}
 Since  \eqref{eq:c333} shows that
\beas 
 P (\fA^c ) \= P \{ \tau_3 \> {\wh{\tau}_3} \+    \sqrt{\l}  \}
  \ls P(A_{R})  \+ P \big( A^c_{R} \cap  \big\{ \tau_3 \> {\wh{\tau}_3} \+ \sqrt{\l} \, \big\}\big) \< \l_o \< \e_o  ,
\eeas
   \eqref{eq:c363}, \eqref{eq:c373}, \eqref{eq:c353},    \eqref{eq:c357} imply that
   \bea
   E \bigg[ \Big| \int_0^{\wh{\tau}_3 } \n  f  (t\+r, X^{t,x}_r  )  dr \- \int_0^{\tau_3 } \n f  (t\+r,X^{t,x}_r) dr  \Big| \bigg]
    & \tn  \ls   & \tn    E \bigg[   \big(  2 \+ (X^{t,x}_*)^p \big)  \Big( \b1_{\fA^c} \n \int_0^\infty \n c(t\+r) dr
  \+  \b1_\fA \|c(\cd)\|  ( \tau_3  \-  \wh{\tau}_3 ) \Big)    \bigg] \qq  \nonumber \\
      & \tn   \<  & \tn   \fC \big( 2  P(\fA^c) \+   \e_o \big)   \+  \sqrt{\l} (2\+ \fM) \|c(\cd)\|
     \< (1\+3\fC) \e_o \, , \label{eq:c381}
  \eea
 and that
\bea
E \big[  \b1_{\fA^c}  \big| \pi \big( \wh{\tau}_3  , X^{t,x}_{\wh{\tau}_3 } \big)
 \- \pi \big( \tau_3  , X^{t,x}_{\tau_3 } \big) \big| \big] \ls
   2 \fC E \big[  \b1_{\fA^c}  \big(2 \+   (X^{t, x}_*   )^p \big)  \big]
  \< 2  \fC \big( 2  P(\fA^c) \+   \e_o \big)  \< 6 \fC   \e_o.  \label{eq:c367}
\eea

Also, we can deduce from    \eqref{eq:c321}, \eqref{eq:c357}, H\"older's inequality,
   \eqref{eq:esti_X_2}  and \eqref{eq:c361}   that
 \bea
  && \hspace{-1cm} E \big[  \b1_\fA   \big| \pi \big( \wh{\tau}_3  , X^{t,x}_{\wh{\tau}_3 } \big)
 \- \pi \big( \tau_3  , X^{t,x}_{\tau_3 } \big) \big| \big]
   \ls      E \big[  \b1_\fA  \rho \big( \tau_3 \- \wh{\tau}_3  \big) \big]
\+ \fC E \Big[ \b1_\fA \Big( \big| X^{t,x}_{\wh{\tau}_3 } \n \-    X^{t,x}_{\tau_3 } \big| \+
 \big| X^{t,x}_{\wh{\tau}_3 }  \n \-    X^{t,x}_{\tau_3 } \big|^p \Big)  \Big] \nonumber \\
  & &  \ls     \rho \big(\sqrt{\l}\big)
 \+  \fC  \bigg\{ E \bigg[ \b1_\fA  \underset{r \in (0,\sqrt{\l}\,]}{\sup} \big| X^{t, x}_{\wh{\tau}_3 +r} \n \-   X^{t,x}_{\wh{\tau}_3 } \big|^p \bigg] \bigg\}^{\frac{1}{p}}
 \+ \fC    E \bigg[ \b1_\fA  \underset{r \in (0,\sqrt{\l}\,]}{\sup} \big| X^{t, x}_{\wh{\tau}_3 +r}  \n \-   X^{t,x}_{\wh{\tau}_3 } \big|^p \bigg]   \nonumber  \\
   & &  \ls    \e_o     \+  \fC(C_p)^{\frac{1}{p}} \big(1\+|x| \big)
 \big( \|c(\cd)\| \l^{\frac{1}{2}} \+ \|c(\cd)\|^{\frac{1}{2}}  \l^{\frac{1}{4}} \big)
  \+    \fC  C_p \big(1\+|x|^p \big) \big( \|c(\cd)\|^p \l^{\frac{p}{2}} \+ \|c(\cd)\|^{\frac{p}{2}}  \l^{\frac{p}{4}} \big) \ls (1\+\fC) \e_o \, . \q  \qq  \label{eq:c365}
\eea
\fi
 Then applying \eqref{eq:c467} with $(t,x,x',\tau) \= \big(\ft,\fx,x,\wh{\tau}_3 \big)$ and
  applying \eqref{eq:c617} with $\z \= \wh{\tau}_3$,    we see from   \eqref{eq:c351b}  that
  \beas
E \big[ \big|  R(\ft,\fx,  \wh{\tau}_3   )\-  R(t,x, \tau_3  ) \big| \big]
 & \tn \ls & \tn  E \big[ \big|  R(\ft,\fx,  \wh{\tau}_3   )\-  R(\ft,x, \wh{\tau}_3  ) \big| \big]
\+ E \big[ \big|  R(\ft,x,  \wh{\tau}_3   )\-  R(t,x, \wh{\tau}_3  ) \big| \big]
\+ E \big[ \big|  R(t,x,  \wh{\tau}_3   )\-  R(t,x, \tau_3  ) \big| \big] \\
 & \tn \ls & \tn  2 \fC    (C_p)^{\frac{1}{p}} \d' \+ 2 \fC C_p (\d')^p  \+ 2 \fC   C^{1/p}_{p,\varpi} (1\+|x|^\varpi) \rho (\d')
 \+ 2 \fC C_{p,\varpi} (1\+|x|^{p \varpi}) \big(\rho (\d')\big)^p  \\
  & \tn & \dn + \rho (\d')
 \+ \fC \rho (\d') \big( 1\+ C_\varpi  (1\+|x|^\varpi ) \big) \+ (2   \+ 10 \fC) \e_o
   \<    (4   \+ 10 \fC) \e_o \=  \e \- \e_o  .
\eeas
  It follows from \eqref{eq:c405} and \eqref{eq:c387b} that
  for any $(\ft,\fx, \fy) \ins [(t\-\d')^+,t\+\d' ] \ti  \ol{O}_{\d'} (x) \ti [(y\-\d')^+,\infty)$,
\bea \label{eq:c395b}
V(\ft,\fx,\fy) \gs V\big( \ft,\fx,(y\-\d')^+  \big) \gs   E \big[ R(\ft,\fx, {\wh{\tau}_3} ) \big]
 \> E \big[  R(t,x, \tau_3  )   \big] \- \e \+ \e_o    \gs    V(t,x,y) \- \e    .
\eea

\no {\bf 2b)} {\it We next show that  $V(\ft,\fx,\fy) \ls   V(t,x,y) \+ \e $,
$ \fa (\ft,\fx,\fy) \ins [(t\-\d')^+,t\+\d' ] \ti \ol{O}_{\d'} (x) \ti [0,y\+\d']  $.   }

 Let   $(\ft,\fx) \ins [(t\-\d')^+,t\+\d' ] \ti \ol{O}_{\d'} (x)$.
 There exists   $\tau_4 \= \tau_4 (t,\ft,x,\fx,y, \e) \ins \cT_{\ft,\fx} (y\+\d')$ such that
 \bea \label{eq:c399b}
 E \big[ R(\ft,\fx,   \tau_4  ) \big]  \gs  V\big(\ft,\fx,y\+\d' \big)  \- \e_o   .
 \eea
 If $ E \big[ \int_0^{\tau_4} \n g (t\+r, X^{t,x}_r) dr \big] \ls y $, 
 we directly set $\wh{\tau}_4 \df \tau_4$.
  Otherwise, we define $  \wh{\tau}_4  \=  \wh{\tau}_4  (t,\ft,x,\fx, y, \e)
 \df \inf\big\{s \ins [0,\infty) \n : E \big[ \int_0^{\tau_4} \n g(t\+r, X^{t,x}_r) dr \big|\cF_s \big]
 \- \int_0^s g(t\+r, X^{t,x}_r) dr  \ls  \fb'     \big\} $ with
     $ \fb' \df E \big[ \int_0^{\tau_4} \n g (t\+r, X^{t,x}_r) dr \big] \- y \> 0  $.
 \if{0}
 Similar to \eqref{eq:c335}, $  \wh{\tau}_4  \=  \wh{\tau}_4  (t,x,\fx, y, \e)
 \df \inf\big\{s \in [0,\infty): E \big[ \int_0^{\tau_4} \n g(t\+r, X^{t,x}_r) dr \big|\cF_s \big]
 \- \int_0^s g(t\+r, X^{t,x}_r) dr  \ls  \fb'     \big\} $  is   a $\bF-$stopping time satisfying
$  E \big[ \int_0^{\tau_4} \n g(t\+r, X^{t,x }_r) dr \big| \cF_{\wh{\tau}_4} \big]
 \- \int_0^{\wh{\tau}_4} \n g(t\+r, X^{t,x }_r) dr  \=   \fb'  $.
Taking expectation $E[\cd]$ yields that
 \bea \label{eq:c337c}
  E \Big[ \int_0^{\wh{\tau}_4} \n g(t\+r, X^{t,x }_r) dr \Big]
 \=   E \Big[ \int_0^{\tau_4} \n g(t\+r, X^{t,x }_r) dr \Big] \-  \fb'
   \= y ,   \q \hb{so } \, {\wh{\tau}_4} \ins \cT_{t,x} (y) .
 \eea
 As $E \big[\int_0^{\tau_4} \n g (t\+r, X^{t,x}_r) dr \big| \cF_{\tau_4} \big]
 \- \int_0^{\tau_4} \n g (t\+r, X^{t,x}_r) dr \= 0 \<  \fb' $,
 we also see that $   \wh{\tau}_4 \ls \tau_4 $.

     Since   $ E \big[ \int_0^{\tau_4} \n g (\ft\+r, X^{\ft,\fx}_r)    dr   \big] \ls y \+ \d' \< y \+ \l$,
 we can deduce from \eqref{eq:c331b} and  \eqref{eq:c337c} that
\beas
2\l & \tn \gs  & \tn  E \Big[ \int_0^\infty \n \big| g(\ft\+r, X^{\ft,\fx}_r) \- g(t\+r, X^{t,x}_r)   \big| dr  \Big]
\+\l\gs E \Big[    \int_0^{\tau_4} \n \big(   g(t\+r, X^{t,x}_r) \- g(\ft\+r, X^{\ft,\fx}_r)  \big) dr \Big] \+ \l \\
& \tn \> & \tn  E \Big[ \int_0^{\tau_4} \n g (t\+r, X^{t,x}_r) dr \Big] \- y \= \fb'
\=  E \Big[ \int_{\wh{\tau}_4}^{\tau_4} \n g(t\+r, X^{t,x}_r) dr  \Big] \\
  & \tn  \gs  & \tn    E \bigg[    \b1_{A^c_{R} \cap \{ \tau_4 >  {\wh{\tau}_4}  + \sqrt{\l} \}}
   \int_{\wh{\tau}_4}^{\tau_4} \n g(t\+r, X^{t,x}_r) dr \bigg]
  \gs  \k_{\n \overset{}{R}} \, \sqrt{\l} \, P \big( A^c_{R} \Cp  \{ \tau_4 \> {\wh{\tau}_4} \+ \sqrt{\l} \, \}\big) .
\eeas
 By \eqref{eq:c357},  $ P \big( A^c_{R} \Cp  \{ \tau_4 \> {\wh{\tau}_4} \+ \sqrt{\l} \, \}\big)
 \< 
 \frac{2 \sqrt{\l}  }{ \k_{\n \overset{}{R}}} \< \l_o / 2  $, proving the claim \eqref{eq:c339c}.

 \fi
    Analogous to \eqref{eq:c339b}, we can deduce from \eqref{eq:c331b} that
  $  \wh{\tau}_4 $ is a $ \cT_{t,x} (y) -$stopping time satisfying
 \beas 
 \wh{\tau}_4 \ls   \tau_4   \q \hb{and} \q
  P \big( A^c_{R} \Cp  \{ \tau_4 \> {\wh{\tau}_4} \+ \sqrt{\l} \, \}\big) \< \l_o / 2 .
 \eeas

  Since  an analogy to \eqref{eq:c381}$-$\eqref{eq:c365} gives that
   $ E \big[ \big|  R(t,x,  \wh{\tau}_4   )\-  R(t,x, \tau_4  ) \big| \big]  \< (2   \+ 10 \fC) \e_o $,
  applying \eqref{eq:c617} with $\z \= \tau_4$ and
  applying \eqref{eq:c467} with $(t,x,x',\tau) \= \big(\ft,x,\fx, \tau_4 \big)$,
  we see from   \eqref{eq:c351b}  that
  \beas
E \big[ \big|  R(t,x,  \wh{\tau}_4  )\-  R(\ft,\fx, \tau_4 ) \big| \big]
 & \tn \ls & \tn E \big[ \big|  R(t,x,  \wh{\tau}_4   )\-  R(t,x, \tau_4  ) \big| \big]
 \+ E \big[ \big| R(t,x, \tau_4  ) \- R(\ft,x,  \tau_4   )   \big| \big]
 \+ E \big[ \big|  R(\ft,x,  \tau_4  ) \-  R(\ft,\fx, \tau_4  ) \big| \big]  \\
 & \tn \ls & \tn (2   \+ 10 \fC) \e_o
   \+ 2 \fC   C^{1/p}_{p,\varpi} (1\+|x|^\varpi) \rho (\d')
 \+ 2 \fC C_{p,\varpi} (1\+|x|^{p \varpi}) \big(\rho (\d')\big)^p \+ \rho (\d') \\
  & \tn & \dn
  + \fC \rho (\d') \big( 1\+ C_\varpi  (1\+|x|^\varpi ) \big)  \+   2 \fC    (C_p)^{\frac{1}{p}} \d' \+ 2 \fC C_p (\d')^p
   \<    (4   \+ 10 \fC) \e_o \=  \e \- \e_o  .
\eeas
It then follows from \eqref{eq:c405} and  \eqref{eq:c399b}    that
for any $(\ft, \fx, \fy) \ins [(t\-\d')^+,t\+\d' ] \ti  \ol{O}_{\d'} (x) \ti [0,y\+\d']$
\beas
 V(\ft,\fx,\fy) \ls V\big(\ft,\fx,y\+\d' \big)   \ls     E \big[ R(\ft,\fx,  \tau_4 ) \big]    \+ \e_o
 \<  E \big[ R(t,x,  \wh{\tau}_4 ) \big] \+ \e  \ls V(t,x,y) \+ \e ,
\eeas
  which together with \eqref{eq:c395b} yields   $\big|V(\ft,\fx,\fy) \-  V(t,x,y)\big| \ls  \e$,
 $ \fa (\ft,\fx,\fy) \ins [(t\-\d')^+,t\+\d' ] \ti  \ol{O}_{\d'} (x) \ti \big[(y\-\d')^+, y\+\d'\big] $. \qed

\subsection{Proofs of Section \ref{sec:shift_prob}}

  \no {\bf Proof of Lemma \ref{lem_element}:}
     Set  $   \L    \df \Big\{ A \sb \O^t \n :
 A \=  \underset{\o \in A}{\cup} \big(\o \otis \O^s \big)  \Big\} $. Clearly,   $\es, \O^t \ins  \L$.
 For any     $A \ins \L$,    we claim that
  \bea   \label{eq:f221}
    \o \otis \O^s   \sb A^c   \hb{ for any }  \o \ins A^c .
 \eea
  Assume not, there exist  an $\o \ins A^c$ and an $ \wt{\o} \ins \O^s$ such that $\o \otis  \wt{\o}  \ins  A $.
 Then $\big(\o \otis  \wt{\o}\big) \otis \O^s \sb  A$ and it follows that
   $\o \ins \o \otis \O^s \=\big(\o \otis  \wt{\o}\big) \otis \O^s \sb  A$. A contradiction appear.
   So \eqref{eq:f221} holds, which shows that $A^c \ins \L $.

   For any $\{A_n\}_{n \in \hN} \sb  \L$, one can deduce that
  $  \underset{n \in \hN}{\cup} A_n \=  \underset{n \in \hN}{\cup} \Big( \underset{\o \in A_n}{\cup} \big(\o \otis \O^s \big) \Big)  \=  \underset{\o \in \underset{n \in \hN}{\cup} A_n}{\cup} \big(\o \otis \O^s \big)   $,
  namely, $ \underset{n \in \hN}{\cup} A_n \ins  \L$.
  Given $r \ins [t,s ]$ and  $ \cE \ins \sB(\hR^d)$, if $ \o \ins     (W^t_r)^{-1}      ( \cE  )$,
    it holds for any $\wt{\o} \ins \O^s$ that   $ \big(\o \otimes_s \wt{\o}\big)(r) \= \o(r) \ins \cE$
   or $ \o \otimes_s \wt{\o} \ins
    (W^t_r)^{-1}      ( \cE  )$,  
   which implies that $ (W^t_r)^{-1}      ( \cE  ) \ins \L$.
  Hence,   $\L$ is   a sigma$-$field of $\O^t$
  containing all   generating sets of $\cF^t_s$.
   It   follows that $ \cF^t_s   \sb \L $, proving the lemma. \qed

 \no {\bf Proof of Lemma \ref{lem_concatenation}:}
 Let us   regard $\o \otis \cd$ as a mapping $\G$ from $\O^s$ to $\O^t$, i.e.,
 $\G(\wt{\o}) \df    \o \otis \wt{\o} $, $ \fa \wt{\o} \ins \O^s$.
 So $A^{s,\o} \= \G^{-1}(A)$ for any $A \sb \O^t$.

 \no {\bf 1)} Assume first that $r \ins [s,\infty)$.
     Given $t' \ins [t,r]  $ and $ \cE \ins \sB(\hR^d)$, we can deduce that
       \beas
    \q    \G^{-1} \big(   (W^t_{t'})^{-1}   (\cE) \big)
    \= \{\wt{\o} \ins \O^s: W^t_{t'} (\o \otis \wt{\o}) \ins \cE \}  \=
    \begin{cases}
    \O^s, \hspace{4.5cm}  \hb{if  $t' \ins  [t,s)$ and $\o(t') \ins  \cE$} ; \\
    \es,  \hspace{4.75cm}  \hb{if  $t' \ins  [t,s)$ and $\o(t') \n \notin \n \cE$} ; \\
    \big\{ \wt{\o} \ins \O^s  \n :
    \o(s) \+ \wt{\o} (t') \ins  \cE     \big\}
    \=     (W^s_{t'})^{-1}     ( \cE'    ) \ins    \cF^s_r     ,   \q \hb{if  $t' \ins  [s,r]$} ;
   \end{cases}
   \eeas
    where $\cE' \df \cE \- \o(s) \= \{ x \- \o(s) \n : x \ins \cE \} \ins \sB(\hR^d) $.
      So all   generating sets of $\cF^t_r$ belong to  $ \L_r  \df \big\{ A \sb \O^t \n : \G^{-1} (A) \ins \cF^s_r \big\} $, which is clearly
          a sigma$-$field of $\O^t$. It follows that  $ \cF^t_r    \sb \L_r  $,
        or     $A^{s,\o} \= \G^{-1} (A) \ins \cF^s_r$ for any  $ A \ins \cF^t_r$.

  On the other hand,  let  $ \wt{A} \ins \cF^s_r$.
  We know from  Lemma \ref{lem_shift_inverse} (1)   that
  $(\Pi^t_s)^{-1}  \big(\wt{A}\,\big) \ins \cF^t_r $.
  Since the continuity of paths in $\O^t$  shows that
  $    \o \otis \O^s
  \=  \big\{ \o' \ins \O^t  \n : \o'(t')  \=  \o(t'),  \fa t' \ins (t,s) \n \cap \n  \hQ  \big\}
         \=      \underset{t' \in  (t,s) \cap \hQ   }{\cap}   (W^t_{t'})^{-1}   \big( \{ \o(t') \} \big)
          \ins  \cF^t_s \sb \cF^t_r   $, one can deduce that
    $ 
      \o \otis \wt{A} \= (\Pi^t_s)^{-1}  \big(\wt{A}\,\big) \cap \big( \o \otis \O^s \big) \ins  \cF^t_r $.

 \no {\bf 2)} Next, we consider the case of   $r \= \infty $.
  Given $r'  \ins [s,\infty)$, since
$ \G^{-1}  (A)    \ins   \cF^s_{r'}  \sb \cF^s $ for any $A \ins \cF^t_{r'} $, we see that
$  \cF^t_{r'} \sb \L   \df \big\{ A \sb  \O^t \n :  \G^{-1}  (A) \ins  \cF^s  \big\} $,
which is clearly a sigma$-$field of $\O^t$.
It follows from Lemma \ref{lem_countable_generate1} (1)  that $ \cF^t \= \si \Big(   \underset{r' \in [t,\infty)}{\cup}    \cF^t_{r'}  \Big) \= \si \Big(   \underset{r' \in [s,\infty)}{\cup}    \cF^t_{r'}  \Big)   \sb \L  $.
So $ A^{s,\o} \= \G^{-1}  (A) \ins  \cF^s$ for any $A \ins \cF^t$.

  On the other hand, let $r' \ins [s,\infty)$.
  Since  $ \G \big(\wt{A}\,\big) \= \o \otis \wt{A}  \ins   \cF^t_{r'} \sb \cF^t $ for any $\wt{A} \ins \cF^s_{r'}$,
 one has  $  \cF^s_{r'} \sb \wt{\L}  \df \big\{ \wt{A} \sb  \O^s \n :  \G  \big(\wt{A}\,\big) \ins  \cF^t  \big\} $.
  Given $\wt{A} \ins \wt{\L}$,
  it is clear that $\G  \big(\wt{A}\,\big) \cup \G  \big(\wt{A}^c\big)$
  is a disjoint union of   $ \G  (  \O^s ) \=  \o \otis \O^s  \ins  \cF^t_s  \sb \cF^t $.
  It    follows that  $ \G  \big(\wt{A}^c\big) \=  ( \o \otis \O^s  )
 \backslash \G  \big(\wt{A}\,\big) \ins  \cF^t  $.
 Also,  it   holds for any $\{\wt{A}_n\}_{n \in \hN} \sb \wt{\L} $ that
 $ \G  \Big(\underset{n \in \hN}{\cup}\wt{A}_n \Big) \= \underset{n \in \hN}{\cup} \G  \big( \wt{A}_n \big)
 \ins  \cF^t  $.
 So  $ \wt{\L} $ is   a sigma$-$field of $ \O^s$ that contains all $\cF^s_{r'}$, $r' \ins [s,\infty) $.
 Then Lemma \ref{lem_countable_generate1} (1) implies  that  $   \cF^s
 \= \si \Big(   \underset{r' \in [s,\infty)}{\cup}    \cF^s_{r'}  \Big)     \sb \wt{\L} $,
        or     $  \o \otis \wt{A} \= \G  \big(\wt{A}\,\big) \ins  \cF^t  $
        for any  $ \wt{A} \ins   \cF^s  $.       \qed

  \no {\bf Proof of Proposition \ref{prop_shift_meas}: 1)}
  Let $\xi $ be an $\hE-$valued random variable on $\O^t$ that is $\cF^t_r  - $measurable
  for some $r \ins [s,\infty]$.
  For any $ \cE \ins \sB(\hE)$, since $\xi^{-1}(\cE) \ins \cF^t_r $,
   Lemma \ref{lem_concatenation} shows   that
  $  \big(\xi^{s,\o}\big)^{-1}(\cE)
    \= \big\{\wt{\o} \ins \O^s \n :    \xi ( \o \otis  \wt{\o} ) \ins \cE \big\}
    \= \big\{\wt{\o} \ins \O^s \n :      \o   \otis  \wt{\o} \ins \xi^{-1}(\cE) \big\}
    \=\big( \xi^{-1}(\cE) \big)^{s,\o} \ins \cF^s_r  $.
  So $\xi^{s,\o}$ is   $\cF^s_r  - $measurable.

\no {\bf 2)} Let $   \{X_r \}_{r \in  [t,\infty)}$ be an   $\hE-$valued, $\bF^t-$adapted   process.
 For any $r \ins [s,\infty)$  and $ \cE \ins \sB(\hE)$,
 since $ X_r \ins \cF^t_r$,   one can deduce from  Lemma \ref{lem_concatenation}  that
 $    \big(X^{s,\o}_r\big)^{-1}(\cE)  
     \=   \big\{\wt{\o} \ins \O^s \n :    X \big(r, \o \otis \wt{\o} \big) \ins \cE \big\}
       \=   \big\{\wt{\o} \ins \O^s \n :    \o \otis \wt{\o} \ins X^{-1}_r (\cE) \big\}
     \= \big(  X^{-1}_s (\cE) \big)^{s,\o} \ins \cF^s_r  $,
     which shows that    $ \big\{X^{s,\o}_r \big\}_{  r \in  [s,\infty)}$ is   $\bF^s-$adapted. \qed

   \no {\bf Proof of Proposition \ref{prop_shift_integr}:}
      In virtue of Theorem 1.3.4 and (1.3.15) of  \cite{Stroock_Varadhan},
 there exists a family $\{P^\o_s \}_{\o \in  \O^t }$ of probabilities on $(\O^t,\cF^t)$,
  called the  regular conditional probability distribution  of $P_t$ with respect to
  the sigma-field $\cF^t_s$, such that

 \no   (\,\,i)  For any $A \ins \cF^t$, the mapping $\o \to P^\o_s (A)$ is $\cF^t_s -$measurable;  
  \bea    \label{rcpd_1}
  \hb{(\,ii)    For any } \xi \ins  L^1  ( \cF^t  ), ~  E_{P^\o_s}  [\xi]
  \=   E_t  \big[   \xi  \big| \cF^t_s  \big] (\o)    \hb{ for  }   \pas ~ \o \ins \O^t ;\hspace{7  cm} \\
    \hb{(iii)   For any }  \o \ins \O^t  ,  ~    P^\o_s  \big( \o \otis \O^s \big) \=1 . \hspace{11.3 cm}  \label{rcpd_2}
      \eea

   \no {\bf 1)}  Given $\o \ins \O^t$,   Lemma \ref{lem_concatenation}
      shows that   $\o \otis \wt{A}    \ins \cF^t$ for any $\wt{A} \ins \cF^s$.
  Then one  can deduce from    \eqref{rcpd_2}  that
  $  P^{s,\o} \big(\wt{A}\,\big) \df  P^\o_s \big( \o \otis \wt{A} \, \big) $, $ \fa  \wt{A} \ins \cF^s $
  defines   a probability measure  on $ ( \O^s , \cF^s  )$.
  We claim that  for  \pas ~ $ \o  \ins  \O^t$
  \bea      \label{eq:f203}
    P^{s,\o} \big(\wt{A}\,\big) \= P_s \big(\wt{A}\,\big), \q  \fa  \wt{A} \ins \cF^s .
  \eea

   To see this, we let $  \wt{A} \ins  \cF^s $.
   Since $ (\Pi^t_s)^{-1}  \big(\wt{A}\,\big)  \ins  \cF^t  $ by Lemma \ref{lem_shift_inverse} (1),
        \eqref{rcpd_2} and   \eqref{rcpd_1} imply  that  for  \pas ~ $ \o  \ins  \O^t$
    \bea \label{eq:ax023}
  P^{s,\o} \big(\wt{A}\,\big) \=  P^\o_s \big(\o \otis \wt{A}\big)
  \=  P^\o_s \big( (\Pi^t_s)^{-1} \big(\wt{A}\,\big) \cap  ( \o \otis \O^s )    \big)
    \= P^\o_s  \big(   (\Pi^t_s)^{-1} \big(\wt{A}\,\big)  \big)
    \=   E_t  \Big[ \b1_{(\Pi^t_s)^{-1}  (\wt{A} )} \big|\cF^t_s  \Big] (\o) .
    \eea
  We can deduce from  Lemma \ref{lem_countable_generate1} (1)  that
    \beas
     (\Pi^t_s)^{-1}(\cF^s) & \tn \= & \tn  (\Pi^t_s)^{-1} \big( \si \big\{ (W^s_r)^{-1} (\cE) \n :
     r \ins [s,\infty),\; \cE \ins \sB(\hR^d) \big\} \big)
     \= \si \big\{ (\Pi^t_s)^{-1} \big( (W^s_r)^{-1} (\cE) \big)  \n : r \ins [s,\infty),\; \cE \ins \sB(\hR^d) \big\} \\
      & \tn \= & \tn  \si \big\{    (W^t_r \- W^t_s  )^{-1} (\cE)    \n :
      r \ins [s,\infty),\; \cE \ins \sB(\hR^d) \big\}
    \=  \si \big( W^t_r \- W^t_s ; r \ins [s,\infty) \big) ,
    \eeas
    which is independent of $\cF^t_s$ under $P_t $.
  Then \eqref{eq:ax023} and  Lemma \ref{lem_shift_inverse} (2)   show  that for  \pas ~ $ \o  \ins  \O^t$,
      \beas
      P^{s,\o} \big(\wt{A}\,\big)  \=   E_t \Big[ \b1_{(\Pi^t_s)^{-1}  (\wt{A})} \big|\cF^t_s \Big] (\o)
       \=   E_t \Big[\b1_{(\Pi^t_s)^{-1}  (\wt{A})} \Big]  \=  P_t \big( (\Pi^t_s)^{-1}   \big(\wt{A}\,\big) \big)
    \=   P_s \big(\wt{A}\,\big)   .
  \eeas

   As  $\sC^s_\infty \df    \Big\{  \underset{i =1}{\overset{m}{\cap}}   (W^s_{s_i})^{-1}   \big( O_{\d_i} (x_i) \big) \n  :
     \,  m \ins  \hN ,\,    s_i \ins \hQ \cp \{s\}
   \hb{ with } s \ls   s_1 \ls \cds \ls s_m   ,\, x_i \ins \hQ^d , \, \d_i \ins \hQ_+ \Big\} $ is a countable set,
   we can find   a   $ \cN \ins \sN^s   $ such that for any $\o \ins \cN^c$,
 $    P^{s,\o} \big(\wt{A}\,\big)  \= P_s \big(\wt{A}\,\big)  $ holds  for each $\wt{A} \ins \sC^s_\infty   $.
  To wit,
  $  \sC^s_\infty    \sb  \L \df \big\{ \wt{A} \ins \O^s \n :  P^{s,\o} \big(\wt{A}\,\big)  \= P_s \big(\wt{A}\,\big) , \,
    \fa \o \ins \cN^c  \big\} $.
      It is easy to see that  $\L$  is   a  Dynkin system. As $\sC^s_\infty $ is closed under intersection,
   Lemma \ref{lem_countable_generate1} (2) and Dynkin System Theorem    show  that
  $   \cF^s  \= \si (\sC^s_\infty)  \sb \L  $. Namely,    it holds  for any $\o \ins \cN^c$ that
  $   P^{s,\o} \big(\wt{A}\,\big)  \= P_s \big(\wt{A}\,\big) $, $ \fa \wt{A} \ins \cF^s  $, proving \eqref{eq:f203}.

  \no {\bf 2)} Now, let $ \xi \ins  L^1  ( \cF^t ) $.
       Proposition \ref{prop_shift_meas} (1) shows  that
   $\xi^{s,\o}$ is $\cF^s - $measurable for any $\o \ins \O^t$.      Also,
  we can deduce from  \eqref{rcpd_1}$-$\eqref{eq:f203}   that
     for     \pas ~ $ \o \ins \O^t $
        \beas
 E_s \big[ | \xi^{s,\o} |  \big]
  &\tn \=&\tn  \int_{\wt{\o} \in \O^s}  \big| \xi^{s,\o}(\wt{\o}) \big| \, d   P^{s,\o}(\wt{\o})
  \=  \int_{\wt{\o} \in \O^s} \big| \xi \big(\o \otis \wt{\o}\big) \big| \,
   d   P^\o_s  \big( \o \otis  \wt{\o}\big)
  \= \int_{\o' \in \o \otimes_s  \O^s}  \big| \xi (\o') \big| \,  d P^\o_s (\o') \\
  &\tn\=&\tn     \int_{\o' \in \O^t  }  \big| \xi (\o') \big| \,  d P^\o_s  (\o')
  \= E_{P^\o_s }\big[ | \xi | \big] \= E_t \big[ |\xi|  \big|\cF^t_s \big](\o)   \< \infty     ,
   \eeas
 thus $ \xi^{s,\o} \ins L^1 \big( \cF^s \big) $.
 Similarly, it holds  for     \pas ~ $ \o \ins \O^t $ that
 $E_s \big[   \xi^{s,\o}   \big]
 \= E_t \big[  \xi  \big|\cF^t_s \big](\o) \ins \hR $.  \qed

    \no {\bf Proof of Proposition \ref{prop_null_set}: 1)} Let $  \cN  $ be     a $P_t-$null set,
    so there exists an $A \ins \cF^t$ with $P_t(A)\=0$ such that $\cN \sb A$.
    For any   $\o \ins \O^t$,    Lemma \ref{lem_concatenation} shows  that
      $  \cN^{s,\o} \= \{\wt{\o} \ins \O^s \n : \o \otis \wt{\o} \ins \cN\}  \sb
      \{\wt{\o} \ins \O^s \n : \o \otis \wt{\o} \ins A \} \=      A^{s,\o} \ins \cF$,
      and we see that
      $ (\b1_A)^{s,\o} (\wt{\o}) \= \b1_{\{ \o \otimes_s \wt{\o} \in A   \}} \= \b1_{\{ \wt{\o} \in A^{s,\o} \} }
      \= \b1_{A^{s,\o}} (\wt{\o}) $, $\fa \wt{\o} \ins \O^s$.
      Then  \eqref{eq:f475} implies   that for       \pas ~ $\o \ins \O^t  $
     \bea \label{eq:ax025}
   P_s   \big( A^{s,\o} \big) \= E_s \big[\b1_{A^{s,\o}}\big]
   \= E_s \big[ (\b1_A)^{s,\o} \big]   \=        E_t  \big[ \b1_A \big| \cF^t_s \big] (\o)       \=      0 ,
        \q \hb{ and thus } \q  \cN^{s,\o} \ins \sN^s  .
        \eea

   Next, let $\xi_1$ and $ \xi_2$ be two real-valued random variables with  $\xi_1 \ls \xi_2$,    \pas ~ Since  $   \cN  \df \{ \o \ins \O^t: \xi_1 (\o) \> \xi_2 (\o) \} \ins \sN^t $,
             \eqref{eq:ax025} leads to that   for   \pas ~ $\o \ins \O^t  $,
  \beas
     0 \= P_s \big( \cN^{s,\o} \big) \= P_s \big\{ \wt{\o} \ins \O^s \n :
          \xi_1 \big(\o \otis  \wt{\o} \big) \> \xi_2 \big(\o \otis  \wt{\o} \big)  \big\}
               \=  P_s  \big\{ \wt{\o} \ins \O^s \n : \xi^{s,\o}_1 (\wt{\o})
               \> \xi^{s,\o}_2 (\wt{\o})    \big\}   .
               \eeas

   \no  {\bf 2)}  Let $\tau \ins \ol{\cT}^t $ with $\tau \gs s$   and let $ r \ins [s,\infty) $.
 As $ A_r \df \{\tau \ls r \} \ins \ol{\cF}^t_r $,
  there exists an $\wt{A}_r  \ins \cF^t_r$   such that $ \cN_r \df A_r \, \D \,  \wt{A}_r  \ins \sN^t   $
    (see e.g. Problem 2.7.3 of \cite{Kara_Shr_BMSC}).
  By Part (1), it holds for all $ \o \ins \O^t$ except on a
   $P_t-$null set $\wh{\cN}_r$ that $\cN^{s,\o}_r \ins \sN^s $.
   Given $ \o \ins  \wh{\cN}_r^c $,
   since $   A^{s,\o}_r \, \D \,  \wt{A}^{s,\o}_r \= \big( A_r \, \D \,  \wt{A}_r \big)^{s,\o}
   \= \cN^{s,\o}_r \in \sN^s $
   and since $ \wt{A}^{s,\o}_r \ins  \cF^s_r $ by Lemma \ref{lem_concatenation},
   we can   deduce that $ A^{s,\o}_r \in \ol{\cF}^s_r $ and it follows that
  \bea    \label{eq:d221}
 \{ \tau^{s,\o} \n  \le \n  r \} \n  = \n \{\wt{\o} \n \in  \n  \O^s \n : \tau^{s,\o} (\wt{\o}) \n \le \n r \}
  \n = \n   \{\wt{\o}  \n \in \n  \O^s \n : \tau (\o  \n \otimes_s \n  \wt{\o})  \n \le \n  r \}
  \n = \n   \{\wt{\o}  \n \in \n  \O^s \n :  \o  \n \otimes_s \n  \wt{\o}   \n \in \n  A_r  \}
   \n = \n  A^{s,\o}_r  \n \in \n  \ol{\cF}^s_r .
 \eea
  Let $\o \ins  \underset{r \in (s,\infty) \cap \hQ  }{\cap} \wh{\cN}^c_r $.
  For any $r  \n \in \n  [s,\infty)$, there exists a sequence $\{r_n\}_{n \in \hN}$
  in $(s,\infty)   \Cp \hQ  $ such that
  $\lmtd{n \to \infty} r_n  \n = \n  r$.
  Then \eqref{eq:d221} and the right-continuity of Brownian filtration $ \ol{\bF}^s$
  (under $P_s$) imply that
  $
 \{ \tau^{s,\o}   \ls r \} \= \underset{n \in \hN}{\cap} \{ \tau^{s,\o}   \ls r_n \} \ins \ol{\cF}^s_{r+} \= \ol{\cF}^s_r $.
  Hence $\tau^{s,\o} \ins  \ol{\cT}^s $. \qed

 \no {\bf Proof of Proposition \ref{prop_shift_FP}:  1)}
   Let $r \ins [s,\infty]$ and
  $\xi$ be an $ \ol{\cF}^t_r -$measurable random variable.
 By Lemma \ref{lem_F_version} (2), there exists
 an $\cF^t_r-$measurable random variable $\wt{\xi}$   that  equals to $   \xi$ except on a $\cN \ins \sN^t$.
   Proposition \ref{prop_shift_meas} (1) shows that
 $\wt{\xi}^{s,\o}  $ is   $\cF^s_r -$measurable for any $\o \ins \O^t$.
 Also,  we see from Proposition \ref{prop_null_set} (1)   that for  \pas ~ $\o \ins \O^t$,
  \bea \label{eq:d137}
   \big\{\wt{\o} \ins \O^s \n : \wt{\xi}^{s,\o} (\wt{\o}) \n \ne \n \xi^{s,\o} (\wt{\o}) \big\}
 \= \big\{\wt{\o} \ins \O^s \n :    \o \otis \wt{\o} \ins \cN \big\} \= \cN^{s,\o}  \ins \sN^s
 \eea
 and thus $ \xi^{s,\o} \ins \ol{\cF}^s_r  $. In particular,
 if $\xi$ is an $ \ol{\cF}^t_s -$measurable and $\wt{\xi}$ is $\cF^t_s-$measurable, then \eqref{eq:d137}
 and \eqref{eq:bb421}  imply that \pas ~ $\o \ins \O^t$,
 $\xi^{s,\o} \= \wt{\xi}^{s,\o} \= \wt{\xi} (\o) \=  \xi(\o)$, $P_s-$a.s.

 Suppose next that $\xi$ is integrable \big(so is $\wt{\xi}$\big).   Proposition \ref{prop_shift_integr} and Lemma \ref{lem_F_version} (1)
  show   that  for    \pas ~ $ \o \ins \O^t $, $\wt{\xi}^{s,\o}$ is integrable (so is $\xi^{s,\o}$) and
 $ E_t  \big[\xi\big| \ol{\cF}^t_s \big](\o) \= E_t  \big[\xi\big| \cF^t_s \big](\o)
  \=  E_t  \big[  \, \wt{\xi} \,\big| \cF^t_s \big](\o)
  \=  E_s  \big[ \,  \wt{\xi}^{s,\o} \big] \= E_s    \big[ \xi^{s,\o} \big] \ins \hR $.

 \no {\bf 2a)} Let $ X \= \{X_r\}_{r \in [t,\infty)}   $ be an $\ol{\bF}^t-$adapted process with \pas ~ continuous paths
and set $\cN_1 \df \big\{\o \ins \O^t \n : \, \hb{the path $X_\cd (\o)$ is not continuous}\big\}   \ins \sN^t$.
In light of Lemma \ref{lem_F_version} (3), we can find an $\hE-$valued, $\bF^t-$predictable   process
$\wt{X} \= \big\{ \wt{X}_r \big\}_{r \in [t,\infty)}$ such that
$ \cN_2 \df \{ \o \ins \O^t \n : \wt{X}_r (\o) \n \ne \n X_r (\o) \hb{ for some } r \ins [t,\infty) \} \ins \sN^t$.
 In particular, $\wt{X}$ is an $\bF^t-$adapted   process. 

  Proposition \ref{prop_shift_meas} (2) shows that
 the shifted process  $ \wt{X}^{s,\o} $ is   $\bF^s-$adapted  for any $\o \ins \O^t$,
 and Proposition \ref{prop_null_set} (1) implies that for any $ \o \ins \O^t $
 except on a $P_t-$null set $\cN_3$ that $ \big(\cN_1 \cp \cN_2\big)^{s,\o}   \ins \sN^s$.
 Let $\o \ins \cN^c_3$. Since
 \beas
   \big\{ \wt{\o} \ins \O^s \n : X^{s,\o}_\cd (\wt{\o}) \hb{ is not  continuous}  \big\}
 \cp \big\{ \wt{\o} \ins \O^s \n : \wt{X}^{s,\o}_r (\wt{\o}) \n \ne \n X^{s,\o}_r (\wt{\o})
 \hb{ for some } r \ins [s,\infty) \big\} 
 \sb  \big(\cN_1 \cp \cN_2\big)^{s,\o} \ins \sN^s ,
 \eeas
 one can deduce that  $ X^{s,\o} $ is an $\ol{\bF}^s-$adapted process with $P_s-$a.s.   continuous paths.

 \no {\bf 2b)} Next, let us further assume that $X \ins  \hC^q_t (\hE)$ for some $q \ins [1,\infty)$.
  Define $\xi \df \underset{r \in [t,\infty) \cap \hQ}{\sup} \,  \big| \wt{X}_r \big|^q  \ins  \cF^t  $.
 As $\xi$ equals to   $ X^q_*   $ on $   (\cN_1 \cp \cN_2)^c  $,
 one has   $ X^q_*  \ins \ol{\cF}^t $ and thus   $ E_t [ \xi ] \= E_t \big[ X^q_*  \big] \< \infty$.
  According to Part (1), it holds for all  $\o \ins \O^t$ except on a $P_t-$null set $\cN_4$ that
 $\xi^{s,\o}$ is $\ol{\cF}^s-$measurable and $P_s-$integrable.

 Let $\o \ins  (\cN_3 \cp \cN_4)^c  $.  For any $\wt{\o} \ins \big(  (\cN_1 \cp \cN_2)^{s,\o} \big)^c
 \= \big(  (\cN_1 \cp \cN_2)^c \big)^{s,\o}$,
   the continuity of the path $    X^{s,\o}_\cd (  \wt{\o}) \= X_\cd (\o \otis \wt{\o}) $
   implies that
$\underset{r \in [s,\infty)}{\sup} \big| X^{s,\o}_r (\wt{\o}) \big|^q
\= \underset{r \in   [s,\infty) \cap \hQ   }{\sup} \big| X_r (\o \otis \wt{\o}) \big|^q
\= \underset{r \in   [s,\infty) \cap \hQ   }{\sup} \big| \wt{X}_r (\o \otis \wt{\o}) \big|^q
  \ls  \xi (\o \otis \wt{\o}) $. 
 It follows   that     $   E_s \Big[\,\underset{r \in [s,\infty)}{\sup} |X^{s,\o}_r|^q \Big]
   \ls E_s \big[ \xi^{s,\o} \big]  \< \infty    $.
  Hence, $ X^{s,\o} \ins  \hC^q_s (\hE)  $ for any $\o \ins (\cN_3 \cp \cN_4)^c$.  \qed

  \no {\bf Proof of Proposition \ref{prop_shift_mart}:}
  Let $M  \= \{M_r\}_{r \in [t,\infty)} \ins \hM_t$.
  By Proposition \ref{prop_shift_FP} (3),
  it holds for \pas ~ $\o \ins \O^t$ that $ M^{s,\o} $ is an $\ol{\bF}^s-$adapted process with $P_s-$a.s.  continuous paths.
  So we only need to show that
  $M^{s,\o}$ is a  uniformly integrable martingale with respect to $\big(\ol{\bF}^s,P_s\big)$
  for \pas ~ $\o \ins \O^t$.

  By the uniform integrability of $M$,
  there exists $\xi \ins L^1 \big(\ol{\cF}^t  \big)$ such that for any $r \ins [s, \infty)$,
  \bea \label{eq:ax101}
    M_r \= E_t \big[\xi \big|\ol{\cF}^t_r \big] , \q  \pas
  \eea
   Set   $ \cN  \df \{\o \ins \O^t \n :
    \hb{the path $ M_\cd (\o) $ is not continuous} \}  \ins \sN^t $.
     Proposition \ref{prop_null_set}  (1) and Proposition \ref{prop_shift_FP} (2) imply that
      for all $\o \ins \O^t$ except on a   $\cN_o \ins \sN^t$, one has   $\cN^{s,\o}   \ins \sN^s$ and
     $ \xi^{s,\o}  \ins L^1 \big(\ol{\cF}^s   \big)$.

 Fix $r \ins [s,\infty)$. As $  M_r \ins L^1 \big(\ol{\cF}^t_r  \big)$,
  Proposition \ref{prop_shift_FP} (2) shows that for all   $ \o \ins \O^t$ except on a $P_t-$null set $\cN^1_r$,
     $M^{s,\o}_r \ins L^1 \big(\ol{\cF}^s_r  \big)$.

  Let $ \wt{A} \ins \ol{\cF}^s_r $.   By Lemma \ref{lem_shift_inverse1b} (2),
  the set  $\cA \df (\Pi^t_s)^{-1} \big(\wt{A} \, \big) $ belongs to $ \ol{\cF}^t _r $,
 so   $ \b1_{\cA} M_r \ins L^1 \big(\ol{\cF}^t_r  \big) $ and $ \b1_{\cA} \xi \ins L^1 \big(\ol{\cF}^t   \big)$.
 Since  it holds for any $\o \ins \O^t $ and $ \wt{\o} \ins \O^s$ that
 $  (\b1_\cA)^{s,\o} (\wt{\o})  \= \b1_{\{\o \otimes_s \wt{\o} \in \cA\}}
 \= \b1_{\{  \Pi^t_s  (\o \otimes_s \wt{\o}) \in \wt{A} \}}
 \= \b1_{\{    \wt{\o}  \in \wt{A} \}}    \= \b1_{\wt{A}} (\wt{\o}) $,
 Proposition \ref{prop_shift_FP} (2) and \eqref{eq:ax101} yield  that  for \pas ~ $\o \ins \O^t$
  \beas
  E_s \big[ \b1_{\wt{A}} M^{s,\o}_r \big] \=    E_t \big[\b1_{\cA} M_r | \ol{\cF}^t_s \big] (\o)
  \=  E_t \Big[\b1_{\cA} E_t \big[\xi| \ol{\cF}^t_r\big] \Big| \ol{\cF}^t_s \Big] (\o)
  \=  E_t \Big[ E_t \big[ \b1_{\cA} \xi| \ol{\cF}^t_r \big] \Big| \ol{\cF}^t_s \Big] (\o)
  \=  E_t \big[\b1_{\cA} \xi | \ol{\cF}^t_s \big] (\o)
  \=  E_s \big[ \b1_{\wt{A}}  \xi^{s,\o}   \big] .
  \eeas

   As  $\sC^s_r  \df    \Big\{  \underset{i =1}{\overset{m}{\cap}}   (W^s_{s_i})^{-1}   \big( O_{\d_i} (x_i) \big) \n  :
     \,  m \ins  \hN ,\,    s_i \ins \hQ_+ \cp \{s\}
   \hb{ with } s \ls   s_1 \ls \cds \ls s_m \ls r   ,\, x_i \ins \hQ^d , \, \d_i \ins \hQ_+ \Big\} $ is a countable set,
   there exists   a   $ \cN^2_r \ins \sN^s   $ such that for any $\o \ins (\cN^2_r)^c$,
 $   E_s \big[ \b1_{\wt{A}} M^{s,\o}_r \big]  \= E_s \big[ \b1_{\wt{A}}  \xi^{s,\o}   \big]  $ holds
  for each $\wt{A} \ins \sC^s_r   $.   To wit,
  $  \sC^s_r    \sb  \L_r \df \big\{ \wt{A} \sb \O^s   \n :  E_s \big[ \b1_{\wt{A}} M^{s,\o}_r \big]
   \= E_s \big[ \b1_{\wt{A}}  \xi^{s,\o}   \big] , \,     \fa \o \ins (\cN^2_r)^c  \big\} $.
      It is easy to see that $\sC^s_r $ is closed under intersection and   $\L$  is   a  Dynkin system.
  Then Lemma \ref{lem_countable_generate1} (2) and Dynkin System Theorem    show  that
  $   \cF^s_r  \= \si (\sC^s_r)  \sb \L_r  $.
  Clearly, $\sN^s $ also belongs to $\L_r$, so
  \bea \label{eq:ax103}
  \ol{\cF}^s_r \= \si (\cF^s_r \cp \sN^s )  \sb \L_r  .
  \eea

  Now, let $\o \ins \cN^c_o \cap \Big( \underset{r \in [s,\infty) \cap \hQ}{\cup} (\cN^1_r \cup \cN^2_r) \Big)^c $.
  For any $r \ins  [s,\infty)$,   \eqref{eq:ax103} shows that
   $   E_s \big[ \b1_{\wt{A}} M^{s,\o}_r \big]
   \= E_s \big[ \b1_{\wt{A}}  \xi^{s,\o}   \big] $, $ \fa   \wt{A} \ins \ol{\cF}^s_r  $ and thus
   $    E_s [ \xi^{s,\o}  | \ol{\cF}^s_r]  \= M^{s,\o}_r $, $P_s - $a.s.
   Since
   $ \{\wt{\o} \ins \O^s \n :   \hb{path $ M^{s,\o}_\cd (  \wt{\o}) $ is not continuous} \}  \sb
  \cN^{s,\o} \ins \sN^s $,  
  we can deduce from the continuity of process
  $ \big\{ E_s [ \xi^{s,\o}  | \ol{\cF}^s_r] \big\}_{r \in [s,\infty)} $
  that $P_s \big\{ M^{s,\o}_r \= E_s \big[ \xi^{s,\o}  | \ol{\cF}^s_r \big]   ,
   ~ \fa r \ins  [s,\infty) \big\} \= 1 $.
   Therefore, $M^{s,\o}$ is a  uniformly integrable
   continuous martingale with respect to $\big(\ol{\bF}^s,P_s\big)$.    \qed

  \no {\bf Proof of Lemma \ref{lem_shift_converge_proba}:}
 Let $ \{\xi_i\}_{i \in \hN} $ be a sequence of $ L^1 \big(\ol{\cF}^t\big)  $
  that converges to 0 in probability $P_t$, i.e.
   \bea  \label{eq:p011}
   \lmtd{i \to \infty} E_t  \big[\b1_{ \{ |\xi_i| > 1/n    \}} \big]
     =  \lmtd{i \to \infty} P_t \big( |\xi_i| > 1 / n   \big) = 0 , \q    \fa  n \in \hN .
    \eea
   In particular, $\lmtd{i \to \infty} E_t  \big[\b1_{ \{ |\xi_i| > 1  \}} \big]    =   0 $
   allows us to extract a subsequence $S_1 = \big\{ \xi^1_i \big\}_{i \in \hN}$
  from $ \{\xi_i\}_{i \in \hN}$ such that $ \lmt{i \to \infty} \b1_{\{|\xi^1_i| > 1\}}  = 0$, $P_t-$a.s.
  Clearly, $S_1$ also satisfies \eqref{eq:p011}. Then by $\lmtd{i \to \infty} E_t  \big[\b1_{ \{ |\xi^1_i| > 1/2  \}} \big]
     =   0$, we  can find    a subsequence $S_2 = \big\{ \xi^2_i \big\}_{i \in \hN}$ of $S_1$
    such that $ \lmt{i \to \infty} \b1_{\{|\xi^2_i| > 1/2 \}}  = 0$, $P_t-$a.s.  Inductively, for each $n \in \hN$ we
    can  select a subsequence $S_{n+1} = \{\xi^{n+1}_i\}_{i \in \hN}$ of $ S_n = \{\xi^n_i\}_{i \in \hN}$
     such that $ \lmt{i \to \infty} \b1_{\big\{ |\xi^{n+1}_i| > \frac{1}{n+1} \big\}} = 0$, $P_t-$a.s.

 For any $ i \in \hN $, we set $\wt{\xi}_i \df \xi^i_i$, which belongs to $S_n$ for  $ n =1,\cds, i$. Given $n \in \hN$,
 since $\{\wt{\xi}_i\}^\infty_{i = n} \subset S_n$,  it holds $P_t-$a.s. that $ \lmt{i \to \infty} \b1_{\big\{|\wt{\xi}_i| > \frac{1}{n}\big\}} = 0$. Then a conditional-expectation version of
  the bound convergence theorem and  Proposition \ref{prop_shift_FP} (2) imply that
  for all $\o \in \O^t$ except on a $P_t-$null set $\cN_n$, $\wt{\xi}_i$ is $\ol{\cF}^s-$measurable and
  \bea  \label{eq:p015}
   0= \lmt{i \to \infty} E_t \Big[ \b1_{ \{|\wt{\xi}_i| > 1/n  \}} \big| \ol{\cF}^t_s\Big](\o)=
    \lmt{i \to \infty} E_{s} \Big[ \big( \b1_{ \{|\wt{\xi}_i| > 1/n  \}} \big)^{s,\o} \Big] .
  \eea

  Let $\o \in \Big(\underset{n \in \hN}{\cup} \cN_n\Big)^c$.
  For any $n \in \hN$, one can deduce that
  \beas
   \big( \b1_{ \{|\wt{\xi}_i| > 1/n  \}} \big)^{s,\o} (\wt{\o})
    =   \b1_{ \big\{|\wt{\xi}_i(\o \otimes_s \wt{\o}) | > 1/n  \big\} }
    =     \b1_{ \big\{ \big|\wt{\xi}^{\,s,\o}_{\, i} (\wt{\o}) \big| > 1/n  \big\}}
    =   \Big( \b1_{ \big\{|\wt{\xi}^{\,s,\o}_{\, i}| > 1/n  \big\}} \Big) (\wt{\o}) , \q \fa \wt{\o} \in \O^s ,
   \eeas
   which together with \eqref{eq:p015} leads to that
   $  \lmt{i \to \infty} P_s \Big(   |\wt{\xi}^{\,s,\o}_{\,i}| > 1 / n    \Big)
   = \lmt{i \to \infty} E_{s} \Big[ \big( \b1_{ \{|\wt{\xi}_i| > 1/n  \}} \big)^{s,\o} \Big] = 0 $.   \qed

  \no {\bf Proof of Proposition \ref{prop_FSDE_shift}:}
 As $\fX \ins \hC^2_t (\hR^l)$ by Corollary \ref{cor_X_Lp_estimate},
 we know from  Proposition \ref{prop_shift_FP} (3)   that  for $P_t-$a.s. $\o  \ins  \O^t$,
 $ \big\{   \fX^{s,\o}_r \big\}_{r \in [s, \infty)}   \ins \hC^2_s (\hR^l)$.

 To show that for $P_t-$a.s. $\o \ins \O^t$, $\fX^{s,\o}$ solves \eqref{FSDE2} over $[s,\infty)$
 with   initial state $\fX_s (\o)$,
 we let $\cN_1$ be the $P_t-$null set such that $\fX$ satisfies \eqref{FSDE2} on $\cN^c_1$.
 Define $M_{s'} \df \int_t^{s'} \n \b1_{\{r > s \}} \si(r,\fX_r)   dW^t_r$, $s' \ins [t , \infty)$.

 \no  {\bf 1)}  By Proposition \ref{prop_shift_FP} (1), there  exists a $P_t-$null set $\cN_2$ such that
 for any $\o \ins \cN^c_2$, $  \fX_s (\o \otimes_s \wt{\o}) \=   \fX_s (\o) $ holds for all $\wt{\o} \ins \O^s$
 except on a $\cN_\o \ins \sN^s$.

   Let $\o \ins  \cN_1^c \Cp \cN^c_2   $ and $ \wt{\o} \ins \cN^c_\o$.
   Implementing  \eqref{FSDE2} on the path $\o   \otis   \wt{\o}$ over period $ [s, \infty)$   yields that
    \bea
   \fX^{s,\o}_{s'} \n  ( \wt{\o})  & \tn \= &  \tn  \fX_{s'} (\o \oti_s \wt{\o})
     =   \fX_s (\o \oti_s \wt{\o}) \+ \int_s^{s'} \n  b \big(r, \fX_r (\o \oti_s \wt{\o})\big) dr
    \+ \Big( \int_s^{s'} \n \si  (r, \fX_r) dW^t_r \Big) (\o \oti_s \wt{\o})  \nonumber \\
    & \tn   = &  \tn
      \fX_s (\o) \+  \int_s^{s'}   b \big(r, \fX^{s,\o}_r (  \wt{\o})  \big) dr    +  M^{s,\o}_{s'} (  \wt{\o}),
      \q   s'  \ins     [s , \infty) . \label{eq:c625}
 \eea
   So it remains to show that   for $P_t-$a.s.~$\o \in \O^t$, it holds   $ P_s-$a.s. that
  \bea \label{eq:p515}
     M^{s,\o}_{s'}
         \=       \int_s^{s'} \si \big(r, \fX^{s,\o}_r    \big) \, dW^s_r   , \q s' \in [s,\infty) .
  \eea

\no  {\bf 2)} Since $\{M_{s'} \}_{s' \in [t , \infty)}$
  is a square-integrable martingale with respect to $ \big(  \ol{\bF}^t,P_t \big) $
  by \eqref{si_cond} and  Corollary \ref{cor_X_Lp_estimate},
  we know that (see e.g. Problem 3.2.27 of    \cite{Kara_Shr_BMSC})
    there is a sequence of $ \hR^{l \times d}-$valued, $ \ol{\bF}^t -$simple processes
   $  \Big\{\Phi^n_r \= \sum_{i \in \hN} \eta^n_i \, \b1_{ \big\{r \in (t^n_i, t^n_{i+1}] \big\} } ,
  \,  r \ins  [t,\infty) \Big\}_{n \in \hN}$ \big(where $\{t^n_i\}_{i \in \hN}  $
  is an increasing sequence   in $[t,\infty)$
   and $\eta^n_i  \ins    \ol{\cF}^t_{  t^n_i}$  for $i \ins \hN$\big) such that
      \beas
     P_t \n  - \n  \lmt{n \to \infty} \int_t^\infty     trace\Big\{ \big(  \Phi^n_r \- \si(r,\fX_r)  \big)
     \big(    \Phi^n_r   \- \si(r,\fX_r)  \big)^T  \Big\} dr \= 0  \q
    \hb{and}   \q
       P_t  \n - \n   \lmt{n \to \infty} \, \underset{s' \in [t,\infty)}{\sup} \big|  M^n_{s'}  - M_{s'}  \big| \= 0   ,
    \eeas
   where $  M^n_{s'} \df \int_t^{s'} \n \Phi^n_r dW^t_r = \sum_{i \in \hN}  \eta^n_i
    \big(  W^t_{s' \land t^n_{i+1}} \-  W^t_{s' \land t^n_i} \big)  $. Then it directly follows that
         \beas
     P_t \n  - \n  \lmt{n \to \infty} \int_s^\infty     trace\Big\{ \big(  \Phi^n_r \- \si(r,\fX_r)  \big)
     \big(    \Phi^n_r   \- \si(r,\fX_r)  \big)^T  \Big\} dr \= 0 ~
    \q  \hb{and}  \q
       P_t  \n - \n   \lmt{n \to \infty} \, \underset{s' \in [s,\infty)}{\sup} \big| M^n_{s'} - M_{s'} \big| =0   .
    \eeas

    By Lemma \ref{lem_shift_converge_proba},  $    \{ \Phi^n \}_{n \in \hN}$ has a subsequence  $  \Big\{ \wh{\Phi}^n_r \= \sum_{i \in \hN} \wh{\eta}^{\,n}_{\,i} \b1_{ \big\{r \in (\wh{t}^{\,n}_{\,i}, \wh{t}^{\,n}_{\, i+1}] \big\} } ,
  \,   r \ins [t,\infty) \Big\}_{n \in \hN}$ such that for any $\o \in \O^t$  except on a $P_t-$null set  $\cN_4 $
   \bea
  0 & \tn \= & \tn P_s \n  - \n  \lmt{n \to \infty} \bigg( \int_s^\infty     trace\Big\{ \big(  \wh{\Phi}^n_r \- \si(r,\fX_r)  \big)
     \big(    \wh{\Phi}^n_r   \- \si(r,\fX_r) \big)^T  \Big\} dr
      \bigg)^{s,\o}       \nonumber   \\
       & \tn \= & \tn  P_s \n  - \n  \lmt{n \to \infty}   \int_s^\infty
        trace\Big\{ \Big( \big( \wh{\Phi}^n \big)^{s,\o}_r \- \si(r,\fX^{s,\o}_r)   \Big)
     \Big( \big( \wh{\Phi}^n \big)^{s,\o}_r \- \si(r,\fX^{s,\o}_r)   \Big)^T  \Big\}  dr   \label{eq:p415a} \\
    \hb{and}\q 0 & \tn \= & \tn   P_s  \n - \n   \lmt{n \to \infty} \bigg( \underset{s' \in [s , \infty)}{\sup}
        \big|  \wh{M}^n_{s'} \- \wh{M}^n_s \- M_{s'}   \big| \,  \bigg)^{s,\o}  \nonumber \\
        & \tn \= & \tn     P_s  \n - \n   \lmt{n \to \infty}  \,  \underset{s' \in [s,\infty)}{\sup}
        \Big| \big( \wh{M}^n  \big)^{s, \o}_{s'} - \big( \wh{M}^n  \big)^{s, \o}_s
          - M^{s, \o}_{s'}   \Big|    ,  \qq  \label{eq:p415b}
    \eea
    where    $ \wh{M}^n_{s'} \df \int_t^{s'} \wh{\Phi}^n_r dW^t_r
    \=  \underset{i \in \hN}{\sum} \, \wh{\eta}^{\,n}_{\,i}  \Big(  W^t_{s' \land \wh{t}^{\,n}_{\, i+1}} \-  W^t_{s' \land \wh{t}^{\,n}_{\,i}} \Big)  $.

   Given $n \in \hN$, let $\ell_n$ be the largest integer such that $\wh{t}^{\,n}_{\,\ell_n} \< s$.
  For any $i \= \ell_n, \ell_n \+ 1, \cds$,   we set  $s^n_i \df \wh{t}^{\,n}_{\,i} \ve s$.
  Since $ \wh{\eta}^{\,n}_{\,i} \ins \ol{\cF}^t_{\wh{t}^{\,n}_{\,i}} \sb \ol{\cF}^t_{s^n_i } $.
   Proposition \ref{prop_shift_FP} (2) shows that
   $ \big( \wh{\eta}^{\,n}_{\,i} \big)^{s,\o} \ins \ol{\cF}^s_{s^n_i }$
   holds for any $\o \ins \O^t $ except on a $P_t-$null set $\cN^n_i$.
  Let   $\o \ins \wh{\O} \df \cN_4^c \Cp \Big(\underset{ n \in \hN}{\cap} \underset{i = \ell_n}{\overset{\infty}{\cap}} (\cN^n_i)^c \Big)$. As $ s^n_{\ell_n} \= s$, one has $\big( \wh{\eta}^{\,n}_{\,\ell_n} \big)^{s,\o} \ins \ol{\cF}^s_s $.
  For any $s' \ins [s,\infty)$ and $\wt{\o} \ins \O^s$,
  \beas
   \big( \wh{\Phi}^n \big)^{s,\o}_{s'} (\wt{\o}) \=   \wh{\Phi}^n_{s'} (\o \otis \wt{\o})
   \=  \sum_{i \in \hN} \wh{\eta}^{\,n}_{\,i} (\o \otis \wt{\o}) \, \b1_{ \big\{s' \in (\wh{t}^{\,n}_{\,i}, \wh{t}^{\,n}_{\, i+1}] \big\} } \= \big( \wh{\eta}^{\,n}_{\,\ell_n} \big)^{s,\o} (  \wt{\o}) \b1_{\{s' \in [s, s^n_{\ell_n+1}]\}} \+
    \sum^\infty_{i = \ell_n+1} \big( \wh{\eta}^{\,n}_{\,i} \big)^{s,\o} (  \wt{\o}) \, \b1_{  \{s' \in  ( s^n_i, s^n_{i+1} ]  \}}  .
   \eeas
   So $ \big\{ \big( \wh{\Phi}^n \big)^{s,\o}_{s'} \big\}_{s' \in [s,\infty)} $
   is an $ \hR^{l \times d}-$valued, $ \ol{\bF}^s -$simple process.
   Applying Proposition 3.2.26 of    \cite{Kara_Shr_BMSC} and using  \eqref{eq:p415a} yield that
    \bea   \label{eq:p417}
      0  =     P_s  \n - \n   \lmt{n \to \infty}  \,  \underset{s' \in [s,\infty)}{\sup}
        \Bigg|  \int_s^{s'} \big( \wh{\Phi}^n \big)^{s,\o}_r dW^s_r
         - \int_s^{s'}   \si(r,\fX^{s,\o}_r)   dW^s_r  \Bigg|    .
    \eea
 For any $\wt{\o} \in \O^s$, one can deduce that
   \beas
   \hspace{-3mm}
    \big( \wh{M}^n  \big)^{s, \o}_{s'} (\wt{\o})   \n  - \n  \big( \wh{M}^n  \big)^{s, \o}_s (\wt{\o})
   & \tn \dn =&  \tn  \dn   
      \sum^\infty_{i = \ell}   \wh{\eta}^{\,n}_{\,i} (\o  \otis  \wt{\o})
     \Big((\o  \otis  \wt{\o}) \big(s'  \dn \land \n  s^n_{i+1} \big)
    \n   - \n  (\o  \otis  \wt{\o}) \big(s'  \dn \land \n  s^n_i \big)  \Big)
    \=   \sum^\infty_{i = \ell}   \big( \wh{\eta}^{\,n}_{\,i} \big)^{s,\o} (  \wt{\o})
     \Big(  \wt{\o} \big(s'  \dn \land \n  s^n_{i+1}  \big)
      \n  -  \n   \wt{\o}  \big(s'  \dn \land \n  s^n_i  \big) \Big) \\
      & \tn   \dn  =  & \tn  \dn
        \sum^\infty_{i = \ell} \big( \wh{\eta}^{\,n}_{\,i} \big)^{s,\o} (  \wt{\o}) \Big(W^s_{s'    \land    s^n_{i+1} }
        \- W^s_{s'    \land    s^n_i } \Big) (  \wt{\o})
   \=  \bigg( \int_s^{s'} \big( \wh{\Phi}^n \big)^{s,\o}_r dW^s_r  \bigg) (\wt{\o})   , \q s' \in [s, \infty)  ,
   \eeas
   which together with \eqref{eq:p415b} and \eqref{eq:p417} shows that \eqref{eq:p515} holds  $P_s -$a.s.
   for any $\o \ins \wh{\O} $.
   Eventually,  we see from \eqref{eq:c625}   that
   $ P_s \big\{ \wt{\o} \ins \O^s \n :   \fX_r (\o \otis \wt{\o}) \= \fX^{s,\o}_r (\wt{\o})
   \=  \cX^{s, \fX_s (\o)}_r (\wt{\o}), \; \fa r \ins [s,\infty) \big\} \= 1 $
   for any $\o \ins \wh{\O} $. \qed

 \subsection{Proof of Section \ref{sec:DPP}}

The proof of the first DPP (Theorem \ref{thm_DPP}) is based on the following auxiliary result.

\begin{lemm} \label{lem_DPP}
 Given $(t,x,y) \ins [0,\infty) \ti \hR^l \ti [0,\infty) $,
 let $\tau \ins \cT^t_x (y)$ and   $\z \ins \ol{\cT}^t_\sharp$. Then
 \bea \label{eq:c627}
 E_t \big[\cR(t,x,\tau)\big] \ls  E_t \bigg[ \b1_{\{\tau \le \z   \}} \cR(t,x,\tau) \+ \b1_{\{\tau > \z  \}}
 \Big( \cV(\z,\cX^{t,x}_\z  , \cY^{t,x,\tau}_\z   )  \+ \int_t^\z \n f  (r,\cX^{t,x}_r   ) dr \Big) \bigg]
 \ls \cV(t,x,y)   .
 \eea

\end{lemm}

 \no {\bf Proof: 1)} {\it Let us start with some basic settings.}

 Denote $(\fX,\fY)  \df  (\cX^{t,x},\cY^{t,x,\tau})$
 and let $\z  $ take values in a countable
 subset $\{t_i\}_{i \in \hN}$ of $[t,\infty)$.
 In light of Lemma \ref{lem_F_version} (3), there exists an $\hR^l-$valued, $\bF^t-$predictable   process
$\wt{\fX} \= \big\{ \wt{\fX}_r \big\}_{r \in [t,\infty)}$ such that
$ \cN   \df \{ \o \ins \O^t \n : \wt{\fX}_r (\o) \n \ne \n \fX_r (\o) \hb{ for some } r \ins [t,\infty) \} \ins \sN^t$.

 Let $i \ins \hN$. By Proposition \ref{prop_null_set} (1), we can find a $P_t-$null set $ \cN_i$ such that
 for any $\o \ins  \cN^c_i$,   $  \cN^{t_i,\o}  $ is a $ P_{t_i}  -$null set.
   For any $r \ins [t,t_i]$,  since  $\wt{\fX}_r \ins \cF^t_r \sb \cF^t_{t_i}$, \eqref{eq:bb421} implies that
  \bea \label{eq:c427}
   \fX_r (\o \otii \wt{\o}) \= \wt{\fX}_r (\o \otii \wt{\o})
  \= \wt{\fX}_r (\o) \= \fX_r (\o) , \q \fa \o \ins   \cN^c \Cp  \cN_i^c , ~ \;
  \fa  \wt{\o} \ins   ( \cN^c )^{t_i,\o}  .
  \eea
 Also Proposition \ref{prop_FSDE_shift} shows that
 for all $\o \ins \O^t$ except on a $P_t-$null set $\wt{\cN}_i$,
  \bea
   \cN^i_\o \df   \Big\{ \wt{\o} \ins \O^{t_i} \n :  \fX^{t_i,\o}_r (\wt{\o})
 \n \ne \n  \cX^{t_i, \fX_{t_i} (\o)}_r (\wt{\o})  ,   \hb{ for some } r \ins [t_i,\infty) \Big\} \ins \sN^{t_i} .  \label{eq:c421}
 \eea

 Let $ \tau_i $ be a $\ol{\cT}^t -$stopping time  with $\tau_i \gs t_i$.
 According to  Proposition \ref{prop_null_set} (2)   and Proposition \ref{prop_shift_FP} (2),
  it holds for all $\o \ins \O^t$ except on a $P_t-$null set $ \wh{\cN}_i$ that
  $ \tau^i_\o \df \tau_i^{t_i,\o}  \ins \ol{\cT}^{t_i}$,
  \bea
  E_t \big[   \cR(t,x,\tau_i)   \big| \ol{\cF}^t_{t_i} \big] (\o)
 \= E_{t_i} \Big[  \big( \cR(t,x,\tau_i) \big)^{t_i,\o}   \Big]
   \; \hb{ and } \;   E_t \bigg[  \int_t^{\tau_i} \n g(r,\fX_r) dr  \Big| \ol{\cF}^t_{t_i} \bigg] (\o)
 \= E_{t_i} \bigg[ \Big(  \int_t^{\tau_i} \n g(r,\fX_r) dr    \Big)^{t_i,\o}   \bigg]  . \q \label{eq:c423}
 \eea

 Let $\o \ins \cN^c \Cp  \cN_i^c \Cp \wt{\cN}_i^c \Cp \wh{\cN}_i^c$.
 Given $\wt{\o} \ins \big(  \cN^{t_i,\o}   \cp  \cN^i_\o \big)^c
  \= ( \cN^c )^{t_i,\o}   \Cp \big(\cN^i_\o\big)^c $,
   \eqref{eq:c421} shows    $ \fX_r (\o \otii \wt{\o}) \=    \cX^{t_i,\fX_{t_i} (\o)}_r  (\wt{\o}) $
    for any $r \ins [t_i,\infty)$.
 In particular, taking $r \= \tau^i_\o (\wt{\o}) $ yields that
 $ 
 \fX  \big(\tau_i (\o \otii \wt{\o}) , \o \otii \wt{\o} \big)
 \= \fX  \big(\tau^i_\o (\wt{\o}) , \o \otii \wt{\o} \big)
 \=  \cX^{t_i,\fX_{t_i} (\o)}_{\tau^i_\o}  (\wt{\o})  $, which together with
   \eqref{eq:c427} leads to that
\beas
\big( \cR(t,x,\tau_i) \big)^{t_i,\o} (\wt{\o})
& \tn \=& \tn  \int_t^{\tau_i (\o \otimes_{t_i} \wt{\o})} \n
f \big(r,\fX_r (\o \otii \wt{\o})\big) dr
\+     \pi \big(\tau_i (\o \otii \wt{\o}),\fX_{\tau_i} (\o \otii \wt{\o})\big) \\
& \tn \=& \tn  \int_t^{t_i} \n f \big(r,\fX_r (\o)\big) dr \+
    \int_{t_i}^{\tau^i_\o (\wt{\o}) } \n f \Big(r,\cX^{t_i, \fX_{t_i} (\o)}_r (\wt{\o})\Big) dr
    \+\pi \Big(\tau^i_\o (\wt{\o}),\cX^{t_i,\fX_{t_i} (\o)}_{\tau^i_\o} (\wt{\o})\Big)     \\
  & \tn \=& \tn  \int_t^{t_i} \n f \big(r,\fX_r (\o)\big) dr \+
  \big( \cR( t_i , \fX_{t_i} (\o) , \tau^i_\o   ) \big) (\wt{\o}) ,
 \eeas
and similarly,
$\big(  \int_t^{\tau_i} \n g(r,\fX_r) dr    \big)^{t_i,\o} (\wt{\o})
\= \Big( \int_{t_i}^{\tau^i_\o  } \n g \Big( r , \cX^{t_i, \fX_{t_i} (\o)}_r \Big)  dr \Big) (\wt{\o})
\+ \int_t^{t_i} \n g \big(r,\fX_r (\o) \big) dr $.
 Taking expectation $E_t[\cd]$, we   see from \eqref{eq:c423}   that
 for $P_t-$a.s. $\o \ins \O^t$,  $ \tau^i_\o $ is a    $\ol{\cT}^{t_i}-$stopping time   satisfying
\bea
  E_t \big[   \cR(t,x,\tau_i)   \big| \ol{\cF}^t_{t_i} \big] (\o)
  & \tn \=   & \tn   E_{t_i} \big[ \cR( t_i , \fX_{t_i} (\o) , \tau^i_\o   ) \big]
  \+ \int_t^{t_i} \n f \big(r,\fX_r (\o)\big) dr ,   \label{eq:c473}  \\
   \hb{and} \q   E_t \bigg[   \int_t^{\tau_i} \n g(r,\fX_r) dr    \Big| \ol{\cF}^t_{t_i} \bigg] (\o)
  & \tn \= & \tn  E_{t_i} \bigg[ \int_{t_i}^{\tau^i_\o  } \n g \Big(r,\cX^{t_i, \fX_{t_i} (\o)}_r \Big)   dr \bigg]
  \+    \int_t^{t_i} \n g\big(r,\fX_r (\o)\big) dr . \hspace{1cm}  \label{eq:c475}
 \eea

 \no {\bf 2)} {\it We next show the first inequality in \eqref{eq:c627}.}

  Let   $i \ins \hN$ and set $\tau_i \df \tau \ve t_i \ins \ol{\cT}^t$.
  We can deduce from \eqref{eq:c475}, \eqref{eq:c471}, \eqref{eq:c485} and \eqref{eq:c473} that
  for \pas ~ $\o \ins \O^t$, $ \tau^i_\o \df \tau^{t_i,\o}_i $ is a    $\ol{\cT}^{t_i}-$stopping time  satisfying
  \beas
   E_{t_i} \bigg[ \int_{t_i}^{\tau^i_\o  } \n g \Big(r,\cX^{t_i, \fX_{t_i} (\o)}_r \Big)   dr \bigg]
   \=  E_t \bigg[   \int_t^{\tau \vee t_i} \n g(r,\fX_r) dr    \Big| \ol{\cF}^t_{t_i} \bigg] (\o)
  \-   \int_t^{t_i} \n g\big(r,\fX_r (\o)\big) dr
 \= \fY_{t_i} (\o) \ins [0,\infty)  ,
 \eeas
 and
 \bea \label{eq:c477}
  E_t \big[   \cR(t,x,\tau_i)   \big| \ol{\cF}^t_{t_i} \big] (\o)
    \=     E_{t_i} \big[ R  \big( t_i , \fX_{t_i} (\o) , \tau^i_\o   \big) \big]
  \+ \int_t^{t_i} \n f \big(r,\fX_r (\o)\big) dr
  \ls \cV\big( t_i , \fX_{t_i} (\o), \fY_{t_i} (\o) \big)   \+ \int_t^{t_i} \n f \big(r,\fX_r (\o)\big) dr .    \q
 \eea

 As $\{\tau \> \z\} \ins \ol{\cF}^t_{\tau \land \z} \sb \ol{\cF}^t_\z $ (see e.g. Lemma 1.2.16 of \cite{Kara_Shr_BMSC}),
 one has
 $ \{\tau \> \z  \= t_i\} \= \{\tau \> \z\} \Cp \{ \z  \= t_i\} \ins \ol{\cF}^t_{t_i} $.
 Then \eqref{eq:c477}   shows  that
 \bea
 && \hspace{-1cm} E_t \big[ \b1_{\{\tau > \z = t_i\}} \cR(t,x,\tau) \big]
    \=     E_t \big[ \b1_{\{\tau > \z = t_i\}} \cR(t,x,\tau_i) \big]
    \=   E_t \Big[ \b1_{\{\tau > \z = t_i\}} E_t  \big[  \cR(t,x,\tau_i)  \big| \ol{\cF}^t_{t_i} \big] \Big] \nonumber \\
  & & \le   E_t \bigg[ \b1_{\{\tau > \z = t_i\}} \Big( \cV\big(t_i,\fX_{t_i}  ,\fY_{t_i}  \big)  \+ \int_t^{t_i} \n f \big(r,\fX_r  \big) dr \Big) \bigg]
      \=    E_t \bigg[ \b1_{\{\tau > \z = t_i\}} \Big( \cV\big(\z,\fX_\z  ,\fY_\z  \big)
      \+ \int_t^\z \n f \big(r,\fX_r  \big) dr \Big) \bigg] . \q  \qq \label{eq:c489}
\eea
 Since   \eqref{eq:c433b}, \eqref{eq:c363}, \eqref{eq:c373} and the first inequality in  \eqref{eq:esti_X_1b} imply that
\bea
  E_t \bigg[ \big| \cV(\z,\fX_\z  ,\fY_\z   )\big|  \+ \int_t^\z \n \big| f  (r,\fX_r   ) \big| dr \bigg]
 & \tn \ls  & \tn    E_t \bigg[ 2 \fC \big( 2\+   C_p (1\+|\fX_\z|^p) \big)
\+ \int_t^\infty c(r) \big( 2 \+ |\fX_\z|^p \big) dr \bigg] \qq  \nonumber \\
 & \tn  \ls  & \tn    2 \fC(3\+ C_p) \+ \fC(1\+2C_p) E_t \big[ \fX^p_* \big]
\< \infty , \label{eq:c491}
\eea
 taking summation over $i \ins \hN$ in \eqref{eq:c489}, we can deduce from the first inequality in \eqref{eq:c483b}
 and the dominated convergence theorem that
 \bea
&& \hspace{-1cm} E_t \big[ \b1_{\{\tau > \z  \}} \cR(t,x,\tau) \big]
 \= E_t \Bigg[ \sum_{i \in \hN} \b1_{\{\tau > \z = t_i\}} \cR(t,x,\tau) \Bigg]
 \= \sum_{i \in \hN} E_t \big[  \b1_{\{\tau > \z = t_i\}} \cR(t,x,\tau) \big]  \nonumber \\
&& \ls \sum_{i \in \hN}  E_t \bigg[ \b1_{\{\tau > \z = t_i\}} \Big( \cV\big(\z,\fX_\z  ,\fY_\z  \big)
      \+ \int_t^\z \n f \big(r,\fX_r  \big) dr \Big) \bigg]
  \=   E_t \Bigg[ \sum_{i \in \hN} \b1_{\{\tau > \z = t_i\}} \Big( \cV\big(\z,\fX_\z  ,\fY_\z  \big)
      \+ \int_t^\z \n f \big(r,\fX_r  \big) dr \Big) \Bigg]   \nonumber \\
&&      \= E_t \bigg[ \b1_{\{\tau > \z  \}} \Big( \cV\big(\z,\fX_\z  ,\fY_\z  \big)
      \+ \int_t^\z \n f \big(r,\fX_r  \big) dr \Big) \bigg] .     \label{eq:c493}
 \eea
  It follows that
  $   E_t \big[\cR(t,x,\tau)\big] \ls
  E_t \Big[ \b1_{\{\tau  \le  \z   \}} \cR(t,x,\tau) \+ \b1_{\{\tau > \z  \}}
  \big( \cV(\z,\fX_\z  ,\fY_\z   )  \+ \int_t^\z \n f  (r,\fX_r   ) dr \big) \Big]  $.

\no {\bf 3)} {\it Now, we   demonstrate  the second inequality in \eqref{eq:c627}.}

 Fix $\e \ins (0,1)$ and   let $i \ins \hN$, $ \fx  \ins \hR^l   $.
 In light of \eqref{eq:c629} and Theorem \ref{thm_continuity} (1),
 there exists $\d_i(\fx) 
 \ins (0,\e/2)$ such that
   \bea \label{eq:c451}
   \fC  (C_p)^{\frac{1}{p}} \d_i (\fx)\+ \fC C_p \big( \d_i (\fx) \big)^p \<  \e/4    ,
   \eea
 and that for any $\fy  \ins   [0,\infty)$,
 \bea \label{eq:c453}
    \big| \cV(t_i,\fx', \fy' ) \-  \cV(t_i,\fx, \fy) \big|  \ls  \e/4 ,
    \q \fa ( \fx' , \fy' )   \ins \ol{O}_{\d_i(\fx)} (\fx) \ti
    \big[  (\fy \- \d_i(\fx))^+ , \fy \+  \d_i(\fx) \big]   .
\eea
  Then (g1), H\"older's inequality and the second inequality in \eqref{eq:esti_X_1b} imply that
\bea
   && \hspace{-2cm} E_{t_i} \bigg[ \int_{t_i}^\vs \n
    \big| g \big(r,\cX^{t_i,\fx}_r \big) \- g \big(r,\cX^{t_i,\fx'}_r \big)\big|    dr \bigg]
    \ls    E_{t_i} \bigg[ \int_{t_i}^\infty \n c(r) \big(\big|\cX^{t_i,\fx}_r  \n \- \n \cX^{t_i,\fx'}_r\big|
\+ \big|\cX^{t_i,\fx}_r \n \- \n  \cX^{t_i,\fx'}_r\big|^p \big) dr \bigg]
  \nonumber  \\
& & \ls   \Big(\int_0^\infty c(r) dr \Big)  E_{t_i} \Big[ \big(\cX^{t_i,\fx}  \n \- \n  \cX^{t_i,\fx'}\big)_* \+ \big(\cX^{t_i,\fx}  \n \- \n  \cX^{t_i,\fx'}\big)^p_* \Big]
\ls \fC  (C_p)^{\frac{1}{p}} |\fx\-\fx'|\+ \fC C_p |\fx\-\fx'|^p  \nonumber  \\
& & \ls    \fC  (C_p)^{\frac{1}{p}} \d_i (\fx)\+ \fC C_p \big( \d_i (\fx) \big)^p
 \<  \e / 4    , \q    \fa  \vs \ins \ol{\cT}^{t_i}  , ~
 \fa \fx' \ins \ol{O}_{\d_i ( \fx )} (\fx) \, . \label{eq:c441}
\eea

 We can find  a sequence  $\big\{(x^i_n,y^i_n)\big\}_{n \in \hN} $ in $ \hR^l \ti [0,\infty)$ such that
$ \hR^l \ti [0,\infty) \= \underset{n \in \hN}{\cup} \cO^i_n \ti \cD^i_n $ with
$\cO^i_n \df O_{\d_i (x^i_n)} (x^i_n)   $
and $\cD^i_n \df
\left\{ \ba{ll}
\big( \big( y^i_n  \-    \d_i(x^i_n)\big)^+  , \;  y^i_n \+   \d_i(x^i_n)  \big) , \q & \hb{if }
 y^i_n \> 0 , \ss \\
 \big[ 0  ,    \d_i(x^i_n)  \big) , & \hb{if }  y^i_n \= 0  .
\ea
\right.
$

   Let $  n \ins \hN $. We set
  $  A^i_n \df  \{\tau \> \z   \= t_i \} \Cp  \big\{ (\fX_{t_i}, \fY_{t_i}) \ins  \cO^i_n \ti \cD^i_n \big\}
  \Cp \cN^c_{t,x,\tau}   \ins \ol{\cF}^t_{t_i} $ and
  $ \cA^i_n  \df  A^i_n \Big\backslash \Big( \underset{n' < n}{\cup} A^i_{n'} \Big) \ins \ol{\cF}^t_{t_i}  $.
 There exists a $\tau^i_n \ins \cT^{t_i}_{ x^i_n} \big(  y^i_n   \big)$ such that
 \bea \label{eq:c455}
 E_{t_i} \big[\cR(t_i,x^i_n,\tau^i_n )\big] \gs \cV\big( t_i,x^i_n , y^i_n   \big)  \- \e/4 .
 \eea
  Lemma \ref{lem_shift_inverse1b} shows that $\tau^i_n (\Pi^t_{t_i})$
 is a $\ol{\cT}^t-$stopping time with values in $[t_i,\infty]$ such that
 $ \big( \tau^i_n (\Pi^t_{t_i}) \big)^{t_i,\o} (\wt{\o})
 \=  \tau^i_n \big(\Pi^t_{t_i} (\o \otimes_{t_i} \wt{\o})\big) \= \tau^i_n (\wt{\o})$
 for any $\o \ins \O^t$ and $  \wt{\o} \ins \O^{t_i}$.
 Also, by \eqref{eq:c473} and \eqref{eq:c475},
  it holds for any $\o \ins \O^t$ except on a $P_t-$null set $\cN^{i,n}$ that
  \bea
     E_t \Big[  R \big(t,x,\tau^i_n (\Pi^t_{t_i})\big) \big|\ol{\cF}^t_{t_i} \Big] (\o)
 & \tn  \=  & \tn     E_{t_i} \big[ \cR(t_i, \fX_{t_i}(\o), \tau^i_n ) \big]
       \+ \int_t^{t_i} \n f \big(r,\fX_r(\o)\big) dr  \label{eq:c479} ,  \\
   \hb{and} \q  E_t \bigg[   \int_t^{\tau^i_n (\Pi^t_{t_i})} g (r, \fX_r)dr \Big|\ol{\cF}^t_{t_i} \bigg] (\o)
 & \tn  \= & \tn  E_{t_i} \bigg[   \int_{t_i}^{\tau^i_n } g \Big(r, \cX^{t_i,\fX_{t_i}(\o)}_r \Big) dr     \bigg]
 \+ \int_t^{t_i} g\big(r,\fX_r(\o)\big) dr .  \hspace{1cm}   \label{eq:c481}
  \eea

 Clearly, the disjoint union  $\underset{i , n \in \hN}{\cup}  \cA^i_n $ satisfies that
 \bea
 \underset{i , n \in \hN}{\cup}  \cA^i_n
 & \tn \= & \tn \underset{i , n \in \hN}{\cup}   A^i_n
 \=   \Big( \underset{i \in \hN}{\cup} \, \{\tau \> \z   \= t_i \} \Cp
    \Big\{ (\fX_{t_i}, \fY_{t_i}) \ins \underset{n \in \hN}{\cup} \cO^i_n \ti \cD^i_n \Big\} \Big) \Cp \cN^c_{t,x,\tau}
    \nonumber \\
   & \tn  \= & \tn  \Big( \underset{i \in \hN}{\cup} \, \{\tau \> \z   \= t_i \} \Big) \Cp \cN^c_{t,x,\tau}
     \= \{\tau \> \z \} \Cp \cN^c_{t,x,\tau}  . \label{eq:c495}
 \eea
 We claim that
 \bhe
  \bea \label{eq:d067}
   \ol{\tau} \df  \b1_{\{\tau \le \z\}} \tau \+     \sum_{i,n \in \hN} \b1_{\cA^i_n} \tau^i_n (\Pi^t_{t_i})
   \+ \b1_{\{\tau > \z\} \cap \cN_{t,x,\tau}} t
   \hb{ belongs to }  \cT^t_x (y\+\e)  .
  \eea
 \ehe

 Let $i,n \ins \hN$ and $\o \ins \cA^i_n \Cp \big( \cN^{i,n} \big)^c $.
 As $\fX_{t_i} (\o) \ins \cO^i_n \= O_{\d_i (x^i_n)} (x^i_n) $,
 \eqref{eq:c451} and the second inequality in \eqref{eq:c483b}    imply
 \bea \label{eq:c487}
    E_{t_i} \Big[ \big| \cR(t_i, \fX_{t_i}(\o), \tau^i_n ) \- \cR(t_i, x^i_n, \tau^i_n ) \big| \Big]
   \ls 2 \fC \Big(   (C_p)^{\frac{1}{p}} \d_i (x^i_n)\+   C_p \big( \d_i (x^i_n) \big)^p \Big)
   \<     \e/2  .
 \eea
 Since $\big|\fX_{t_i} (\o) \- x^i_n \big| \ve \big| \fY_{t_i} (\o) \- y^i_n\big|   \< \d_i (x^i_n) $,
 applying \eqref{eq:c453} with $(\fx,\fy) \= (x^i_n,y^i_n) $ and $(\fx',\fy') \= \big( \fX_{t_i} (\o), \fY_{t_i} (\o) \big) $,
 we can deduce from \eqref{eq:c479}, \eqref{eq:c487} and \eqref{eq:c455}    that
 \beas
   E_t \Big[  R \big(t,x,\tau^i_n (\Pi^t_{t_i})\big) \big|\ol{\cF}^t_{t_i} \Big] (\o)
  & \tn  \>  & \tn   E_{t_i} \big[ \cR(t_i, x^i_n , \tau^i_n ) \big]
       \+ \int_t^{t_i} \n f \big(r,\fX_r(\o)\big) dr \-   \e/2
   \gs   \cV\big( t_i,x^i_n , y^i_n   \big) \+ \int_t^{t_i} \n f \big(r,\fX_r(\o)\big) dr   \- \frac{3 }{4} \e \nonumber \\
  & \tn   \gs   & \tn   \cV\big(t_i, \fX_{t_i} (\o), \fY_{t_i} (\o) \big)    \+ \int_t^{t_i} f(r,  \fX_r (\o) ) dr  \-  \e .
 \eeas
  Taking expectation $E_t [\cd]$   over $ \cA^i_n $ yields that
\beas
 E_t \big[ \b1_{\cA^i_n} \cR(t,x,\ol{\tau})   \big]
 & \tn \= & \tn  E_t \Big[ \b1_{\cA^i_n} R \big(t,x,\tau^i_n (\Pi^t_{t_i})\big)  \Big]
 \= E_t \bigg[ \b1_{\cA^i_n} E_t \Big[  R \big(t,x,\tau^i_n (\Pi^t_{t_i})\big) \big|\ol{\cF}^t_{t_i} \Big]  \bigg]  \\
 & \tn  \gs & \tn  E_t \Big[ \b1_{\cA^i_n} \Big(  \cV\big(\z, \fX_\z , \fY_\z  \big)
  \+ \int_t^\z f(r, \fX_r ) dr \- \e \Big) \bigg]  .
\eeas

Similar to \eqref{eq:c493}, taking summation up  over $i,n \ins \hN$, we can deduce from
\eqref{eq:c495}, \eqref{eq:c491}, the first inequality in \eqref{eq:c483b} and the dominated convergence theorem that
\beas
 E_t \big[ \b1_{\{\tau > \z\}} \cR(t,x,\ol{\tau})   \big]
 \gs  E_t \bigg[ \b1_{\{\tau > \z\}} \Big(  \cV\big(\z, \fX_\z , \fY_\z  \big)
  \+ \int_t^\z f(r, \fX_r ) dr \- \e \Big) \bigg]  .
\eeas
It thus follows that
$ \cV(t,x,y\+\e) \gs E_t \big[ \cR(t,x,\ol{\tau}) \big]    \gs  E_t  \Big[ \b1_{\{\tau \le \z \}}  \cR(t,x,\tau)
\+  \b1_{\{\tau > \z\}} \big(  \cV( \z, \fX_\z , \fY_\z     )
\+ \int_t^\z \n f(r,\fX_r) dr \big)   \Big]  \- \e $.
As $\e \n \to \n \infty$,  the second inequality in \eqref{eq:c627} follows from
 the continuity of $\cV$ in $y$ \big(i.e. \eqref{eq:c629} and Theorem \ref{thm_continuity} (1)\big). \qed

 \no {\bf Proof of Theorem \ref{thm_DPP}:} Fix $t \ins [0,\infty)$.

\no {\bf  1)} Let $(x,y) \ins \hR^l \ti [0,\infty)$ and
let  $\{\z (\tau) \}_{\tau \in \cT^t_x (y)} $ be a family of $  \ol{\cT}^t_\sharp -$stopping times.
 For any $ \tau \ins \cT^t_x (y) $, taking $\z \= \z(\tau)$ in \eqref{eq:c627} yields that
  \beas
 E_t \big[\cR(t,x,\tau)\big] \ls  E_t \bigg[ \b1_{\{\tau \le \z (\tau)  \}} \cR(t,x,\tau) \+ \b1_{\{\tau > \z (\tau) \}}
 \Big( \cV(\z (\tau),\cX^{t,x}_{\z (\tau)}  , \cY^{t,x,\tau}_{\z (\tau)}   )
  \+ \int_t^{\z (\tau)} \n f  (r,\cX^{t,x}_r   ) dr \Big) \bigg] \ls \cV(t,x,y)   .
 \eeas
Taking supremum over  $ \tau \ins \cT^t_x (y) $
\big(or taking supremum over  $ \wh{\tau} \ins \cT^t_x (y) $ if $y \> 0$\big),
we can deduce \eqref{eqn_DPP} from \eqref{eq:d187}.

\no {\bf  2)} Next, assume that $\cV(s,x,y)$ is continuous in $(s,x,y) \ins [t,\infty) \ti \hR^l \ti (0,\infty)$.

 We fix $(x,y) \ins \hR^l \ti [0,\infty)$ and
 a family  $\{\z (\tau) \}_{\tau \in \cT^t_x (y)} $ of $  \ol{\cT}^t -$stopping times.
 Let   $ \tau \ins \cT^t_x (y) $, $n \ins \hN$ and define
 \beas
 \z_n \= \z_n(\tau) \df \b1_{\{\z (\tau) = t \}}t \+ \sum_{i \in \hN} \b1_{\{\z (\tau) \in  (t+(i-1)2^{-n},t+i2^{-n}]\}} (t\+i2^{-n}) \ins \ol{\cT}^t .
 \eeas
 Applying \eqref{eq:c627} with $\z \= \z_n$ yields that
 \bea \label{eq:c631}
 E_t \big[\cR(t,x,\tau)\big] \ls   E_t \bigg[ \b1_{\{\tau \le \z_n   \}} \cR(t,x,\tau) \+ \b1_{\{\tau > \z_n  \}}
 \Big( \cV \big(\z_n,\cX^{t,x}_{\z_n}  , \cY^{t,x,\tau}_{\z_n}   \big)
  \+ \int_t^{\z_n} \n f  (r,\cX^{t,x}_r   ) dr \Big) \bigg]
  \ls   \cV(t,x,y )   .
 \eea
 An analogy to \eqref{eq:c491}  shows that
 \bea \label{eq:d155}
  \big| \cV \big(\z_n,\cX^{t,x}_{\z_n}  , \cY^{t,x,\tau}_{\z_n}   \big)\big|
   \+ \int_t^{\z_n} \n \big| f (r,\cX^{t,x}_r   ) \big| dr  \ls    2 \fC(3\+ C_p) \+ \fC(1\+2C_p)   \fX^p_*
   \ins L^1 \big(\ol{\cF}^t\big) .
   \eea

 We claim that $   \cY^{t,x,\tau}_{\z (\tau)} \> 0 $, \pas ~ on $\{\tau \> \z (\tau)  \}$. To see it, we
 set $A \df \{\tau \>\z(\tau)\} \Cp \big\{\cY^{t,x,\tau}_{\z (\tau)} \= 0 \big\} \ins \ol{\cF}^t_{\z (\tau)}$
 and can deduce that
 \beas
  0 & \tn \=  & \tn  E_t \Big[ \b1_A \cY^{t,x,\tau}_{\z (\tau)}  \Big]
  \=  E_t \bigg[ \b1_A \Big( E_t \Big[  \int_t^\tau g(r,\cX^{t,x}_r) dr \Big|\ol{\cF}^t_{\z (\tau)}\Big]
  \- \int_t^{\z (\tau)} g(r,\cX^{t,x}_r) dr \Big) \bigg] \\
   & \tn \=  & \tn  E_t \bigg[   E_t \Big[ \b1_A \int_{\z (\tau)}^\tau g(r,\cX^{t,x}_r) dr \Big| \ol{\cF}^t_{\z (\tau)} \Big]   \bigg]
  \=  E_t \Big[ \b1_A  \int_{\z (\tau)}^\tau g(r,\cX^{t,x}_r) dr \Big] ,
 \eeas
 which implies that $ \b1_A  \int_{\z (\tau)}^\tau g(r,\cX^{t,x}_r) dr \=0 $, ~ \pas ~
 It follows from the strict positivity of function $g$   that $P_t(A) \= 0$, proving the claim.
 As $\lmtd{n \to \hN} \z_n \= \z (\tau)$,
 one has $ \lmtd{n \to \hN} \b1_{\{\tau \le \z_n   \}}   \= \b1_{\{\tau \le \z (\tau)  \}}$.
 The continuity of function $\cV$  in $(s,\fx,\fy) \ins [t,\infty) \ti \hR^l \ti (0,\infty)$ and
 the continuity of processes $\big(\cX^{t,x},\cY^{t,x,\tau}\big)$
 then show that
 $ \lmt{n \to \hN} \b1_{\{\tau \le \z_n   \}} \cV \big(\z_n,\cX^{t,x}_{\z_n}  , \cY^{t,x,\tau}_{\z_n}   \big)
 \= \b1_{\{\tau \le \z (\tau)  \}} \cV \big(\z (\tau),\cX^{t,x}_{\z (\tau)}  , \cY^{t,x,\tau}_{\z (\tau)}   \big) $, \pas

 Letting $n \n \to \n \infty$ in \eqref{eq:c631}, we can deduce from
 \eqref{eq:d155},  the first inequality in \eqref{eq:c483b}
 and the dominated convergence theorem that
 \bea \label{eq:d135}
 E_t \big[\cR(t,x,\tau)\big] \ls   E_t \bigg[ \b1_{\{\tau \le \z (\tau)   \}} \cR(t,x,\tau) \+ \b1_{\{\tau > \z (\tau)  \}}
 \Big( \cV \big(\z (\tau),\cX^{t,x}_{\z (\tau)}   , \cY^{t,x,\tau}_{\z (\tau)}   \big)
  \+ \int_t^{\z (\tau)} \n f  (r,\cX^{t,x}_r   ) dr \Big) \bigg]
  \ls   \cV(t,x,y )   . \q
  \eea
 Taking supremum over  $ \tau \ins \cT^t_x (y) $
 \big(or taking supremum over  $ \wh{\tau} \ins \cT^t_x (y) $ if $y \> 0$\big),
 we   obtain \eqref{eqn_DPP} again from \eqref{eq:d187}.     \qed


\no {\bf Proof of Proposition \ref{prop_control_cont}:}
 Let us simply denote $\tau( t, x,  \a  )$ by $\tau_o$.
 For $n \ins \hN$, an analogy   to \eqref{eq:ax041} shows that
 \beas
 \ul{\tau}_n   \df \inf\big\{ s \ins [t,\infty) \n : Y^{t,x,\a}_s \= 1/n \big\}
 \q \hb{and} \q
 \ol{\tau}_n   \df \inf\big\{ s \ins [t,\infty) \n : Y^{t,x,\a}_s \= -1/n \big\}
 \eeas
 define two $\ol{\cT}^t-$stopping times.

 By definition, $\a \= M \- K$ for some $(M,K) \ins \hM_t \ti \hK_t$.
 It holds for all $\o \ins \O^t$ except on a   $P_t-$null set $\cN$   that
 $M_\cd (\o)$ is a continuous path, that $K_\cd (\o) $ is an  continuous increasing path
 and that $\ol{\tau}_n (\o) \< \infty$ for any $n \ins  \hN$.

 \no {\bf  1)} {\it We first show that }
\bea \label{eq:d035}
 \lmtu{n \to \infty} \ul{\tau}_n \= \lmtd{n \to \infty} \ol{\tau}_n \= \tau_o   \q  \pas
\eea

 Let $\o \ins \cN^c$   and set  $\ul{\tau} (\o) \df \lmtu{n \to \infty} \ul{\tau}_n(\o) \ls \tau_o (\o)$.
  The   continuity of   path $Y^{t,x,\a}_\cd (\o) $ implies that
 $ Y^{t,x,\a} \big(\ul{\tau}_n(\o),\o\big) \= 1/n $,   $\fa n \ins \hN $  and thus
 $ Y^{t,x,\a} \big(\ul{\tau}(\o),\o\big)
 \= \lmt{n \to \infty} Y^{t,x,\a} \big(\ul{\tau}_n(\o),\o\big) \= 0 $.
 It follows that $ \tau_o (\o) \= \ul{\tau} (\o) \= \lmtu{n \to \infty} \ul{\tau}_n(\o)$.

 On the other hand, we define a  $\ol{\cT}^t-$stopping time
 $\ol{\tau}    \df \lmtd{n \to \infty} \ol{\tau}_n   \gs \tau_o   $ and let $\o \ins \cN^c$.
 For any $  n \ins \hN $, as $ \ol{\tau}_n (\o) \< \infty $,
   the   continuity of   path $Y^{t,x,\a}_\cd (\o) $ again gives that
 $ Y^{t,x,\a} \big(\ol{\tau}_n(\o),\o\big) \= - 1/n $. Letting $n \to \infty$ yields that
 $ Y^{t,x,\a} \big(\ol{\tau}(\o),\o\big)
 \= \lmt{n \to \infty} Y^{t,x,\a} \big(\ol{\tau}_n(\o),\o\big) \= 0 $.

 Since $M$ is a uniformly integrable martingale, we know  from   the optional sampling theorem   that
\beas
 E_t \bigg[  K_{\ol{\tau}} \- K_{\tau_o} \+ \int^{\ol{\tau}}_{\tau_o}
 g(r, \cX^{t,x}_r ) dr   \bigg] \= E_t \Big[ M_{\ol{\tau}} \- M_{\tau_o}
 \- Y^{t,x,\a}_{\ol{\tau}}\+ Y^{t,x,\a}_{\tau_o}  \Big]  \= 0 ,
\eeas
which implies $  K_{\ol{\tau} }  - \n K_{\tau_o} \+ \int^{\ol{\tau}}_{\tau_o}
g(r, \cX^{t,x}_r ) dr    \= 0 $, \pas ~
Then one can deduce from the strict positivity of function $g$   that $\tau_o \= \ol{\tau}   \= \lmtd{n \to \infty} \ol{\tau}_n  $,  \pas, proving \eqref{eq:d035}.

\no {\bf 2)} Next, let $\e \ins (0,1)$  and  set $\e_o \df (4 \+  10 \fC  )^{-1} \e $.
  As $\fM \df E_t \big[   (\cX^{t,x}_*)^p \big] \< \infty$
 by the first inequality in \eqref{eq:esti_X_1b}, we can find   $\l_o \= \l_o (t,x, \e)
 \ins  \big(0,\e_o\big) $ such that
\bea \label{eq:d033}
E_t \big[\b1_A \, (\cX^{t,x}_*)^p \big]
\< \e_o \; \hb{ for any $A \ins \ol{\cF}^t $ with $P_t(A) \< \l_o $\,. }
\eea
There    exists   $R \=  R(t,x, \e) \ins (0,\infty)$ such that the set
 $A_R \df \big\{  \cX^{t,x}_*  \> R \big\} \ins \ol{\cF}^t$
 satisfies  $P_t(A_R) \< \l_o / 2 $.

Let   $ \l \= \l(t,x, \e) \ins   (0, 1  )  $  satisfy that
\bea
 && \hspace{3cm}   \l  \ls  \frac{\e_o}{  (2\+ \fM)  \|c(\cd)\|  }   \n \land \rho^{-1}(  \e_o ) \q \hb{and} \q \label{eq:d039} \\
 &&  (C_p)^{\frac{1}{p}} \big(1\+|x| \big)
 \big( \|c(\cd)\| \l  \+ \|c(\cd)\|^{\frac{1}{2}}  \l^{\frac{1}{2}} \big)
    \+    C_p \big(1\+|x|^p \big) \big( \|c(\cd)\|^p \l^p \+ \|c(\cd)\|^{\frac{p}{2}}  \l^{\frac{p}{2}} \big) \ls \e_o .  \qq  \qq   \label{eq:d041}
\eea
 We pick up $ \d \= \d(t,x, \e) \ins \Big(0, \frac{1}{\fC \fn}\big( \frac{\l_o}{2C_p} \big)^{\frac{1}{p}}\Big) $ such that
  \bea \label{eq:d043}
  \fC(C_p)^{\frac{1}{p}} \d \+   \fC  C_p \d^p   \ls \l \ld  \e_o  ,
  \eea
  Set $\O_n \df \{ \tau_o \- \l \ls \ul{\tau}_n \ls \ol{\tau}_n \ls \tau_o \+ \l \}
\ins \ol{\cF}^t$ for any $n \ins \hN $.
As  \eqref{eq:d035} implies that $P_t \big(\underset{n \in \hN}{\cup} \O_n \big) \=1$,
there exists  $\fn \ins \hN$ such that $P_t (\O_\fn) \> 1 \- \l_o / 2 $.

  Now, fix $x' \ins \ol{O}_\d (x)$ and simply denote $\tau(t,x',\a) $ by $\tau'$.
  We define $ A'   \df \big\{ (\cX^{t,x'} \- \cX^{t,x})_* \ls  (   \fC \fn  )^{-1} \big\} \ins \ol{\cF}^t$.
  The second inequality in \eqref{eq:esti_X_1b} shows that
 \beas
P_t \big( (A')^c \big)  \=   \fC^p \fn^p E_t \big[ (\cX^{t,x'} \- \cX^{t,x})^p_* \big]
\ls C_p  \fC^p \fn^p |x'\-x|^p \ls C_p  \fC^p \fn^p \d^p \< \l_o/2 .
\eeas
 So the  set $\cA \df A' \Cp \O_\fn \ins \ol{\cF}^t$ satisfies that
 $P_t (\cA^c) \= P_t \big( (A')^c \cp \O^c_\fn \big) \ls  P_t \big( (A')^c \big) \+ P_t \big( \O^c_\fn \big)
 \< \l_o \< \e_o $.

 Let   $\o \ins \cA  $.  Since it holds for any $s \ins [t,\infty)$ that
  \beas
  \big| Y^{t,x',\a}_s (\o) \- Y^{t,x,\a}_s (\o) \big|
  & \tn \ls & \tn  \int_t^s \n  \big|g \big(r,\cX^{t,x'}_r (\o)\big) \- g \big(r,\cX^{t,x}_r (\o)\big)\big| dr
  \ls  \int_t^s \n  c(r) \Big(  \big|\cX^{t,x'}_r  \- \cX^{t,x}_r   \big|
  \ve \big|\cX^{t,x'}_r   \- \cX^{t,x}_r   \big|^p \Big) (\o) dr \\
   & \tn \ls & \tn  \big(   \fC \fn \big)^{-1} \int_t^\infty \n c(r) dr \ls 1/\fn ,
  \eeas
  we see that
  \beas
  \q  Y^{t,x',\a}_s (\o) \gs Y^{t,x,\a}_s (\o) \- 1/\fn \> 0 , ~ \fa s \ins \big[ t ,  \ul{\tau}_n(\o) \big)
  \q \hb{and} \q
  Y^{t,x,\a}_s (\o) \gs Y^{t,x',\a}_s (\o) \- 1/\fn \> -   1/\fn , ~ \fa s \ins \big[ t ,  \tau' (\o) \big) .
  \eeas
  The former implies that $ \tau'(\o) \gs \ul{\tau}_n(\o) $ while
  the latter means that $\ol{\tau}_n(\o) \gs  \tau' (\o)$. In summary,
\bea \label{eq:d049}
\tau_o   \- \l \ls \ul{\tau}_{\fn}   \ls \tau'    \ls \ol{\tau}_{\fn}   \ls \tau_o   \+ \l \; \hb{ on } \; \cA .
\eea

 By an analogy to  \eqref{eq:c381} and \eqref{eq:c367},
  we can deduce from  \eqref{eq:c363}, \eqref{eq:c373},  \eqref{eq:d033}, \eqref{eq:d049} and \eqref{eq:d039}   that
   \bea
   E_t \bigg[ \int_{\tau_o \land \tau'}^{\tau_o \vee \tau'} \big|f(r,\cX^{t,x}_r)\big| dr \bigg]
    & \tn  \ls   & \tn    E_t \bigg[   \big(  2 \+ (\cX^{t,x}_*)^p \big)  \Big( \b1_{\cA^c} \n \int_t^\infty \n c( r) dr
  \+  \b1_\cA \|c(\cd)\|  | \tau'  \-  \tau_o | \Big)    \bigg] \qq  \nonumber \\
      & \tn   \<  & \tn   \fC \big( 2  P_t(\cA^c) \+   \e_o \big)   \+  \l (2\+ \fM) \|c(\cd)\|
     \< (1\+3\fC) \e_o \, , \label{eq:d045} \\
 \hb{and} \q  E_t \big[  \b1_{\cA^c}  \big| \pi \big( \tau'  , \cX^{t,x}_{\tau' } \big)
 \- \pi \big( \tau_o  , \cX^{t,x}_{\tau_o } \big) \big| \big]
  & \tn \ls  & \tn
   2 \fC E_t \big[  \b1_{\cA^c}  \big(2 \+   (\cX^{t, x}_*   )^p \big)  \big]
  \< 2  \fC \big( 2  P_t(\cA^c) \+   \e_o \big)  \< 6 \fC   \e_o \, . \hspace{1.5cm}  \label{eq:d047}
\eea
 And similar to   \eqref{eq:c365}, H\"older's inequality, \eqref{eq:c321}, \eqref{eq:d039},
 \eqref{eq:d049},     \eqref{eq:esti_X_2b} and  \eqref{eq:d041} imply  that
 \bea
  && \hspace{-1cm} E_t \big[  \b1_\cA   \big| \pi \big( \tau'  , \cX^{t,x}_{\tau' } \big)
 \- \pi \big( \tau_o  , \cX^{t,x}_{\tau_o } \big) \big| \big]
   \ls      E_t \big[  \b1_\cA  \rho \big( |\tau' \- \tau_o |    \big) \big]
\+ \fC E_t \Big[ \b1_\cA \Big( \big| \cX^{t,x}_{\tau' } \n \-    \cX^{t,x}_{\tau_o } \big| \+
 \big| \cX^{t,x}_{\tau' }  \n \-    \cX^{t,x}_{\tau_o } \big|^p \Big)  \Big] \nonumber \\
  & &  \ls     \rho  (\l)
 \+  \fC  \bigg\{ E_t \bigg[ \b1_\cA  \underset{r \in (0,\l ]}{\sup}
 \big| \cX^{t, x}_{\tau' \land \tau_o +r} \n \-   \cX^{t,x}_{\tau' \land \tau_o} \big|^p \bigg] \bigg\}^{\frac{1}{p}}
 \+ \fC    E_t \bigg[ \b1_\cA  \underset{r \in (0,\l ]}{\sup} \big| \cX^{t, x}_{\tau' \land \tau_o + r}  \n \-   \cX^{t,x}_{\tau' \land \tau_o} \big|^p \bigg]   \nonumber  \\
   & &  \ls    \e_o     \+  \fC(C_p)^{\frac{1}{p}} \big(1\+|x| \big)
 \big( \|c(\cd)\| \l  \+ \|c(\cd)\|^{\frac{1}{2}}  \l^{\frac{1}{2}} \big)
  \+    \fC  C_p \big(1\+|x|^p \big) \big( \|c(\cd)\|^p \l^p \+ \|c(\cd)\|^{\frac{p}{2}}  \l^{\frac{p}{2}} \big) \ls (1\+\fC) \e_o \, . \q  \qq  \label{eq:d051}
\eea
 Combining \eqref{eq:d045}, \eqref{eq:d047} and \eqref{eq:d051}   yields that
 \beas
 E_t \big[   \big|\cR  (t, x, \tau'  ) \- \cR  (t, x, \tau_o  ) \big| \big]
 \ls E_t \bigg[ \int_{\tau_o \land \tau'}^{\tau_o \vee \tau'} |f(r,\cX^{t,x}_r)| dr \+
 \Big| \pi \big(\tau_o, \cX^{t,x}_{\tau_o} \big) \- \pi \big(\tau', \cX^{t,x}_{\tau'} \big) \Big| \bigg] \< (2   \+ 10 \fC) \e_o \, ,
 \eeas
  which together with   \eqref{eq:c467} and \eqref{eq:d043} leads to  that
 \beas
 \;  \q  E_t \big[ \big|  \cR(t,x',  \tau'   )\-  \cR(t,x, \tau_o  ) \big| \big]
\ls E_t \big[ \big|  \cR(t,x',  \tau'   )\-  \cR(t,x, \tau'  ) \big| \big]
\+ E_t \big[ \big|  \cR(t,x,  \tau'   )\-  \cR(t,x, \tau_o  ) \big| \big]
\< (4   \+ 10 \fC) \e_o \=  \e  \, . \qq \hb{\qed}
\eeas

   \no  {\bf Proof of Proposition \ref{prop_surm}:}
   Let $(t,x,y ) \ins [0,\infty) \ti  \hR^l \ti  (0,\infty)   $.

 \no {\bf 1)}  Let $\a \ins \fA_t (y)$. Since $\tau (t,x,\a) \< \infty $, \pas ~ by \eqref{eq:ax041},
    the  continuity of process  $Y^{t,x,\a}$ implies   that
 \bea  \label{eq:ax043}
 \a_{\tau (t,x,\a)  }  \=  \int_t^{\tau (t,x,\a)  } \n g    (r, \cX^{t,x}_r  ) dr     , \q  \pas
 \eea
  One can then deduce from  the   uniform integrability of the   $\big(\ol{\bF}^t,P_t\big)-$supermartingale $\a$
  and the optional sampling theorem   that
  $ E_t \Big[ \int_t^{\tau (t,x,\a)} g    (r, \cX^{t,x}_r  ) dr \Big] \=
  E_t \big[ \a_{\tau (t,x,\a)}   \big] \ls E_t [\a_t]  \= y $,
  namely,  $ \tau (t,x,\a) \ins \cT^t_x(y) $.
  As $   Y^{t,x,\a}_t \= \a_t \= y \> 0$, \pas,
  we also derive from the  continuity of process  $Y^{t,x,\a}$   that $\tau (t,x,\a) \> t $, \pas ~
  Thus  $\a   \n \to \n \tau (t,x,\a)$   is   a   mapping
  from $  \fA_t (y) $ to    $ \wh{\cT}^t_x (y) $.

 \no {\bf 2)}    Next, let $\tau \ins \wh{\cT}^t_x (y) $
 and set $\d \df y \- E_t \big[   \int_t^\tau  g(r, \cX^{t,x}_r) dr   \big] \gs 0 $. Clearly,
 $  M_s \df \d \+ E_t \Big[  \int_t^\tau  g(r, \cX^{t,x}_r) dr  \Big|\ol{\cF}^t_s \Big] \gs 0  $, $  s \ins [t,\infty) $
 is a   uniformly integrable  continuous martingale with respect to $\big(\ol{\bF}^t,P_t\big)$, i.e., $M \ins \hM_t$.

 Define   $J_s \df \underset{s' \in [t,s]}{\inf} E_t \big[  \tau \- t   |\ol{\cF}^t_{s'}\big] $,  $s \ins [t,\infty)$
 and let $\cN$ be the $P-$null set such that for any  $\o \ins \cN^c$, the path $E_t \big[  \tau \- t   |\ol{\cF}^t_\cd\big] (\o)$ is continuous and $E_t \big[  \tau \- t   |\ol{\cF}^t_s \big] (\o)   \gs 0$, $\fa s \ins [t,s] \cap \hQ$.
 For any $ \o \ins \cN^c$, we can deduce that $ J_\cd(\o)$ is a nonnegative, continuous decreasing process.
 Given $s \ins [t,\infty)$, set $\xi_s \df \underset{s' \in [t,s] \cap \hQ}{\inf} E_t \big[  \tau \- t   |\ol{\cF}^t_{s'}\big]$, which is $\ol{\cF}^t_s-$measurable random variable.
 The continuity of process $ E_t \big[  \tau  \- t   |\ol{\cF}^t_s\big] $, $s \ins [t,\infty)$ shows that
  $J_s  \= \xi_s  $ on $ \cN^c$,
 so $J_s $ is also $\ol{\cF}^t_s-$measurable.   It follows that
  \bea \label{eq:d225}
    K_s \df \d \bigg[ 1 \ld  \Big(   \frac{s\-t}{J_s}   \Big)^+ \bigg]  \ins [0,\d] , \q s \ins [t, \infty)
    \eea
  is an $\ol{\bF}^t-$adapted continuous increasing process.
 Since $\tau \> t$, \pas, one has  $J_t \= E_t \big[  \tau \- t |\ol{\cF}^t_t \big] \= E_t  [ \tau \- t ] \>  0 $, \pas ~
 and thus $   K_t \= 0 $, \pas ~ To wit, $K \ins \hK_t$.

      Set  $ \a \df M \- K  $. It is clear that
      \beas
      \a_s \= M_s \- K_s \gs \d \+ E_t \bigg[  \int_t^\tau  g(r, \cX^{t,x}_r) dr  \Big|\ol{\cF}^t_s \bigg] \- \d
      \= E_t \bigg[  \int_t^\tau  g(r, \cX^{t,x}_r) dr  \Big|\ol{\cF}^t_s \bigg], \q \fa s \in [t,\infty) .
      \eeas
       As $ \a_t \= M_t \- K_t \= \d \+
      E_t \big[   \int_t^\tau  g(r, \cX^{t,x}_r) dr  \big] \+ 0  \= y $, \pas,
      we see that $\a \ins \fA_t (y)$.
  Since 
      $J_\tau \ls E_t \big[ \(\tau \- t) |\ol{\cF}^t_\tau\big]
      \=  \tau \- t  $, $P_t-$a.s.,
      one has   $  K_\tau \= \d $, \pas ~ and thus
    \bea \label{eq:ax045}
  \a_\tau   \= M_\tau - \d  \=  \int_t^\tau  g(r, \cX^{t,x}_r) dr , \q \pas
   \eea
 This shows   $ \tau(t,x,\a) \ls \tau $, \pas ~
 On the other hand,   subtracting   \eqref{eq:ax043}   from \eqref{eq:ax045} and applying
   the optional sampling theorem to $\a$ again yield   that
   $  0 \ls E_t \big[ \int_{\tau(t,x,\a)}^\tau   g    (r, \cX^{t,x}_r  ) dr \big]
   \= E_t \big[  \a_\tau \- \a_{\tau(t,x,\a)}  \big]
     \ls 0 $.
   The strict positivity of function $g$ then implies that $ \tau(t,x,\a) \= \tau  $, \pas \qed

 Similar to Lemma \ref{lem_DPP}, the following auxiliary result is crucial for
 proving the second DPP of $\cV$ (Theorem \ref{thm_DPP2}).

\begin{lemm} \label{lem_DPP2}
 Given $(t,x,y) \ins [0,\infty) \ti \hR^l \ti (0,\infty) $,
 let $\a \ins \fA_t (y)$ and let $\z \ins \ol{\cT}^t_\sharp$. Then
 \bea
  E_t \big[   \cR \big(t,x,\tau(t,x,\a)\big) \big]
  & \tn \ls  & \tn  E_t \bigg[ \b1_{\{\tau(t,x,\a) \le \z   \}} \cR \big(t,x,\tau(t,x,\a)\big)
 \+ \b1_{\{\tau(t,x,\a) > \z  \}}
 \Big( \cV \big(\z,\cX^{t,x}_\z  ,  Y^{t,x,\a}_\z   \big)  \+ \int_t^\z \n f  (r,\cX^{t,x}_r   ) dr \Big) \bigg] \nonumber \\
  & \tn  \ls  & \tn   \cV(t,x,y )    .   \label{eq:d021}
 \eea

\end{lemm}

 \no {\bf Proof:}
 Suppose that  $\a \= M\-K$ for some $(M,K) \ins \hM_t \ti \hK_t$.
 We denote   $(\fX,\fY,\wh{\tau})  \df  \big(\cX^{t,x}, Y^{t,x,\a},\tau(t,x,\a)\big)$
 and let $\z  $ take values in a countable subset $\{t_i\}_{i \in \hN}$ of $[t,\infty)$.

 \no {\bf 1)} {\it Let us  start with the first inequality in \eqref{eq:d021}.  }

 Since $\a$ is a uniformly integrable continuous supermartingales
  with respect to $\big(\ol{\bF}^t,P_t\big)$, one has $ \a_{\wh{\tau}} \= \int_t^{\wh{\tau}} \n  (r,\fX_r ) dr $
  and   the optional sampling theorem implies that
 \bea
 \fY_{\wh{\tau} \land \z}  & \tn  \=  & \tn    \a_{\wh{\tau} \land \z}
 \- \int_{t}^{\wh{\tau} \land \z} \n g (r,\fX_r) dr
 \gs E_t \big[ \a_{\wh{\tau}}  \big| \ol{\cF}^{t}_{\wh{\tau} \land \z} \big]
 \- \int_t^{\wh{\tau} \land \z} \n g (r,\fX_r) dr \nonumber \\
   & \tn  \=  & \tn    E_t \Big[ \int_t^{\wh{\tau}} \n g (r,\fX_r) dr
  \Big| \ol{\cF}^{t}_{\wh{\tau} \land \z}\Big]
 \- \int_t^{\wh{\tau} \land \z} \n g (r,\fX_r) dr
 \= \cY^{t,x,\wh{\tau}}_{\wh{\tau} \land \z}  , \q \pas \label{eq:d214}
 \eea
 As $\wh{\tau} \ins \wh{\cT}^t_x(y)$ by Proposition \ref{prop_surm}, we see from \eqref{eq:c627} that
 \beas
 \q  && \hspace{-2cm} E_t \bigg[ \b1_{\{\wh{\tau} \le \z   \}} \cR \big(t,x,\wh{\tau}\big) \+ \b1_{\{\wh{\tau} \, > \z  \}}
 \Big( \cV(\z,\fX_\z  , \fY_\z  )  \+ \int_t^\z \n f  (r,\fX_r   ) dr \Big) \bigg] \\
 && \gs  E_t \bigg[ \b1_{\{\wh{\tau} \le \z   \}} \cR \big(t,x,\wh{\tau}\big) \+ \b1_{\{\wh{\tau} \, > \z  \}}
 \Big( \cV(\z,\fX_\z  , \cY^{t,x,\wh{\tau}}_\z   )  \+ \int_t^\z \n f  (r,\fX_r   ) dr \Big) \bigg]
 \gs E_t \big[\cR(t,x,\wh{\tau})\big]  ,
 \eeas
 proving the first  inequality in \eqref{eq:d021}.

 \no {\bf 2)} {\it The proof of the second inequality in \eqref{eq:d021} is relatively lengthy,
 we split it into several steps. }

 By an analogy to \eqref{eq:d063}, we must have either $ P_t \{\z \= t \} \= 1 $  or $ P_t \{\z \> t \} \= 1 $.
 If $ P_t \{\z \= t \} \= 1 $, as $\fY_t \= \a_t \= y \> 0$, \pas, one has
 $\wh{\tau} \= \tau(t,x,\a) \>  t \= \z$, \pas ~ Then
 \beas
   E_t \bigg[ \b1_{\{\wh{\tau} \le \z   \}} \cR \big(t,x,\wh{\tau}\big)
 \+ \b1_{\{\wh{\tau} \, > \z  \}}
 \Big( \cV \big(\z,\fX_\z  ,  \fY_\z   \big)  \+ \int_t^\z \n f  (r,\fX_r   ) dr \Big) \bigg]
  \= E_t \big[    \cV  (t,\fX_t  ,  \fY_t   )    \big]
 \=  E_t \big[    \cV  (t, x ,  y   )    \big] \=  \cV  (t, x ,  y   )   .
 \eeas

 So  let us suppose   that $t_1 \> t$ in the rest of this proof.
 There exists a $P_t-$null set $\cN$ such that for any $\o \ins \cN^c$,
 $M_\cd (\o)$ is a continuous path and  $K_\cd (\o)$ is an  continuous increasing path.
 By the uniform integrability of $M $,
 there exists $\xi  \ins L^1 \big(\ol{\cF}^t \big)$ such that \pas
 \bea \label{eq:d107}
   M_s \= E_t \big[\xi  \big| \ol{\cF}^t_s \big]  , \q \fa s \ins [t, \infty) .
 \eea

 For any $i \ins \hN$, similar  to \eqref{eq:c427} and \eqref{eq:c421},
  it holds for all $\o \ins \O^t$ except on a $P_t-$null set $\cN_i$   that
  \bea \label{eq:d007}
  \hspace{-5mm}
   \cN^i_\o  \df \Big\{ \wt{\o} \ins \O^{t_i} \dn : \fX_s   (\o \otii \wt{\o}) \n \ne \n \fX_s   (\o)
   \hb{ for some }  s \ins [t,t_i] \hb{ or } \fX_r (\o \otii \wt{\o})
 \n \ne \n  \cX^{t_i, \fX_{t_i} (\o)}_r (\wt{\o})    \hb{ for some } r \ins [t_i,\infty) \Big\} \ins \sN^{t_i}.
  \eea

 \no {\bf 2a)} Fix $\e \ins (0,1)$.
 The first inequality in \eqref{eq:c483b} and an analogy to \eqref{eq:c491} show that
 \beas
  E_t \Big[ \big| \cR(t,x,\wh{\tau}) \big| \+ \big| \cV(\z,\fX_\z  ,\fY_\z   )\big|
   \+ \int_t^\z \n \big| f  (r,\fX_r   ) \big| dr \Big]
 \ls \Psi(x)\+ 2 \fC(3\+ C_p) \+ \fC(1\+2C_p) E_t \big[ \fX^p_* \big]  \< \infty   .
 \eeas
 So there exists  $\l  \= \l  (t,x, \a, \e) \ins   (0,1) $ such that
\bea \label{eq:d123}
E_t \bigg[\b1_A \Big( \big| \cR(t,x,\wh{\tau}) \big| \+ \big| \cV(\z,\fX_\z  ,\fY_\z   )\big|  \+ \int_t^\z \n \big| f  (r,\fX_r   ) \big| dr \Big) \bigg]
\< \e / 5 \; \hb{ for any $A \ins \ol{\cF}^t $ with $P_t(A) \< \l  $\,. }
\eea
 We can find $ \cI_o \ins \hN$ such that  $   P_t  \{\z \> t_{\cI_o} \}  \< \l/2  $.

  Let $i \= 1,\cds,\cI_o$ and $(\fx,\fy) \ins \hR^l \ti (0,\infty)$.
  In light of \eqref{eq:c629} and Theorem \ref{thm_continuity} (1),
  there exists  $\d_i(\fx,\fy) 
  \ins \big(0, 1   \ld  \fy  \ld \e  \big)$ such that
   \bea
  \big| \cV(t_i,\fx', \fy' ) \-  \cV(t_i,\fx, \fy) \big|  \ls  \e/5 ,
    \q \fa ( \fx' , \fy' )   \ins \ol{O}_{\d_i(\fx,\fy)} (\fx) \ti
    \big[  \fy \- \d_i(\fx,\fy)  , \fy \+  \d_i(\fx,\fy) \big] .  \label{eq:d053}
\eea
 By \eqref{eq:d189}, there exists  $    \a (t_i,\fx,\fy)    \ins \fA_{t_i} \big(\fy\-\d_i(\fx,\fy)\big) $  such that
 \bea \label{eq:d089}
 \cV\big(t_i,\fx, \fy\-\d_i(\fx,\fy)\big)
 \= \underset{\wt{\a} \in \fA_{t_i} (\fy\-\d_i(\fx,\fy))}{\sup} \,
 E_{t_i} \Big[  \cR \big(t_i, \fx, \tau( t_i, \fx, \wt{\a} )  \big) \Big] \ls
 E_{t_i} \Big[  \cR \big(t_i, \fx, \tau( t_i, \fx,  \a (t_i,\fx,\fy) ) \big) \Big] \+ \e/5  ,
\eea
and Proposition \ref{prop_control_cont} shows that
for some $ \wh{\d}_i(\fx,\fy) 
\ins \big(0, \d_i(\fx,\fy)  \big] $
\bea
  E_{t_i} \Big[ \big| \cR \big(t_i, \fx', \tau( t_i, \fx', \a (t_i,\fx,\fy)  ) \big)
  \- \cR \big(t_i, \fx, \tau( t_i, \fx, \a (t_i,\fx,\fy) ) \big)   \big| \Big]
 \ls \e/5 , \q \fa \fx' \ins  \ol{O}_{\wh{\d}_i(\fx,\fy)} (\fx) .   \label{eq:d055}
\eea
 Let us simply write $\cO_i (\fx,\fy)$ for the open set $ O_{\wh{\d}_i(\fx,\fy)} (\fx) \ti
    \big( \fy \- \wh{\d}_i(\fx,\fy)  , \fy \+  \wh{\d}_i(\fx,\fy) \big) $.

 \ss   Since  \eqref{eq:d057} implies that $\fY_{t_i} (\o) \> 0$ for any $\o \ins   \{\wh{\tau} \> t_i\}$, one has
 \beas
 P_t \{ \wh{\tau} \>   t_i \} \= P_t \big( \{ \wh{\tau} \>   t_i \} \Cp \big\{ (\fX_{t_i},\fY_{t_i}) \ins \hR^l \ti (0,\infty) \big\} \big) \= \lmtu{R \to \infty} P_t \big( \{ \wh{\tau} \>    t_i \} \Cp \big\{ (\fX_{t_i},\fY_{t_i})
 \ins \ol{O}_R (0) \ti [R^{-1},R] \big\} \big)  .
 \eeas
 So there exists $ R_i \ins (0,\infty)$ such that
    \bea \label{eq:d077}
    P_t \big( \{ \wh{\tau} \>    t_i \} \Cp  \big\{ (\fX_{t_i},\fY_{t_i})
    \n \notin \n \ol{O}_{R_i} (0) \ti [R^{-1}_i,R_i] \big\} \big)    \ls   \frac{\l}{2^{i+1}}    ,
    \eea
  and we can find a finite   subset   $\big\{   ( x^i_n, y^i_n ) \big\}^{\fn_i}_{n = 1}$ of
  $ \ol{O}_{R_i} (0) \ti [R^{-1}_i,R_i] $ such that
  $ \underset{n = 1}{\overset{\fn_i}{\cup}} \, \cO_i ( x^i_n, y^i_n ) \n \supset \n \ol{O}_{R_i} (0) \ti [R^{-1}_i,R_i] $.

 Let $ n \= 1 , \cds , \fn_i $ and define  $A^i_n \df \{ \wh{\tau} \> \z \= t_i \} \Cp \big\{ (\fX_{t_i},\fY_{t_i})
 \ins \cO_i ( x^i_n, y^i_n ) \big\} \ins \ol{\cF}^t_{t_i} $. Clearly,
 \beas
    \fY_{t_i} (\o)  \- y^i_n  \ins \big( \n   - \wh{\d}_i(x^i_n,y^i_n),   \wh{\d}_i(x^i_n,y^i_n) \big)
 \sb  \big( \n  - \d_i(x^i_n,y^i_n)  , \d_i(x^i_n,y^i_n) \big) , \q \fa \o \ins   A^i_n .
 \eeas
 We also set $ \cA^i_n \df A^i_n \Big\backslash \Big( \underset{n' < n}{\cup} A^i_{n'} \Big) \ins \ol{\cF}^t_{t_i} $
 and define a $\ol{\cF}^t_{t_i}-$measurable random variable
 $\eta^i_n \df \b1_{\cA^i_n} \big( \fY_{t_i}   \- y^i_n \+ \d_i(x^i_n,y^i_n)   \big) \ins \big[0, 2 \d_i(x^i_n,y^i_n) \big) $.
 Suppose that $\a^{i,n} \df \a(t_i,x^i_n, y^i_n)$ equals to  $M^{i,n} \- K^{i,n}$
 for some   $\big(M^{i,n} , K^{i,n}\big) \ins \hM_{t_i} \ti \hK_{t_i}$.
 By the uniform integrability of $M^{i,n}$,
 there exists $\xi^{i,n} \ins L^1 \big(\ol{\cF}^{t_i}  \big)$ such that $P_{t_i}-$a.s.
 \bea \label{eq:d079}
   M^{i,n}_s \= E_{t_i} \big[\xi^{i,n} \big| \ol{\cF}^{t_i}_s \big]  , \q \fa s \ins [t_i, \infty) .
 \eea
 Let $\cN^{i,n}$ be the $P_{t_i}-$null set such that for any $\wt{\o} \ins (\cN^{i,n})^c$,
 $M^{i,n}_\cd (\wt{\o})$ is a continuous path;   $K^{i,n}_\cd (\wt{\o})$ is an  continuous increasing path;
 and
 \bea \label{eq:d075}
 \a^{i,n}_{t_i} (\wt{\o}) \= M^{i,n}_{t_i} (\wt{\o}) \= y^i_n \-   \d_i(x^i_n,y^i_n) \> 0  .
 \eea
 As   $(\Pi^t_{t_i})^{-1}(\cN^{i,n})$ is a $P_t-$null set by Lemma \ref{lem_shift_inverse1b} (1),
 one can deduce from Lemma \ref{lem_shift_inverse1b} (2) that
 \bea \label{eq:d141}
 \hspace{-3mm}
 \ba{r}
 \hb{$ M^{i,n}_s (\Pi^t_{t_i})$, $s \ins  [t_i,\infty)$ is
 an $\ol{\bF}^t-$adapted continuous  process with $ M^{i,n}_{t_i} (\Pi^t_{t_i}) \= y^i_n \-   \d_i(x^i_n,y^i_n) $, \pas ~
 and} \q \ss \\
 \hb{$ K^{i,n}_s (\Pi^t_{t_i})$, $s \ins  [t_i,\infty)$ is
 an $\ol{\bF}^t-$adapted,  continuous  increasing process   with    $K^{i,n}_{t_i} (\Pi^t_{t_i}) \= 0$, \pas} \q
 \ea
 \eea
 An analogy to \eqref{eq:d057} and \eqref{eq:ax041} shows that
  $\nu^i_n \df  \inf\big\{s \ins [t_i,\infty) \n : \a^{i,n}_s  ( \Pi^t_{t_i}    )
  \- \n \int_{t_i}^s \n g  (r, \fX_r   ) dr  \= 0 \big\}$
 defines  a $\ol{\cT}^t-$stopping time.
 Since $ \a^{i,n}_{t_i} (\Pi^t_{t_i}) \> 0 $, \pas ~ by \eqref{eq:d075},
 we see that   $\nu^i_n \> t_i $, \pas  ~ and thus $E_t \big[ \nu^i_n \- t_i \big|\ol{\cF}^t_{t_i}\big] \> 0$, \pas ~
 \if{0}
 Set $A \df \big\{ E_t \big[ \nu^i_n \- t_i \big|\ol{\cF}^t_{t_i}\big] \= 0 \big\} \ins \ol{\cF}^t_{t_i}$.
 \beas
  0 \=  E_t \Big[ \b1_A  E_t \big[ \nu^i_n \- t_i \big|\ol{\cF}^t_{t_i}\big] \Big]
  \=  E_t \Big[  E_t \big[ \b1_A  (\nu^i_n \- t_i) \big|\ol{\cF}^t_{t_i}\big] \Big]
  \=  E_t \big[ \b1_A  (\nu^i_n \- t_i) \big] .
 \eeas
 which implies that $ \b1_A  (\nu^i_n \- t_i) \=0 $, ~ \pas ~ Thus $P_t(A) \= 0$.
 \fi

 Similar to the proof of Proposition \ref{prop_surm},
 $ J^{i,n}_s \df \underset{s' \in  [t_i,s]}{\inf} E_t \big[ \nu^i_n \- t_i \big| \ol{\cF}^t_{s'}\big] $,
 $ s \ins [t_i,\infty) $ is an $\ol{\bF}^t-$adapted, non-negative, continuous decreasing process
 such that $ J^{i,n}_{t_i} \= E_t \big[ \nu^i_n \- t_i \big|\ol{\cF}^t_{t_i}\big] \> 0$, \pas ~
 and that $ J^{i,n}_{\nu^i_n} \ls E_t \Big[ \nu^i_n \- t_i \big| \ol{\cF}^t_{\nu^i_n}\Big]
 \= \nu^i_n \- t_i $, \pas ~
 Then   $ \cK^{i,n}_s   \df   \eta^i_n
 \Big[ 1 \ld \big( \frac{s - t_i}{J^{i,n}_s} \big)^{\n +} \, \Big] \gs 0 $, $ s \ins [t_i,\infty) $
 defines an $\ol{\bF}^t-$adapted,   continuous increasing process over period $[t_i,\infty)$ such that
 \bea \label{eq:d143}
 \cK^{i,n}_{t_i}   \= 0 \q \hb{and} \q \cK^{i,n}_{\nu^i_n} \= \eta^i_n \, \hb{ holds except on a $P_t-$null set } \cN^{i,n}_K .
 \eea

 Set $\cA_\sharp \df  \underset{i = 1}{\overset{\cI_o}{\cup}} \,
   \underset{n = 1}{\overset{\fn_i}{\cup}} \cA^i_n  \ins \ol{\cF}^t_{t_{\cI_o}}$
 and $\cN_\sharp \df  \cN \cp \Big( \underset{i = 1}{\overset{\cI_o}{\cup}} \,
   \underset{n = 1}{\overset{\fn_i}{\cup}} (\Pi^t_{t_i})^{-1}(\cN^{i,n}) \Big)   \ins \sN^t$.
 We claim that
  \bhe
  \bea \label{eq:d071}
  \ba{rll}
  \ol{M}_s & \tn  \df  & \tn  M_s \+
    \sum^{\cI_o}_{i=1} \, \sum^{\fn_i}_{n=1} \,   \b1_{\{s \ge t_i  \} \cap \cA^i_n }
    \big( M^{i,n}_s  ( \Pi^t_{t_i}    ) \- M_s \+ M_{t_i}
     \- y^i_n \+ \d_i(x^i_n,y^i_n)       \big)  , ~ s \ins [t,\infty)
  \, \hb{ is of } \,  \hM_t , \qq  \ss  \\
  \hb{and } \;   \ol{K}_s  & \tn  \df  & \tn  K_s \+
  \sum^{\cI_o}_{i=1}   \sum^{\fn_i}_{n=1} \,  \b1_{\{s \ge t_i  \} \cap \cA^i_n }
 \big( K^{i,n}_s  ( \Pi^t_{t_i}    )  \- K_s \+ K_{t_i}   \+  \cK^{i,n}_s \big)    ,
 ~ s \ins [t,\infty)    \, \hb{ is of } \,  \hK_t  .
 \ea
 \eea
 \ehe

 As $t_1\>t$ by assumption, it holds \pas ~ that $\ol{M}_t   \= M_t \= y$.
 So $\ol{\a}   \df \ol{M}   \- \ol{K} \ins \fA_t(y)$.

 \no {\bf 2b)} {\it Setting  $\ol{\tau} \df \tau(t,x,\ol{\a})$, we next show that $\ol{\tau} \= \wh{\tau} $, \pas ~
on $ \{\wh{\tau} \ls  \z \}     \cp  \big( \{\wh{\tau} \> \z \} \Cp \cA_\sharp^c \big) $.}

Since \eqref{eq:d071}   shows that
\bea \label{eq:d103}
\hb{$\big(\ol{M}_s(\o), \ol{K}_s(\o)\big) \= \big(M_s(\o),K_s(\o)\big) $,   $\fa (s,\o) \ins  \big( [t,\infty) \ti \cA_\sharp^c \big) \cp \[ t,\z \[ $\,,}
\eea
we obtain that
\bea \label{eq:d099}
 \ol{\a}_s(\o) \= \a_s(\o), \q \fa (s,\o) \ins  \big( [t,\infty) \ti  \cA_\sharp^c   \big)
 \cp \big[ \dn \big[ t, \z   \big[ \dn \big[ \,.
\eea
So for any $\o \ins \cA_\sharp^c$, one has
$ \wh{\tau} (\o)   \=    \big( \tau(t,x,\a) \big) (\o) \= \inf\Big\{ s \ins [t ,\infty) \n :  \a_s ( \o   )
\- \int_t^s \n g \big(r, \fX_r ( \o   ) \big) dr \= 0 \Big\}
     \=   \inf\Big\{ s \ins [t ,\infty) \n :  \ol{\a}_s ( \o   )
\- \int_t^s \n g \big(r, \fX_r ( \o   ) \big) dr \= 0 \Big\}
\= \big( \tau(t,x,\ol{\a}) \big) (\o) \= \ol{\tau} ( \o   ) $.

 Let   $\o \ins \{\wh{\tau} \ls  \z   \ls \cI_o\}  \Cp \cN_\sharp^c  $.
By   \eqref{eq:d099},
\bea \label{eq:d101}
\ol{\a} (s,\o) \= \a (s,\o) , \q \fa s \ins   [t,\z(\o)  )   .
\eea
If $ \wh{\tau} (\o) \<  \z (\o) $, one can deduce from \eqref{eq:d101} that
\beas
\wh{\tau} (\o) & \tn \= & \tn 
\inf\Big\{ s \ins [t ,\infty) \n :  \a_s ( \o   )
\- \int_t^s \n g \big(r, \fX_r ( \o   ) \big) dr \= 0 \Big\}
\= \inf\Big\{ s \ins [t ,\z(\o)) \n : \a_s ( \o   )
\- \int_t^s \n g \big(r, \fX_r ( \o   ) \big) dr \= 0 \Big\} \\
 & \tn \= & \tn  \inf\Big\{ s \ins [t ,\z(\o)) \n : \ol{\a}_s ( \o   )
\- \int_t^s \n g \big(r, \fX_r ( \o   ) \big) dr \= 0 \Big\} ,
\eeas
which implies that
$\ol{\tau} ( \o   ) \= 
\inf\big\{ s \ins [t ,\infty) \n : \ol{\a}_s ( \o   )
\- \int_t^s \n g \big(r, \fX_r ( \o   ) \big) dr \= 0 \big\}
\= \wh{\tau} (\o) $.

 Otherwise, suppose that $\wh{\tau} (\o) \=  \z (\o)$. The definition of $\tau(t,x,\a)$ and \eqref{eq:d101} show that
 \bea \label{eq:d105}
\ol{\a} (s,\o) \= \a (s,\o) \> \int_t^s \n g \big(r, \fX_r(\o)\big) dr, ~ \fa s \ins   [t,\z(\o)  )
\q \hb{and} \q \a \big(\z(\o),\o\big) \=  \int_t^{\z(\o)} \n g \big(r, \fX_r(\o)\big) dr .
 \eea
 As $\ol{M}_\cd (\o) ,  M_\cd (\o),\ol{K}_\cd (\o) ,  K_\cd (\o)$ are all continuous paths by the proof of \eqref{eq:d071},
 we see from  \eqref{eq:d103}   and \eqref{eq:d105}    that
 \beas
 \ol{\a} \big(\z (\o),\o\big) \= \big( \ol{M} \- \ol{K}\big) \big(\z (\o),\o\big)
 \= ( M   \- K ) \big(\z(\o),\o\big) \= \a \big(\z(\o),\o\big) \=
 \int_t^{\z(\o)} \n g \big(r, \fX_r(\o)\big) dr ,
 \eeas
 which means that   $ \ol{\tau} ( \o   ) \= \big( \tau(t,x,\ol{\a}) \big) (\o) \= \z (\o) \= \wh{\tau} (\o)  $.
 Hence, we have  verified that
 \bea \label{eq:d125}
 \hb{ $\ol{\tau} \= \wh{\tau} $, \; \pas ~
 on }   \cA_\sharp^c   \cp   \{\wh{\tau} \ls  \z   \ls  \cI_o\}
 \=  \{\wh{\tau} \ls  \z \}     \cp  \big( \{\wh{\tau} \> \z \} \Cp \cA_\sharp^c \big)  .
 \eea

  \no {\bf 2c)} {\it
  Let $i \= 1, \cds ,\cI_o$ and $n \ins 1,\cds,\fn_i$.
   In this step, we demonstrate that }
  \beas
  E_t \big[ \b1_{  \cA^i_n  }   \cR  ( t , x ,  \ol{\tau}  )    \big]
  \gs E_t \Big[ \b1_{  \cA^i_n  } \Big(\cV \big(t_i, \fX_{t_i}   , \fY_{t_i}   \big)    \+ \int_t^{t_i} f(r,  \fX_r   ) dr \- 4\e/5 \Big) \Big] .
  \eeas

 Set $\wh{\cN}^{i,n} \df \{\o \ins  \O^t \n : \nu^i_n (\o) \= \infty  \} \ins \sN^t $ and
 $ G^i_n \df   \cA^i_n \Cp  \big(  \cN_i  \cp \wh{\cN}^{i,n} \cp \cN^{i,n}_K \cp (\Pi^t_{t_i})^{-1}(\cN^{i,n}) \big)^c
 \ins \ol{\cF}^t_{t_i} $.
  Let $\o \ins  G^i_n $.  The definition  of $\nu^i_n$ shows that
   \bea \label{eq:d145}
    \a^{i,n}_s  ( \Pi^t_{t_i} (\o)   ) \>   \int_{t_i}^s \n g \big(r, \fX_r (\o)  \big) dr     ,
    \q \fa s \ins \big[t_i, \nu^i_n(\o)\big) \q \hb{and} \q  \a^{i,n}   \big( \nu^i_n (\o), \Pi^t_{t_i} (\o)   \big)
    \= \int_{t_i}^{\nu^i_n (\o)} \n g \big(r, \fX_r (\o)  \big) dr    .
  \eea
 Since $\o \ins \cA^i_n \sb \{\wh{\tau} \>   t_i  \}$ and since
\beas
 \ol{\a}_s(\o)  & \tn  \=  & \tn   \b1_{\{s < t_i  \}} \a_s (\o) \+ \b1_{\{s \ge t_i  \}}
 \Big( \a^{i,n}_s  \big( \Pi^t_{t_i} (\o)   \big) \+  \int_t^{t_i} g \big(r,\fX_r (\o)\big) dr
 \+  \eta^i_n (\o) \-  \cK^{i,n}_s (\o)  \Big)    , \q  s \ins [t,\infty) ,  
\eeas
we can deduce from \eqref{eq:d143} and \eqref{eq:d145} that
\beas
 \ol{\a}_s (\o  ) & \tn \= & \tn  \a_s (\o  ) \> \int_t^s g \big(r,\fX_r (\o  )\big) dr , \q \fa s \ins [t,t_i) ,   \\
   \ol{\a}_s(\o)  & \tn \gs  & \tn     \a^{i,n}_s  \big( \Pi^t_{t_i} (\o)   \big)
   \+  \int_t^{t_i} g \big(r,\fX_r (\o)\big) dr
    \>   \int_t^s \n g \big(r, \fX_r (\o)  \big) dr  , \q \fa s \ins \big[t_i, \nu^i_n(\o)\big)  ,   \\
 \hb{and} \q     \ol{\a} \big(\nu^i_n (\o),\o\big)  & \tn \=  & \tn  \a^{i,n}   \big( \nu^i_n (\o), \Pi^t_{t_i} (\o)   \big)
  \+  \int_t^{t_i} g \big(r,\fX_r (\o)\big) dr
    \= \int_t^{\nu^i_n (\o)} \n g \big(r, \fX_r (\o)  \big) dr  ,
    \eeas
     which   implies that
     \bea \label{eq:d147}
     \ol{\tau} (\o) \= \big(\tau (t,x,\ol{\a})\big)(\o)\= \nu^i_n (\o) , \q \fa \o \ins G^i_n .
     \eea

  Similar to Problem 2.7.3 of \cite{Kara_Shr_BMSC},
 there exists $  \wt{G}^i_n \ins \cF^t_{t_i}$   such that $ \cN^{i,n}_G \df  G^i_n  \D  \wt{G}^i_n  \ins  \sN^t  $.
 By Proposition \ref{prop_null_set} (1), it holds for all $\o \ins \O^t$ except on a $P_t-$null set
 $\wh{\cN}^{i,n}_G$ that $  \big( \cN^{i,n}_G \big)^{t_i,\o} \ins \sN^{t_i}$.

  Now, let $\o \ins G^i_n \Cp  \wt{G}^i_n  \Cp \big( \wh{\cN}^{i,n}_G \big)^c$
  and $\wt{\o} \ins \big(  (  \cN^{i,n}_G  )^{t_i,\o}  \cup \cN^i_\o \big)^c   $.
  As $\o \ins  \wt{G}^i_n $
  and $\wt{\o} \ins \big(  ( \cN^{i,n}_G  )^{t_i,\o} \big)^c \= \big(  ( \cN^{i,n}_G  )^c \big)^{t_i,\o}$, Lemma \ref{lem_element} shows that   $\o \otii \wt{\o} \ins \wt{G}^i_n$
  and thus $  \o \otii \wt{\o} \ins  \wt{G}^i_n \Cp \big( \cN^{i,n}_G  \big)^c \= G^i_n \Cp  \wt{G}^i_n
  \sb G^i_n    $. Applying \eqref{eq:d147} with $\o \= \o \otii \wt{\o}$, we see from \eqref{eq:d007} that
  \beas
  \ol{\tau} (\o \otii \wt{\o}) & \tn \= & \tn  \nu^i_n (\o \otii \wt{\o})
  \= \inf\big\{s \ins [t_i,\infty) \n : \a^{i,n}_s  ( \wt{\o}   )
  \- \n \int_{t_i}^s \n g \big(r, \fX_r (\o \otii \wt{\o})  \big) dr  \= 0 \big\} \\
   & \tn   \=  & \tn   \inf\Big\{s \ins [t_i,\infty) \n : \a^{i,n}_s  ( \wt{\o}  )
  \- \int_{t_i}^s \n g \Big(r, \cX^{t_i,\fX_{t_i}(\o)}_r (   \wt{\o} ) \Big) dr \= 0 \Big\}
  \= \big( \tau \big(t_i, \fX_{t_i}(\o), \a^{i,n} \big) \big) (\wt{\o}) \n =: \n \tau^{i,n}_\o (\wt{\o}) .
  \eeas
  Then \eqref{eq:d007}   again shows  that
 \beas
 \big( \cR  ( t , x , \ol{\tau}  ) \big)^{t_i,\o} (\wt{\o})
 & \tn \=& \tn  \big( \cR  ( t , x , \ol{\tau}  ) \big) (\o \otii \wt{\o})
 \= \int_t^{\ol{\tau} (\o \otimes_{t_i} \wt{\o})} \n
f \big(r,\fX_r (\o \otii \wt{\o})\big) dr
\+  \pi \Big(\ol{\tau} (\o \otii \wt{\o}),\fX  \big( \ol{\tau} (\o \otii \wt{\o}), \o \otii \wt{\o} \big)\Big)   \\
& \tn \=& \tn  \int_t^{t_i} \n f \big(r,\fX_r (\o)\big) dr \+
\int_{t_i}^{\tau^{i,n}_\o (\wt{\o}) } \n f \Big(r,\cX^{t_i, \fX_{t_i} (\o)}_r (\wt{\o})\Big) dr
\+ \pi \Big(\tau^{i,n}_\o (\wt{\o}),\cX^{t_i,\fX_{t_i} (\o)}  \big( \tau^{i,n}_\o (\wt{\o}),\wt{\o}\big)\Big)     \\
  & \tn \=& \tn  \int_t^{t_i} \n f \big(r,\fX_r (\o)\big) dr \+
  \big( \cR( t_i , \fX_{t_i} (\o) , \tau^{i,n}_\o   ) \big) (\wt{\o}) .
 \eeas
 Taking expectation $E_{t_i}[\cd]$ over $\wt{\o} \ins \O^{t_i}$ except the $P_{t_i}-$null set
 $  (  \cN^{i,n}_G  )^{t_i,\o}  \cup \cN^i_\o $ yields that
 \beas
 E_{t_i} \Big[ \big( \cR  ( t , x , \ol{\tau}  ) \big)^{t_i,\o} \Big]
 \= E_{t_i} \big[ \cR \big( t_i , \fX_{t_i} (\o) , \tau  (t_i, \fX_{t_i}(\o), \a^{i,n} )  \big) \big]
 \+  \int_t^{t_i} \n f \big(r,\fX_r (\o)\big) dr .
 \eeas
 Since $\big( \fX_{t_i} (\o), \fY_{t_i} (\o)\big) \ins O_{\wh{\d}_i(x^i_n,y^i_n)} (x^i_n) \ti
 \big( y^i_n \- \wh{\d}_i(x^i_n,y^i_n)  , y^i_n \+  \wh{\d}_i(x^i_n,y^i_n) \big)$,
 using \eqref{eq:d055} with $(\fx,\fy,\fx') \= \big(x^i_n,y^i_n,\fX_{t_i} (\o) \big)$
 and applying \eqref{eq:d053} with $(\fx,\fy,\fx',\fy') \= \big(x^i_n,y^i_n,x^i_n,y^i_n \- \d_i(x^i_n, y^i_n) \big)$
 and $(\fx,\fy,\fx',\fy') \= \big(x^i_n,y^i_n,\fX_{t_i} (\o), \fY_{t_i} (\o) \big)$ respectively,
 we can deduce from \eqref{eq:d089}   that
 \bea
 \q && \hspace{-2cm}  E_{t_i} \Big[ \big( \cR  ( t , x , \ol{\tau}  ) \big)^{t_i,\o} \Big]
   \- \int_t^{t_i} f(r,  \fX_r (\o) ) dr
   \gs     E_{t_i} \Big[ \cR \big(t_i, x^i_n,  \tau(t_i, x^i_n, \a^{i,n} ) \big)
 \Big]   \- \e/5
    \gs   \cV \big(t_i, x^i_n, y^i_n \- \d_i(x^i_n, y^i_n) \big)   \- 2\e / 5 \nonumber \\
 &   &   \gs \cV \big(t_i, x^i_n, y^i_n   \big)   \- 3 \e/5
   \gs   \cV \big(t_i, \fX_{t_i} (\o), \fY_{t_i} (\o) \big)   \- 4\e/5 ,
  \q \fa \o \ins G^i_n \Cp  \wt{G}^i_n  \Cp \big( \wh{\cN}^{i,n}_G \big)^c  . \label{eq:d121}
 \eea

 The first inequality in \eqref{eq:c483b} and   Proposition \ref{prop_shift_FP} (2) imply that
  $  E_t \big[ \cR  ( t , x , \ol{\tau}  ) \big|\ol{\cF}^t_{t_i} \big] (\o)
  \= E_{t_i} \Big[ \big( \cR  ( t , x , \ol{\tau}  ) \big)^{t_i,\o} \Big]  $ for \pas ~ $\o \ins \O^t$.
 As   $ \b1_{ G^i_n \cap  \wt{G}^i_n   } \= \b1_{ G^i_n } \b1_{ \wt{G}^i_n   } \=
  \b1_{ G^i_n }    \= \b1_{  \cA^i_n }  $,   \pas,
 we can derive from \eqref{eq:d121} that
\beas
\q && \hspace{-1cm} E_t \big[ \b1_{  \cA^i_n } \cR  ( t , x , \ol{\tau}  ) \big]
\= E_t \Big[ \b1_{  \cA^i_n } E_t \big[  \cR  ( t , x , \ol{\tau}  ) \big| \ol{\cF}^t_{t_i} \big] \Big]
\= E_t \Big[ \b1_{ G^i_n \cap  \wt{G}^i_n   } E_{t_i} \big[  ( \cR  ( t , x , \ol{\tau}  )  )^{t_i,\o} \big] \Big] \\
&& \gs E_t \bigg[ \b1_{ G^i_n \cap  \wt{G}^i_n   } \Big( \cV \big(t_i, \fX_{t_i}  , \fY_{t_i}   \big) \+ \int_t^{t_i} f(r,  \fX_r   ) dr \- 4\e/5 \Big) \bigg]
\= E_t \bigg[ \b1_{  \cA^i_n }
\Big( \cV \big(\z, \fX_\z  , \fY_\z   \big) \+ \int_t^\z f(r,  \fX_r   ) dr \- 4\e/5 \Big) \bigg] .
\eeas
Taking summation over $n \ins 1,\cds,\fn_i$ and $i \= 1, \cds ,\cI_o$ and using the conclusion of Part 2 yield that
\bea \label{eq:d131}
 E_t \big[  \cR  ( t , x , \ol{\tau}  ) \big]
 \gs E_t \bigg[ \b1_{ \{\wh{\tau} \le  \z \}     \cup  ( \{\wh{\tau} \, > \z \} \cap \cA_\sharp^c  )}
 \cR  ( t , x , \wh{\tau}) \+
 \b1_{  \cA_\sharp  } \Big( \cV \big(\z, \fX_\z  , \fY_\z   \big) \+ \int_t^\z f(r,  \fX_r   ) dr \Big) \bigg] \- 4\e/5 .
\eea

  \no {\bf 2d)} Since    $ \underset{n = 1}{\overset{\fn_i}{\cup}} \cA^i_n
   \= \underset{n = 1}{\overset{\fn_i}{\cup}}  A^i_n
 \= \{ \wh{\tau} \> \z \= t_i \} \Cp \Big\{ (\fX_{t_i},\fY_{t_i})
 \ins \underset{n = 1}{\overset{\fn_i}{\cup}} \cO_i ( x^i_n, y^i_n ) \Big\} $ for $i\=1,\cds,\cI_o$,
 one can deduce   that
  \beas
 \{\wh{\tau}  \> \z \} \Cp \cA_\sharp^c \= \{\wh{\tau}  \> \z \> \cI_o \} \cp
 \bigg( \underset{i = 1}{\overset{\cI_o}{\cup}} \Big( \{ \wh{\tau} \> \z \= t_i \} \Cp \Big\{ (\fX_{t_i},\fY_{t_i})
 \n \notin \n \underset{n = 1}{\overset{\fn_i}{\cup}} \cO_i ( x^i_n, y^i_n ) \Big\} \Big) \bigg)  ,
 \eeas
 and \eqref{eq:d077} implies that  $P_t \big( \{\wh{\tau}  \> \z \} \Cp \cA_\sharp^c \big)
 \ls P_t\{\z \> \cI_o\} \+  \sum_{i=1}^{\cI_o} P_t \big( \{ \wh{\tau} \>   t_i \} \Cp \big\{ (\fX_{t_i},\fY_{t_i})
 \n \notin \n  \ol{O}_{R_i} (0) \ti [R^{-1}_i,R_i] \big\} \big)  \< \l$. It then follows from \eqref{eq:d123} that
 \beas
 \qq && \hspace{-1.5cm} \bigg| E_t \bigg[ \b1_{\{\wh{\tau} \, > \z \} \cap \cA_\sharp^c}
 \Big( \cR  ( t , x , \wh{\tau}) \- \cV \big(\z, \fX_\z  , \fY_\z   \big) \- \int_t^\z f(r,  \fX_r   ) dr \Big) \bigg]\bigg| \\
 && \ls E_t \bigg[ \b1_{\{\wh{\tau} \, > \z \} \cap \cA_\sharp^c}  \Big( \big| \cR(t,x,\wh{\tau}) \big|
 \+ \big| \cV(\z,\fX_\z  ,\fY_\z   )\big|  \+ \int_t^\z \n \big| f  (r,\fX_r   ) \big| dr \Big) \bigg] \< \e / 5 ,
 \eeas
 which together with \eqref{eq:d131} and \eqref{eq:d189} leads to that
 \beas
   \cV(t,x,y) \gs  E_t \big[  \cR  ( t , x , \ol{\tau}  ) \big]
 \gs E_t \bigg[ \b1_{ \{\wh{\tau} \le  \z \}     }
 \cR  ( t , x , \wh{\tau}) \+
 \b1_{ \{\wh{\tau} \, > \z \} } \Big( \cV \big(\z, \fX_\z  , \fY_\z   \big)
 \+ \int_t^\z f(r,  \fX_r   ) dr \Big) \bigg] \- \e .
 \eeas
 Letting $\e \to \infty$ yields the second inequality in \eqref{eq:d021}. \qed

 \no {\bf Proof of Theorem \ref{thm_DPP2}:} Fix $t \ins [0,\infty)$.

\no {\bf  1)} Let $(x,y) \ins \hR^l \ti (0,\infty)$ and  $\{\z (\a) \}_{\a   \in \fA_t (y)} $ be a family of $  \ol{\cT}^t_\sharp -$stopping times.
For any $ \a \ins \fA_t (y) $, taking $\z \= \z (\a)$ in \eqref{eq:d021} yields that
 \beas
  E_t \big[   \cR \big(t,x,\tau(t,x,\a)\big) \big]
  & \tn \ls  & \tn  E_t \bigg[ \b1_{\{\tau(t,x,\a) \le \z (\a)  \}} \cR \big(t,x,\tau(t,x,\a)\big)   \\
  & \tn   &    + \b1_{\{\tau(t,x,\a) > \z (\a) \}}
 \Big( \cV \big(\z (\a),\cX^{t,x}_{\z (\a)} ,  Y^{t,x,\a}_{\z (\a)}  \big)
  \+ \int_t^{\z (\a)} \n f  (r,\cX^{t,x}_r   ) dr \Big) \bigg] \ls  \cV(t,x,y )    .
 \eeas
 Taking supremum over  $ \a \ins \fA_t (y) $, we obtain \eqref{eqn_DPP2} from \eqref{eq:d189}.

\no {\bf  2)} Next, suppose that $\cV(s,x,y)$ is continuous in $(s,x,y) \ins [t,\infty) \ti \hR^l \ti (0,\infty)$.

 We fix $(x,y) \ins \hR^l \ti (0,\infty)$ and  a family  $\{\z (\a) \}_{\a   \in \fA_t (y)} $
 of $  \ol{\cT}^t -$stopping times.
 Let   $ \a \ins \fA_t (y) $, $n \ins \hN$ and define
 \beas
  \z_n \=  \z_n (\a) \df \b1_{\{\z (\a) = t \}}t \+ \sum_{i \in \hN}
 \b1_{\{\z (\a) \in  (t+(i-1)2^{-n},t+i2^{-n}]\}} (t\+i2^{-n}) \ins \ol{\cT}^t .
 \eeas
 Applying \eqref{eq:d021} with $\z \= \z_n$ yields that
 \beas
 E_t \big[   \cR \big(t,x,\tau(t,x,\a)\big) \big]
 & \tn \ls  & \tn  E_t \bigg[ \b1_{\{\tau(t,x,\a) \le \z_n   \}} \cR \big(t,x,\tau(t,x,\a)\big) \\
   & \tn   &   + \b1_{\{\tau(t,x,\a) > \z_n  \}}
 \Big( \cV \big(\z_n,\cX^{t,x}_{\z_n}  , Y^{t,x,\a}_{\z_n}   \big)
  \+ \int_t^{\z_n} \n f  (r,\cX^{t,x}_r   ) dr \Big) \bigg]  \ls     \cV(t,x,y )   .
 \eeas
\if{0}
 An analogy to \eqref{eq:c491}  shows that
 $ \big| \cV \big(\z_n,\cX^{t,x}_{\z_n}  , Y^{t,x,\a}_{\z_n}   \big)\big|
   \+ \int_t^{\z_n} \n \big| f (r,\cX^{t,x}_r   ) \big| dr  \ls    2 \fC(3\+ C_p) \+ \fC(1\+2C_p)   \fX^p_* $,
 which is $P_t-$integrable. Since $\lmtd{n \to \hN} \z_n \= \z$, one has
 $ \lmtd{n \to \hN} \b1_{\{\tau(t,x,\a) \le \z_n   \}} \= \b1_{\{\tau(t,x,\a) \le \z   \}} $.
 \fi
 As   $n \to \infty$, using similar arguments to those that lead to \eqref{eq:d135}
  we can deduce from the continuity of function $\cV$
 in $(s,\fx,\fy) \ins [t,\infty) \ti \hR^l \ti (0,\infty)$,
 the continuity of processes $\big(\cX^{t,x},Y^{t,x,\a}\big)$, 
 and the dominated convergence theorem that
 \beas
 E_t \big[   \cR \big(t,x,\tau(t,x,\a)\big) \big]
 & \tn \ls  & \tn  E_t \bigg[ \b1_{\{\tau(t,x,\a) \le \z(\a)   \}} \cR \big(t,x,\tau(t,x,\a)\big) \\
  & \tn  &  + \b1_{\{\tau(t,x,\a) > \z(\a)  \}}
 \Big( \cV \big(\z(\a),\cX^{t,x}_{\z(\a)}  , Y^{t,x,\a}_{\z(\a)}   \big)
  \+ \int_t^{\z(\a)} \n f  (r,\cX^{t,x}_r   ) dr \Big) \bigg]
     \ls     \cV(t,x,y )   .
 \eeas
 Taking supremum over  $ \a \ins \fA_t (y) $ and using \eqref{eq:d189} yield \eqref{eqn_DPP2} again.     \qed

 \subsection{Proof of Section \ref{sec:PDE}}

 \ss \no {\bf Proof of Theorem \ref{thm_visc_exist}:} Under \eqref{eq:c601} and   \eqref{eq:c611},
 Theorem \ref{thm_continuity} (2) and \eqref{eq:c629} show  that
 $\cV$  is  continuous in $(t, x, y ) \ins [0,\infty) \ti \hR^l \ti  [0,\infty)  $.
 By \eqref{eq:c433b}, $\cV (  t,x,0  ) \= \pi (t,x) $ for any $ (t,x) \ins [0,\infty) \ti \hR^l  $.

 \no {\bf 1)}  We first show that  $\cV$  is a viscosity supersolution   of \eqref{eq:PDE}.

 Let $(t_o,x_o,y_o ) \ins (0,\infty) \ti \hR^l \ti (0,\infty)   $
 and  let   $ \phi \ins C^{1,2,2}\big([0,\infty) \ti \hR^l \ti [0,\infty)\big)$
 such that    $\cV \- \phi$ attains a strict local  minimum $0$  at $(t_o,x_o,y_o)$. So there exists a
    $\d_o \ins  \big(0,    t_o  \ld y_o \big) $ such that for any $(t,x,y) \in   O_{\d_o} (t_o,x_o,y_o)
     \big\backslash \big\{ (t_o,x_o,y_o) \big\}$
   \bea \label{eq:d167}
     (\cV \- \phi) (t,x,y) \>  (\cV \- \phi) (t_o,x_o,y_o) \= 0   \hb{ and }
     \big|D_x \phi (t,x,y) \- D_x \phi (t_o,x_o,y_o)  \big| \ve
     \big| \pa_y \phi (t,x,y) \- \pa_y \phi (t_o,x_o,y_o)\big| \<  1 . \q
   \eea

  According to \eqref{eq:c601} and   \eqref{eq:c611},   the functions $b,\si,f,g$ are   continuous in $(t, x)$.  Then
 \bea \label{eq:d193}
  \wh{\phi} \,  (t,x,y) \df  -  \pa_t \phi (t,x,y) \- \cL_x   \phi (t,x,y)  \+ g(t,x)    \pa_y \phi (t,x,y)
 \- f(t,x)  , \q  \fa (t,x,y) \ins [0,\infty) \ti \hR^l \ti [0,\infty)
 \eea
 is also a continuous function.

 To show $ \wh{\phi} \,  (t_o,x_o,y_o) \- \cH \phi  (t_o,x_o,y_o) \gs 0 $,
 it suffices to verify that for any $a \ins \hR^d$
  \beas
 \wh{\phi} \,  (t_o,x_o,y_o) \-   \hb{$\frac12$} |a|^2   \pa^2_y \phi (t_o,x_o,y_o)
\- \big(D_x (\pa_y \phi (t_o,x_o,y_o))  \big)^T \dn \cd \n \si (t_o,x_o) \n \cd \n  a \gs  0  .
  \eeas
   Assume not, i.e. there exists   an $\fra \ins \hR^d$ such that
  \beas
  \e \df  \frac12 |\fra|^2   \pa^2_y \phi (t_o,x_o,y_o)
\+ \big(D_x (\pa_y \phi (t_o,x_o,y_o))  \big)^T \dn \cd \n \si (t_o,x_o) \n \cd \n  \fra
\- \wh{\phi} \,  (t_o,x_o,y_o)   \> 0  .
  \eeas
  Using the continuity of $ \si, \phi$ and $\wh{\phi}$, we can find   some $ \d  \ins (0,\d_o)$ such that
  \bea \label{eq:d169}
  \frac12 |\fra|^2   \pa^2_y \phi (t,x,y) \+ \big(D_x (\pa_y \phi (t,x,y))  \big)^T
  \dn \cd \n \si (t,x) \n \cd \n  \fra  \- \wh{\phi} \,  (t,x,y)   \gs \frac12 \e \> 0 ,
  \q \fa (t,x,y) \ins \ol{O}_\d (t_o,x_o,y_o) .
  \eea

 Clearly, $M_s \df y_o \+  \fra^T  \n \cd    W^{t_o}_s $, $  s \ins [t_o,\infty) $ is a uniformly integrable
 continuous martingale with respect to $\big( \ol{\bF}^{t_o}, P_{t_o} \big)$.
 By taking $K \equiv 0$, we have $\a^o \df M   \ins \fA_{t_o} (y_o)$.
 As $ \Th_s \df \big(s,\cX^{t_o,x_o}_s  , Y^{t_o,x_o,\a^o}_s\big) $, $s \ins [t_o,\infty)$
 are $\ol{\bF}^{t_o}-$adapted continuous processes with $\Th_{t_o} \= (t_o,x_o,y_o)$, \pas,
 $ \z 
  \df   \inf\big\{ s \ins [t_o, \infty) \n : \Th_s   \n  \notin \n  \ol{O}_\d (t_o,x_o,y_o)   \big\} $
 defines  an $\ol{\bF}^{t_o}$-stopping time with $t_o \< \z   \ls t_o \+ \d$, $P_{t_o}-$a.s.
 Since
 \bea \label{eq:d171}
 \hb{$\Th_s  \ins \ol{O}_\d (t_o,x_o,y_o)$ on the stochastic interval $\[t_o,\z \[$\,,}
 \eea
 \eqref{eq:d169}, \eqref{eq:d167}, \eqref{si_cond} and \eqref{eq:c601} imply that
 \bea
  & \tn  & \tn \hspace{-2.5cm} \frac12 |\fra|^2   \pa^2_y \phi (\Th_r) \+ \big(D_x (\pa_y \phi (\Th_r))  \big)^T
  \dn \cd \n \si (r,\cX^{t_o,x_o}_r) \n \cd \n  \fra  \- \wh{\phi} \,  (\Th_r)   \gs \frac12 \e \> 0  \nonumber   \\
   \hb{and} \q
  \big| D_x \phi (\Th_r)\big| \big| \si(r,\cX^{t_o,x_o}_r) \big|   \+ \big|\pa_y \phi (\Th_r)\big| |\fra|
  & \tn \ls  & \tn   \big(1\+ |D_x \phi (t_o,x_o,y_o)|\big)
 \big( | \si(t_o,x_o)| \+ \sqrt{\|c(\cd)\|} \,  \d \+ \sqrt{\|c(\cd)\|} \rho(\d)(1\+|x_o|^\varpi) \big)  \nonumber  \\
  & \tn  & \tn  \;   + \big(1\+ | \pa_y \phi (t_o,x_o,y_o)|\big) |\fra| \< \infty   \label{eq:d195}
   \eea
 holds on $\[t_o,\z \[$\,.
 Applying It\^o's formula to process $\big\{ \phi (\Th_s) \big\}_{s \in [t_o,\infty)}$ then  yields that
 \bea
 \hspace{-0.7cm}
 \phi  ( \Th_\z ) \- \phi (t_o,x_o,y_o) & \dn \dn \=  &  \dn \dn
   \int_{t_o}^\z \n \Big(  \pa_t \phi (\Th_r)   \-  g (r, \cX^{t_o,x_o}_r) \pa_y \phi (\Th_r)    \+   \cL_x \phi (\Th_r)
   \+ \frac12   |\fra|^2  \pa^2_y \phi (\Th_r)
   \+     \big(D_x (\pa_y \phi (\Th_r))\big)^T \n \cd \n \si  (r, \cX^{t_o,x_o}_r) \n  \cd \n \fra \Big) dr \nonumber \\
    &  \dn \dn  &  \dn \dn  \+ \int_{t_o}^\z \n \big( (D_x \phi (\Th_r))^T \n \cd \n \si(r,\cX^{t_o,x_o}_r)
   \+ \pa_y \phi (\Th_r) \n \cd \n \fra^T \big) dW^{t_o}_r , \nonumber \\
 &  \dn \dn  \gs &  \dn \dn - \n \int_{t_o}^\z  \n   f (r, \cX^{t_o,x_o}_r)  dr
 \+ \int_{t_o}^\z  \n  \big( (D_x \phi (\Th_r))^T \n \cd \n \si(r,\cX^{t_o,x_o}_r)
   \+ \pa_y \phi (\Th_r) \n \cd \n \fra^T \big) dW^{t_o}_r , \q P_{t_o} \n -a.s. \label{eq:d173}
 \eea

 Set $ \fm_1 \df \underset{(t,x,y) \in \pa O_\d (t_o,x_o,y_o)}{\min} (\cV \- \phi) (t,x,y) \> 0 $ by \eqref{eq:d167}.
 The continuity of process $ \Th $  and \eqref{eq:d171}   show  that
   \bea  \label{eq:d179}
   P_{t_o} \big\{  \Th_\z \ins \pa O_\d (t_o,x_o,y_o)  \big\} \=
   P_{t_o} \big\{  Y^{t_o,x_o,\a^o}_s \gs y_o \- \d 
   \> 0 , \, \fa  s \ins [t_o,\z] \big\} \= 1  ,
   \eea
  the latter of which implies that
 \bea \label{eq:d181}
 \tau(t_o,x_o,\a^o)  \> \z \> t_o , \q P_{t_o}- a.s.
 \eea
 Taking expectation $E_{t_o} [\cd]$ in \eqref{eq:d173}
 and applying Theorem \ref{thm_DPP2} (2) with $\z (\a) \equiv \z$,
 we can derive from \eqref{eq:d195},  \eqref{eq:d179}      that
     \bea
     && \hspace{-0.7cm}  \phi (t_o,x_o,y_o)  \+ \fm_1   \ls    E_{t_o}  \bigg[    \phi  ( \Th_\z     )
\+ \int_t^\z \n f(r,\cX^{t_o,x_o}_r) dr    \bigg] \+ \fm_1
  \ls    E_{t_o}  \bigg[    \cV  ( \Th_\z     )
\+ \int_t^\z \n f(r,\cX^{t_o,x_o}_r) dr    \bigg] \nonumber \\
 &&   \=  E_{t_o} \n \bigg[ \b1_{\{\tau (t_o,x_o,\a^o) \le \z \}}  \cR \big(t_o,x_o, \tau (t_o,x_o,\a^o) \big)
  \+  \b1_{\{\tau (t_o,x_o,\a^o) > \z\}} \Big(  \cV  ( \Th_\z     )
\+ \n  \int_t^\z \n f(r,\cX^{t_o,x_o}_r) dr \Big)   \bigg] \label{eq:d175} \\
 && \ls  \underset{\a    \in \fA_{t_o} (y_o)   }{\sup} \,
 E_{t_o}  \bigg[ \b1_{\{\tau (t_o,x_o,\a ) \le \z \}}  \cR \big(t_o,x_o, \tau (t_o,x_o,\a ) \big)
   \+  \b1_{\{\tau (t_o,x_o,\a ) > \z\}} \Big(  \cV  \big( \z, \cX^{t_o,x_o}_\z ,  Y^{t_o,x_o,\a }_\z     \big)
 \+ \int_t^\z \n f(r,\cX^{t_o,x_o}_r) dr \Big)   \bigg]  \nonumber  \\
 && \n =  \cV (t_o,x_o,y_o)      \= \phi (t_o,x_o,y_o)  .  \nonumber
 \eea
 A contradiction appears.

 We can also employ   the first DPP (Theorem \ref{thm_DPP})
 to induce the incongruity: Denote $ \tau_o \df \tau (t_o,x_o,\a^o)$.
 By the continuity of process $Y^{t_o,x_o,\a^o}$,
 \bea \label{eq:d177}
  y_o \+ \fra^T \n \cd \n W^{t_o}_{\tau_o}
  \= \int_{t_o}^{\tau_o} \n g \big(r,\cX^{t_o,x_o}_r\big) dr  , \q   P_{t_o}   - a.s.
 \eea
 So $E_{t_o} \big[ \int_{t_o}^{\tau_o} \n g \big(r,\cX^{t_o,x_o}_r\big) dr \big] \= y_o$,
 which together with \eqref{eq:d181} shows $\tau_o \ins \wh{\cT}^{t_o}_{x_o} (y_o)$.
 On the other hand, taking conditional expectation $  E_{t_o} \big[ \cd  \big| \ol{\cF}^{t_o}_\z \big] $
 in \eqref{eq:d177}, one can deduce from \eqref{eq:d181} and the optional sampling theorem  that
 \beas
 Y^{t_o,x_o,\a^o}_\z & \tn \= & \tn  y_o \+ \fra^T \n \cd \n W^{t_o}_\z
 \- \int_{t_o}^\z \n g \big(r,\cX^{t_o,x_o}_r\big) dr
 \= E_{t_o} \big[ y_o \+ \fra^T \n \cd \n W^{t_o}_{\tau_o}  \big| \ol{\cF}^{t_o}_\z \big]
 \- \int_{t_o}^\z \n g \big(r,\cX^{t_o,x_o}_r\big) dr \\
  & \tn  \=  & \tn   E_{t_o} \Big[ \int_{t_o}^{\tau_o} \n g \big(r,\cX^{t_o,x_o}_r\big) dr \Big| \ol{\cF}^{t_o}_\z\Big]
 \- \int_{t_o}^\z \n g \big(r,\cX^{t_o,x_o}_r\big) dr
 \= \cY^{t_o,x_o,\tau_o}_\z , \q P_{t_o}  \n -   a.s.
 \eeas
 Then we can apply Theorem \ref{thm_DPP} (2) with $\z (\a) \equiv \z$ to continue the deduction  in \eqref{eq:d175}
 \beas
     \phi (t_o,x_o,y_o) \+ \fm_1   & \tn  \ls  & \tn
      E_{t_o} \n \bigg[ \b1_{\{\tau_o \le \z \}}  \cR \big(t_o,x_o, \tau_o \big)
  \+  \b1_{\{\tau_o > \z\}} \Big(  \cV  ( \z, \cX^{t_o,x_o}_\z ,  \cY^{t_o,x_o,\tau_o}_\z     )
\+ \n  \int_t^\z \n f(r,\cX^{t_o,x_o}_r) dr \Big)   \bigg]   \\
  & \tn  \ls  & \tn  \underset{\tau   \in \wh{\cT}^{t_o}_{x_o} (y_o)   }{\sup} \,
 E_{t_o}  \bigg[ \b1_{\{\tau  \le \z \}}  \cR \big(t_o,x_o, \tau_o \big)
   \+  \b1_{\{\tau  > \z\}} \Big(  \cV  \big( \z, \cX^{t_o,x_o}_\z ,  \cY^{t_o,x_o,\tau}_\z     \big)
 \+ \int_t^\z \n f(r,\cX^{t_o,x_o}_r) dr \Big)   \bigg]     \\
  & \tn  =  & \tn  \cV (t_o,x_o,y_o)      \= \phi (t_o,x_o,y_o)  .
 \eeas
 The contradiction recurs.
 Therefore,  $\cV$  is a viscosity supersolution   of \eqref{eq:PDE}.

  \no {\bf 2)} Next, we demonstrate that  $\cV$  is also a viscosity subsolution   of \eqref{eq:PDE2}.

 Let $(t_o,x_o,y_o ) \ins (0,\infty) \ti \hR^l \ti (0,\infty)   $
 and  let   $ \vf \ins C^{1,2,2}\big([0,\infty) \ti \hR^l \ti [0,\infty)\big)$
 such that    $\cV \- \vf$ attains a strict local  maximum $0$  at $(t_o,x_o,y_o)$.
 So there exists a
    $\l_o \ins  \big(0,    t_o  \ld y_o \big) $ such that for any $(t,x,y) \in   O_{\l_o} (t_o,x_o,y_o)
     \big\backslash \big\{ (t_o,x_o,y_o) \big\}$
   \bea \label{eq:d167b}
   (\cV \- \vf) (t,x,y) \<  (\cV \- \vf) (t_o,x_o,y_o) \= 0   \hb{ and }
     \big|D_x \vf (t,x,y) \- D_x \vf (t_o,x_o,y_o)  \big| \ve
     \big| \pa_y \vf (t,x,y) \- \pa_y \vf (t_o,x_o,y_o)\big| \< 1 . \q
   \eea

  Similar to \eqref{eq:d193},  $ \wh{\vf} \,  (t,x,y) \df  -  \pa_t \vf (t,x,y)  - \cL_x   \vf (t,x,y)
   + g(t,x)    \pa_y \vf (t,x,y)   - f(t,x) $, $ \fa (t,x,y) \ins [0,\infty) \ti \hR^l \ti [0,\infty) $
  defines a continuous function.
 If $ \ol{\cH} \vf (t_o,  x_o,   y_o)   \= \infty$, then
 $\wh{\vf} \,  (t_o,x_o,y_o) \- \ol{\cH} \vf (t_o,x_o,y_o)  \ls 0$
 holds automatically.

 So let us just consider the case $ \ol{\cH} \vf (t_o,   x_o,   y_o)   \< \infty$.
 By \eqref{eq:d191}, there exists   $\wt{\l}_o \ins (0,\l_o)$ such that
 $ \cH  \vf (t,x,y) \ls \ol{\cH} \vf (t_o,   x_o,   y_o)\+1 \< \infty $
 and thus $ \pa^2_y \vf (t,x,y) \ls 0 $, $ \fa  (t,x,y) \ins O_{\wt{\l}_o} (t_o,   x_o,   y_o)$.
 If one had  $\pa_y \vf (t_o,x_o,y_o) \ls 0 $, \eqref{eq:c405} and \eqref{eq:d167b}   would imply that
 \beas
 \vf (t_o,x_o,y ) & \tn \= & \tn  \vf (t_o,x_o,y_o ) \+ \int_{y_o}^y \vf_y (t_o,x_o,s) ds
 \=  \vf (t_o,x_o,y_o ) \+ (y \- y_o) \n \cd \n \pa_y \vf (t_o,x_o,y_o)
 \+ \int_{y_o}^y \n \int_{y_o}^s \vf^2_y (t_o,x_o,r) dr ds \\
  & \tn \ls & \tn  \vf (t_o,x_o,y_o ) \=  \cV (t_o,x_o,y_o ) \ls \cV (t_o,x_o,y ),
  \q  \fa y \ins \big( y_o, y_o\+\wt{\l}_o \big) .
 \eeas
 which contradicts with the strict local  maximum of $\cV \- \vf$  at $(t_o,x_o,y_o)$. Hence we must have
 \bea \label{eq:d209}
  \pa_y \vf (t_o,x_o,y_o) \> 0  .
 \eea

 To draw a contradiction, we assume  that
 \beas
 \eps \df \wh{\vf} \,  (t_o,x_o,y_o) \- \ol{\cH} \vf (t_o,x_o,y_o) \> 0 .
 \eeas
 According to \eqref{eq:d209} and the continuity of $\wh{\vf}$, there exists $\l \ins (0,\l_o)$ such that
 for any $(t,x,y) \ins \ol{O}_\l (t_o,x_o,y_o)$
 \bea \label{eq:d169b}
 \pa_y \vf (t,x,y) \gs 0 \q \hb{and} \q
 \cH \vf  (t,x,y) \ls \ol{\cH} \vf (t_o,x_o,y_o) \+ \eps/2  \= \wh{\vf} \,  (t_o,x_o,y_o) \-  \eps/2
 \ls   \wh{\vf} (t,x,y)    .
 \eea

  Fix $\a   \ins \fA_{t_o} (y_o)$, so $\a \= M^\a\-K^\a$ for some $(M^\a,K^\a) \ins \hM_{t_o} \ti \hK_{t_o}$.
 In light of the martingale representation theorem, one can find
   $ \fq^\a \ins \hH^{2,{\rm loc}}_{t_o}$ such that
 \bea \label{eq:d197}
 P_{t_o} \Big\{ \int_{t_o}^s |\fq^\a_r|^2 dr \< \infty, \,  \fa s \ins [t_o,\infty) \Big\}
 \= P_{t_o} \Big\{ M^\a_s \= \int_{t_o}^s (\fq^\a_r)^T d W^{t_o}_r  , \, \fa s \ins [t_o,\infty) \Big\} \= 1 .
 \eea
 As $ \Th^\a_s \df \big(s, \cX^{t_o,x_o}_s  , Y^{t_o,x_o,\a}_s   \big) $   $s \ins [t_o,\infty)$
 are $\ol{\bF}^{t_o}-$adapted continuous processes   with $ \Th^\a_{t_o} \= (t_o,x_o,y_o)$, $P_{t_o}-$a.s.,
 \bea \label{eq:d217}
 \z^\a   \df   \inf\big\{ s \ins [t_o, \infty) \n : \Th^\a_s \n \notin \n  \ol{O}_\l (t_o,x_o,y_o )  \big\}
    \eea
 defines  an $\ol{\bF}^{t_o}-$stopping time with $t_o \< \z^\a    \ls t_o \+ \l$, $P_{t_o}-$a.s.
 The continuity of processes $\Th^\a$  implies that  $P_{t_o}-$a.s.
 \bea \label{eq:d171b}
  \Th^\a_s \ins   \ol{O}_\l (t_o,x_o,y_o), \q  \fa s \ins  [t_o,\z^\a   ]  .
 \eea

 Similar to \eqref{eq:d195}, we can deduce from
 \eqref{eq:d171b},   \eqref{eq:d169b} and \eqref{eq:d167b}   that for $P_{t_o}-$a.s.
 \bea
   &    &  \hspace{-7mm} \pa_y \vf (\Th^\a_r) \gs 0, ~
   \frac12 \big|\fq^\a_r\big|^2   \pa^2_y \vf (\Th^\a_r) \+ \big(D_x (\pa_y \vf (\Th^\a_r))  \big)^T
  \dn \cd \n \si (r,\cX^{t_o,x_o}_r) \n \cd \n \fq^\a_r  \- \wh{\vf} \,  (\Th^\a_r)
  \ls \cH \vf (\Th^\a_r)  \-  \wh{\vf}  \,  (\Th^\a_r)  \ls 0   \q \hb{and}  \label{eq:d199}  \\
  &&  \hspace{-7mm} \big| D_x \vf (\Th^\a_r)\big| \big| \si(r,\cX^{t_o,x_o}_r) \big|   \+ \big|\pa_y \vf (\Th^\a_r)\big| \big|\fq^\a_r\big|
  \ls     \big(1\+ |D_x \vf (t_o,x_o,y_o)|\big)
 \big( | \si(t_o,x_o)| \+ \sqrt{\|c(\cd)\|} \,  \d \+ \sqrt{\|c(\cd)\|} \rho(\d)(1\+|x_o|^\varpi) \big) \nonumber \\
   &    &  \hspace{6.6cm}      + \big(1\+ | \pa_y \vf (t_o,x_o,y_o)|\big) \big|\fq^\a_r\big| \< \infty ,
     \q \fa s \ins  [t_o,\z^\a   ] . \label{eq:d201}
   \eea

 Let $n \ins \hN$ and define $\z^\a_n \df  \inf\big\{ s \ins [t_o, \infty) \n : \int_{t_o}^s |\fq^\a_r|^2 dr \> n    \big\} \ld \z^\a \ins \ol{\cT}^{t_o}$.
 Applying It\^o's formula to process $\big\{ \vf (\Th^\a_s) \big\}_{s \in [t_o,\infty)}$,
  and using \eqref{eq:d199} yield     that
 \beas
 \q  && \hspace{-1.2cm}
  \vf  \big( \Th^\a_{\z^\a_n} \big) \- \vf (t_o,x_o,y_o) \nonumber \\
 & & \=
   \int_{t_o}^{\z^\a_n} \n \Big(  \pa_t \vf (\Th^\a_r)   \-  g (r, \cX^{t_o,x_o}_r) \pa_y \vf (\Th^\a_r)    \+   \cL_x \vf (\Th^\a_r)
   \+ \frac12   \big|\fq^\a_r\big|^2  \pa^2_y \vf (\Th^\a_r)
   \+  \big(D_x (\pa_y \vf (\Th^\a_r))\big)^T \n \cd \n \si  (r, \cX^{t_o,x_o}_r) \n  \cd \n \fq^\a_r \Big) dr \nonumber \\
    &  &    \q \- \int_{t_o}^{\z^\a_n} \n \pa_y \vf (\Th^\a_r) d K^\a_r
     \+ \int_{t_o}^{\z^\a_n} \n \big( (D_x \vf (\Th^\a_r))^T \n \cd \n \si(r,\cX^{t_o,x_o}_r)
   \+ \pa_y \vf (\Th^\a_r) \n \cd \n (\fq^\a_r)^T \big) dW^{t_o}_r , \nonumber \\
 &  & \ls   - \int_{t_o}^{\z^\a_n}  \n   f (r, \cX^{t_o,x_o}_r)  dr
 \+ \int_{t_o}^{\z^\a_n}  \n  \big( (D_x \vf (\Th^\a_r))^T \n \cd \n \si(r,\cX^{t_o,x_o}_r)
   \+ \pa_y \vf (\Th^\a_r) \n \cd \n (\fq^\a_r)^T \big) dW^{t_o}_r , \q P_{t_o} \n -a.s.
 \eeas
 Taking expectation $E_{t_o} [\cd]$, we see from \eqref{eq:d201} that
 $ \vf (t_o,x_o,y_o) \gs E_{t_o} \Big[ \vf  \big( \Th^\a_{\z^\a_n} \big)
 \+ \int_{t_o}^{\z^\a_n}   f (r, \cX^{t_o,x_o}_r)  dr \Big] $.
 Since \eqref{eq:d171b} and the continuity of $f$ show  that
 $  \big| \vf  \big( \Th^\a_{\z^\a_n} \big)  \big|
   \+ \int_{t_o}^{\z^\a_n}  \n \big|  f (r, \cX^{t_o,x_o}_r)  \big| dr
    \ls  \underset{(t,x,y) \in  \ol{O}_\l (t_o,x_o,y_o)}{\max}| \vf (t,x,y)|
    \+ \l \underset{(t,x ) \in  \ol{O}_\l (t_o,x_o )}{\max}|  f (t,x )| $
   and since $\lmtu{n \to \infty} \z^\a_n \= \z^\a$, $P_{t_o}-$a.s. by \eqref{eq:d197},
   we can derive from the dominated convergence theorem   that
   \bea \label{eq:d211}
     \vf (t_o,x_o,y_o) \gs \lmt{n \to \infty} E_{t_o} \n \bigg[ \vf  \big( \Th^\a_{\z^\a_n} \big)
     \+ \int_{t_o}^{\z^\a_n} \n  f \big(r, \cX^{t_o,x_o}_r\big)  dr \bigg]
   \= E_{t_o} \n \bigg[ \vf  \big( \Th^\a_{\z^\a } \big)
   \+ \int_{t_o}^{\z^\a} \n   f \big(r, \cX^{t_o,x_o}_r\big)  dr \bigg] .
   \eea

 Set $ \fm_2 \df \underset{(t,x,y) \in \pa O_\l (t_o,x_o,y_o)}{\min} (\vf \- \cV  ) (t,x,y) \> 0 $ by \eqref{eq:d167b}.
 The continuity of $\Th^\a$ and \eqref{eq:d171b} show   that
 \beas
 P_{t_o} \big\{   \Th^\a_{\z^\a}  \ins \pa O_\l (t_o,x_o,y_o) \big\}
 \=  P_{t_o} \big\{  Y^{t_o,x_o,\a }_s \gs y_o \- \l 
   \> 0 , \, \fa  s \ins [t_o,\z^\a] \big\} \= 1 .
 \eeas
  The latter of which implies that  $ \tau(t_o,x_o,\a )  \> \z^\a \> t_o $,  $ P_{t_o}- $a.s.,
 which together with  \eqref{eq:d211} leads to that
     \bea
 \q  &&  \hspace{-2cm}   \vf (t_o,x_o,y_o) \- \fm_2
  \gs    E_{t_o}  \bigg[    \vf  ( \Th^\a_{\z^\a}     ) \- \fm_2
\+ \int_t^{\z^\a} \n f(r,\cX^{t_o,x_o}_r) dr    \bigg]
\gs   E_{t_o}  \bigg[    \cV  ( \Th^\a_{\z^\a}     )
\+ \int_t^{\z^\a} \n f(r,\cX^{t_o,x_o}_r) dr    \bigg] \nonumber \\
  &&  \gs   E_{t_o} \n \bigg[ \b1_{\{\tau (t_o,x_o,\a) \le \z^\a \}}  \cR \big(t_o,x_o, \tau (t_o,x_o,\a) \big)
  \+  \b1_{\{\tau (t_o,x_o,\a) > \z^\a\}} \bigg(  \cV  \big( \Th^\a_{\z^\a}     \big)
\+ \n  \int_t^{\z^\a} \n f(r,\cX^{t_o,x_o}_r) dr \bigg)   \bigg] .   \label{eq:d219}
 \eea
 Taking supremum over $\a   \in \fA_{t_o} (y_o)$, we can deduce from Theorem \ref{thm_DPP2} that
 \beas
 \hspace{-3mm}
  \vf (t_o,x_o,y_o)  \- \fm_2
 & \dn \dn \gs  & \dn \dn   \underset{\a   \in \fA_{t_o} (y_o)   }{\sup} \,
 E_{t_o}  \bigg[ \b1_{\{\tau (t_o,x_o,\a) \le \z^\a \}}  \cR \big(t_o,x_o, \tau (t_o,x_o,\a) \big)
   \+  \b1_{\{\tau (t_o,x_o,\a) > \z^\a\}} \bigg(   \cV  \big( \Th^\a_{\z^\a}     \big)
 \+ \int_t^{\z^\a} \n f(r,\cX^{t_o,x_o}_r) dr \bigg)   \bigg]   \\
  & \dn \dn  \=  & \dn \dn  \cV (t_o,x_o,y_o)      \= \vf (t_o,x_o,y_o) .
 \eeas
 A contradiction appears.

 We can also use  the first DPP (Theorem \ref{thm_DPP})
 to get the incongruity: Let $\tau \ins \wh{\cT}^{t_o}_{x_o} (y_o)$.
 By Proposition \ref{prop_surm},    $\tau \= \tau(t_o,x_o,\a)$ for some $\a \ins \fA_{t_o} (y_o)$.
 Let $\z^\a$ be the $\ol{\bF}^{t_o}-$stopping time defined in \eqref{eq:d217}.
 Similar to  \eqref{eq:d214}, one has
 $   Y^{t_o,x_o,\a}_{\tau \land \z^\a}    \gs   \cY^{t_o,x_o,\tau}_{\tau \land \z^\a} $, $P_o-$a.s.
 It follows from \eqref{eq:d219} that
     \beas
      \vf (t_o,x_o,y_o) \- \fm_2
 & \tn \gs & \tn  E_{t_o} \n \bigg[ \b1_{\{\tau   \le \z^\a \}}  \cR  (t_o,x_o, \tau   )
  \+  \b1_{\{\tau   > \z^\a\}} \bigg(  \cV  \big( \z^\a, \cX^{t_o,x_o}_{\z^\a} ,   Y^{t_o,x_o,\a}_{\z^\a}    \big)
\+ \n  \int_t^{\z^\a} \n f(r,\cX^{t_o,x_o}_r) dr \bigg)   \bigg] \\
& \tn \gs & \tn E_{t_o} \n \bigg[ \b1_{\{\tau   \le \z^\a \}}  \cR  (t_o,x_o, \tau   )
  \+  \b1_{\{\tau   > \z^\a\}} \bigg(  \cV  \big( \z^\a, \cX^{t_o,x_o}_{\z^\a} ,   \cY^{t_o,x_o,\tau}_{\z^\a}    \big)
\+ \n  \int_t^{\z^\a} \n f(r,\cX^{t_o,x_o}_r) dr \bigg)   \bigg].
 \eeas
 Taking supremum over $\tau \ins \wh{\cT}^{t_o}_{x_o} (y_o)$ and applying
   Theorem \ref{thm_DPP} (2)  yield that
 \beas
    \vf (t_o,x_o,y_o) \- \fm_2   & \tn  \gs  & \tn
 \underset{\tau \ins \wh{\cT}^{t_o}_{x_o} (y_o)}{\sup}  E_{t_o} \n \bigg[ \b1_{\{\tau   \le \z^\a \}}  \cR  (t_o,x_o, \tau   )
  \+  \b1_{\{\tau   > \z^\a\}} \bigg(  \cV  \big( \z^\a, \cX^{t_o,x_o}_{\z^\a} ,   \cY^{t_o,x_o,\tau}_{\z^\a}    \big)
\+ \n  \int_t^{\z^\a} \n f(r,\cX^{t_o,x_o}_r) dr \bigg)   \bigg]    \\
  & \tn  =  & \tn  \cV (t_o,x_o,y_o)      \= \vf (t_o,x_o,y_o)  .
 \eeas
 The contradiction appears   again.  \qed

\appendix
\renewcommand{\thesection}{A}
\refstepcounter{section}
\makeatletter
\renewcommand{\theequation}{\thesection.\@arabic\c@equation}
\makeatother

\section{Appendix}

  \begin{lemm} \label{lem_countable_generate1}

  Let $t \ins [0,\infty)$.
  \(1\) The sigma$-$field $\cF^t$ satisfies
 $ \sB (\O^t) \=  \si \big(W^t_s; s \ins [t,\infty) \big) \= \si \Big(   \underset{s \in [t,\infty)}{\cup}  \cF^t_s  \Big)  $.

\no   \(2\) For any $s \ins [t,\infty]$, the sigma$-$field $\cF^t_s$ can be countably generated by
 $   \sC^t_s \df    \Big\{  \underset{i =1}{\overset{m}{\cap}}  ( W^t_{t_i} )^{-1}  \big( O_{\d_i} (x_i) \big) \n  :
     \,  m \ins  \hN ,\,    t_i \ins \hQ_+ \cp \{t\}
   \hb{ with } t \ls  t_1 \ls \cds \ls t_m \ls s  ,\, x_i \ins \hQ^d , \, \d_i \ins \hQ_+ \Big\}  $.

   \end{lemm}

       \no {\bf Proof:  1a)} Let $\o \ins \O^t $ and $\d \ins (0,\infty)$.
 For any $ n \ins \hN$ with $n \> 1/\d$, since all paths in $\O^t$ are continuous,
  we can deduce   that
 \beas
 \ol{O}_{\d-1/n} (\o) 
 & \tn \= & \tn  \{\o' \ins \O^t \n :  |\o'(s)\-\o (s)| \ls  \d \- 1/n , \, \fa s \ins [t,\infty) \}
 \= \{\o' \ins \O^t \n :  |\o'(s)\-\o (s)| \ls  \d \- 1/n , \;\fa s \ins [t,\infty) \Cp \hQ  \} \\
 & \tn  \= & \tn  \underset{s \in [t,\infty) \cap \hQ}{\cap} \{\o' \ins \O^t \n : W^t_s(\o') \ins \ol{O}_{\d - 1/n}(\o(s)) \}
 \= \underset{s \in [t,\infty) \cap \hQ}{\cap} (W^t_s)^{-1}  \big( \ol{O}_{\d - 1/n}(\o(s)) \big)
 \ins \si \big(W^t_s; s \ins [t,\infty) \big)  .
 \eeas
 It follows that $O_\d (\o) = \underset{n \in \hN}{\cup} \ol{O}_{\d-1/n} (\o)
 \ins \si \big(W^t_s; s \ins [t,\infty) \big)  $.
 As $\sB (\O^t)$ is generated by   open sets $\{ O_\d (\o) \n : \o \ins \O^t , \d \ins (0,\infty) \} $,
 one thus has $ \sB (\O^t) \sb \si \big(W^t_s; s \ins [t,\infty) \big) $.

  Next,  let $ s \ins [t,\infty) $,      $x \ins \hR^d$ and $\d \ins (0,\infty)$.
   Given     $\o \ins (W^t_s)^{-1} (O_\d (x))$,   set $\l \= \l(s,x,\o) \df \d \- |W^t_s(\o) \- x| \> 0$.
  Since
   \beas
 |W^t_s(\o') \-x| \ls   |\o'(s) \- \o(s)| \+ |\o(s) \- x|
 \ls \|\o'\-\o\|_t \+  |\o(s) \- x|  \< \l \+ |\o(s) \- x| \= \d , \q   \fa \o' \ins O_\l(\o) ,
   \eeas
   we see that   $ O_\l(\o) \sb (W^t_s)^{-1} \big( O_\d (x) \big) $
   and thus    $(W^t_s)^{-1} \big(O_\d (x)\big)$ is an open set under the uniform norm $\| \cd \|_t $.
   Then
 $  O_\d (x)  \ins \L_s \df \{\cE \sb \hR^d  \n : (W^t_s)^{-1} (\cE) \ins \sB (\O^t) \}$,
    which is a sigma$-$field of $\hR^d$.
   As $\sB(\hR^d)$ is generated by   open sets $\{O_\d (x) \n : x \ins \hR^d, \d \ins (0,\infty) \}$,
   one has $ \sB(\hR^d) \sb \L_s $, which implies that
   $ \si \big(W^t_s; s \ins [t,\infty) \big) \=
   \si \big\{ (W^t_s)^{-1} (\cE) \n : s \ins [t,\infty),\; \cE \ins \sB(\hR^d) \big\} \sb \sB (\O^t) $.

  \no {\bf 1b)} Clearly,  $\cF^t_s \= \si \big(W^t_r; r  \ins  [t,s] \big) \sb \si \big(W^t_r; r  \ins  [t,\infty) \big) $, $\fa s \ins [t,\infty)$. It follows that
     $ \si \Big(\underset{s \in [t,\infty)}{\cup} \cF^t_s \Big) \sb \si \big(W^t_s; s  \ins  [t,\infty) \big) $.
   On the other hand, since
   $ (W^t_s)^{-1} (\cE) \ins \cF^t_s \sb \si \Big(\underset{r \in [t,\infty)}{\cup} \cF^t_r \Big) $
   for any $ s \ins [t,\infty) $ and $ \cE \ins \sB(\hR^d) $,
   we have    $ \si \big(W^t_s; s \ins [t,\infty) \big) \=
   \si \big\{ (W^t_s)^{-1} (\cE) \n : s \ins [t,\infty),\; \cE \ins \sB(\hR^d) \big\}
   \sb \si \Big(\underset{s \in [t,\infty)}{\cup} \cF^t_s \Big) $.

  \if{0}

  \no {\bf 1c)} Clearly,
  $  \si^t \df \si \Big(   \underset{s \in [t,\infty)}{\cup}
    \big( \sB([t,s]) \oti  \cF^t_s \big) \Big) \sb \sB [t,\infty) \oti \cF^t $.
    To see the inverse inclusion, we
    let $s \ins [t,\infty)$ and set
  $\wt{\L}_s \df \big\{A \sb \O^t \n : \cE  \ti A \ins  \si^t  , \fa \cE  \ins \sB([t,s]) \big\}$,
  which is clearly a sigma$-$field of $\O^t$.

  Given $r \ins [t,\infty) $,
  since $\cE \ti A \ins \sB([t,s \ve r]) \oti \cF^t_{s \vee r} \sb \si^t $
  for any $\cE  \ins \sB([t,s])$ and $A  \ins \cF^t_r $, we see that
$\cF^t_r \sb \wt{\L}_s $. Then Part 1a and Part 1b show that
$ \cF^t \= \si \Big( \underset{r \in [t,\infty)}{\cup}    \cF^t_r  \Big) \sb \wt{\L}_s $.
Namely, $ \si^t $ contains all generating sets of $ \sB([t,s]) \oti \cF $. It follows that
  $\sB([t,s]) \oti \cF \sb \si^t$.

  Let $\cE \ins [t,\infty) $ and $ A \ins \cF $. For any $n \ins \hN$,
    as  $\cE \Cp [t,t\+n] \ins \sB([t,t\+n])$,  one has
 $(\cE \ti A) \Cp ([t,t\+n] \ti \O^t) \= \big(\cE \Cp [t,t\+n]\big) \ti A \ins \sB([t,t\+n]) \oti \cF \sb \si^t$.
 Taking union over $n \ins  \hN$ yields that
 $ \cE \ti A \= \underset{n \in \hN}{\cup} \big( (\cE \ti A) \Cp ([t,t\+n] \ti \O^t) \big)
 \ins \si^t $, which shows that $\si^t$ contains all generating sets of $ \sB [t,\infty) \oti \cF^t $.
 It follows that $\sB [t,\infty) \oti \cF^t \sb \si^t$.

 \fi

 \no {\bf 2)}
  Fix $s \ins [t,\infty]$. Define $[t,s\ran \df [t,s] $ if $s \< \infty$ and $ [t,s\ran \df [t,\infty) $ if $s \= \infty$.
  Let   $r \ins [t,s\ran   $ and let
  $\{r_i\}_{i \in  \hN} \sb \{t\} \cup \big( (t,r) \cap  \hQ_+ \big)  $  with $\lmtu{i \to \infty} r_i \= r $.
    For any $x \ins \hQ^d$ and $ \d \ins \hQ_+$,  the continuity of paths in $\O^t$ implies that
 $   
  (W^t_r)^{-1} \big(O_\d(x)\big) \= \underset{n \= \lceil 2/\d   \rceil   }{\overset{\infty}{\cup}}
 \underset{m \in  \hN }{\cup}  \underset{i  > m}{\cap} \Big(  (W^t_{r_i})^{-1}  \big(\ol{O}_{\d-\frac{1}{n}}(x)\big) \Big)
 \ins \si  ( \sC^t_s   ) $.
 Thus $O_\d (x) \ins \wh{\L}_r  \df   \big\{ \cE  \sb  \hR^d \n :  (W^t_r)^{-1}  (\cE  )
 \ins  \si (\sC^t_s) \big\} $, which is a sigma$-$field of $ \hR^d $.
 Since $\sB(\hR^d)$ can also be generated by   $\{O_\d (x) \n : x \ins \hQ^d , \,  \d \ins \hQ_+ \}$,
 we see that      $\sB(\hR^d)  \sb \wh{\L}_r$.
 It follows that  $ \cF^t_s \= \si \big(W^t_r; r \ins [t, s \ran \big)
 \=  \si \big\{ (W^t_r)^{-1} (\cE )  \n : r  \ins  [t, s \ran , \cE  \ins  \sB(\hR^d) \big\}
   \sb  \si  ( \sC^t_s    )$.
     On the other hand, it is clear that    $   \si (\sC^t_s)  \sb  \si   \big\{  (W^t_r)^{-1}  (\cE ) :
 r \ins [t, s \ran , \cE \ins \sB(\hR^d) \big\}
 \= \si \big(W^t_r; r \ins [t,s \ran  \big) \= \cF^t_s $. \qed

    \begin{lemm}  \label{lem_shift_inverse} 
 Let $0 \ls t \ls s \< \infty  $.

 \no   \(1\) The  mapping  $\Pi^t_s$  is $\cF^t_r \big/ \cF^s_r-$measurable for any $r \ins [s,\infty]$.
 Then for each $\bF^s-$stopping time $\tau  $,
 $\tau (\Pi^t_s) $ is a $ \bF^t-$stopping time with values in $[s,\infty]$.

 \no   \(2\) The law of $ \Pi^t_s $ under $P_t$ is    $ P_s $: i.e.,
  $ P_t   \circ    (\Pi^t_s)^{-1} \big(\wt{A}\,\big)   \=  P_s \big(\wt{A}\,\big) $, $  \fa   \wt{A} \ins  \cF^s  $.

\end{lemm}

     \no {\bf Proof:}
   {\bf 1)}
    \if{0}
     Let $r \ins [s,\infty)$.   For any $s' \ins  [s,r] $ and $ \cE \ins  \sB(\hR^d)$,
   \bea    \label{eq:g313}
       (\Pi^t_s)^{-1}   \big(   (W^s_{s'})^{-1}   (\cE) \big)   \=
    \big\{\o \ins \O^t \n : W^s_{s'} \big((\Pi^t_s) (\o) \big)  \ins  \cE  \big\}
     \=  \big\{\o \ins \O^t \n : \o(s')  \- \o(s) \ins   \cE \big\}
     \=   (W^t_{s'}    \-  W^t_s)^{-1} (\cE)  \ins   \cF^t_r . ~
    \eea
  So all   generating sets of $\cF^s_r$ belong to $\L_r \df \big\{ \wt{A} \sb \O^s \n : (\Pi^t_s)^{-1}  \big(\wt{A}\,\big) \ins \cF^t_r \big\}$,
  which is clearly a sigma$-$field of $\O^s$. It follows that $\cF^s_r \sb  \L_r$,
  or  $(\Pi^t_s)^{-1}  \big(\wt{A}\,\big) \ins \cF^t_r$ for  any $\wt{A} \ins \cF^s_r$.
     \fi
  For any $r \ins [s,\infty)$, an analogy to Lemma A.1 of \cite{ROSVU} shows that
  $(\Pi^t_s)^{-1}  \big(\wt{A}\,\big) \ins \cF^t_r \sb \cF^t $ for  any $\wt{A} \ins \cF^s_r$.
  Set $ \wt{\L}  \df \big\{ \wt{A} \sb \O^s \n : (\Pi^t_s)^{-1}  \big(\wt{A}\,\big) \ins \cF^t  \big\}$,
  which is 
  a sigma$-$field of $\O^s$. As $\cF^s_r   \sb \wt{\L} $ for any $r \ins [s,\infty)$,
  we see from Lemma \ref{lem_countable_generate1} (1) that
   $\cF^s \= \si \Big(   \underset{r \in [s,\infty)}{\cup}    \cF^s_r  \Big)   \sb \wt{\L} $, i.e.,
   $(\Pi^t_s)^{-1}  \big(\wt{A}\,\big) \ins \cF^t $ for  any $\wt{A} \ins \cF^s $.

 Let $\tau  $ be an $\bF^s-$stopping time.  For any $r \ins [s,\infty)$, since $\{\tau \ls r \} \ins \cF^s_r$,
 we see   that $\big\{ \tau (\Pi^t_s) \ls r \big\} \= (\Pi^t_s)^{-1} \big( \{\tau \ls r \} \big) \ins \cF^t_r $\,.
 Thus  $\tau (\Pi^t_s) $ is a $\bF^t-$stopping time with values in $[s,\infty]$.

 \no {\bf 2)}   Next, let us demonstrate that  the induced probability $\wt{P} \df P_t  \n \circ \n  (\Pi^t_s)^{-1}$
  equals to $  P_s  $ on $\cF^s$.
   Since the Wiener measure $P_s$ on $ ( \O^s, \cF^s )$ is unique (see e.g. Proposition I.3.3 of \cite{revuz_yor}),
    it suffices to show that  the canonical process $W^s $
 is a Brownian motion on $\O^s$ under $\wt{P}$:
  Let $ s \ls r \< r' \< \infty$. For any $\cE \ins  \sB(\hR^d)$,
 one can deduce that
 \bea
   (\Pi^t_s)^{-1} \big(\big(W^s_{r'} \- W^s_r\big)^{-1}(\cE)\big)
  & \tn \=   & \tn  \big\{\o \ins \O^t:  W^s_{r'} \big((\Pi^t_s) (\o) \big) \- W^s_r \big((\Pi^t_s) (\o) \big)  \ins  \cE  \big\} \nonumber  \\
   & \tn  \=   & \tn  \big\{\o \ins \O^t: \o(r')  \-  \o(s)  \- \big( \o(r)  \-  \o(s) \big) \ins   \cE \big\}
        \=     (W^t_{r'} \-  W^t_r)^{-1} (\cE)  .  \label{eq:g319}
        \eea
         So $  \wt{P}\big(  \big(W^s_{r'} \- W^s_r\big)^{-1}(\cE) \big)
 \= P_t \big((W^t_{r'} \-  W^t_r)^{-1} (\cE) \big)  $,
   which shows that   the distribution of $W^s_{r'} \- W^s_r$ under $ \wt{P}$ is the same as that of
   $W^t_{r'} \-  W^t_r$ under $P_t $ (a $d-$dimensional normal distribution with mean $0$ and variance matrix
    $  (r' \- r)I_{d \times d} $).

     On the other hand,  for any $\wt{A} \ins \cF^s_r$,
     since   $ (\Pi^t_s)^{-1} \big(\wt{A}\,\big) \ins \cF^t_r$ is independent of
       $W^t_{r'} \-  W^t_r$ under $P_t$,
    \eqref{eq:g319}     implies  that
           \beas
 && \hspace{-1cm} \wt{P}\big(\wt{A} \cap \big(W^s_{r'} \- W^s_r\big)^{-1}(\cE) \big)
   \=     P_t \Big((\Pi^t_s)^{-1} \big(\wt{A}\,\big) \cap  (\Pi^t_s)^{-1}
 \big(\big(W^s_{r'} \- W^s_r\big)^{-1}(\cE)\big)  \Big) \\
 &  &     \= P_t \big((\Pi^t_s)^{-1} \big(\wt{A}\,\big) \big) \cd  P_t \Big(   (\Pi^t_s)^{-1} \big(\big(W^s_{r'} \- W^s_r\big)^{-1}(\cE)\big) \Big)
   \= \wt{P} \big(\wt{A}\,\big) \cd \wt{P} \Big( \big(W^s_{r'} \- W^s_r\big)^{-1}(\cE) \Big) , \q \fa \cE \ins \sB(\hR^d) ,
 \eeas
  which shows that $ W^s_{r'} \- W^s_r$ is independent of $ \cF^s_r $ under $\wt{P}$.
  Hence,  $W^s $  is a Brownian motion on $\O^s$ under $\wt{P}$. \qed

We   have the following extension of Lemma \ref{lem_shift_inverse}.

\begin{lemm} \label{lem_shift_inverse1b}
  Let $0 \ls t \ls s \< \infty  $.

\no \(1\)    For  any $P_s-$null set $\wt{\cN}$, $ (\Pi^t_s)^{-1} \big(\wt{\cN}\big)$ is   a $P_t-$null set.

\no \(2\) For any $r \ins [s,\infty]$, the  mapping  $ \Pi^t_s $  is $\ol{\cF}^t_r \big/ \ol{\cF}^s_r-$measurable.
Then for each $\tau \ins \ol{\cT}^s$, $\tau (\Pi^t_s) $ is a $\ol{\cT}^t-$stopping time with values in $[s,\infty]$.

 \no \(3\)    $ P_t   \circ    (\Pi^t_s)^{-1} \big(\wt{A}\,\big)   \=  P_s \big(\wt{A}\,\big)  $
 holds for any $    \wt{A} \ins  \ol{\cF}^s  $.

\end{lemm}

    \no {\bf Proof: 1)}      Let     $ \wt{\cN} \ins \sN^s   $, so there exists
    an $\wt{A} \ins \cF^s$ such that $\wt{\cN} \sb \wt{A}$ and  $   P_s    \big(\wt{A}\,\big) \= 0$.
      Lemma \ref{lem_shift_inverse} implies that
    $ (\Pi^t_s)^{-1} \big(\wt{A}\,\big) \ins \cF^t    $
    and that $  P_t   \big( (\Pi^t_s)^{-1} \big(\wt{A}\,\big) \big)
    \= P_s \big(\wt{A}\,\big) \= 0$. As $ (\Pi^t_s)^{-1} \big(\wt{\cN}\big) \sb
     (\Pi^t_s)^{-1} \big(\wt{A}\,\big) $, we see that $  (\Pi^t_s)^{-1} \big(\wt{\cN}\big) \ins \sN^t  $.

\no {\bf 2)} Given $r \ins [s,\infty]$,    Lemma \ref{lem_shift_inverse} (1) shows that
 $\cF^s_r \sb  \L_r \df \big\{ \wt{A} \sb \O^s \n : (\Pi^t_s)^{-1} \big(\wt{A}\,\big) \ins \cF^t_r \big\}
 \sb \ol{\L}_r  \df \big\{ \wt{A} \sb \O^s \n : (\Pi^t_s)^{-1} \big(\wt{A}\,\big) \ins \ol{\cF}^t_r \big\}$,
  which is clearly a sigma$-$field of $\O^s $. Since    $\sN^s  \sb \ol{\L}_r $ by Part (1),
  it follows that   $ \ol{\cF}^s_r \= \si (\cF^s_r \cp \sN^s) \sb \ol{\L}_t $,
  i.e. $ (\Pi^t_s)^{-1}  \big(\wt{A}\,\big)  \ins \ol{\cF}^t_r $ for any $\wt{A} \ins \ol{\cF}^s_r $.

Let $\tau \ins \ol{\cT}^s$.  For any $r \ins [s,\infty)$, since $\{\tau \ls r \} \ins \ol{\cF}^s_r$,
 we see   that $\big\{ \tau (\Pi^t_s) \ls r \big\} \= (\Pi^t_s)^{-1} \big( \{\tau \ls r \} \big) \ins \ol{\cF}^t_r $\,.
 Thus  $\tau (\Pi^t_s) $ is a $\ol{\cT}^t-$stopping time with values in $[s,\infty]$.

\no {\bf 3)}  Let $\wt{A} \ins \ol{\cF}^s $. Similar to Problem 2.7.3 of \cite{Kara_Shr_BMSC},
   there exists an $\cA  \ins \cF^s $   such that $ \wt{A}   \D    \cA  \ins  \sN^s   $.
   Since
   \beas
   (\Pi^t_s)^{-1}   \big(  \wt{A}   \D    \cA \big)
   & \tn \= & \tn  (\Pi^t_s)^{-1}   \big(  (\wt{A}   \Cp    \cA^c) \cp (\cA   \Cp    \wt{A}^c)  \big)
   \= (\Pi^t_s)^{-1}      (\wt{A}   \Cp    \cA^c) \cp (\Pi^t_s)^{-1} (\cA   \Cp    \wt{A}^c) \\
   & \tn \= & \tn  \Big( (\Pi^t_s)^{-1} \big(\wt{A}\,\big) \Cp \big( (\Pi^t_s)^{-1}  (\cA) \big)^c \Big) \cp
   \Big( (\Pi^t_s)^{-1}  (\cA ) \Cp \big( (\Pi^t_s)^{-1} \big(\wt{A}\,\big) \big)^c \Big)
     \=   (\Pi^t_s)^{-1} \big(\wt{A}\,\big) \D   (\Pi^t_s)^{-1} \big(\cA\big) ,
   \eeas
   we know from Part (1) that $(\Pi^t_s)^{-1} \big(\wt{A}\,\big) \D   (\Pi^t_s)^{-1} \big(\cA\big)$ is a    $P_t-$null set.
   So by Part (2) and Lemma \ref{lem_shift_inverse} (1),
   the    $\ol{\cF}^t-$measurable random variable    $ (\Pi^t_s)^{-1} \big(\wt{A}\,\big)   $ equals to
   the    $\cF^t-$measurable random variable    $ (\Pi^t_s)^{-1} \big(\cA\big)   $, \pas ~ Then Lemma \ref{lem_shift_inverse} (2) yields that    $    P_t \big( (\Pi^t_s)^{-1} \big(\wt{A}\,\big) \big)
    \= P_t  \big( (\Pi^t_s)^{-1}  (\cA ) \big)    \= P_s  ( \cA ) \= P_s  ( \wt{A} )   $.    \qed

  \begin{lemm}  \label{lem_F_version}

    Let   $t \ins [0,\infty)$.

   \no    \(1\)  For any $\xi   \ins    L^1 \big( \ol{\cF}^t , \hE\big)$ and $s   \ins   [t,\infty]$,
  $  
    E_t   \big[ \xi \big| \ol{\cF}^t_s   \big]    \n  =  \n   E_t    [ \xi  |  \cF^t_s  ] $,      \pas ~ Consequently,   an  $\hE-$valued martingale
      \(resp.\;local martingale or semi-martingale\)  with respect to $(\bF^t, P_t )$
 is also a martingale \(resp.\;local martingale or semi-martingale\) with respect to $\big(\ol{\bF}^t, P_t \big)$.

   \no  \(2\)  For any $s  \ins  [t,\infty]$ and any  $\hE-$valued, $\ol{\cF}^t_s- $measurable random variable $\xi$,
   there exists an $\hE-$valued, $\cF^t_s- $measurable random variable $ \wt{\xi} $ such that
 $ \wt{\xi} \= \xi $,    \pas ~

    \no  \(3\)    For any $\hE-$valued, $\ol{\bF}^t-$adapted  process $X \= \{X_s\}_{s \in [t,\infty)}$
   with  \pas ~ left-continuous paths, there exists an   $\hE-$valued, $\bF^t-$predictable  process
   $\wt{X} \= \big\{ \wt{X}_s \big\}_{s \in [t,\infty)}$
   such that $ \big\{ \o \ins \O^t \n : \wt{X}_s  (\o)  \n \ne \n X_s (\o) $
   for some $s \ins [t,\infty)  \big\}  \ins  \sN^t$.

\if{0}
    \no  \(4\)  For any $\hE-$valued, $\ol{\bF}^t-$progressively measurable  process $X \= \{X_s\}_{s \in [t,\infty)}$,
     there exists an   $\hE-$valued, $\bF^t-$predictable  process
   $\wt{X} \= \{\wt{X}_s\}_{s \in [t,\infty)}$
   such that $   \wt{X}_s  (\o)  \=  X_s (\o) $ for   $ds \ti P_t-$a.s. $(s,\o) \ins [t,\infty) \ti \O^t$.
\fi

\end{lemm}

     \no {\bf Proof:  1)}
  Let  $\xi   \ins    L^1 \big( \ol{\cF}^t , \hE\big)$ and $s   \ins   [t,\infty]$.
 For any $A  \ins   \ol{\cF}^t_s \= \si \big( \cF^t_s \cup \sN^t \big)$,
 similar to Problem 2.7.3 of \cite{Kara_Shr_BMSC},
   there exists an $\wt{A}  \ins \cF^t_s$   such that $ A   \D    \wt{A}  \ins  \sN^t  $.
 Then   we can deduce that
     $
       \int_A   \xi  d P_t    \n  =  \n    \int_{\wt{A}}  \xi  d P_t
 \n =  \n  \int_{\wt{A}}  E_t \big[ \xi \big|  \cF^t_s  \big]  d P_t
     \n    =  \n  \int_A   E_t \big[ \xi \big|  \cF^t_s  \big]   d P_t  
     $,
     which implies that $E_t  \big[ \xi \big| \ol{\cF}^t_s \big]
      \n  =  \n   E_t  \big[ \xi \big|  \cF^t_s  \big] $,  \pas ~ 

   \no  {\bf 2)}  Let   $s   \ins  [t,\infty]$ and let $\xi$ be   an $\hE-$valued, $\ol{\cF}^t_s-$measurable
 random variable. We first assume  $ \hE \= \hR $.
 For any $ n \ins \hN $, we set $  \xi_n \df ( \xi  \ld n ) \ve (-n) \ins  \ol{\cF}^t_s $
 and see from Part (1) that
    $  \wt{\xi}_n \df E_t \big[ \xi_n   \big|  \cF^t_s  \big]
    \= E_t \big[ \xi_n \big| \ol{\cF}^t_s  \big] \=  \xi_n  $,     \pas ~ Clearly, the random variable
    $ \wt{\xi} \df \Big( \lsup{n \to \infty} \wt{\xi}_n \Big)
    \b1_{\big\{ \lsup{n \to \infty} \wt{\xi}_n < \infty \big\}} $ is
      $\cF^t_s -$measurable  and   satisfies
    \beas
    \wt{\xi} \= \Big( \lsup{n \to \infty}  \xi_n \Big)
    \b1_{\big\{ \lsup{n \to \infty}  \xi_n < \infty \big\}} \= \xi \b1_{\{\xi < \infty\}} \= \xi , \q   \pas ~ \eeas
    When $\hE \= \hR^k $ for some $k \> 1$,  let $\xi^i$ be the i-th component of $\xi$, $i \= 1,\cds,k$. We denote by $\wt{\xi}^i$ the
 real-valued, $ \cF^t_s -$measurable random variable such that $\wt{\xi}^i \= \xi^i$,  \pas ~ Then $\wt{\xi} \= (\wt{\xi}^1,\cds,\wt{\xi}^k)$ is
 an $\hE-$valued, $\cF^t_s-$measurable   random variable such that
  $ \wt{\xi} \= \xi $,  \pas ~

        \no {\bf 3)} 
  Let  $X \= \{X_s\}_{s \in [t,\infty)}$   be    an $\hE-$valued, $\ol{\bF}^t-$adapted    process
  with   \pas ~ left-continuous paths. Like Part (2), it suffices to discuss the case of   $ \hE \= \hR $.
   For any $s \ins [t,\infty) \Cp \hQ   $,
   Part (2) shows that there exists a real-valued, $\cF^t_s-$measurable random variable $\sX_s $
   such that       $      \sX_s \=   X_s   $,     \pas ~ Define  $\cN \df  \big\{\o  \ins  \O^t  \n : $
   the path $  X_\cd (\o)$ is not left-continuous\big\}$\, \cup \,
  \Big( \underset{s \in  [t,\infty) \cap \hQ  }{\cup } \{ X_s  \n \ne \n  \sX_s \} \Big) \in \sN^t  $.

  For any $ n \ins  \hN$, set $t^n_i \= t\+ i/n$, $\fa i \ins \hN \cp \{0\}$.
   Since
  $     X^n_s \df  \sX_t \b1_{\{s=t\}} \+  \sum^{  n^2 }_{i=1}
   \sX_{ t^n_{i-1}  } \b1_{\{s \in ( t^n_{i-1} ,  t^n_i    ]\}} $, $s \ins [t,\infty)$
    is a real$-$valued, $\bF^t-$predictable process,  we see that
  $  \wt{X}_s \df \Big( \lsup{n \to \infty} X^n_s \Big) \b1_{\big\{ \lsup{n \to \infty} X^n_s < \infty\big\}}  $,
  $   s \ins [t,\infty)$   also  defines a real$-$valued, $\bF^t-$predictable process.

  Let $ \o \ins \cN^c$ and $s \ins (t,\infty)$. For any $n \ins \hN$ with $n \gs s \- t $,  since
   $s \ins \big(s_n \-\frac{1}{n}, s_n \big] $ with $ s_n \df t \+ \frac{\lceil n (s \- t) \rceil}{n}  
  $, one has  $  X^n_s (\o) \= \sX_{s_n - \frac{1}{n}   }(\o) = X_{s_n - \frac{1}{n}  }(\o)$.
  Clearly, $\lmt{n \to \infty} (s_n \- \frac{1}{n}) \= s$.
    As $ n \n \to \n  \infty$,   the left-continuity of $X$ shows that
 $   \lmt{n \to \infty} X^n_s (\o) =  \lmt{n \to \infty} X_{s_n \- \frac{1}{n}}(\o) = X_s(\o) $,
 which implies that  $\cN^c \sb \big\{\o \ins \O^t \n :  \wt{X}_s(\o) \= X_s(\o) $, $\fa s \ins [t,\infty)  \big\} $.
 \qed

 \if{0}

 \ss \no {\bf 4)} Let $X \= \{X_s\}_{s \in [t,\infty)}$  be  $\hE-$valued, $\ol{\bF}^t-$progressively measurable  process.
 Like Part (2), we only need to discuss the case   $ \hE \= \hR $.

   Suppose first that the process $\{X_s\}_{s \in [t,\infty)}$ is bounded.
   Since $ \sX_s  \df  \int_t^s X_r dr $, $s \ins [t,\infty)$ defines a  real-valued, $\ol{\bF}^t-$adapted
   process with   continuous paths,
   using  Part (3) we can find a real-valued, $\bF^t-$predictable process  $\wt{\sX}$
   such that $ \wh{\cN} \df  \big\{\o \ins \O^t \n : \wt{\sX}_s (\o) \n \ne \n \sX_s (\o) \hb{ for some }s \in [t, \infty) \big\} \ins  \sN^t $. 

   For any $n \ins \hN$,
   $ X^n_s  \df  n \big(\wt{\sX}_s \- \wt{\sX}_{(s-1/n) \vee t} \big) $, $s \ins  [t,\infty)$
   is also a   real-valued,  $\bF^t-$predictable   process. 
   It follows that $  \wt{X}_s \df  \Big( \lsup{n \to \infty} X^n_s \Big) \b1_{\big\{ \lsup{n \to \infty} X^n_s < \infty\big\}}  $,
     $s \ins [t,\infty)$    again  defines a real-valued, $\bF^t-$predictable  process.
   Set $\cD
      \df    \big\{(s,\o) \in [t,\infty) \times \O^t : \wt{X}_s (\o) \n \ne \n X_s (\o)  \big\}
      \ins \sB\big([t,\infty)\big) \oti \ol{\cF}^t  $.
     For any $\o \ins  \wh{\cN}^c$,  the Lebesgue differential theorem implies that
     \beas
     \lmt{n \to \infty} X^n_s (\o) \=  \lmt{n \to \infty} n \big( \sX_s (\o) \- \sX_{(s-1/n) \vee t} (\o) \big)
       \= \lmt{n \to \infty} n \int_{(s-1/n) \vee t}^s X_r (\o) dr \= X_s (\o) , \q  \hb{for a.e. }s \in [t,\infty) .
     \eeas
      So $ \wt{X}_s (\o) \= X_s (\o)$   or $\b1_\cD (s,\o) \= 0$ for a.e. $ s \ins [t,\infty)$.
      Then we know from     Fubini's Theorem   that
     $ ( ds \ti dP_t ) (\cD)  \= \int_{\O^t} \big(  \int_t^\infty  \b1_{\cD} (s,\o)  ds \big) dP_t(\o)
     \= \int_{ \wh{\cN}^c} \big(  \int_t^\infty   \b1_{\cD} (s,\o)  ds \big)  dP_t(\o) \= 0 $.

  \ss Next, let $\{X_s\}_{s \in [t,\infty)}$ be a general  real-valued, $\ol{\bF}^t-$progressively measurable process.
   For any   $m \ins  \hN$, the above argument shows that for process
   $ \big\{  (-m ) \ve \big( X_s \ld m \big) \big\}_{s \in [t,\infty)}$,
   there exists a real-valued, $\bF^t-$predictable process
   $\big\{\wt{X}^m_s\big\}_{s \in [t,\infty)}$ such that
   the set $\cD_m \df \big\{(s,\o) \ins [t,\infty) \times  \O^t  \n   :
   \wt{X}^m_s (\o) \n \ne \n (-m ) \ve  ( X_t (\o) \ld m  )  \big\}   \ins \sB\big([t,\infty)\big) \oti \ol{\cF}^t $
    satisfies  $(ds \ti dP_t)(\cD_m) \= 0 $.
   It follows that
     $  \wt{X}_s \df  \Big( \lsup{m \to \infty} \wt{X}^m_s \Big) \b1_{\big\{ \lsup{m \to \infty} \wt{X}^m_s < \infty\big\}}  $,
     $s \ins [t,\infty)$    is  a real-valued, $\bF^t-$predictable   process,
     and   the monotone convergence theorem implies that
     $(ds \ti dP_t) \big( \underset{m \in \hN }{\cup } \cD_m \big) \=
     \lmtu{m \to \infty} \, (ds \ti dP_t) \big( \underset{ i=1 }{\overset{m}{\cup} } \cD_i \big)\= 0 $.
     Given  $ (s,\o) \ins   \underset{m \in \hN }{\cap } \cD^c_m $, we can deduce that
      $ X_s (\o) \= \lmt{m \to \infty} (-m ) \ve  ( X_s (\o) \ld m  )
    \= \lmt{m \to \infty} \wt{X}^m_s (\o) $
    and thus obtain that    $ \wt{X}_s (\o) \= X_s (\o)  $. \qed

 \fi

\begin{eg} \label{eg_slm}
Suppose that $d \= 1$ and $\fg \df \underset{(t,x) \in (0,\infty) \times \hR^l}{\sup} g(t,x) \< \infty $.
Given   $(t,x ) \ins [0,\infty) \ti  \hR^l   $,
there exist $y \ins (0,\infty)$ and  $\fq \ins \hH^{2,{\rm loc}}_t (\hR)$ such that
$\a_s \df y \+ \int_t^s \fq_r dW^t_r$, $s \ins [t,\infty)$
is a positive strict local martingale with respect to $\big(\ol{\bF}^t,P_t\big)$
that satisfies $ E_t \big[\int_t^{\tau (t,x,\a)} g(r,\cX^{t,x}_r) dr\big] \< y$.
\end{eg}

 \no {\bf Proof:} Let $q \ins (1, \infty)  $. In light of \cite{Delbaen_Shirakawa_2002},
 the solution $\{\U_s\}_{s \in [t,\infty)}$ to
\beas
\U_s \= 1 \+ \int_t^s (\U_r)^q dW^t_r , \q s \ins [t,\infty)
\eeas
is positive strict local martingale with respect to $(\ol{\bF}^t,P_t)$,
So there exists a $\fs \in (0, \infty)$ such that $E_t [\U_{\fs}] \< 1 $.

 Let  $y \ins [1\+ \fg(\fs\-t),\infty)  $   and set $\fq^o_s \df (\U_s)^q \> 0$, $s \ins [t,\infty)$.
 For any $  n \ins  \hN$,
  the $\ol{\bF}^t-$stopping times  $\z_n \df \inf\{s \ins [t,\infty) \n : |\U_s \- 1| \> n \}$ satisfies that
  $E_t \big[ \int_t^{\z_n} (\fq^o_r)^2 dr \big] \= E_t \big[ |\U_{\z_n} \- 1|^2 \big] \ls n^2$.
  So it holds except on a $P_t-$null set $\cN_n$ that
 $ \int_t^{\z_n} (\fq^o_r)^2 dr \< \infty $.
 Since $\U$ is also a supermartingale such that $\U_\infty \df \lmt{s \to \infty} \U_s$ exists in $[0,\infty)$, \pas,
  the continuity of process $\U$ implies that for all $\o \ins \O^t$ except on a $P_t-$null set $\wt{\cN}$,
 $ \z_\fn (\o) \= \infty $ for some $\fn \= \fn (\o) \ins \hN$.
 Given $\o \ins   \Big( \underset{n \in \hN}{\cap} \cN^c_n \Big) \cap \wt{\cN}^c $, one has
 $ \int_t^\infty |\fq^o_r (\o)|^2 dr \= \int_t^{\z_\fn (\o)} |\fq^o_r (\o)|^2 dr \< \infty $.
  Thus,   $\fq^o \ins  L^{2,{\rm loc}}_t(\hR) $.

 Set $\a^o_s \df y \+ \int_t^s \fq^o_r dW^t_r
 \= \U_s \+ y \- 1 \> 0 $, $s \ins [t,\infty)$.  As it holds \pas ~ that
 \beas
  \int_t^s g(r,\cX^{t,x}_r) dr  \ls  \fg(\fs\-t) \< y \-  1 \+ \U_s 
  \= \a^o_s  ,    \q \fa s \ins [t,\fs] ,
  \eeas
 we see that $ \fs \< \tau_o \df \tau (t,x,\a^o)$, \pas ~

 Next, let us define $\fq_s \df \b1_{\{s \le \tau_o\}} \fq^o_s$, $ s \ins [t,\infty) $, which is clearly of
 $ L^{2,{\rm loc}}_t(\hR) $. Then    $\a_s \df y \+ \int_t^s \fq_r dW^t_r \= \a^o_{\tau_o \land s}$, $\fa s \ins [t,\infty)$
 and it follows that
 $\tau(t,x,\a) \= \tau_o   \> \fs$, \pas ~
 Since  $ \a^o \= \U\+y\-1$ is a positive continuous supermartingale,
 we can deduce from  the continuity of $\a$ and the optional sampling theorem that
   \beas
  \hspace{2.7cm}  E_t \Big[ \int_t^{\tau (t,x,\a)} g(r,\cX^{t,x}_r) dr \Big]
  \= E_t \big[ \a_{\tau (t,x,\a)} \big]  \= E_t \big[ \a^o_{\tau_o} \big]
  \ls E_t \big[\a^o_{\fs}\big] 
  \= E [ \U_{  \fs} \+  y \- 1    ] \< y . \hspace{2.7cm} \hb{\qed}
 \eeas

\subsection{Proofs of Starred Statements in Section \ref{sec:proofs}}

 \no {\bf Proof of \eqref{eq:d067}:}
   Given $s \ins [t,\infty)$, let $i_o$ be the largest integer such that $t_{i_o} \ls s  $.
  For any $i \= 1, \cds, i_o$ and $n \ins \hN$,
  since   $\{\tau \ls \z\} \ins \ol{\cF}^t_{\tau \land \z} \sb \ol{\cF}^t_\tau$,   one can deduce that
  $ \{ \ol{\tau} \ls s \} \= \big( \{\tau \ls \z\} \Cp \{\tau \ls s\} \big) \cp \Big(  \underset{i \le i_o}{\cup} \,
   \underset{n \in \hN}{\cup} \,  \cA^i_n \Cp \big\{\tau^i_n (\Pi^t_{t_i}) \ls s \big\} \Big)
   \cp \big( \{\tau > \z\} \Cp \cN_{t,x,\tau} \big) \ins \ol{\cF}^t_s $.
  So $\ol{\tau} \ins  \ol{\cT}^t$.

   For  $i,n \ins \hN$ and $\o \ins \cA^i_n \Cp \big( \cN^{i,n} \big)^c \sb \{\tau \> \z \= t_i\} $,
 since $\fX_{t_i} (\o) \ins \cO^i_n \= O_{\d_i (x^i_n)} (x^i_n) $
 and $\fY_{t_i} (\o) \ins \cD^i_n   \sb \big(   y^i_n  \-   \e/2   , \infty \big) $,
 applying \eqref{eq:c441} with $(\fx,\fx',\vs) \= \big(x^i_n, \fX_{t_i} (\o),  \tau^i_n \big)$,
   we see from \eqref{eq:c481} that
  \bea
   E_t \bigg[   \int_t^{\tau^i_n (\Pi^t_{t_i})} g (r, \fX_r)dr \Big|\ol{\cF}^t_{t_i} \bigg] (\o)
 & \tn \< & \tn  E_{t_i} \bigg[  \int_{t_i}^{\tau^i_n } g \big(r, \cX^{t_i,x^i_n}_r \big) dr  \bigg] \+ \e/2
 \+ \int_t^{t_i} g\big(r,\fX_r(\o)\big) dr
  \ls y^i_n   \+ \e/2 \+ \int_t^{t_i} g\big(r,\fX_r(\o)\big) dr   \nonumber   \\
 & \tn  \< & \tn  \fY_{t_i} (\o)  \+ \e \+ \int_t^{t_i} g\big(r,\fX_r(\o)\big) dr
 \= E_t \bigg[   \int_t^\tau  g (r, \fX_r)dr \Big|\ol{\cF}^t_{t_i} \bigg] (\o) \+ \e  .   \label{eq:d229}
 \eea
  Taking summation over $i,n \ins \hN $,
  one can deduce from \eqref{eq:c495} and the  monotone convergence theorem that
 \beas
 \q && \hspace{-1.2cm} E_t \bigg[ \b1_{\{\tau > \z \}} \int_t^{\ol{\tau}} g (r, \fX_r)dr  \bigg]
   \=    E_t \bigg[ \sum_{i,n \in \hN} \b1_{\cA^i_n} \int_t^{\tau^i_n (\Pi^t_{t_i})} g (r, \fX_r)dr  \bigg]
 \= \sum_{i,n \in \hN}  E_t \bigg[ \b1_{\cA^i_n} \int_t^{\tau^i_n (\Pi^t_{t_i})} g (r, \fX_r)dr  \bigg] \\
  && \=   \sum_{i,n \in \hN}  E_t \Bigg[ \b1_{\cA^i_n} E_t \bigg[  \int_t^{\tau^i_n (\Pi^t_{t_i})} g (r, \fX_r)dr \Big|\ol{\cF}^t_{t_i} \bigg] \Bigg]
 \ls \sum_{i,n \in \hN}  E_t \Bigg[ \b1_{\cA^i_n} E_t \bigg[   \int_t^\tau  g (r, \fX_r)dr \Big|\ol{\cF}^t_{t_i} \bigg] \Bigg]
 \+  \e    \\
  &&  \=    \sum_{i,n \in \hN}   E_t \bigg[  \b1_{\cA^i_n} \int_t^\tau  g (r, \fX_r)dr   \bigg] \+  \e
 \=     E_t \bigg[  \sum_{i,n \in \hN}  \b1_{\cA^i_n} \int_t^\tau  g (r, \fX_r)dr   \bigg] \+  \e
 \=  E_t \bigg[  \b1_{\{\tau > \z\}} \int_t^\tau  g (r, \fX_r)dr   \bigg]  \+  \e    .
 \eeas
 It   follows that
 $  E_t \Big[  \int_t^{\ol{\tau}} g (r, \fX_r)dr  \Big]
 \ls  E_t \big[     \int_t^\tau  g (r, \fX_r)dr   \big]  \+  \e   \ls y  \+  \e  $.
 Thus  $\ol{\tau} \ins \cT^t_x (y\+  \e)$. \qed

  \no {\bf Proof of \eqref{eq:d071}: 1)}
  Given $i\=1,\cds,\cI_o$, if $A \ins \ol{\cF}^t_{t_i}$ and if
  $ \{ \U_s \}_{ s \in [t_i,\infty) }$ is an $\ol{\bF}^t-$adapted continuous  process
  over period $[t_i,\infty)$   with $\U_{t_i} \= 0 $, \pas,
  one can easily deduce that $ \big\{ \b1_{\{s \ge t_i \} \cap A} \U_s \big\}_{s \in [t,\infty)}$
  is an $\ol{\bF}^t-$adapted continuous  process starting from $0$.
  Then we see from \eqref{eq:d141}   that
  $\ol{M}$ is an $\ol{\bF}^t-$adapted continuous   process 
  and $\ol{K}$ is an $\ol{\bF}^t-$adapted continuous  increasing  process with $\ol{K}_t \= K_t \= 0 $, \pas ~

  For any $\o \ins \cA^c_\sharp \Cp \cN^c$, $\ol{K}_\cd (\o) \= K_\cd (\o)$ is an increasing path;
  for $i\=1,\cds,\cI_o$ and $n \= 1, \cds, \fn_i$,
  it holds for any $\o \ins \cA^i_n \Cp (\Pi^t_{t_i})^{-1} \big( \big(\wt{\cN}^{i,n}\big)^c \big)$ that
  \beas
   \ol{K}_s (\o) \=   \b1_{\{s < t_i  \}} K_s (\o) \+
  \b1_{\{s \ge t_i  \}} \big( K^{i,n}_s  ( \Pi^t_{t_i} (\o)   )    \+ K_{t_i} (\o)  \+  \cK^{i,n}_s (\o) \big)    ,
  \q  s \ins [t,\infty)
  \eeas
  is an also increasing path. Thus, $\ol{K}$ has   increasing paths except the $P_t-$null set $\cN_\sharp$.

\no {\bf 2)} {\it To show   $ \ol{M} $ is a uniformly integrable martingale with respect to $(\ol{\bF}^t,P_t)$, we define }
\beas
 \q   \ol{\xi} \df   \sum^{\cI_o}_{i=1} \, \sum^{\fn_i}_{n=1} \,
    \b1_{ \cA^i_n }  \big( \xi^{i,n}   ( \Pi^t_{t_i}    ) \- \xi \+ M_{t_i} \- y^i_n \+ \d_i(x^i_n,y^i_n)    \big)
    \=   \b1_{\cA_\sharp}  ( M_\z \- \xi) \+ \sum^{\cI_o}_{i=1} \, \sum^{\fn_i}_{n=1} \,
    \b1_{ \cA^i_n }  \big( \xi^{i,n}   ( \Pi^t_{t_i}    )  \- y^i_n \+ \d_i(x^i_n,y^i_n)    \big) .
\eeas
Since $M$ is a uniformly integrable continuous $\big(\ol{\bF}^t, P_t\big)-$martingale,
we know from the optional sampling theorem (e.g. Theorem II.3.2 of \cite{revuz_yor})   that
$M_\z \= E_t \big[ \xi \big|\ol{\cF}^t_\z \big]$, \pas ~ It follows that
\bea \label{eq:d109}
E_t \big[|M_\z|\big] \= E_t \Big[ \big| E_t \big[ \xi \big|\ol{\cF}^t_\z \big] \big| \Big]
\ls E_t \Big[  E_t \big[ |\xi| \big|\ol{\cF}^t_\z \big]   \Big] \= E_t \big[ |\xi| \big] \< \infty .
\eea

 Given $i\=1,\cds,\cI_o$ and $n\=1,\cds,\fn_i$, as $\xi^{i,n}   $ is $\ol{\cF}^{t_i}-$measurable,
 Lemma \ref{lem_shift_inverse1b} (2) implies that $ \xi^{i,n}   (\Pi^t_{t_i}) $ is $\ol{\cF}^t-$measurable.
 By Proposition \ref{prop_shift_FP} (2), it holds for \pas ~ $\o \ins \O^t$ that
 \bea \label{eq:d081}
 E_t \big[| \xi^{i,n}   (\Pi^t_{t_i}) | \big| \ol{\cF}^t_{t_i} \big]  (\o)
 \= E_{t_i} \Big[ \big( \big| \xi^{i,n}   (\Pi^t_{t_i}) \big| \big)^{t_i,\o} \Big]
 \= E_{t_i} \big[   | \xi^{i,n}     |   \big] \< \infty  .
 \eea
 Taking expectation $E_t[\cd]$ and using \eqref{eq:d109}, one can deduce that
 \beas
 \q E_t \big[  \, | \ol{\xi} \,  | \big]  \ls E_t \big[ |\xi|\+|M_\z|\big] \+ \sum^{\cI_o}_{i=1} \, \sum^{\fn_i}_{n=1}
 \big(  E_t \big[   | \xi^{i,n}   (\Pi^t_{t_i}) | \big] \+ y^i_n \- \d_i(x^i_n,y^i_n)   \big)
 \= 2 E_t \big[ |\xi| \big] \+ \sum^{\cI_o}_{i=1} \, \sum^{\fn_i}_{n=1}
 \big( E_{t_i} \big[   | \xi^{i,n}     |   \big] \+ y^i_n    \big)   \< \infty ,
 \eeas
 which shows that $\ol{\xi} \ins L^1 \big(\ol{\cF}^t\big)$.

  Fix $s  \in [t , \infty)$. We    denote by $i_o$ the largest integer such that $t_{i_o} \ls s$.
  Let $i\=1,\cds,i_o$ and $n\=1,\cds,\fn_i$.
  In light of Proposition \ref{prop_shift_FP} (2),   there exists $\fN^{i,n} \ins \sN^t$ such that
\bea \label{eq:d083}
E_t \big[    \xi^{i,n}  ( \Pi^t_{t_i}  )     \big| \ol{\cF}^t_s \big](\o)
\= E_s \Big[  \big(  \xi^{i,n}  ( \Pi^t_{t_i}  )    \big)^{s,\o} \Big]
\= E_s \Big[  ( \xi^{i,n} )^{s,\Pi^t_{t_i}(\o)} \Big] , \q \fa \o \ins (\fN^{i,n})^c .
\eea
The last equality   uses  the fact that $ \Pi^t_{t_i} (\o \otis \wh{\o}) \= \Pi^t_{t_i} (\o) \otis \wh{\o}  $
for all $\wh{\o} \ins \O^s$. Applying Proposition \ref{prop_shift_FP} (2) again
and using \eqref{eq:d079}, we can find   $\wt{\cN}^{i,n} \ins \sN^{t_i}$ such that
\bea \label{eq:d085}
M^{i,n}_s (\wt{\o})  \= E_{t_i} \big[ \xi^{i,n} \big| \ol{\cF}^{t_i}_s \big] (\wt{\o})
\= E_s \big[ (\xi^{i,n})^{s,\wt{\o}} \big] , \q \fa \wt{\o} \ins \big(\wt{\cN}^{i,n}\big)^c .
\eea

  By  Lemma \ref{lem_shift_inverse1b} (1), $ \ol{\fN}^{i,n} \df \fN^{i,n}
  \cup \big((\Pi^t_{t_i})^{-1} \big(\wt{\cN}^{i,n}\big)\big) $ is a $P_t-$null set.
  For any $\o \ins \big(\ol{\fN}^{i,n}\big)^c \= (\fN^{i,n})^c \cap \big((\Pi^t_{t_i})^{-1} \big(\wt{\cN}^{i,n}\big)\big)^c
  \=  (\fN^{i,n})^c \cap \Big( (\Pi^t_{t_i})^{-1} \big( \big(\wt{\cN}^{i,n}\big)^c \big) \Big)  $,
  using \eqref{eq:d083} and taking $\wt{\o} \= \Pi^t_{t_i} (\o) $ in \eqref{eq:d085} yield that
  \beas
   E_t \big[    \xi^{i,n}  ( \Pi^t_{t_i}  )     \big| \ol{\cF}^t_s \big](\o)
  \= E_s \Big[  (\xi^{i,n})^{s,\Pi^t_{t_i}(\o)} \Big] \= M^{i,n}_s \big( \Pi^t_{t_i} (\o) \big)  .
  \eeas
  Then  we see from \eqref{eq:d107} that
  \bea
  \qq \qq && \hspace{-2cm} E_t \bigg[  \sum^{i_o}_{i=1} \, \sum^{\fn_i}_{n=1} \,
    \b1_{ \cA^i_n }  \big( \xi^{i,n}   ( \Pi^t_{t_i}    ) \- \xi \+ M_{t_i} \- y^i_n \+ \d_i(x^i_n,y^i_n)    \big) \Big| \ol{\cF}^t_s \bigg]
    \= \sum^{i_o}_{i=1} \, \sum^{\fn_i}_{n=1} \b1_{ \cA^i_n } \Big( E_t \big[    \xi^{i,n}   ( \Pi^t_{t_i}    )  \- \xi  \big| \ol{\cF}^t_s \big]   \+ M_{t_i} \- y^i_n \+ \d_i(x^i_n,y^i_n) \Big) \nonumber  \\
 &&    \=   \sum^{i_o}_{i=1} \, \sum^{\fn_i}_{n=1} \b1_{ \cA^i_n } \big( M^{i,n}_s \big( \Pi^t_{t_i}   \big)
 \- M_s \+ M_{t_i}
 \- y^i_n \+ \d_i(x^i_n,y^i_n) \big) \= \ol{M}_s \- M_s  , \q \pas  \label{eq:d087}
  \eea

  If $i_o \= \cI_o$, \eqref{eq:d087} just shows that $   E_t \big[ \, \ol{\xi} \, \big| \ol{\cF}^t_s \big]
  \= \ol{M}_s \- M_s  $,   \pas ~ and thus $ \ol{M}_s \= E_t \big[ \, \ol{\xi} \+ \xi   \big| \ol{\cF}^t_s \big] $.
   Suppose next that $i_o \< \cI_o$.
  For $i\=i_o\+1,\cds,\cI_o$ and $n\=1,\cds,\fn_i$,
   using an analogy to \eqref{eq:d081}, we can deduce from \eqref{eq:d079} and \eqref{eq:d075} that
   for \pas ~ $\o \ins \O^t$,
   $  E_t \big[  \xi^{i,n}   (\Pi^t_{t_i})  \big| \ol{\cF}^t_{t_i} \big]  (\o)
 \= E_{t_i} \big[    \xi^{i,n}       \big] \= y^i_n \-   \d_i(x^i_n,y^i_n) $.
  It follows that
  \beas
 \q && \hspace{-1.5cm} E_t \bigg[  \sum^{\cI_o}_{i=i_o+1} \, \sum^{\fn_i}_{n=1} \,
    \b1_{ \cA^i_n }  \big( \xi^{i,n}   ( \Pi^t_{t_i}    ) \- \xi \+ M_{t_i} \- y^i_n \+ \d_i(x^i_n,y^i_n)    \big) \Big| \ol{\cF}^t_s \bigg]   \\
 &&    \=  \sum^{\cI_o}_{i=i_o+1} \, \sum^{\fn_i}_{n=1}
    E_t \Big[ \b1_{ \cA^i_n } \Big( E_t \big[     \xi^{i,n}   ( \Pi^t_{t_i}    ) \- \xi \big| \ol{\cF}^t_{t_i} \big]
    \+ M_{t_i} \- y^i_n \+ \d_i(x^i_n,y^i_n)    \Big) \Big| \ol{\cF}^t_s \Big] \= 0 , \q \pas ,
  \eeas
  which together with \eqref{eq:d087} yields $ \ol{M}_s \= E_t \big[ \, \ol{\xi} + \xi  \big| \ol{\cF}^t_s \big]   $,
  \pas ~ again. Therefore,  $ \ol{M} $ is a uniformly integrable martingale
  with respect to $(\ol{\bF}^t,P_t)$.

 \no {\bf 3)} {\it We now prove that $E_t \big[ \, \ol{K}_* \big] \< \infty $.}

 Let $i\=1,\cds,\cI_o$ and $n\=1,\cds,\fn_i$.
 It is clear that $ \underset{s \in [t_i,\infty) \cap \hQ}{\sup} \,  K^{i,n}_s  $ is $\cF^{t_i}-$measurable.
 Since  the  continuity of $K^{i,n}$ implies that
  $ K^{i,n}_* (\wt{\o}) \= \underset{s \in [t_i,\infty) \cap \hQ}{\sup} \,  K^{i,n}_s (\wt{\o}) $
  for any $\wt{\o} \ins (\cN^{i,n})^c$, the random variable $K^{i,n}_*$ is $\ol{\cF}^{t_i}-$measurable
 and we thus know from Lemma \ref{lem_shift_inverse1b} (2)   that $ K^{i,n}_*   (\Pi^t_{t_i}) $ is $\ol{\cF}^t-$measurable.
 An analogy to \eqref{eq:d081} then shows that   for \pas ~ $\o \ins \O^t$
 \bea \label{eq:d095}
 E_t \big[  K^{i,n}_*  ( \Pi^t_{t_i} ) \big| \ol{\cF}^t_{t_i} \big]  (\o)
 \= E_{t_i} \big[  \big(  K^{i,n}_*  ( \Pi^t_{t_i} )   \big)^{t_i,\o} \big]
 \= E_{t_i} \big[   K^{i,n}_*    \big] \< \infty  .
 \eea

 As   $ \eta^i_n   \ls 2 \d_i(x^i_n,y^i_n) $ on $\cA^i_n$, it holds for any $s \ins [t,\infty)$ that
\beas
  \ol{K}_s    \=   \b1_{\cA_\sharp^c} K_s  \+
  \sum^{\cI_o}_{i=1}   \sum^{\fn_i}_{n=1} \b1_{   \cA^i_n  } \Big(  \b1_{ \{s < t_i  \}} K_s  \+
  \b1_{ \{s \ge t_i  \}}  \big(   K^{i,n}_s  ( \Pi^t_{t_i}    )   \+ K_{t_i}   \+  \cK^{i,n}_s    \big) \Big)
  \ls   K_* \+ \sum^{\cI_o}_{i=1}   \sum^{\fn_i}_{n=1}   \big(   K^{i,n}_*  ( \Pi^t_{t_i}    )
  \+ 2 \d_i(x^i_n,y^i_n)  \big)   .
\eeas
 Taking supremum over $s \ins [t,\infty) $ and  taking expectation $ E_t [\cd]$, we see from   \eqref{eq:d095} that
\beas
 \hspace{0.7cm}  E_t \big[ \, \ol{K}_* \big]
   \ls    E_t [K_*] \+ \sum^{\cI_o}_{i=1} \, \sum^{\fn_i}_{n=1}
  \big( E_t \big[  K^{i,n}_*  ( \Pi^t_{t_i} ) \big]  \+ 2 \d_i(x^i_n,y^i_n)  \big)
     \=     E_t [K_*] \+  \sum^{\cI_o}_{i=1} \, \sum^{\fn_i}_{n=1}
  \big( E_{t_i} \big[   K^{i,n}_*    \big]  \+ 2 \d_i(x^i_n,y^i_n)  \big) \< \infty . \hspace{0.7cm} \hb{\qed}
\eeas

\bibliographystyle{siam}
\bibliography{OSMC_bib}

\begin{thebibliography}{10}

\bibitem{AKKK_2017}
{\sc S.~Ankirchner, N.~Kazi-Tani, M.~Klein, and T.~Kruse}, {\em Stopping with
  expectation constraints: 3 points suffice},  (2017).
\newblock Available on \url{https://hal.archives-ouvertes.fr/hal-01525439}.

\bibitem{AKK_2015}
{\sc S.~Ankirchner, M.~Klein, and T.~Kruse}, {\em A verification theorem for
  optimal stopping problems with expectation constraints}, Applied Mathematics
  \& Optimization,  (2017), pp.~1--33.

\bibitem{ABG_1949}
{\sc K.~J. Arrow, D.~Blackwell, and M.~A. Girshick}, {\em Bayes and minimax
  solutions of sequential decision problems}, Econometrica, 17 (1949),
  pp.~213--244.

\bibitem{Balzer_Jansen_2002}
{\sc T.~Balzer and K.~Janssen}, {\em A duality approach to problems of combined
  stopping and deciding under constraints}, Math. Methods Oper. Res., 55
  (2002), pp.~431--446.

\bibitem{Bayraktar_Huang_2013}
{\sc E.~Bayraktar and Y.~Huang}, {\em On the multi-dimensional controller and
  stopper games}, SIAM J. Control Optim., 51 (2013), pp.~1263--1297.

\bibitem{OS_CRM}
{\sc E.~Bayraktar, I.~Karatzas, and S.~Yao}, {\em Optimal stopping for dynamic
  convex risk measures}, Illinois J. Math., 54 (2010), pp.~1025--1067.

\bibitem{Bayraktar_Miller_2016}
{\sc E.~Bayraktar and C.~Miller}, {\em Distribution-constrained optimal
  stopping}, to appear in Mathematical Finance,  (2016).
\newblock Available on \url{https://arxiv.org/abs/1604.03042}.

\bibitem{OSNE1}
{\sc E.~Bayraktar and S.~Yao}, {\em Optimal stopping for non-linear
  expectations---{P}art {I}}, Stochastic Process. Appl., 121 (2011),
  pp.~185--211.

\bibitem{OSNE2}
\leavevmode\vrule height 2pt depth -1.6pt width 23pt, {\em Optimal stopping for
  non-linear expectations---{P}art {II}}, Stochastic Process. Appl., 121
  (2011), pp.~212--264.

\bibitem{ROSVU}
\leavevmode\vrule height 2pt depth -1.6pt width 23pt, {\em On the robust
  optimal stopping problem}, SIAM J. Control Optim., 52 (2014), pp.~3135--3175.

\bibitem{RDG}
\leavevmode\vrule height 2pt depth -1.6pt width 23pt, {\em On the robust
  {D}ynkin game}, Ann. Appl. Probab., 27 (2017), pp.~1702--1755.

\bibitem{RDOSRT}
\leavevmode\vrule height 2pt depth -1.6pt width 23pt, {\em Optimal stopping
  with random maturity under nonlinear expectations}, Stochastic Process.
  Appl., 127 (2017), pp.~2586--2629.

\bibitem{BEES_2016}
{\sc M.~Beiglboeck, M.~Eder, C.~Elgert, and U.~Schmock}, {\em Geometry of
  distribution-constrained optimal stopping problems},  (2016).
\newblock Available on \url{https://arxiv.org/abs/1612.01488}.

\bibitem{BPZ_2016}
{\sc O.~Bokanowski, A.~Picarelli, and H.~Zidani}, {\em State-constrained
  stochastic optimal control problems via reachability approach}, SIAM J.
  Control Optim., 54 (2016), pp.~2568--2593.

\bibitem{BEI_2009}
{\sc B.~Bouchard, R.~Elie, and C.~Imbert}, {\em Optimal control under
  stochastic target constraints}, SIAM J. Control Optim., 48 (2009/10),
  pp.~3501--3531.

\bibitem{BER_2015}
{\sc B.~Bouchard, R.~Elie, and A.~R\'eveillac}, {\em B{SDE}s with weak terminal
  condition}, Ann. Probab., 43 (2015), pp.~572--604.

\bibitem{BET_2009}
{\sc B.~Bouchard, R.~Elie, and N.~Touzi}, {\em Stochastic target problems with
  controlled loss}, SIAM J. Control Optim., 48 (2009/10), pp.~3123--3150.

\bibitem{BMN_2014}
{\sc B.~Bouchard, L.~Moreau, and M.~Nutz}, {\em Stochastic target games with
  controlled loss}, Ann. Appl. Probab., 24 (2014), pp.~899--934.

\bibitem{riedel2012}
{\sc X.~Cheng and F.~Riedel}, {\em Optimal stopping under ambiguity in
  continuous time}, Mathematics and Financial Economics,  (2012), pp.~1--40.

\bibitem{CDK-2006}
{\sc P.~Cheridito, F.~Delbaen, and M.~Kupper}, {\em Dynamic monetary risk
  measures for bounded discrete-time processes}, Electron. J. Probab., 11
  (2006), pp.~no. 3, 57--106 (electronic).

\bibitem{CRS_1971}
{\sc Y.~S. Chow, H.~Robbins, and D.~Siegmund}, {\em Great expectations: the
  theory of optimal stopping}, Houghton Mifflin Co., Boston, Mass., 1971.

\bibitem{Delbaen_2006}
{\sc F.~Delbaen}, {\em The structure of m-stable sets and in particular of the
  set of risk neutral measures}, in In memoriam {P}aul-{A}ndr\'e {M}eyer:
  {S}\'eminaire de {P}robabilit\'es {XXXIX}, vol.~1874 of Lecture Notes in
  Math., Springer, Berlin, 2006, pp.~215--258.

\bibitem{Delbaen_Shirakawa_2002}
{\sc F.~Delbaen and H.~Shirakawa}, {\em No arbitrage condition for positive
  diffusion price processes}, Asia-Pac. Financ. Mark., 9 (2002).

\bibitem{DTX_2012}
{\sc J.~Detemple, W.~Tian, and J.~Xiong}, {\em An optimal stopping problem with
  a reward constraint}, Finance Stoch., 16 (2012), pp.~423--448.

\bibitem{Dokuchaev_1996}
{\sc N.~Dokuchaev}, {\em Random process optimal stopping in the problem with
  constraints}, Theory of probability and its applications, 41 (1996),
  pp.~761--768.

\bibitem{ETZ_2014}
{\sc I.~Ekren, N.~Touzi, and J.~Zhang}, {\em Optimal stopping under nonlinear
  expectation}, Stochastic Process. Appl., 124 (2014), pp.~3277--3311.

\bibitem{El_Karoui_1981}
{\sc N.~El~Karoui}, {\em Les aspects probabilistes du contr\^ole stochastique},
  in Ninth {S}aint {F}lour {P}robability {S}ummer {S}chool---1979 ({S}aint
  {F}lour, 1979), vol.~876 of Lecture Notes in Math., Springer, Berlin, 1981,
  pp.~73--238.

\bibitem{Follmer_Schied_2004}
{\sc H.~F{\"o}llmer and A.~Schied}, {\em Stochastic finance}, vol.~27 of de
  Gruyter Studies in Mathematics, Walter de Gruyter \& Co., Berlin,
  extended~ed., 2004.
\newblock An introduction in discrete time.

\bibitem{Horiguchi_2001b}
{\sc M.~Horiguchi}, {\em Stopped {M}arkov decision processes with multiple
  constraints}, Math. Methods Oper. Res., 54 (2001), pp.~455--469 (2002).

\bibitem{Horiguchi_2001c}
{\sc M.~Horiguchi, M.~Kurano, and M.~Yasuda}, {\em Markov decision processes
  with constrained stopping times}, Proceedings of the IEEE Conference on
  Decision and Control, 1 (2000), pp.~706--710.

\bibitem{Kallblad_2017}
{\sc S.~K\"allblad}, {\em A dynamic programming principle for
  distribution-constrained optimal stopping},  (2017).
\newblock Available on \url{https://arxiv.org/abs/1703.08534}.

\bibitem{Kara_Shr_BMSC}
{\sc I.~Karatzas and S.~E. Shreve}, {\em Brownian motion and stochastic
  calculus}, vol.~113 of Graduate Texts in Mathematics, Springer-Verlag, New
  York, second~ed., 1991.

\bibitem{Kara_Shr_MF}
{\sc I.~Karatzas and S.~E. Shreve}, {\em Methods of mathematical finance},
  vol.~39 of Applications of Mathematics (New York), Springer-Verlag, New York,
  1998.

\bibitem{Karatzas_Sudderth_2001}
{\sc I.~Karatzas and W.~D. Sudderth}, {\em The controller-and-stopper game for
  a linear diffusion}, Ann. Probab., 29 (2001), pp.~1111--1127.

\bibitem{Kara_Zam_2005}
{\sc I.~Karatzas and I.~M. Zamfirescu}, {\em Game approach to the optimal
  stopping problem}, Stochastics, 77 (2005), pp.~401--435.

\bibitem{Kara_Zam_2008}
\leavevmode\vrule height 2pt depth -1.6pt width 23pt, {\em Martingale approach
  to stochastic differential games of control and stopping}, Ann. Probab., 36
  (2008), pp.~1495--1527.

\bibitem{Kennedy_1982}
{\sc D.~P. Kennedy}, {\em On a constrained optimal stopping problem}, J. Appl.
  Probab., 19 (1982), pp.~631--641.

\bibitem{Lempa_J_2012}
{\sc J.~Lempa}, {\em Optimal stopping with information constraint}, Appl. Math.
  Optim., 66 (2012), pp.~147--173.

\bibitem{Liang_G_2015}
{\sc G.~Liang}, {\em Stochastic control representations for penalized backward
  stochastic differential equations}, SIAM J. Control Optim., 53 (2015),
  pp.~1440--1463.

\bibitem{Liang_Wei_2016}
{\sc G.~Liang and W.~Wei}, {\em Optimal switching at {P}oisson random
  intervention times}, Discrete Contin. Dyn. Syst. Ser. B, 21 (2016),
  pp.~1483--1505.

\bibitem{LSMS_1995}
{\sc F.~J. L\'opez, M.~San~Miguel, and G.~Sanz}, {\em Lagrangean methods and
  optimal stopping}, Optimization, 34 (1995), pp.~317--327.

\bibitem{Makasu_2009}
{\sc C.~Makasu}, {\em Bounds for a constrained optimal stopping problem},
  Optim. Lett., 3 (2009), pp.~499--505.

\bibitem{Menaldi_Robin_2016}
{\sc J.~L. Menaldi and M.~Robin}, {\em On some optimal stopping problems with
  constraint}, SIAM J. Control Optim., 54 (2016), pp.~2650--2671.

\bibitem{Miller_C_2017a}
{\sc C.~Miller}, {\em Nonlinear {PDE} approach to time-inconsistent optimal
  stopping}, SIAM J. Control Optim., 55 (2017), pp.~557--573.

\bibitem{Moustakides_1986}
{\sc G.~V. Moustakides}, {\em Optimal stopping times for detecting changes in
  distributions}, Ann. Statist., 14 (1986), pp.~1379--1387.

\bibitem{Neveu_1975}
{\sc J.~Neveu}, {\em Discrete-parameter martingales}, North-Holland Publishing
  Co., Amsterdam, revised~ed., 1975.
\newblock Translated from the French by T. P. Speed, North-Holland Mathematical
  Library, Vol. 10.

\bibitem{NZ_2015}
{\sc M.~Nutz and J.~Zhang}, {\em Optimal stopping under adverse nonlinear
  expectation and related games}, Ann. Appl. Probab., 25 (2015),
  pp.~2503--2534.

\bibitem{Pedersen_Peskir_2016}
{\sc J.~L. Pedersen and G.~Peskir}, {\em Optimal mean-variance selling
  strategies}, Math. Financ. Econ., 10 (2016), pp.~203--220.

\bibitem{Pedersen_Peskir_2017}
{\sc J.~L. Pedersen and G.~Peskir}, {\em Optimal mean-variance portfolio
  selection}, Math. Financ. Econ., 11 (2017), pp.~137--160.

\bibitem{Peskir_2012}
{\sc G.~Peskir}, {\em Optimal detection of a hidden target: the median rule},
  Stochastic Process. Appl., 122 (2012), pp.~2249--2263.

\bibitem{Peskir_Shiryaev_2006}
{\sc G.~Peskir and A.~Shiryaev}, {\em Optimal stopping and free-boundary
  problems}, Lectures in Mathematics ETH Z\"urich, Birkh\"auser Verlag, Basel,
  2006.

\bibitem{Pontier_Szpirglas_1984}
{\sc M.~Pontier and J.~Szpirglas}, {\em Optimal stopping with constraint}, in
  Analysis and optimization of systems, {P}art 2 ({N}ice, 1984), vol.~63 of
  Lect. Notes Control Inf. Sci., Springer, Berlin, 1984, pp.~82--91.

\bibitem{revuz_yor}
{\sc D.~Revuz and M.~Yor}, {\em Continuous martingales and {B}rownian motion},
  vol.~293 of Grundlehren der Mathematischen Wissenschaften [Fundamental
  Principles of Mathematical Sciences], Springer-Verlag, Berlin, third~ed.,
  1999.

\bibitem{Riedel_2009}
{\sc F.~Riedel}, {\em Optimal stopping with multiple priors}, Econometrica, 77
  (2009), pp.~857--908.

\bibitem{Shiryayev_1978}
{\sc A.~N. Shiryayev}, {\em Optimal stopping rules}, Springer-Verlag, New
  York-Heidelberg, 1978.
\newblock Translated from the Russian by A. B. Aries, Applications of
  Mathematics, Vol. 8.

\bibitem{Snell_1952}
{\sc J.~L. Snell}, {\em Applications of martingale system theorems}, Trans.
  Amer. Math. Soc., 73 (1952), pp.~293--312.

\bibitem{Stroock_Varadhan}
{\sc D.~W. Stroock and S.~R.~S. Varadhan}, {\em Multidimensional diffusion
  processes}, Classics in Mathematics, Springer-Verlag, Berlin, 2006.
\newblock Reprint of the 1997 edition.

\bibitem{Urusov_2005}
{\sc M.~A. Urusov}, {\em On a property of the moment at which brownian motion
  attains its maximum and some optimal stopping problems}, Theory of
  Probability \& Its Applications, 49 (2005), pp.~169--176.

\end{thebibliography}

\end{document}